\documentclass[hidelinks]{amsart}

\usepackage{macros}

\begin{document}

\title{A dendroidal approach to operadic right modules and manifold calculus}

\author{Miguel Barata}
\address{Purdue University, 150 North University Street, IN 47907, West Lafayette}
\email{barata@purdue.edu}

\maketitle

\vspace{-2.5em}

\begin{abstract}
\leftskip4em
\rightskip4em
   In this work we study the homotopy theory of the category $\RMod{\mathcal{P}}$ of right modules over a simplicial operad $\mathcal{P}$ via the formalism of forest spaces $\fSpaces$, as introduced by Heuts, Hinich and Moerdijk in \cite{HeutsHinnichMoerdijk}. In particular, we show that, for $\mathcal{P}$ a simplicial closed and $\Sigma$-free operad, there exists a Quillen equivalence between the projective model structure on $\RMod{\mathcal{P}}$, and the contravariant model structure on the slice category $\fSpaces_{/N\mathcal{P}}$ over the dendroidal nerve of $\mathcal{P}$. 
   
   As an application, we comment on how this result can be used to simplify the computation of derived mapping spaces between operadic right modules, and use this formalism to analyse the components and layers of the Goodwillie--Weiss tower coming from embedding calculus.
\end{abstract}

\tableofcontents


\section{Introduction}

The purpose of the present work is to provide an analysis of the theory of \textit{simplicial right $\mathcal{P}$-modules} over a simplicial operad $\mathcal{P}$ via the framework of dendroidal sets and dendroidal spaces. In this way, we will show in a precise sense that the theory of simplicial operadic right modules corresponds to the theory of right fibrations on a presheaf category modelled on forests, which we will discuss in more detail below. In particular, we will see how this result will allow us to compute derived mapping spaces of operadic right modules by instead computing a derived mapping space in this presheaf category of forests, which, as we will explain, tends to be an easier task. For any notation that we use, see the last section of the Introduction.

The classical theory of operads has had a great influence in the recent development and understanding of homotopy theory. The term was first coined in \cite{MayGeometry} by May in his seminal work on the structure of iterated loop spaces where the current definition of an uncoloured operad was given. However, the idea of behind operads as a way of studying homotopy coherent algebraic structures already appeared in previous work of Stasheff \cite{StasheffHomotopy} on the homotopical associaitivity of operations in topological spaces. We remark that the effectiveness of the theory of operads has also extended itself beyond homotopy theory: just to name two examples, the study of operads governing topological fields theories \cite{GetzlerKapranov, GetzlerOperads, GetzlerBV}, and the study of deformation problems in algebra, together with an extension of Koszul duality for more general algebraic structures \cite{GinzburgKapranov, KontsevichSoibelamn, LodayVallette}.

As of the last few years, the theory of operads has also surfaced in certain areas of differential topology and low-dimensional topology \cite{AroneTurchinRational, SalvatoreTurchin}, in particular in the study of the homotopy type of embedding spaces and its connection with the $\mathbb{E}_d$-operad. Since this will be our main motivating example for this text, let us briefly discuss it: let $M$ be a $d$-dimensional smooth manifold with no boundary, and consider the following associated collection of spaces of smooth embeddings 
$$ \mathbb{E}_M(k) = \mathsf{Emb}\left( \underline{k} \times D^d, M\right)$$
for each $k \geq 0$, as well as the spaces $\mathbb{E}^{\mathsf{fr}}_d(k) = \mathsf{Emb}( \underline{k} \times D^d, D^d)$. Associated with this data we have a "right action" of $\mathbb{E}^{\mathsf{fr}}_d$ on $\mathbb{E}_M$
\begin{equation}
\mathbb{E}_M(k) \times \mathbb{E}^{\mathsf{fr}}_d(j_1) \times \cdots \times \mathbb{E}^{\mathsf{fr}}_d(j_k) \longrightarrow \mathbb{E}_M(j_1 + \cdots + j_k)
\label{right_action}
\end{equation}
given by replacing the $i^{th}$ embedded disk $D^d \to M$ by the $j_i$ disks coming from $\mathbb{E}_d^{\mathsf{fr}}(j_i)$. If we additionally note that $\mathbb{E}_d^{\mathsf{fr}}$ is what is usually known as the \textit{framed little $d$-disks operad}, we can contextualize the existence of this action within the operadic world by saying that $\mathbb{E}_M$ is an operadic right module over the operad $\mathbb{E}^{\mathsf{fr}}_d$. As this construction is clearly natural on the manifold component with respect to smooth embeddings, we get a map
$$ \alpha_{M,N} : \mathsf{Emb}(M,N) \longrightarrow \Map_{\mathbb{E}^{\mathsf{fr}}_d}(\mathbb{E}_M, \mathbb{E}_N)$$
where $N$ is a smooth manifold of dimension at least $d$, which can still be seen as a right $\mathbb{E}^{\mathsf{fr}}_d$-module. Here the right hand side is a derived mapping space for right $\mathbb{E}^{\mathsf{fr}}_d$-modules, and from the perspective of embedding calculus, this term coincides with the limit term $T_{\infty} \mathsf{Emb}(M,N)$ of the Taylor tower coming from the Goodwillie--Weiss embedding calculus, and $\alpha_{M,N}$ is the canonical comparison map. This operadic approach to manifold calculus was introduced by Weiss and Boavida de Brito in \cite{WeissBritoHomotopySheaves} via homotopy sheaves, and is also present in the work of Arone--Turchin \cite{AroneTurchinRational} on spaces of long knots and, more recently, in the work of Krannich and Kupers \cite{KrannichKupersDisc}. One of the main questions in the subject concerns whether $\alpha_{M,N}$ is a weak homotopy equivalence, for which certain definitive answers are known: 
\begin{itemize}
    \item By the famous disjunction results of Goodwillie and Klein \cite{GoodwillieMultiple, GoodwillieKleinMultipleSmooth, GoodwillieKleinPoincareMultiple}, we know that $\alpha_{M,N}$ is a weak homotopy equivalence if, for instance, the codimension $\dim(N) - \dim(M) $ is at least 3.
    \item For manifolds of dimension at most 2 it is known that $\alpha_{M, N}$ is a weak homotopy equivalence by the results of Krannich and Kupers in \cite{KrannichKupersSurfaces}.
    \item If $M$ is $\mathbb{R}^d$, then $\alpha_{M,N}$ is a weak homotopy equivalence, as explained after Definition 4.2 of \cite{WeissBritoHomotopySheaves} in the language of $\mathcal{J}_{\infty}$-sheaves.
\end{itemize}

 Going back to the world of operads and in order to put things in more concrete terms, suppose we are given a coloured operad $\mathcal{P}$ with set of colours $C$ and with values in a symmetric monoidal category $\mathcal{V}$, which will always either mean $\Sets$ or $\sSets$ with the cartesian monoidal structure. Then a \textit{right $\mathcal{P}$-module} is a collection of objects 
\begin{equation}
 M(c_1, \ldots, c_n) \in \mathcal{V}
 \label{description}
\end{equation}
where $n \geq 1$ and $c_i \in C$ are colours of the operad $\mathcal{P}$, which admit right action maps very similar to the ones in \eqref{right_action}, compatible with the unit operations and operadic composition of $\mathcal{P}$ is a sensible way. Although operadic right modules have not been as studied as operadic algebras, this form of operadic action already has appeared in the literature: for instance, the work by Kapranov and Manin on Morita equivalence of categories of modules \cite{KapranovManin}, the analysis of spaces of long knots by Arone--Turchin in \cite{AroneTurchinRational} and the study of configuration spaces of finite products of manifolds by Dwyer--Hess--Knudsen \cite{DwyerHessKnudsen}. We also want to mention the book by Fresse \cite{FresseModules}, which develops the theory of operadic right modules for algebraic operads in a different way than the one we will consider.

As we mentioned in the initial paragraph, our intent with this work is to understand the theory of simplicial operadic right modules, but we will also be interested in realizing this analysis within the more general context of $\infty$-operads. With this in mind, we will from now on use dendroidal sets and dendroidal spaces as models for $\infty$-operads; by \cite{HeutsHinnichMoerdijk, BarwickOperator, ChuHaugsengHeuts}, it is know that both of these models are equivalent to Lurie's definition of homotoyp $\infty$-operads.

The category of \textit{dendroidal sets} is a presheaf category $\dSets$ first considered by Moerdijk and Weiss \cite{MoerdijkWeissDendroidal}, based on a certain diagram category $\Tree$, which has objects given by certain finite, non-planar, rooted trees. One can construct a functor 
$$N : \mathsf{Op} \longrightarrow \dSets$$
which associates to an operad $\mathcal{P}$ in $\Sets$ a dendroidal set $N\mathcal{P} \in \dSets$ which captures the operadic structure. Therefore, one should think of the dendroidal presheaf categories as playing a similar role to that of simplicial objects in the study of categories. The literature on dendroidal sets and dendroidal spaces is quite extensive, and the recent textbook \cite{HeutsMoerdijkDendroidal} gathers some of the essential aspects of this theory.

Before stating our first results, let us make a small comment on our set up using dendroidal sets. In our attempt to relate right $\mathcal{P}$-modules to presheaves on $\Tree$, we will see that the right action by $\mathcal{P}$ is, informally speaking, related to what happens to a tree near its (unique) root edge. On the other hand, the collection of objects \eqref{description} defining an operadic module involves multiple inputs $(c_1, \ldots, c_n)$ instead of just a single colour. These two observations together will lead us to notice that trees are insufficient for our purposes since they have a single root edge: instead, we can enlarge $\Tree$  to a more general category of \textit{forests} $\For$, which has as objects finite collections of trees. This category was first considered by Heuts--Hinich--Moerdijk in \cite{HeutsHinnichMoerdijk} when trying to relate dendroidal $\infty$-operads to Lurie's $\infty$-operads. We will write from now  $\fSets$ and $\fSpaces$ for the respective presheaf categories of \textit{forest sets} and \textit{forest spaces}. We remark that although these are not the same thing as the dendroidal presheaf categories, the morphisms in $\For$ are defined in terms of the morphisms in $\Tree$ and therefore we don't deviate ourselves too much from the dendroidal formalism. In particular, every dendroidal set defines a forest set, by seeing each forest as a coproduct of trees.

\vspace{0.5em}

\noindent\textbf{Operadic right modules via dendroidal objects}

Our first step in trying to understand operadic right modules will be to show that the slice categories $\fSets_{/V}$ admit a version of the contravariant model structure, as we make precise now.

\vspace{0.5em}

\noindent\textbf{\cref{ctv_character}.}\textit{
    Let $V$ be a forest set. The slice category $\fSets_{/V}$ admits a left proper, cofibrantly generated model structure such that the following conditions hold:
    \begin{enumerate}[label=(\alph*)]
        \item The cofibrations are the normal monomorphisms over $V$.
        \item The weak equivalences between normal forest sets are the maps $A \to B$ over $V$ such that, for any right fibration $X \to V$, the map
        $$ \Hom_V(B,X) \longrightarrow \Hom_V(A,X)$$
        is a categorical equivalence of $\infty$-categories.
        \item The fibrant objects are the forest right fibrations over $V$.
        \item The fibrations between fibrant objects are the forest right fibrations.
    \end{enumerate}}

\vspace{0.5em}

We will write $\left( \fSets_{/V} \right)_{\ctv}$ for the contravariant model structure. For the reader who is less knowledgeable on dendroidal sets, we note that although a version of the covariant model structure for dendroidal sets exists (and is easily checked to also exist for the forest context), it is not clear that this implies the existence of the contravariant model structure, contrary to what happens with simplicial sets. The problem is that both of the diagram categories $\Tree$ and $\For$ don't admit an involution which interchanges left fibrations and right fibrations; this is also evident by the very different behaviour of dendroidal and forest left and right fibrations. Therefore the proof of \cref{ctv_character} will involve some technical work with forest sets and right fibrations, to which we fully dedicate Section 3.3.

We also define a version of the contravariant model structure for the overcategories $\fSpaces_{/V}$ with $V$ now being a forest space. In this case, depending on whether we consider $\fSpaces_{/V}$ to carry the projective or the Reedy model structure, we get two versions of the contravariant model structure by left Bousfield localization at the morphisms
$$ \rho(F) \longrightarrow F$$
over $V$ for every forest $F$, where $\rho(F) \subseteq F$ is the forest of the root edges of $F$. We will call these mode structures the \textit{projective/Reedy contravariant model structures}. We will show that the fibrant objects of these model categories are the maps $p \colon X \to V$ which are projective/Reedy fibrations and such that the diagram below, induced by $\rho(F) \to F$ for any forest $F$,
$$
\begin{tikzcd}
X(F) \arrow[d] \arrow[r] & X(\rho(F)) \arrow[d] \\
V(F) \arrow[r]           & V(\rho(F))    
\end{tikzcd}
$$
is a homotopy pullback square. We will say that $p$ is a \textit{right fibration} if it satisfies this last condition.

One of our main results is the following Quillen equivalence between the homotopy theory of simplicial right $\mathcal{P}$-modules and the homotopy theory of right fibrations in $\fSpaces_{/N\mathcal{P}}$, under some very mild restrictions on the simplicial operad $\mathcal{P}$.

\vspace{0.5em}

\noindent\textbf{\cref{eqv_ctv_RMod}.}\textit{
    Let $\mathcal{P}$ be a closed $\Sigma$-free simplicial operad. Then there exists a Quillen equivalence}
    $$
    \begin{tikzcd}
	{\RMod{\mathcal{P}} } & {\left( \fSpaces_{/N\mathcal{P}} \right)_{\mathsf{P}, \ctv}}
	\arrow["{N_{\mathcal{P}}}"', shift right, from=1-1, to=1-2]
	\arrow["{\tau_{\mathcal{P}}}"', shift right, from=1-2, to=1-1]
\end{tikzcd}$$
    \textit{between the projective model structure on the category of right $\mathcal{P}$-modules $\RMod{\mathcal{P}}$, and the projective contravariant model structure on $\fSpaces_{/N\mathcal{P}}$.}

\vspace{0.5em}

In fact, we also show in \cref{monoidal_eqv} that we can upograde the equivalence above to one between monoidal model categories, where the categories in question are equipped with the \textit{concatenation product}. On the side of right $\mathbb{P}$-modules, this will correspond to the coproduct of manifolds when restricted to the modules $\mathbb{E}_M$, and on the forest side this monoidal structure is induced by the operations of concatenating forests.

Let us make a few remarks concerning the statement above. Firstly, one should think of the adjunction $(\tau_{\mathcal{P}}, N_{\mathcal{P}})$ as being a version of the nerve adjunction but for operadic right modules. Informally speaking, the forest space $N_{\mathcal{P}} M$ simultaneously keeps track of the elements $m$ of the module $M$, and of the the operations of $\mathcal{P}$ that can act on said $m$.

Secondly, a result of this kind has already been proven for operadic algebras and left fibrations in dendroidal sets and spaces \cite{HeutsMoerdijkLeft, BoavidaMoerdijk}. In that sense, the comparison above gives a different perspective on the different behaviour of operadic algebras and operadic right modules, since left fibrations and right fibrations behave quit differently in the dendroidal and forest context, as we already mentioned. Of course, this was already known: an operadic algebra is an example of an operadic left module, and these are rather different from operadic right modules.

Our second remark concerns the proof of the actual statement. The difficult step in showing that the adjunction $(\tau_{\mathcal{P}}, N_{\mathcal{P}})$ is a Quillen equivalence will be in proving that the unit map
$$ X \longrightarrow N_\mathcal{P} \tau_{\mathcal{P}} X$$
is a contravariant weak equivalence, whenever $X$ is a projectively cofibrant forest space over $V$. In order to do so, we will use the method explained in \cite{HeutsMoerdijkDendroidal} concerning \textit{simplicial diagrams of model categories}. Using this, we will show that, for $V$ a fixed forest space, one can upgrade contravariant weak equivalences in $\left( \fSets_{/V_i}\right)_\ctv$ for each simplicial degree $i \geq 0$ to contravariant weak equivalences in the contravariant model structure in $\fSpaces_{/V}$. The problem then becomes one of understanding contravariant weak equivalences in $\fSets$, for which we give a useful criteria in \cref{root_anodyne_ctv_we} that will be enough for our purposes, whereas in $\fSpaces$ it tends to be more complicated to check whether given maps are contravariant weak equivalences.

We finalize this introduction by presenting some consequences of \cref{eqv_ctv_RMod}. As we showed in our motivating example on manifold calculus, the computation of derived mapping spaces for operadic right modules is an interesting question. The next result follows immediately from \cref{eqv_ctv_RMod}, and shows that this kind of computation can be carried in the world of forest spaces.

\vspace{0.5em}

\noindent\textbf{\cref{compute_map_spaces}.}\textit{
    Let $\mathcal{P}$ be a closed simplicial $\Sigma$-free operad, and suppose $M$ and $L$ are simplicial right $\mathcal{P}$-modules which are projectively fibrant. Then the simplicial map between the derived mapping spaces
    $$ \Map_{\mathcal{P}}(M,L) \longrightarrow \Map_{N \mathcal{P}}( N_\mathcal{P} M, N_{\mathcal{P}} L)$$
    is a weak homotopy equivalence. Here the left hand side is computed using the projective model structure on $\RMod{\mathcal{P}}$, and the right hand side using the projective contravariant model structure on $\fSpaces_{/N\mathcal{P}}$.}

\vspace{0.5em}
    
To put into perspective the usefulness of such a result, we observe that being cofibrant in $\RMod{\mathcal{P}}$ is usually a difficult task, as it is already the case when one wants to find a cofibrant simplicial operad. However, due to the wealth of model structures of forest spaces, the mapping space on the right hand side above can be computed more easily. For instance, one can, up to replacing $N\mathcal{P}$ by a Reedy fibrant object, compute the right hand side using  the Reedy contravariant model structure, where the cofibrant objects are much more plentiful and easier to work with (and under the restrictions on the operad $\mathcal{P}$ in the statement of the theorem, every object in $\fSpaces_{/N\mathcal{P}}$ is automatically cofibrant). We will follow up on this question on computing derived mapping spaces of operadic right modules in our next work.

The next result gives a homotopy invariance of the category of operadic right modules on the operadic input up to Dwyer--Kan equivalence.

\vspace{0.5em}

\noindent\textbf{\cref{DK_eqv}.}\textit{
    Let $\Phi \colon \mathcal{P} \to \mathcal{Q}$ be an operad map between closed simplicial $\Sigma$-free operads. If $\Phi$ is a Dwyer-Kan equivalence, then there is a Quillen equivalence}
    \[\begin{tikzcd}
	{\RMod{\mathcal{P}}} & {\RMod{\mathcal{Q}}}
	\arrow["{\Phi^*}"', shift right, from=1-1, to=1-2]
	\arrow["{\Phi_!}"', shift right, from=1-2, to=1-1]
\end{tikzcd}\]
\textit{between the respective projective model structures.}

\vspace{0.5em}

Finally, we also prove that the category of right modules over a normal $\infty$-operad $V$ (note that by \cref{eqv_ctv_RMod} one can think of these as right fibrations over $V$) is Quillen equivalent to the category of operadic right modules over an actual operad $w_!(V)$.

\vspace{0.5em}

\noindent\textbf{\cref{strictification}.}\textit{
    Let $V$ be a normal closed $\infty$-operad. Then there exists a simplicial operad $w_!(V)$ and a natural zig-zag of Quillen equivalences}
\[\begin{tikzcd}[cramped]
	{\left( \fSets_{/V}\right)_{\ctv}} & \cdots & {\mathsf{RMod}_{w!(V)}.}
	\arrow[shift right, from=1-1, to=1-2]
	\arrow[shift right, from=1-2, to=1-1]
	\arrow[shift right, from=1-2, to=1-3]
	\arrow[shift right, from=1-3, to=1-2]
\end{tikzcd}\]

\vspace{0.5em}

\noindent\textbf{A dendroidal approach to Goodwillie--Weiss manifold calculus}

We now move on to explaining the main results of Chapter 5. This section of the text concerns with applying our results on right fibrations of forest spaces in order to give a description  of the Goodwillie--Weiss tower for studying embedding spaces of smooth manifolds.

Our first result concerns the existence of the analogue of the Goodwillie--Weiss tower for any unital simplicial operad $\mathcal{P}$. More explicitly, we show that, given any right $\mathcal{P}$-modules $M$ and $L$, there exist towers converging to the derived mapping spaces 
$$\Map_{\mathcal{P}}(M,L) \hspace{2.5em} \mathrm{and}\hspace{2.5em} \Map_{N_d \mathcal{P}}\left(N_{\mathcal{P}} M, N_{\mathcal{P}} L\right),$$
where the former is computed using the projective model structure on $\mathsf{RMod}_{{\mathcal{P}}}$, and the latter with the projective contravariant model structure on $\fSpaces_{/N_d \mathcal{P}}$. We will respectively refer to these towers as the \textit{operadic tower} and the \textit{forest tower} from now on.

Before stating our result, let us explain how the operadic and forest towers arise via respectively filtering $M$ and $N_{\mathcal{P}} M$. As we want to be more general, we will replace $N_{\mathcal{P}} M$ with an arbitrary forest space $X$ over some fixed $V$.

\begin{itemize}
    \item \textit{Operadic filtration}: taking inspiration from the embedding tower, we define $M^{\smallleq k}$ to be the truncated collection $\{ M(j) : 1 \leq j \leq k \}$, which also inherits a truncated right $\mathcal{P}$-module action. We set 
    $$\mathsf{Map}^{\smallleq k}_{\mathcal{P}}(M, L)$$
    to be the derived mapping space of natural transformations $M^{\smallleq k} \to L^{\smallleq k}$ respecting the operadic action.
    \item \textit{Forest filtration}: let $\FilFor{k} \subseteq \For$ be the full subcategory spanned by the forest composed of at most $k$ trees. For any forest space $X \in \fSpaces_{/V}$, we can consider the restriction of the presheaf morphisms $X^{\smallleq k} \to V^{\smallleq k}$ to $\FilFor{k}$, and form the derived mapping space
    $$ \Map^{\smallleq k}_V(X,Y)$$
    of truncated maps $X^{\smallleq k} \to Y^{\smallleq k}$ over $V^{\smallleq k}$.
\end{itemize}

The following theorem summarizes the main results of Sections 4.1. and 4.2. We write below $\mathbb{R}N_{\mathcal{P}}$ for the right derived nerve functor.

\begin{statement}{\cref{main_thm}}
    Let $M$ and $L$ be right modules over a unital simplicial operad $\mathcal{P}$. Then the operadic and forest filtrations lead to towers at the level of derived mapping spaces
$$
	\left\{ \Map^{\smallleq k}_\mathcal{P}(M, L)\right\}_{k \geq 1} \xlongrightarrow{N_{\mathcal{P}}} \left\{ \Map^{\smallleq k}_{N_d \mathcal{P}}(\mathbb{R}N_{\mathcal{P}} M, \mathbb{R}N_{\mathcal{P}} L)\right\}_{k \geq 1}
	$$
with a morphism between them induced from the nerve functor $N_{\mathcal{P}}$. Moreover, both towers converge, in the sense that there are weak homotopy equivalences
\[\begin{tikzcd}[cramped]
	{\Map_{\mathcal{P}}(M, L)} & {\Map_{N_d \mathcal{P}}(\mathbb{R}N_{\mathcal{P}} M, \mathbb{R}N_{\mathcal{P}} L)} \\
	{\holim_k \Map^{\smallleq k}_\mathcal{P}(M, L)} & { \holim_k \Map^{\smallleq k}_{N_d \mathcal{P}}(\mathbb{R}N_{\mathcal{P}} M, \mathbb{R}N_{\mathcal{P}} L).}
	\arrow["{N_{\mathcal{P}}}", from=1-1, to=1-2]
	\arrow["\simeq"', from=1-1, to=2-1]
	\arrow["\simeq"', from=1-2, to=2-2]
	\arrow["{N_{\mathcal{P}}}", from=2-1, to=2-2]
\end{tikzcd}\]
Finally, if $\mathcal{P}$ is also assumed to be $\Sigma$-free, then all morphisms above coming from applying $N_{\mathcal{P}}$ are weak homotopy equivalences.
\end{statement}

The remaining sections of Chapter 4 are dedicated to analysing the forest tower. As its first stage can already be computed using Theorem C, we instead focus on understanding the \textit{layers}, that is, the fibers
$$ \mathsf{Fib}_k(X,Y) \coloneqq \mathrm{hofiber} \left( \Map_V^{\smallleq k}(X,Y) \longrightarrow \Map_{V}^{\smallleq k-1}(X,Y) \right)$$
of the connecting maps of the tower for $k \geq 2$. We fix from now on $X, Y$ forest spaces over a dendroidal spaces $V$.

Following the approach of Göppl and Weiss in \cite{GopplWeiss}, we show that the layers of the tower can be described via a latching-matching decomposition 
$$ \mathsf{Latch}_k(X) \longrightarrow X(k) \longrightarrow \mathsf{Match}_k(X),$$
where $X(k)$ is the space obtained by evaluating $X$ at the forest of $k$ edges. One should think of the spaces $\mathsf{Latch}_k(X)$ and $\mathsf{Match}_k(X)$ as respectively forming the initial and terminal way of extending $X^{\smallleq k-1}$ to a presheaf on $\FilFor{k}$. In particular, they are completely determined by the truncation $X^{\smallleq k-1}$, and moreover keep track of the the compatibility conditions that any extension of $X^{\smallleq k-1}$ needs to satisfy in order to properly define a presheaf. The next result makes this precise.

\begin{statement}{\cref{better_layer_description}}
Let $V$ be a reduced dendroidal Segal space, and suppose $X, Y$ are forest spaces over $V$. For $k \geq 2$, the $k^{th}$ layer of the forest tower for $\Map_V(X,Y)$ can be described as the total homotopy fiber of the square
$$
\begin{tikzcd}[cramped]
	{\Map^{\Sigma_k}_V(X(k), Y(k))} & {\Map^{\Sigma_k}_V(\mathsf{Latch}_k(X), Y(k))} \\
	{\Map^{\Sigma_k}_V(X(k), \mathsf{Match}_k(Y))} & {\Map^{\Sigma_k}_V(\mathsf{Latch}_k(X), \mathsf{Match}_k(Y))}
	\arrow[from=1-1, to=1-2]
	\arrow[from=1-1, to=2-1]
	\arrow[from=1-2, to=2-2]
	\arrow[from=2-1, to=2-2]
\end{tikzcd}
$$
    where the arrows are determined by restriction along the latching map $\mathrm{Latch}_k(X) \to X(k)$, and by composition along the matching map $Y(k) \to \mathrm{Match}_k(Y)$. 
\end{statement}

Stated differently, \cref{better_layer_description} says that $\mathsf{Fib}_k(X,Y)$ is the space of $\Sigma_k$-equivariant maps $X(k) \to Y(k)$ over $V(k)$ which make the diagram
\[\begin{tikzcd}[cramped]
	{\mathsf{Latch}_k(X)} & {X(k)} & {\mathsf{Match}_k(X)} \\
	{\mathsf{Latch}_k(Y)} & {Y(k)} & {\mathsf{Match}_k(Y)}
	\arrow[from=1-1, to=1-2]
	\arrow[from=1-1, to=2-1]
	\arrow[from=1-2, to=1-3]
	\arrow[dashed, from=1-2, to=2-2]
	\arrow[from=1-3, to=2-3]
	\arrow[from=2-1, to=2-2]
	\arrow[from=2-2, to=2-3]
\end{tikzcd}\]
commute up to homotopy. We note that Krannich and Kupers in \cite{KrannichKupersRightModules} provide the same sort of result via a model-independent approach, and hence our methods using the language of forests are different and independent from theirs.

Let us give an idea of the main step of the proof \cref{better_layer_description}, namely by explaining how one can reduce the computation of $\mathsf{Fib}_k(X,Y)$ to a question about the spaces $X(k)$ and $Y(k)$. The remaining needed ingredients for deducing \cref{better_layer_description} are obtained in Section 4.3.2 via an excision result for presheaf categories first proved by Göppl--Weiss \cite{GopplWeiss}. This was later generalized by Krannich--Kupers \citep[Thm. 3.8]{KrannichKupersRightModules} -- with a proof attributed to Ayala--Mazel-Gee--Rozenblyum -- and also independently shown by Haine--Ramzi--Steinebrunner \citep[Thm. 5.12]{HaineRamziSteinebrunner}. 

Consider the inclusion of full subcategories
    $$ \FilFor{k-1} \subseteq \mathsf{Bas}^{\smallleq k} \subseteq \FilFor{k}$$
    where $\mathsf{Bas}^{\smallleq k}$ is spanned by all forests of filtration at most $k-1$, and the forest $k \cdot \eta$ of $k$ edges. \cref{reduction_basic_forests} states that restriction of presheaves along the inclusion $\mathsf{bas} \colon \mathsf{Bas}^{\smallleq k} \to \FilFor{k}$ defines a weak homotopy equivalence 
    \begin{equation}
    \Map_V(X, Y) \xlongrightarrow{\simeq} \Map_V(\mathsf{bas}^* X, \mathsf{bas}^* Y)
    \label{mapping_intro}
    \end{equation}
    at the level of mapping spaces. The full proof of this result is the content of Section 4.3.1., and it consists in interpolating $\mathsf{Bas}^{\smallleq k} \subseteq \FilFor{k}$ via intermediate full subcategories successively removing the unwanted forests. 
    
    The main observation that is used in showing the equivalence \eqref{mapping_intro} is that, supposing $X \to V$ to be a right fibration, $X$ satisfies a Segal condition, imported from the Segal condition of $V$ which we assume from the start. This condition on $X$ roughly says that, for any $F \in \FilFor{k}$, the space $X(F)$ is completely determined by the value of $X$ on all corollas $C_\ell$, as well as all forests of edges $n \cdot \eta$ with $n \leq k$. As all the forests we have just mentioned are contained in $\FilFor{k-1}$ \textit{except} for $k \cdot \eta$, we conclude that the only essential part one needs to know when extending $X^{\smallleq k-1}$ to $\FilFor{k}$ must be contained in $k \cdot \eta$. This is exactly the data carried by the category $\mathsf{Bas}^{\smallleq k}$.

    In Section 4.4 we apply \cref{better_layer_description} to obtain connectivity estimates for the embedding tower. These are already known from the work of Goodwillie and Klein, but our methods for proving them are new.
    \begin{statement}{\cref{layers_estimates_connectivity}}
        Let $M^d$ and $N^{d+n}$ be smooth manifolds. For $k \geq 2$, the $k^{th}$ layer of the Goodwillie--Weiss tower on $\mathrm{Emb}(M,N)$, that is, the homotopy fiber of 
    $$ T_{k} \mathrm{Emb}(M,N) \longrightarrow T_{k-1}\mathrm{Emb}(M,N)$$
    has connectivity $k \left( \dim N - \hdim M -2 \right)-(\dim N -2)$. If one additionally assumes that $\dim N - \hdim M \geq 3$, then the map
    $$ T_{\infty} \mathrm{Emb}(M,N) \longrightarrow T_{k}\mathrm{Emb}(M,N)$$
    has connectivity given by at least $k \left( \dim N - \hdim M -2 \right)-(\dim N -2)$.
    \end{statement}

    By appealing to the Goodwillie--Weiss Convergence theorem, together with the Hirsch--Smale theorem for identifying $T_1 \mathrm{Emb}(M,N)$ with the space of immersions $\mathrm{Imm}(M,N)$, the following connectivity result follows at once.

\begin{statement}{\cref{conn_imm_emb}}
    Let $M^d$ and $N^{d+n}$ be smooth manifolds such that $\dim N - \hdim M \geq 3$ holds. Then the inclusion of immersions into embeddings
    $$ \mathrm{Emb}(M,N) \longrightarrow \mathrm{Imm}(M,N)$$
    is $(\dim N - 2 \cdot \hdim M -2)$-connected.
\end{statement}

\subsection*{Notation}

\leavevmode

Here we record some of the notation that we will use throughout.

\begin{itemize}
    \item We write $\underline{n}$ for the finite set $\left\{ 1, \ldots, n \right\}$.
    \item We write $\Sigma_n$ for the symmetric group on $n$ letters.
    \item We write $\Fin$ for the category of finite sets.
    \item We write $\Fin^{\cong}$ for the wide subcategory of $\Fin$ on the bijections.
    \item We write $J$ for the nerve of the groupoid with object set $\{ 0, 1\}$ and an isomorphism between them.
    \item For $\mathcal{C}$ a small (simplicial) category, we write $\mathcal{Y} \colon \mathcal{C} \to \mathsf{Psh}(\mathcal{C})$ for the Yoneda embedding, and also for its simplicially-enriched version.
    \item For $\mathcal{C}$ a small category, we will write $\mathcal{C}(-,-)$ for the sets of morphisms. 
    \item For $\mathcal{M}$ a model category, we will write $\mathsf{Map}(-,-)$ for the derived mapping space with respect to the model structure of $\mathcal{M}$. We will usually omit the reference to the model category since it will be clear throughout.
    \item Given a rooted tree $T$, we write $E(T), \: V(T)$ and $L(T)$ for the sets of edges, vertices and leaves, respectively. These naturally extend to the context of forests $F$, for which we also define $R(F)$ to be the set of roots of $F$.
    \item The word \textit{spaces} will, unless otherwise stated, always mean \textit{simplicial sets}.
    \item By \textit{$\infty$-operad} we will always mean a dendroidal $\infty$-operad, that is, a dendroidal set having the inner horn filling condition.
\end{itemize}

\vspace{1em}

\noindent\textbf{Acknowledgments.} 

The author would like to thank Ieke Moerdijk for carefully reading previous versions of this work and pointing out many places where it could be improved, and also Pedro Boavida de Brito for his comments and corrections on earlier drafts. I would also like to thank Gijs Heuts for his interest in this work, and the author's PhD committee for reading the thesis version of this article and their comments on it. The author was supported by the Dutch Research Council (NWO) through the grant 613.009.147.

\section{The homotopy theory of right modules over an operad}

In this section we want to recall some concepts coming from operad theory, namely that of a right module over a coloured operad $\mathcal{P}$, which is the main object we wish to study in this text. Before doing so, we will begin by defining what a coloured operad is, mostly just to fix some needed notation. From here on out, our operads will be enriched in $\mathcal{V}$, which will always be assumed to be either $(\Sets, \times)$ or $(\sSets, \times)$; in particular, we will always use the product symbol for the tensor product. However, most of what we will say also applies if more generally $\mathcal{V}$ is a closed symmetric monoidal category.

For the next definition, we will say that a finite rooted tree is a \textit{corolla} if it has exactly one vertex, and a map of corollas is just a bijection of their sets of leaves, which are the edges of the corolla which are not the root. We note that maps of corollas can be composed in the obvious way. For any finite set $I$, we write $C_I$ for a corolla with leaf set given by $I$; if $I = \underline{n}$, we will abbreviate the notation to $C_n$. For more on the language of trees, we refer the reader to Section 3.

\begin{definition}
    Let $C$ be a set. We define a \textit{$C$-corolla} to be a pair $(C_I, \alpha)$, where $C_I$ is a corolla with leaf set given by a finite set $I$, and  $\alpha \colon E(C_I) \to C$ is a function.
    
    We set a morphism between such objects $(C_I, \alpha) \to (C_J, \beta)$ to be a corolla map $f \colon C_I \to C_J$ such that the diagram
$$
\begin{tikzcd}
E(C_I) \arrow[rd, "\alpha"'] \arrow[rr, "E(f)"] &   & E(C_J) \arrow[ld, "\beta"] \\
                                                & C &                           
\end{tikzcd}$$
    commutes. The $C$-corollas, together with their respective morphisms as in the definition above, form a groupoid $\mathsf{Cor}(C)$. 
\end{definition}

Given any sets $I$ and $J$, and a choice of element $a \in I$, we define a new set $ I \circ_a J = \left( I- \{ a \} \right) \amalg J.$ This construction can be extended to the groupoid of $C$-corollas $\mathsf{Cor}(C)$: to any two objects $(C_I, \alpha)$ and $(C_J, \beta)$, we can associate a new one $(C_{I \circ_a J}, \alpha \circ_a \beta)$, where
$$(\alpha \circ_a \beta)(i) = 
\begin{cases}
   \alpha(i) \: &\mathrm{if} \: i \in A - a \\ 
   \beta(i) \: &\mathrm{if} \: i \in B
\end{cases}
$$
\begin{definition}
    A \textit{coloured $\mathcal{V}$-operad} is given by a set of colours $C$ and a contravariant functor
    $$ \mathcal{P} : \mathsf{Cor}(C)^{\mathsf{op}} \longrightarrow \mathcal{V},$$
    together, for each $a \in I$, with structure maps called \textit{partial compositions}
    $$ \circ_a \colon \mathcal{P}(C_I, \alpha) \times \mathcal{P}(C_J, \beta) \longrightarrow \mathcal{P}(C_{I \circ_a J}, \alpha \circ_a \beta),$$
    which are natural with respect to $(C_I, \alpha)$, $(C_J, \beta)$ and $a$, and satisfy the following conditions:
    \begin{itemize}
        \item For each colour $c \in C$, there is an operation $\mathrm{id}_c \in \mathcal{P}(C_1, 1 \mapsto c)$ which acts as the identity with respect to the partial composition operations.

        \item The partial composition maps are associative in the following sense: for $a \in I$ and $b \in J$, the following diagram
\[\begin{tikzcd}[cramped]
	{\mathcal{P}(C_I, \alpha) \times \mathcal{P}(C_J, \beta) \times \mathcal{P}(C_K, \gamma)} & {} & {\mathcal{P}(C_{I \circ_a J}, \alpha \circ_a \beta) \times \mathcal{P}(C_K, \gamma)} \\
	{\mathcal{P}(C_I, \alpha) \times \mathcal{P}(C_{J \circ_b K},\gamma)} & {} & {\mathcal{P}(C_{I \circ_a J \circ_b K}, \alpha \circ_a \beta \circ_b \gamma)}
	\arrow["{\circ_{a} \times \mathcal{P}(C_K, \gamma)}", from=1-1, to=1-3]
	\arrow["{\mathcal{P}(C_I, \alpha) \times \circ_b }"', from=1-1, to=2-1]
	\arrow["{\circ_b}", from=1-3, to=2-3]
	\arrow["{\circ_a}", from=2-1, to=2-3]
\end{tikzcd}\]
commutes.
        \item The partial composition along different elements commute: if $a, a' \in I$ are distinct elements, then the following diagram
        \[\begin{tikzcd}[cramped]
	{\mathcal{P}(C_I, \alpha) \times \mathcal{P}(C_J, \beta) \times \mathcal{P}(C_K, \gamma)} & {} & {\mathcal{P}(C_{I \circ_{a'} K}, \alpha \circ_{a'} \gamma) \times \mathcal{P}(C_J, \beta)} \\
	{\mathcal{P}(C_{I \circ_a J}, \alpha \circ_a \beta) \times \mathcal{P}(C_{K},\gamma)} & {} & {\mathcal{P}(C_{I \circ_a J \circ_{a'} K}, \alpha \circ_a \beta \circ_{a'} \gamma)}
	\arrow["{\circ_{a'} \times \mathcal{P}(C_J, \beta)}", from=1-1, to=1-3]
	\arrow["{\circ_a \times \mathcal{P}(C_K, \gamma)}"', from=1-1, to=2-1]
	\arrow["{\circ_a}", from=1-3, to=2-3]
	\arrow["{\circ_{a'}}", from=2-1, to=2-3]
\end{tikzcd}\]
commutes, where in the top horizontal arrow one should first precompose with the map that switches the second and third factors.
    \end{itemize}

    A morphism of operads $\mathcal{P} \to \mathcal{Q}$ is then a natural transformation of functors $\mathsf{Cor}(C)^{\mathsf{op}} \to \mathcal{V}$ such that they are compatible with the partial composition maps of $\mathcal{P}$ and $\mathcal{Q}$ in the obvious way. 

    \label{definition_operad}
\end{definition}

\begin{remark}
The usual way of defining a coloured operad as a set of colours $C$ together with objects of operations $\mathcal{P}(c_1, \ldots, c_n ; c_0)$ can be recovered from this one by evaluating the contravariant functor $\mathcal{P}$ at the corolla $(C_n, c_{\bullet})$ which evaluates the leaf $i$ to $c_i$, and the root to $c_0$. Then the equivariance condition corresponds to $\mathcal{P}$ being a functor, and the unitality and associativity will of course be the usual unitality and associativity conditions for an operad. In fact, one can go back and forth between these two different definitions of a coloured operad, see Section 1.8 of \cite{HeutsMoerdijkDendroidal} for more details on this. 

The reason for us adopting this "coordinate-free" approach is that it makes the definitions and proofs more readable: for instance, stating the associativity axioms for the partial composition maps is certainly an easier task.
\label{alternative_operad_description}
\end{remark}

Before defining what a right module over a coloured operad is, we will need to define a category which is a slight variation of $\mathsf{Cor}(C)$. For more details on why we need to do this, see \cref{edges_vs_corollas}.

\begin{definition}
    Let $C$ be a finite set. We define the category of \textit{$C$-edges} $\Edges(C)$ to be the core of the slice category $\Fin_{/C}$, that is:
    \begin{itemize}
        \item The objects of $\Edges(C)$ are the pairs $(I, \alpha)$ of a finite set $I$ and a function $\alpha \colon I \to C$.
        \item A morphism $(I, \alpha) \to (J, \beta)$ is given by a bijection $f \colon I \to J$ such that $\beta f= \alpha$.
    \end{itemize}
\end{definition}

\begin{definition}
    Let $\mathcal{P}$ be a coloured $\mathcal{V}$-operad with set of colours $C$. We say that a contravariant functor
    $$ M : \Edges(C)^{\mathsf{op}} \longrightarrow \mathcal{V}$$
    is a \textit{right $\mathcal{P}$-module} if, for each finite set $I$ and element $a \in I$, there are partial action maps
$$ \cdot_a : M(I, \alpha) \times \mathcal{P}(C_J, \beta) \longrightarrow M(I \circ_a J, \alpha \circ_a \beta)$$
    which are natural on $(I, \alpha)$, $(C_J, \beta)$ and $a$, and such the conditions below hold:
    \begin{itemize}
        \item For each colour $c \in C$, the identity operation $\mathrm{id}_c$ acts as the identity via the partial action maps.
         \item The right action should be compatible with the operad structure of $\mathcal{P}$: for $a \in I$ and $b \in J$, the following diagram
\[\begin{tikzcd}[cramped]
	{M(I, \alpha) \times \mathcal{P}(C_J, \beta) \times \mathcal{P}(C_K, \gamma)} & {} & {M(I \circ_a J, \alpha \circ_a \beta) \times \mathcal{P}(C_K, \gamma)} \\
	{M(I, \alpha) \times \mathcal{P}(C_{J \circ_b K},\gamma)} & {} & {M(I \circ_a J \circ_b K, \alpha \circ_a \beta \circ_b \gamma)}
	\arrow["{\cdot_{a} \times \mathcal{P}(C_K, \gamma)}", from=1-1, to=1-3]
	\arrow["{M(I, \alpha) \times \circ_b }"', from=1-1, to=2-1]
	\arrow["{\cdot_b}", from=1-3, to=2-3]
	\arrow["{\cdot_a}", from=2-1, to=2-3]
\end{tikzcd}\]
commutes.
        \item The actions along different elements commute: if $a, a' \in I$ are distinct elements, then the following diagram
        \[\begin{tikzcd}[cramped]
	{M(I, \alpha) \times \mathcal{P}(C_J, \beta) \times \mathcal{P}(C_K, \gamma)} & {} & {M(I \circ_{a'} K, \alpha \circ_{a'} \gamma) \times \mathcal{P}(C_J, \beta)} \\
	{M(I \circ_a J, \alpha \circ_a \beta) \times \mathcal{P}(C_{K},\gamma)} & {} & {M(I \circ_a J \circ_{a'} K, \alpha \circ_a \beta \circ_{a'} \gamma)}
	\arrow["{\cdot_{a'} \times \mathcal{P}(C_J, \beta)}", from=1-1, to=1-3]
	\arrow["{\cdot_a\times \mathcal{P}(C_K, \gamma)}"', from=1-1, to=2-1]
	\arrow["{\cdot_{a'}}", from=1-3, to=2-3]
	\arrow["{\cdot_{a'}}", from=2-1, to=2-3]
\end{tikzcd}\]
commutes, where again in the top horizontal arrow one should first precompose with the map that switches the second and third factors.
    \end{itemize}

A morphism of right $\mathcal{P}$-modules $M \to L$ is a natural transformation of functors $\Edges(C)^{\mathsf{op}} \to \mathcal{V}$ compatible with the partial composition maps of $M$ and $L$ in the obvious way. We will write $\RMod{\mathcal{P}}$ for the category determined by this data.
\label{definition_right_module}
\end{definition}

\begin{remark}
Whereas in \cite{AroneTurchinRational} one defines uncoloured operads and their respective right modules as certain presheaves parameterized by the category of finite sets $\Fin$, for the coloured setting we had to use the category $\Edges(C)$ instead of $\mathsf{Cor}(C)$. This is mainly due to the output variable not being necessary when defining $M$, and this subtlety doesn't appear in the uncoloured setting, since in that case there is only one possible choice of colour and therefore $\Edges(\ast) = \mathsf{Cor}(\ast) = \Fin^{\cong}$.
\label{edges_vs_corollas}
\end{remark}

\begin{remark}
    A similar observation to the one made about operads in \cref{alternative_operad_description} equally applies to the case of right modules over an operad. Namely, instead of specifying the objects $M(I, \alpha)$ for each $(I, \alpha) \in \Edges(C)$, we could instead provide objects
    $$M(c_1, \ldots, c_n)$$
    for each $n$-tuple of colours $(c_1, \ldots, c_n)$, and it is easily checked that these two possible descriptions of a right module are indeed equivalent.

    We also note here that we could have defined right modules using total action maps, instead of just partial ones: given a $C$-edges $(I, \alpha)$ and $(C_J, \beta)$, together with a function of finite sets $f \colon I \to J$, the choice of all partial compositions $\circ_a$ with $a \in I$ determines a map
    \begin{equation}
    M(J, \beta) \times \prod_{j \in J} \mathcal{P}\left( C_{f^{-1}(j)}, \alpha_j \right) \longrightarrow M(I, \alpha),
        \label{total_action_maps}
    \end{equation}
    which is independent of the way we apply the partial compositions by the axioms of a right module. Here $\alpha_j$ is the restriction of $\alpha$ to $f^{-1}(j)$, and this map also sends the root of the corolla to $\beta(j)$. We call morphism \eqref{total_action_maps} a \textit{total action map} for the module structure.
    \label{remark_total_action}
\end{remark}

\begin{remark}
    Let us make a few brief remarks on the theory of left $\mathcal{P}$-modules. It is not complicated to adapt \cref{definition_right_module} to also work for operadic left modules. We will not make this precise, and in fact will not really delve into theory of left $\mathcal{P}$-modules in this work for two main reasons:
    \begin{enumerate}[label=(\roman*)]
        \item Firstly, the operadic theory of left $\mathcal{P}$-modules is substantially different from that of right $\mathcal{P}$-modules. Therefore, it is not true that our arguments involving dendroidal objects and right modules have an immediate dual for their left counterpart, and consequently left $\mathcal{P}$-modules would require a different analysis from the one we propose here.
        \item A more crucial reason is that the theory of left $\mathcal{P}$-modules is interesting to the extent that, as far as we are aware, the main relevant examples of these types of modules are $\mathcal{P}$-algebras. However, the theory of $\mathcal{P}$-algebras via dendroidal sets and dendroidal spaces has already been thoroughly studied: for instance, the work of Heuts \cite{HeutsAlgebrasInfinityOperads} on the covariant model structure for dendroidal sets, and the approach of Boavida de Brito--Moerdijk in \cite{BoavidaMoerdijk} for dendroidal spaces. 
    \end{enumerate}
\end{remark}

Before providing some examples of right modules over an operad, we will first give a standard characterization of right $\mathcal{P}$-modules as a certain category of enriched functors on the \textit{symmetric monoidal envelope} of $\mathcal{P}$.

\begin{definition}
    Let $\mathcal{P}$ be a coloured $\mathcal{V}$-operad with set of colours $C$. We define the \textit{symmetric monoidal envelope} of $\mathcal{P}$ to be the symmetric monoidal $\mathcal{V}$-category $\left( \Prop(\mathcal{P}), \boxplus, \emptyset \right)$ given by the following data:
    \begin{itemize}[label=$\diamond$]
        \item The objects of $\Prop(\mathcal{P})$ are the $C$-edges $(I, \alpha)$.
        \item The object of morphisms in $\Prop(\mathcal{P})$ with source $(C_I, \alpha)$ and target $(C_J, \beta)$ is given by
        $$ \Prop(\mathcal{P})\left( (I, \alpha), (J, \beta) \right) = \coprod_{f \colon I \to J} \prod_{j \in J} \mathcal{P}\left( C_{f^{-1}(j)}, \alpha_j\right).$$
        Here $f\colon I \to J$ is a function of finite sets, and $(C_{f^{-1}(j)}, \alpha_j) \in \mathsf{Cor}(C)$ is defined as in \cref{remark_total_action}. The composition in $\Prop(\mathcal{P})$ comes from the operadic composition in $\mathcal{P}$, see \citep[Def. 4.1]{AroneTurchinRational} for a detailed explanation of this.
        \item The monoidal structure is given on objects by
        $$ (I, \alpha) \boxplus (J, \beta) = (I \amalg J, \alpha \amalg \beta)$$
        and the monoidal unit is $\emptyset = (\emptyset, \alpha)$.
    \end{itemize}
    \label{Prop_Definition}
\end{definition}

\begin{proposition}
    Let $\mathcal{P}$ be a coloured operad. Then there is an isomorphism of categories
    $$ \RMod{\mathcal{P}}(\mathcal{V}) \cong \Fun^{\mathcal{V}}(\Prop(\mathcal{P})^{\mathsf{op}}, \mathcal{V})$$
    between the category of right $\mathcal{P}$-modules with values in $\mathcal{V}$, and the category of $\mathcal{V}$-enriched contravariant functors $\Prop(\mathcal{P})^{\mathsf{op}} \to \mathcal{V}$. In particular, the categories $\RMod{\mathcal{P}}(\Sets)$ and $\RMod{\mathcal{P}}(\sSets)$ over a fixed operad $\mathcal{P}$ are (co)complete.
    \label{right_module_envelope}
\end{proposition}
\begin{proof}
This is a well-established result on right modules over an \textit{uncoloured} operad, it is for instance mentioned \citep[Prop. 1.2.6]{KapranovManin} or \citep[Lem. 2.3]{DwyerHessKnudsen}. As we found no proof of the coloured version of this statement, we provide it here.

Suppose first that we are given a right $\mathcal{P}$-module $M$ and want to associate to it an enriched functor $F_M \colon \Prop(\mathcal{P})^{\mathsf{op}} \to \mathcal{V}$, and fix $C$ to the set of colours. Since both $\mathsf{Edges}(C)$ and $\Prop(\mathcal{P})$ share the same set of objects, we define $F_M$ to be the same as $M$ on objects. As for the morphisms, we recall from \eqref{total_action_maps} that for any $(I, \alpha), (J, \beta) \in \Prop(\mathcal{P})$ we have maps
$$
    M(J, \beta) \times \prod_{j \in J} \mathcal{P}\left( C_{f^{-1}(j)}, \alpha_j \right) \longrightarrow M(I, \alpha) $$
    for each choice of function $f \colon I \to J$. Taking the coproduct over all such maps leads to a morphism
    $$ M(J, \beta) \times \Prop(\mathcal{P})\left( (I, \alpha), (J, \beta) \right) \longrightarrow M(I, \alpha)$$
    which defines the action of $F_M$ at the level of morphisms. The extra compatibility conditions needed for $F_M$ to define a functor are now encoded in the extra restrictions we have put on how the action maps of $M$. In the other direction, the module associated to an enriched functor $F \in \Fun^{\mathcal{V}}(\Prop(\mathcal{P})^{\mathsf{op}}, \mathcal{V})$ can be obtained by reversing the previous argument.
\end{proof}

\begin{remark}
    In the identification above we don't put any extra condition which accounts for the monoidal structure of $\Prop(\mathcal{P})$. However, we note that if $\mathcal{P}$ is a closed operad (that is, the space of constants $\mathcal{P}(-;c)$ is a singleton for each colour $c \in C$), then for any right $\mathcal{P}$-module $M$ we have maps
    $$ M(c_1, \ldots, c_n, c_{n+1}, \ldots, c_m) \longrightarrow M(c_1, \ldots, c_n) \times M(c_{n+1}, \ldots, c_m)$$
    by using the action by the constants to "remove" input colours. In fact, since these maps come from the module structure, they actually give $M$ the structure of a colax monoidal functor.
    \label{colax_monoidal}
\end{remark}

\begin{example}
Let us give some examples of operadic right modules.

    \begin{enumerate}[label=(\alph*)]
        \item  A comparison between \cref{definition_operad} and \cref{definition_right_module} shows that any operad $\mathcal{P}$ is a right $\mathcal{P}$-module in a natural way.

        \item Let $\mathbf{Comm}$ be the commutative operad in $\Sets$. Then one easily checks that its monoidal envelope is isomorphic to the category of finite sets $\mathsf{Fin}$. Therefore, a right $\mathsf{Comm}$-module is the same as a presheaf on $\Fin$, by \cref{right_module_envelope}.

        \item Given any object $(I, \alpha)$ of $\mathsf{Env}(\mathcal{P})$, the representable presheaf $\mathcal{Y}_{(I, \alpha)}$ will define a right $\mathcal{P}$-module by \cref{right_module_envelope}, which corresponds to the \textit{free right $\mathcal{P}$-module} generated by the colours defined by $(I, \alpha)$. We will come back to these modules at the end of Section 3.
        
        \item For the next example we will consider operads with values in the category of topological spaces $\mathsf{Top}$ instead of simplicial sets, since it is more convenient within the context of the example.
        
        Consider a smooth manifold $M$ of dimension $d$ without boundary. One can associate to $M$ a presheaf $\mathbb{E}_M \colon \Fin^{\mathsf{op}} \to \mathsf{Top}$, determined on objects by
        $$ \mathbb{E}_M(I) = \mathsf{Emb}(I \times D^d, M)$$
        for any finite set $I$. Here $\mathsf{Emb}(I \times D^d, M) \subseteq C^{\infty}(I \times D^d, M)$ is the subspace of smooth embeddings $I \times D^d \to M$ equipped with the weak smooth topology. Consider also the presheaf $\mathbb{E}^{\mathsf{fr}}_d \colon \Fin^{\mathsf{op}} \to \mathsf{Top}$ given by 
        $$ \mathbb{E}^{\mathsf{fr}}_d(I) = \mathsf{Emb}(I \times D^d, D^d),$$
        which actually has the structure of an operad via the function composition of the embeddings. We will call $\mathbb{E}^{\mathsf{fr}}_d$ the \textit{framed little d-disks operad}\footnote{The nomenclature in the literature for this operad is unfortunate, since this is usually know as the \textit{framed} version of the little $d$-disks operad, and the version where the spaces of operations have the homotopy type of \textit{framed} embeddings is the \textit{non-framed} version of this operad.}. Then, for each $a \in I$, there exist maps
        $$ \mathsf{Emb}(I \times D^d, M) \times \mathsf{Emb}(J \times D^d, D^d) \longrightarrow \mathsf{Emb}\left( (I \circ_a J) \times D^d, M \right)$$
        which sends a pair $\left( \{ \varphi_i : i \in I\}, \{ \psi_j : j \in J \}\right)$ in the domain to an embedding $\{ \tilde{\varphi}_i : i \in I \circ_a J\}$ given by
        $$
        \tilde{\varphi}_i =
        \begin{cases}
           \varphi_i \: &\mathrm{if} \; i \not\in J \\
           \varphi_a \psi_i &\mathrm{if} \; i \in J.
        \end{cases}
        $$

        Then one easily check that these maps endow $\mathbb{E}_M$ with the structure of a right $\mathbb{E}^{\mathsf{fr}}_d$-module which is of central importance in the study of the manifold calculus of Goodwillie and Weiss, as we have already mentioned in the introduction.
        
        \item \label{example_FM_manifolds} The Fulton--MacPherson operad $\mathsf{FM}_d$ is obtained via the Fulton--MacPherson compactification of the configuration spaces $\mathsf{Conf}_k(\mathbb{R}^d)$. Inspired by this, one might suspect that more generally the compactification of $\mathsf{Conf}_k(M)$ for a general manifold $M^d$ will lead to a right $\mathsf{FM}_d$-module; in this example we will explain how this is the case if $M$ is parallelizable by describing how the general Fulton--MacPherson works for a general smooth manifold. For the details on this, see \cite{SinhaFM} and also \cite{IdrissiRealHomotopy}. We assume knowledge of the case for euclidean spaces, see \cite{HeutsMoerdijkDendroidal}.

        Throughout this example, we always suppose that $M$ is embedded in some Euclidean space $\mathbb{R}^N$ and identify $M$ with the image of this embedding; by \citep[Thm. 4.7]{SinhaFM}, the diffeomorphism type of the compactification will not depend on this choice. Using the Euclidean metric on $\mathbb{R}^N$, we can consider the map
        $$ \mathsf{Conf}_k(M) \longrightarrow M^{\times k} \times (S^{N-1})^{\binom{k}{2}} \times [0, + \infty]^{\binom{k}{3}},$$
        where the first component is just the usual inclusion, and the other two are the relative direction and distance maps respectively. We define the \textit{Fulton--MacPherson compactification} $\mathsf{FM}_M(k)$ to be the closure of the image of the map above, and according to \citep[Thm. 4.4]{SinhaFM} it is a manifold with corners of dimension $k \cdot \mathrm{dim}(M)$. One should be warned that $\mathsf{FM}_M(k)$ is not a compact manifold unless $M$ is also compact!

        One can describe the elements of $\mathsf{FM}_M(k)$ as virtual configurations in $M$ in the same way as we had before, but now one needs to be careful that collisions in $M$ at a point $x \in M$ lead to configurations in the tangent space $T_x M$ (this subtlety didn't appear before since we used that $\mathbb{R}^d$ is parallelizable). This is formalized in Sections 3.2 and 3.3 of \cite{SinhaFM}. 
        
        Sinha describes the stratification of $\mathsf{FM}_M(k)$ using a certain category $\Psi_{\underline{k}}$ of labelled trees. However, it will be more helpful in further sections that we reformulate this category in terms of \textit{forests}, that is, in terms of finite collections of trees $\bigoplus_{i \in I} T_i \coloneqq \{ T_i : i \in I\}$. Let $\mathsf{fBin}_k$ be the set of pairs $(F, \lambda)$, where $F = \bigoplus_{i \in I} T_i$ is a forest with each $T_i$ having only vertices of degree at least 2, and $\lambda \colon I \to \coprod_{i \in I} L(T_i)$ is a bijection. Then we claim that to any element of $\mathsf{FM}_M(k)$ one can associate a forest of $\mathsf{fBin}_k$: for instance, 
        to the virtual configuration in $M$ below
        \vspace{0.1em}

\begin{center}
    \begin{tikzpicture}[scale=0.6]
       \draw[thick] (0,0)--(10,0)--(15,2)--(5,2)--cycle;
       \filldraw[black] (3,0.5) circle (1.5pt);
       \node[scale=0.8] at (3.4,0.5) {$x$};
       \draw[thick] (3,2.75) ellipse (2.5 and 0.5);
       \draw[thin, dashed] (0.55,2.65)--(3,0.5)--(5.45,2.65);
       \filldraw[black] (2.5,3) circle (1.5pt);
       \node[scale=0.8] at (2.9,3) {$x_1$};
       \filldraw[black] (3.7,2.6) circle (1.5pt);
       \node[scale=0.8] at (4.1,2.6) {$x_3$};

       \filldraw[black] (10,1.5) circle (1.5pt);
       \node[scale=0.8] at (10.4,1.5) {$y$};
       \draw[thick] (10,3.75) ellipse (2.5 and 0.5);
       \draw[thin, dashed] (7.55,3.65)--(10,1.5)--(12.45,3.65);
       \filldraw[black] (8,3.8) circle (1.5pt);
       \node[scale=0.8] at (8.4,3.8) {$x_4$};
       \filldraw[black] (9.5,4) circle (1.5pt);
       \node[scale=0.8] at (9.9,4) {$x_5$};
       \filldraw[black] (12,3.75) circle (1.5pt);
       \draw[thick] (12,6) ellipse (2.5 and 0.5);
       \draw[thin, dashed] (9.55,5.9)--(12,3.75)--(14.45,5.9);
       \filldraw[black] (11,6) circle (1.5pt);
       \node[scale=0.8] at (11.4,6) {$x_2$};
       \filldraw[black] (12.5,6.2) circle (1.5pt);
       \node[scale=0.8] at (12.9,6.2) {$x_6$};

       \node[scale=1] at (14,1) {$M$};
    \end{tikzpicture}
\end{center}
\vspace{1em}
we associate the following forest with two trees:

        \begin{center}
            \begin{tikzpicture}[scale=0.75]
                \filldraw[black] (-1.5,0) circle (1.5pt);
                \draw (-1.5,0)--(-1.5,-1);
                \draw (-2,1)--(-1.5,0)--(-1,1);
                \node[scale=0.75] at (-2,1.2) {1};
                \node[scale=0.75] at (-1,1.2) {3};

                \node at (0,0) {$\oplus$};

                \filldraw[black] (1.5,0) circle (1.5pt);
                \draw (1.5,1)--(1.5,0)--(1.5,-1);
                \draw (0.75,1)--(1.5,0)--(2.25,1);
                \node[scale=0.75] at (1.5,1.2) {5};
                \node[scale=0.75] at (0.75,1.2) {4};
                \filldraw[black] (2.25,1) circle (1.5pt);
                \draw (1.75,2)--(2.25,1)--(2.75,2);
                \node[scale=0.75] at (1.75,2.2) {2};
                \node[scale=0.75] at (2.75,2.2) {6};
            \end{tikzpicture}
        \end{center}

    We can endow $\mathsf{fBin}_k$ with a poset structure defined as follows. Firstly, contracting the inner edge of one of component trees leads to a forest which is larger in this ordering. Secondly, the removal of the root vertex $r$ of a tree $T$ results in a forest where we replace $T$ with the trees above $r$, and this again produces a larger forest for this poset structure. Below we exemplify how this poset structure works for $k = 3$, and one should compare this to Figure 2.5 of \cite{SinhaFM} to see the correspondence between Sinha's trees and our forests:
    \begin{center}
        \begin{tikzpicture}[scale=0.5]
            \filldraw[black] (-6,-1) circle (1.5pt);
            \draw (-6,-1)--(-6,-2);
            \draw (-6.75,0)--(-6,-1)--(-5.25,0);
            \filldraw[black] (-5.25,0) circle (1.5pt);
            \draw (-6,1)--(-5.25,0)--(-4.5,1);
            \node[scale=0.6] at (-6.75,0.2) {$2$};
            \node[scale=0.6] at (-6,1.2) {$3$};
            \node[scale=0.6] at (-4.5,1.2) {$1$};

            \filldraw[black] (-6,3) circle (1.5pt);
            \draw (-6,3)--(-6,2);
            \draw (-6.75,4)--(-6,3)--(-5.25,4);
            \filldraw[black] (-5.25,4) circle (1.5pt);
            \draw (-6,5)--(-5.25,4)--(-4.5,5);
            \node[scale=0.6] at (-6.75,4.2) {$1$};
            \node[scale=0.6] at (-6,5.2) {$2$};
            \node[scale=0.6] at (-4.5,5.2) {$3$};

            \filldraw[black] (-6,-5) circle (1.5pt);
            \draw (-6,-5)--(-6,-6);
            \draw (-6.75,-4)--(-6,-5)--(-5.25,-4);
            \filldraw[black] (-5.25,-4) circle (1.5pt);
            \draw (-6,-3)--(-5.25,-4)--(-4.5,-3);
            \node[scale=0.6] at (-6.75,-3.8) {$3$};
            \node[scale=0.6] at (-6,-2.8) {$1$};
            \node[scale=0.6] at (-4.5,-2.8) {$2$};

            \filldraw[black] (0,4) circle (1.5pt);
            \draw (0,5)--(0,4)--(0,3);
            \draw (-0.75,5)--(0,4)--(0.75,5);
            \node[scale=0.6] at (-0.75,5.2) {$1$};
            \node[scale=0.6] at (0,5.2) {$2$};
            \node[scale=0.6] at (0.75,5.2) {$3$};
            
            \draw (-1,1)--(-1,0);
            \node at (0,0.5) {$\oplus$};
            \draw (1,1)--(1,0);
            \filldraw[black] (1,1) circle (1.5pt);
            \draw (0.5,1.75)--(1,1)--(1.5,1.75);
            \node[scale=0.6] at (-1,1.2) {$1$};
            \node[scale=0.6] at (0.5,1.95) {$2$};
            \node[scale=0.6] at (1.5,1.95) {$3$};

            \draw (-1,-2)--(-1,-3);
            \node at (0,-2.5) {$\oplus$};
            \draw (1,-2)--(1,-3);
            \filldraw[black] (1,-2) circle (1.5pt);
            \draw (0.5,-1.25)--(1,-2)--(1.5,-1.25);
            \node[scale=0.6] at (-1,-1.8) {$2$};
            \node[scale=0.6] at (0.5,-1.05) {$3$};
            \node[scale=0.6] at (1.5,-1.05) {$1$};

            \draw (-1,-5)--(-1,-6);
            \node at (0,-5.5) {$\oplus$};
            \draw (1,-5)--(1,-6);
            \filldraw[black] (1,-5) circle (1.5pt);
            \draw (0.5,-4.25)--(1,-5)--(1.5,-4.25);
            \node[scale=0.6] at (-1,-4.8) {$3$};
            \node[scale=0.6] at (0.5,-4.05) {$1$};
            \node[scale=0.6] at (1.5,-4.05) {$2$};

            \draw (5,-1)--(5,-2);
            \node at (5.5,-1.5) {$\oplus$};
            \draw (6,-1)--(6,-2);
            \node at (6.5,-1.5) {$\oplus$};
            \draw (7,-1)--(7,-2);
            \node[scale=0.6] at (5,-0.8) {$1$};
            \node[scale=0.6] at (6,-0.8) {$2$};
            \node[scale=0.6] at (7,-0.8) {$3$};

            \draw[<-] (4.5,0)--(2,4);
            \draw[<-] (4.2,-1)--(2,1); 
            \draw[<-] (4.2,-1.8)--(2,-2);
            \draw[<-] (4.5,-2.8)--(2,-5);
            \draw[<-] (-1.5,4.2)--(-4,4.2);
            \draw[<-] (-1.5,1)--(-4,3.5);
            \draw[<-] (-1.5,-2)--(-4,-0.5);
            \draw[<-] (-1.5,3.8)--(-2.6,2.348);
            \draw (-2.8,2.084)--(-4,0.5);
            \draw[<-] (-1.5,3.4)--(-2,1.8);
            \draw (-2.1,1.48)--(-2.85,-0.92);
            \draw (-2.95,-1.24)--(-4,-4.6);
            \draw[<-] (-1.5,-5)--(-4,-5);
        \end{tikzpicture}
    \end{center}

    For each $F \in \mathsf{fBin}_k$ we have a corresponding topological space $\mathsf{FM}_M(F)$ -- in Sinha's notation these correspond up to diffeomorphism to the spaces $C_F(M)$ and $D_F(M)$ -- which, using $F$ as a blueprint, provides the needed configurations of points on $M$ and on the Fulton--MacPherson compactification of the tangent spaces at these points. One should be warned that $\mathsf{FM}_M(F)$ isn't simply a product of configuration spaces, unless $M$ is a parallelizable manifold \citep[Prop. 3.10]{SinhaFM}.

We can assemble this forest description into a decomposition of the compactification as a set
$$ \mathsf{FM}_M(k) = \bigcup_{F \in \mathsf{fBin}_k} \mathsf{FM}_M(F),$$
and the topological considerations translate into the following statements:
\begin{itemize}
    \item The open stratum is parameterized by the forest $k \cdot \eta$ consisting of only $k$ edges, and $\mathsf{FM}_M(k \cdot \eta)$ is the configuration space $\mathsf{Conf}_k(M)$. Moreover, the compactification map is the inclusion of the open stratum.
    \item The boundary $\partial \mathsf{FM}_M(k)$ is stratified by the forests $F \in \mathsf{fBin}_k$ with at least one vertex. The different strata are glued to each other via the poset ordering of $\mathsf{fBin}_k$ -- that is, via the removal of the root vertex of a tree or the contraction of an inner edge -- which corresponds to the collision of points in $M$, or of points in a configuration in some tangent space, respectively.
\end{itemize}

As we mentioned above, if $M$ is parallelizable then there is a diffeomorphism
$$ \mathsf{FM}_M \left( \bigoplus_{i = 1}^n T_i \right) \cong \mathsf{Conf}_n(M) \times \prod_{i=1}^n \mathsf{FM}_d(T_i),$$
which, in exactly the same way as we had for $\mathsf{FM}_d$, will lead to a continuous map
$$ \mathsf{FM}_M(n) \times \prod_{i=1} \mathsf{FM}_d(k_i) \longrightarrow \mathsf{FM}_d(k_1+ \cdots + k_n)$$
by grafting trees to the top of the components of each forest indexing the stratification of $\mathsf{FM}_M(k)$. This gives the symmetric sequence defined by $\mathsf{FM}_M(k)$ the structure of a right $\mathsf{FM}_d$-module given by stacking configurations together.

If $M$ is not parallelizable, then the action we just discussed will not make sense in general: indeed, in order to "attach" a configuration of points to the tangent space of $x \in M$ one needs an identification with $\mathbb{R}^d$ which is moreover consistent in the rest of the manifolds. In other words, we would need a trivialization of $TM$. In spite of that, we stress that $\mathsf{FM}_M(k)$ still exists as a manifold with boundary stratified by $\mathsf{fBin}_k$, which is the only aspect we will need further on at the end of the article.
    \end{enumerate}

\label{examples_right_modules}
\end{example}

\begin{remark}
    The interplay between the Fulton--MacPherson compactification of a manifold and operad theory has already been explored in the literature \cite{CamposFramedConfigurations, CamposIdrissi, CamposSurfaces}, namely in connection with the study of the real homotopy type of configuration spaces of manifolds and in finding algebraic models for it. We also remark that the fact that $M$ is smooth is essential in the example above: in \cite{KupersFulton} it is argued that the Fulton--MacPherson compactification cannot be extended to topological manifolds if one asks for certain basic properties to hold.
\end{remark}

\begin{remark}
    Notice that we have built two compactifications of $\mathsf{Conf}_k(\mathbb{R}^d)$, namely $\mathsf{FM}_d(k)$ and $\mathsf{FM}_{\mathbb{R}^d}(k)$ (this last one is not actually compact, as we have observed). According to \cite[Lem. 4.12]{SinhaFM}, the affine group $\mathsf{Aff}(d)$ acts on the non-compact version $\mathsf{FM}_{\mathbb{R}^d}(k)$, and the quotient by this action leads to a weak homotopy equivalences
    $$ \mathsf{FM}_{\mathbb{R}^d}(k) \xlongrightarrow{\simeq} \mathsf{FM}_{\mathbb{R}^d}(k)/\mathsf{Aff}(d) \simeq \mathsf{FM}_d(k).$$
as this is the affine group is contractible. However, we note that
    $$ \dim \mathsf{FM}_{\mathbb{R}^d}(k) = k d \hspace{2em} \mathrm{but} \hspace{2em} \dim \mathsf{FM}_d(k) = kd -(d+1).$$
    \label{FM_Rd}
\end{remark}

As an immediate consequence of \cref{right_module_envelope}, the category of simplicial operadic right modules $\RMod{\mathcal{P}}$ carries a model structure coming from the usual projective model structure for simplicial presheaves, which we will call the \textit{projective model structure} on $\RMod{\mathcal{P}}$. We register this in the proposition below.

\begin{proposition}
    Let $\mathcal{P}$ be a coloured simplicial operad. The category of simplicial right $\mathcal{P}$-modules $\RMod{\mathcal{P}}$ admits a model structure such that a map of right $\mathcal{P}$-modules $M \to N$ is a weak equivalence/fibration if, for each $C$-edge $(I, \alpha)$, the map
        $$M(I, \alpha) \longrightarrow N(I, \alpha)$$
    is a  weak homotopy equivalence/Kan fibration of simplicial sets. Moreover, this model structure is cofibrantly generated and left proper.
    \end{proposition}

\begin{remark}
    Some of the results we just mentioned still hold in more generality by allowing for other types of enrichment $(\mathcal{V}, \otimes, \mathbf{1})$. The identification in \cref{right_module_envelope} if $\mathcal{V}$ is a closed symmetric monoidal category, for instance. As for the existence of the projective model structure, one needs to be a bit more restrictive about what kind of enrichment one takes, see \citep[Thm. 4.4]{MoserEnrichedProjective} for more on this.
    \end{remark}

    Another immediate consequence of \cref{right_module_envelope} is that $\RMod{\mathcal{P}}$ inherits a closed symmetric monoidal structure coming from the Day convolution of enriched functors, which in turn exists due to the monoidal structure on $\mathsf{Env}(\mathcal{P})$. This can be compactly described as a certain coend, or equivalently we can write $(M \boxplus L)(I, \alpha)$ as a quotient of the coproduct
    $$ \coprod_{(J_1, \beta_1), (J_2, \beta_2) \in \mathsf{Env}(\mathcal{P})} \mathsf{Env}(\mathcal{P})((I, \alpha), (J_1, \beta_1) \boxplus (J_2, \beta_2)) \times M(J_1, \beta_1) \times L(J_2, \beta_2),$$
with the quotient defined as follows: for $\varphi \colon (I, \alpha) \to (J_1, \beta_1) \boxplus (J_2, \beta_2)$ and a pair $(m', \ell') \in  M(J'_1, \beta_1')  \times L(J'_2, \beta_2')$, then for any maps $g_i \colon (J_i, \beta_i) \to (J'_i, \beta'_i)$ we identify $((g_1 \boxplus g_2) \varphi, m', \ell')$ with $(\varphi, g_1^* m', g_2^* \ell')$. We can further simplify this coproduct, which we record in the next proposition. 

\begin{proposition}
    Let $\mathcal{P}$ be a closed simplicial operad. If $M$ and $L$ are right $\mathcal{P}$-modules, then:
    \begin{enumerate}[label=(\alph*)]
        \item For every object $(I, \alpha) \in \Prop(\mathcal{P})$ there is a natural isomorphism
        $$ (M \boxplus L)(I, \alpha) \cong \coprod_{(I_1, \alpha_1), (I_2, \alpha_2)} M(I_1, \alpha_1) \times L(I_2, \alpha_2),$$
        where the coproduct is taken over subsets $I_1, I_2 \subseteq I$ inducing an isomorphism $I_1 \amalg I_2 \xrightarrow{\cong} I$. Restricting $\alpha \colon I \to C$ to these subsets defines objects $(I_1, \alpha_1)$ and $(I_2, \alpha_2)$ of $\Prop(\mathcal{P})$.

        \item If $M$ and $L$ are fibrant modules, then $M \boxplus L$ is also a fibrant module.
    \end{enumerate}
    \label{coend_simplify_modules}
\end{proposition}

\begin{proof}
    From the definition of the symmetric monoidal envelope, a morphism 
    $$ (f, p_1, p_2) : (I, \alpha) \longrightarrow (J_1, \alpha_1) \boxplus (J_2, \alpha_2)$$
    in $\mathsf{Env}(\mathcal{P})$ is determined by a function of finite sets $f \colon I \to J_1 \amalg J_2$, together with a certain tuple $(p_1, p_2)$ of operations in $\mathcal{P}$, with in turn $p_i$ a tuple of operations with output colours in $(J_i, \alpha_i)$. 
    
    Setting $I_i = f^{-1}(J_i)$ and $f_i$ to be the restriction of $f$ to $I_i$, one can easily check that there is a factorization in $\mathsf{Env}(\mathcal{P})$
\[\begin{tikzcd}[cramped]
	{(I, \alpha)} & {(J_1, \beta_1) \boxplus(J_2, \beta_2)} \\
	{(I_1, \alpha_1) \boxplus (I_2, \alpha_2)}
	\arrow["{(f, p_1,p_2)}", from=1-1, to=1-2]
	\arrow["\cong", from=2-1, to=1-1]
	\arrow["{(f_1,p_1) \boxplus (f_2, p_2)}"', from=2-1, to=1-2]
\end{tikzcd}\]
where the right hand morphism is the isomorphism induced by the inclusion $I_1 \amalg I_2 \xrightarrow{\cong} I$ and such that all the associated operations of $\mathcal{P}$ are identities. Comparing this with the coend formula for the tensor product gives us (a). Statement (b) is an immediate consequence of (a) now.
\end{proof}

    For the next result, given $F \in \Fun^{\mathcal{V}}(\Prop(\mathcal{P})^{\mathsf{op}}, \mathcal{V})$ and $K \in \sSets$, we will write $F \boxtimes K \in \Fun^{\mathcal{V}}(\Prop(\mathcal{P})^{\mathsf{op}}, \mathcal{V})$ for the functor given by
    $$ (F \boxtimes K)(I, \alpha) = F(I, \alpha) \times K.$$

    \begin{proposition}
    Let $\mathcal{P}$ be a coloured simplicial operad. Then the monoidal category $(\RMod{\mathcal{P}}, \boxplus, \mathcal{Y}_\emptyset)$ carries the structure of a monoidal model category with respect to the projective model structure.
    \label{monoidality_projective}
    \end{proposition}
    \begin{proof}
        Since the unit $\mathcal{Y}_{\emptyset}$ of the monoidal structure is projectively cofibrant, we only need to the check that the tensor product defines a left Quillen bifunctor.
        
        The cofibrations and trivial cofibrations in $\RMod{\mathcal{P}}$ are generated by 
        $$ \mathcal{Y}_{(I, \alpha)} \boxtimes A \longrightarrow \mathcal{Y}_{(I, \alpha)} \boxtimes B,$$
        where $A \to B$ stands for horn inclusions $\Lambda^i [n] \to \Delta^n$ in the case of cofibrations, and boundary inclusions $\partial \Delta^n \to \Delta^n$ for the trivial cofibrations. The pushout-product with respect to $\boxplus$ of such a map with $\mathcal{Y}_{(I', \alpha')} \boxtimes A' \longrightarrow \mathcal{Y}_{(I', \alpha')} \boxtimes B'$ is given by
        $$ \mathcal{Y_{(I, \alpha) \boxplus (I', \alpha)}} \boxtimes \left( (A \times B')\cup_{A \times A'} (B \times A ')\right) \longrightarrow  \mathcal{Y_{(I, \alpha) \boxplus (I', \alpha)}} \boxtimes (B \times B'),$$
        where we used that the Yoneda embedding is strong monoidal with respect to Day convolution. It is now immediate from this formula that $\boxplus \colon \RMod{\mathcal{P}} \times \RMod{\mathcal{P}} \to \RMod{\mathcal{P}}$ is a left Quillen bifunctor with respect to the projective model structure.
    \end{proof}

\section{The contravariant model structure on $\fSets_{/V}$}

We will start this section by giving a brief overview of the tools we will need from the theory of dendroidal sets, as first introduced by Moerdijk and Weiss in \cite{MoerdijkWeissDendroidal}, followed by the framework of forest sets, as first studied in \cite{HeutsHinnichMoerdijk}, that will be the one we will actually make use of in this work. Our preference for the forest formalism instead of the dendroidal one is a consequence from the way we defined operadic right modules $M$. Namely, the data of such an algebraic structure is determined by sets
$$M(c_1, \ldots, c_n)$$
for each $n$-tuple of colours $(c_1,\ldots, c_n)$ of $\mathcal{P}$, and this influences our choice of framework as it will make it more convenient to have multiple root edges instead of a single one.

After a giving a quick account of some aspect of the homotopy theory of forest sets in Section 3.2 and 3.3, we will construct a version of the contravariant model structure on forest sets over a fixed presheaf $V$. We will explain how this model category on forest sets relates to the study of right fibrations in this category, and in the final section we will shed some light on how the theory of operadic right modules appears within the forest framework.

\subsection{The categories $\Tree$ and $\dSets$}

\leavevmode

The tree category $\Tree$ was first introduced in \cite{MoerdijkWeissDendroidal} as a diagram category that parameterizes an object equipped with multiple algebraic operations with finite inputs. In this sense, one should think of $\Tree$ as playing a similar role to the one the simplex category $\mathsf{\Delta}$ with respect to category theory, but now for the theory of operads.

Let us start by specifying what the objects of $\Tree$ are. We will be very brief in our definitions but the interested reader should consult \cite{HeutsMoerdijkDendroidal} for a more thorough reference.

\begin{definition}
    The objects of $\Tree$ are trees which are finite, non-planar and rooted.
\end{definition}

To explain what rooted means above, we mention that the edges of a tree can be of two possible types: \textit{inner edges}, which are adjacent to two vertices, and \textit{outer edges}, which are adjacent to only one vertex. We then say that a tree $T$ is \textit{rooted} if an outer edge has been singled out (which we call the \textit{root edge} of $T$).

Let us explain some more important vocabulary concerning trees $T$. There is a natural ordering on the set of edges $E(T)$, namely we will say $e_1 \leq e_2$ if the unique path in $T$ from $e_2$ to the root contains $e_1$. It follows from this that either one of two things can happen at each vertex $v$:
\begin{itemize}
    \item The vertex is adjacent to one and only one edge, in which case we say $v$ is a \textit{stump}.
    \item The vertex is adjacent to an edge $e_0$ which is lower in the ordering above to all other incoming edges $e_1, \ldots, e_n$. In this case we say $e_0$ is the \textit{outcoming edge} of $v$, and the remaining edges form the \textit{incoming edges}.
\end{itemize}

There is also a unique vertex which contains the root edge as its outcoming edge, which we will refer to as the \textit{root vertex}. Besides this special vertex, we will also single out the \textit{leaf vertices}, which are the vertices of $T$ which are either stumps or all their incoming edges are outer edges of $T$. Associated with this, we also have the notion of a \textit{leaf}, which is any outer edge which is not the root edge.

Let us also single out some special examples of trees that have a prominent role in the dendroidal framework: we will write $\eta$ for the tree which is just a single edge; for each $k \geq 0$ we have a \textit{corolla} $C_k$, which is a tree with exactly one vertex and $k$ incoming edges; finally, for each $k \geq 0$ we have the \textit{linear tree} $\Delta^k$ which has $k$ vertices, all of which are adjacent to exactly 2 edges (this completely describes such trees).

\begin{example}
Below we have an example of an object $T$ of $\Tree$, with some of the vertices and edges labelled:

    
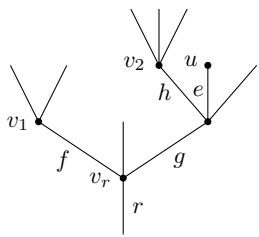
\begin{center}
\begin{tikzpicture}[scale=0.75]
    \draw (0,1)--(0,0)--(0,-1);
     \filldraw[black] (0,0) circle (1.5pt);
    \draw (-1.5, 1)--(0,0);
    \filldraw[black] (-1.5,1) circle (1.5pt);
    \draw (-2, 2) -- (-1.5, 1) -- (-1, 2);
    \draw (1.5, 1)--(0,0);
    \filldraw[black] (1.5,1) circle (1.5pt);
    \draw ({1.5-sin(60)}, 2)--(1.5, 1)--({1.5+sin(60)}, 2);
    \draw (1.5,2)--(1.5, 1);
    \filldraw[black] (1.5,2) circle (1.5pt);
    \filldraw[black] ({1.5-sin(60)},2) circle (1.5pt);
    \draw ({1.5-sin(60)-sin(30)},3)--({1.5-sin(60)},2)--({1.5-sin(60)+sin(30)},3);
    \draw ({1.5-sin(60)},3)--({1.5-sin(60)},2);
     \node [label={[xshift=0.2cm, yshift=-0.7cm, scale=0.9]$r$}] {};
    \node [label={[xshift=-1.4cm, yshift=0.4cm, scale=0.9]$v_1$}] {};
    \node [label={[xshift=0.15cm, yshift=1.2cm, scale=0.9]$v_2$}] {};
    \node [label={[xshift=0.9cm, yshift=1.25cm, scale=0.9]$u$}] {};
    \node [label={[xshift=-0.8cm, yshift=-0.15cm, scale=0.9]$f$}] {};
    \node [label={[xshift=0.75cm, yshift=-0.1cm, scale=0.9]$g$}] {};
    \node [label={[xshift=1cm, yshift=0.85cm, scale=0.9]$e$}] {};
    \node [label={[xshift=0.55cm, yshift=0.8cm, scale=0.9]$h$}] {};
   \node [label={[xshift=-0.3cm, yshift=-0.35cm, scale=0.9]$v_r$}] {};
\end{tikzpicture}
\end{center}

The inner edges of $T$ are exactly the edges $e$, $f$, $g$ and $h$, with $r$ being the root edge and all the other unnamed remaining ones being the leaves. The leaf vertices are $u$, $v_1$ and $v_2$, the root vertex is $v_r$, and the (only) stump is $u$. For this tree, we have that $r$ is the minimal edge (as it is in any tree), $g \leq h$, and $f$ is not comparable with both $e, g$ and $h$.

Each vertex defines a corolla, with the leaves being the incoming edges and the root the outcoming edge; and for instance the path from any leaf of $v_2$ to the root corresponds to a copy of the linear tree $\Delta^3$.    
\label{example_trees}
\end{example}

The morphisms of $\Tree$ are a bit more involved: briefly, every tree $T \in \Tree$ defines a coloured operad $\Tree(T)$, and we set the morphisms $T_1 \to T_2$ to be exactly the operad maps $\Tree(T_1) \to \Tree(T_2)$. For the sake of brevity, we will not follow this route, but see Section 3 of \cite{MoerdijkWeissDendroidal} for more details on the definition of $\mathsf{\Omega}(T)$. Instead, we will provide a more convenient and practical description in terms of generators and relations, where the generating morphisms can be of the following types:
\begin{itemize}
    \item If $e$ is an inner edge of $T$, then there is an \textit{inner face map} $\partial_e T \to T$, where $\partial_e T$ is the tree obtained from $T$ by contracting the edge $e$.
    \item If $v$ is a leaf vertex, then there is a \textit{leaf face map} $\partial_v T \to T$, where $\partial_v T$ is the tree obtained from $T$ by removing the vertex $v$ and all of its incoming edges.
    \item If $T$ is a tree with a root vertex such that has a unique  incoming edge $e$ which is also an inner edge, then there is a \textit{root face map} $\partial_r T \to T$, where $\partial_r T$ is the tree obtained from $T$ by removing the root vertex and keeping the part of $T$ not containing this vertex.
    \item If $e$ is any edge of $T$, then there is a \textit{degeneracy map} $\sigma_e T \to T$, where $\sigma_e T$ is obtained from $T$ by dividing in half the edge $e$ via adding an extra vertex.
    \item The isomorphisms of trees $S \xrightarrow{\cong} T$.
\end{itemize}

These generating operations are related to each other via the \textit{dendroidal identities}, which are a more elaborate version of the usual simplicial relations and are listed in Section 3.3.4 of \cite{HeutsMoerdijkDendroidal}. It can then be shown that any morphism $T_1 \to T_2$ in $\Tree$ can be written, uniquely up to isomorphism, as the composite of
$$ T_1 \xrightarrow{\sigma} T'_1 \xrightarrow{\cong} T'_2 \xrightarrow{\partial} T_2,$$
where $\sigma$ is a degeneracy (that is, a composite of elementary degeneracies) and $\partial$ is a face map (that is, a composite of elementary face maps). For more details on the structure of the morphisms of the tree category and the proofs for the facts mentioned, see Section 3.3 of \cite{HeutsMoerdijkDendroidal}.

From the point of operad theory, the main interest in the tree category $\Tree$ stems from the existence of the \textit{dendroidal nerve}, which is a fully faithful functor \citep[Ex. 4.2]{MoerdijkWeissDendroidal}
$$ N : \mathsf{Op} \longrightarrow \dSets \hspace{2em} \text{given by} \hspace{2em} N\mathcal{P}_T = \mathsf{Op}\left( \Omega(T), \mathcal{P} \right)$$
from the category of operads $\mathsf{Op}$ to the category of \textit{dendroidal sets} $\dSets$, the latter being defined as the presheaf category $\Fun(\Tree^{\mathsf{op}}, \Sets)$. Here one should think of an element $x \in N\mathcal{P}_T$ as a copy of the tree $T$ decorated with operations of $\mathcal{P}$ at each vertex, with the incoming edges and outcoming edge of said vertices being labelled with the corresponding input and output colours.

The existence of the linear trees $\Delta^k$ defines a fully faithful functor $\iota: \mathsf{\Delta} \to \Tree$. After left Kan extending along this functor, one defines a fully faithful colimit-preserving functor at the presheaf level $\iota_! \colon \sSets \to \dSets$, explicitly given by
$$
(\iota_! X)_T = 
\begin{cases}
X_k  &\mathrm{if} \: T = \Delta^k \: \text{for some} \: k,\\
    \emptyset &\mathrm{otherwise.}  
    \end{cases}
$$

The functor allows us to relate the categorical and the dendroidal nerves, since it takes part in the following commutative diagram:
$$
\begin{tikzcd}
\mathsf{Cat} \arrow[d, hook] \arrow[r, "N"] & \sSets \arrow[d, hook, "\iota_!"] \\
\mathsf{Op} \arrow[r, "N"]                  & \dSets.                
\end{tikzcd}
$$

Here the left inclusion sends a category $\mathcal{C}$ to the operad with colours given by the objects of $\mathcal{C}$, operations of arity one $c \to d$ coinciding with the original categorical morphisms with source $c$ and target $d$, and all the sets of operations of higher arity are empty.

In order to study the category of dendroidal sets, the notion of a horn of a representable $T$ plays a major role. These can exist in three types:
\begin{itemize}
    \item If $e$ is an inner edge of $T$, then the \textit{inner horn} $\Lambda^e T \subseteq T$ is the union of all elementary faces of $T$, except for $\partial_e T$.
    \item If $v$ is a leaf vertex of $T$, then the \textit{leaf horn} $\Lambda^v T \subseteq T$ is the union of all elementary faces of $T$, except for $\partial_v T$.
    \item If $T$ admits a root face, then the \textit{root face} $\Lambda^r T \subseteq T$ is the union of all elementary faces, except for the root face $\partial_r T$.
\end{itemize}

We can also define the \textit{boundary} of a tree $T$, which is the dendroidal subset $\partial T \subseteq T$ given by the union of \textit{all} faces of $T$. We observe that, via the inclusion $\iota_! \colon \sSets \to \dSets$, the inner, leaf and root horns will correspond to the simplicial inner, left and right horns, respectively, and the simplicial boundary to the dendroidal boundary.

The existence of horns gives us a clear way of constructing a model for an $\infty$-operad by mimicking Joyal's definition of a quasicategory.

\begin{definition}
    Let $X$ be a dendroidal set. We say $X$ is a \textit{(dendroidal) $\infty$-operad} if, for any tree $T \in \Tree$ and inner edge $e \in E(T)$, the solid diagram
    $$
    \begin{tikzcd}
\Lambda^e T \arrow[r] \arrow[d, hook] & X \\
T \arrow[ru, dashed]                  &  
\end{tikzcd}
    $$
    admits a dashed extension.
    \label{definition_infinty_operad}
\end{definition}

In the category of dendroidal sets, one also has a notion of anodyne morphisms and its different versions, together with their fibration counterparts. Since these will play an essential role in the sections to follow, let us briefly go over their definitions here.

\begin{definition}
    Let $\mathcal{A}$ be a class of morphisms in a small category $\mathcal{C}$. We say $\mathcal{A}$ is \textit{saturated} if $\mathcal{A}$ is closed under pushouts, transfinite compositions and retracts.
\end{definition}

\begin{definition}
    Let $f \colon X \to Y$ be a dendroidal map.
    \begin{itemize}
        \item We say $f$ is a \textit{normal monomorphism} if $f$ is a monomorphism and, for each tree $T$, the group $\mathsf{Aut}_{\Tree}(T)$  acts freely on the set $Y_T - f\left( X_T \right)$. We say $Y$ is a \textit{normal dendroidal set} if $\emptyset \to Y$ is a normal monomorphism, that is, $\mathsf{Aut}_{\Tree}(T)$ acts freely on $Y_T$ for every tree $T$.
        \item We say that $f$ is \textit{inner anodyne} if it belongs to the saturated class generated by the inner horn inclusions $\Lambda^e T \to T$. A dendroidal map is an \textit{inner fibration} if it has the left lifting property with respect to all inner anodynes.
        \item We say that $f$ is \textit{leaf anodyne} if it belongs to the saturated class generated by the leaf horn inclusions $\Lambda^v T \to T$ and the inner horn inclusions. A dendroidal map is an \textit{left fibration} if it has the left lifting property with respect to all leaf anodynes.
        \item We say that $f$ is \textit{root anodyne} if it belongs to the saturated class generated by the root horn inclusions $\Lambda^r T \to T$, whenever $T$ admits a root face map, and the inner horn inclusions. A dendroidal map is an \textit{right fibration} if it has the left lifting property with respect to all root anodynes.
    \end{itemize}
\end{definition}

One last aspect of dendroidal sets we will want to mention is the existence of a version of the Boardman--Vogt tensor product. This was first defined by Moerdijk--Weiss
in \cite{MoerdijkWeissDendroidal} and is given on trees $T_1$ and $T_2$ by the formula
$$ T_1 \otimes T_2 = N \left( \Tree(T_1) \otimes \Tree(T_2) \right),$$
where $\Tree(T_1) \otimes \Tree(T_2)$ denotes the Boardman--Vogt tensor product of operads. We can then extend it to $\dSets$ by taking left Kan extensions on each variable, so that in particular the dendroidal tensor product will preserve colimits in each variable separately.

As the Boardman--Vogt product of categories (seen as operads with only unary operations) reduces to the cartesian product of small categories, it follows that the dendroidal tensor product restricted to $\sSets$ coincides with the cartesian product. However, we note that on $\dSets$ this bifunctor is quite more complicated and has been thoroughly analysed in the literature, see Chapters 4 and 6 of \cite{HeutsMoerdijkDendroidal} for a survey of some these results. Although we will need some of these statements, it will actually be more purposeful for us to state them within the context of forest sets.

\subsection{The categories $\For$ and $\fSets$}

\leavevmode

As we have previously mentioned, our study of operadic right modules will require us to extend the category $\Tree$ to allow for \textit{forests}, that is, finite collections of trees. We will do this by considering a new category $\For$, first studied by Heuts, Hinnich and Moerdijk in \cite{HeutsHinnichMoerdijk} as an intermediate step in comparing Lurie's notion of an $\infty$-operad, and the dendroidal model we presented in \cref{definition_infinty_operad}.

Recall that, for any tree $T$, the set of edges $E(T)$ admits a poset structure, where $e_1 \leq e_2$ if the unique path in $T$ from $e_2$ to the root edge contains $e_1$. We say that two non-empty subsets $A, B \subseteq E(T)$ are \textit{independent} if the elements of $A$ can't be compared to the elements of $B$ in this ordering.

\begin{definition}
    The category of forests $\For$ is the category given by the following data:
    \begin{itemize}[label=$\diamond$]
        \item The objects are finite collections $\{ T_i : i \in I \}$ of objects of $\Tree$, where $I$ is a finite set. We will usually write $\bigoplus_{i \in I} T_i$ for this object, and refer to trees $T_i$ as the \textit{components} of the forest.
        \item A morphism in $\For$
        $$\bigoplus_{i \in I} T_i \longrightarrow \bigoplus_{j \in J} S_j $$ 
        is determined by a pair $(f, \{ \alpha_i : i \in I \})$, where $f\colon I \to J$ is a function of finite sets and $\alpha_i \colon T_i \to S_{f(i)}$ is a morphism in $\Tree$, for each $i \in I$. We additionally ask for the following condition to be satisfied: for any $j \in J$ and distinct elements $a, b \in f^{-1}(j)$, the subsets $\alpha_a(E(T_a))$ and $\alpha_b(E(T_b))$ of the edge set $E(S_j)$ are independent. 
        \item The composition is determined by the composition of functions of finite sets and the composition in $\Tree$.
    \end{itemize}
    \label{forest_definition}
\end{definition}

\begin{notation}
    If $F$ is a forest, we will write $R(F)$ for the set of root edges of $F$, which are the roots edges of each component tree of $F$. The roots of $F$ naturally form a forest $\rho(F)$, together with a forest monomorphism $\rho(F) \to F$.
\end{notation}

The forest category comes equipped with a bifunctor
$$ - \oplus - : \For \times \For \longrightarrow \For$$
which sends a pair of forests $(F,G)$ to the forest $F \oplus G$ formed by both $F$ and $G$. As already remarked in \cite{HeutsHinnichMoerdijk}, we observe that $\oplus$ doesn't define a coproduct for $\For$: indeed, the extra independence restriction that appears on the definition of the morphisms of $\For$ implies that a codiagonal map $\mathrm{id} \oplus \mathrm{id} \colon T \oplus T \to T$ does not exist for any tree $T$. 

Moreover, seeing any tree as a forest with a single component gives rise to a fully faithful inclusion $u \colon \Tree \hookrightarrow \For$. For the sake of notation, whenever $T \in \Tree$ is a tree, we will usually write the forest $u(T)$ simply as $T$. An important change that comes from adopting the forest framework instead of the dendroidal one is that any tree $T$ now admits a root face: indeed, any tree $T \neq \eta$ can be constructed by grafting to the leaves of its root corolla $C_k$ a certain finite family of trees $T_1, \ldots, 
 T_k$. As this is a convenient inductive procedure for constructing trees, let us fix the following notation
$$ T = C_k \circ (T_1, \ldots, T_k)$$
for this grafting procedure.  We set $\partial_r T = T_1 \oplus \cdots \oplus T_k \to T$ to be the \textit{root face map} of $T$; note that this new face in general will not coincide with the root face of $T$ when it exists! This happens only when the root vertex has exactly one incoming edge, as in the otherwise $\partial_r T$ will have more than one component.

\begin{example}

The tree $T$ that we considered in \cref{example_trees} doesn't admit a root face in the category $\Tree$. However, it admits one in $\For$ (as does any other tree), which is represented in the figure below, together with the labels we already used in the original example.

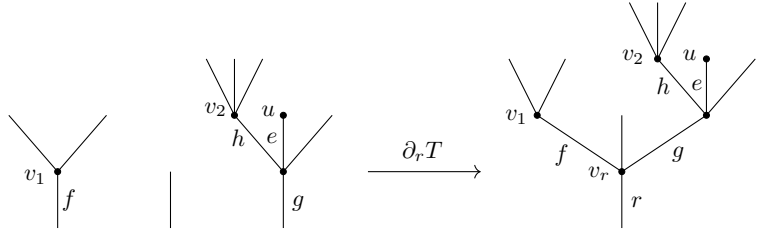
\begin{center}
\begin{tikzpicture}[scale=0.75]
    \draw (-10,0) -- (-10,-1);
    \filldraw[black] (-10,0) circle (1.5pt);
    \draw ({-10-sin(60)}, 1)--(-10, 0)--({-10+sin(60)}, 1);
    \draw (-8,0) -- (-8,-1);
    \draw (-6,1) -- (-6,0) -- (-6,-1);
    \filldraw[black] (-6,0) circle (1.5pt);
    \draw ({-6-sin(60)}, 1)--(-6, 0)--({-6+sin(60)}, 1);
    \draw ({-6-sin(60)-sin(30)},2)--({-6-sin(60)},1)--({-6-sin(60)+sin(30)},2);
    \draw ({-6-sin(60)},2)--({-6-sin(60)},1);
     \filldraw[black] (-6,1) circle (1.5pt);
     \filldraw[black] ({-6-sin(60)},1) circle (1.5pt);
     \draw[->]        (-4.5,0)   -- (-2.5,0);
    \node [label={[xshift=-7.35cm, yshift=-0.75cm, scale=0.9]$f$}] {};
    \node [label={[xshift=-4.3cm, yshift=-0.75cm, scale=0.9]$g$}] {};
    \node [label={[xshift=-7.8cm, yshift=-0.4cm, scale=0.9]$v_1$}] {};
    \node [label={[xshift=-4.65cm, yshift=0.15cm, scale=0.9]$e$}] {};
    \node [label={[xshift=-5.1cm, yshift=0.1cm, scale=0.9]$h$}] {};
    \node [label={[xshift=-4.7cm, yshift=0.5cm, scale=0.9]$u$}] {};
    \node [label={[xshift=-5.4cm, yshift=0.5cm, scale=0.9]$v_2$}] {};
    \node [label={[xshift=-2.65cm, yshift=-0.1cm, scale=0.9]$\partial_r T$}] {};
    
    \draw (0,1)--(0,0)--(0,-1);
     \filldraw[black] (0,0) circle (1.5pt);
    \draw (-1.5, 1)--(0,0);
    \filldraw[black] (-1.5,1) circle (1.5pt);
    \draw (-2, 2) -- (-1.5, 1) -- (-1, 2);
    \draw (1.5, 1)--(0,0);
    \filldraw[black] (1.5,1) circle (1.5pt);
    \draw ({1.5-sin(60)}, 2)--(1.5, 1)--({1.5+sin(60)}, 2);
    \draw (1.5,2)--(1.5, 1);
    \filldraw[black] (1.5,2) circle (1.5pt);
    \filldraw[black] ({1.5-sin(60)},2) circle (1.5pt);
    \draw ({1.5-sin(60)-sin(30)},3)--({1.5-sin(60)},2)--({1.5-sin(60)+sin(30)},3);
    \draw ({1.5-sin(60)},3)--({1.5-sin(60)},2);
     \node [label={[xshift=0.2cm, yshift=-0.7cm, scale=0.9]$r$}] {};
    \node [label={[xshift=-1.4cm, yshift=0.4cm, scale=0.9]$v_1$}] {};
    \node [label={[xshift=0.15cm, yshift=1.2cm, scale=0.9]$v_2$}] {};
    \node [label={[xshift=0.9cm, yshift=1.25cm, scale=0.9]$u$}] {};
    \node [label={[xshift=-0.8cm, yshift=-0.15cm, scale=0.9]$f$}] {};
    \node [label={[xshift=0.75cm, yshift=-0.1cm, scale=0.9]$g$}] {};
    \node [label={[xshift=1cm, yshift=0.85cm, scale=0.9]$e$}] {};
    \node [label={[xshift=0.55cm, yshift=0.8cm, scale=0.9]$h$}] {};
   \node [label={[xshift=-0.3cm, yshift=-0.35cm, scale=0.9]$v_r$}] {};
\end{tikzpicture}
\end{center}

\end{example}

The following easy proposition shows that one can single out a simple set of morphisms generating $\For$.

\begin{lemma}
    The morphisms of $\For$ are generated under composition by the following types of morphisms:
    \begin{enumerate}[label=(\alph*)]
        \item For any trees $T, S \in  \Tree$, any morphism $T \to S$ coming from $\Tree$.
        \item For any forest $F, G \in \For$, all summand inclusions $F \to F \oplus G$.
        \item For any tree $T \in \Tree$, the root face $\partial_{\mathrm{root}}T \to T$.
    \end{enumerate}
    together with the condition that $- \oplus -\colon \For \times \For \to \For$ is a bifunctor satisfying symmetry and associativity.
    \label{arrows_For}
\end{lemma}

\begin{proof}
    See \citep[Lem. 3.1.1]{HeutsHinnichMoerdijk} for a proof.
\end{proof}

The forest category also admits a factorization of its morphisms which is very similar to the one we mentioned for dendroidal sets. In order to explain how this works, let us first define what are the elementary faces and degeneracies of a forest.

\begin{definition}
    Consider a forest $F$ of the form $F= F' \oplus T \oplus F''$, where $T$ is a tree, and $F'$ and  $F''$ are (possibly empty) forests.
    \begin{enumerate}[label=(\alph*)]
        \item The \textit{elementary degeneracies} of $F$ are the maps of the form
        $$ F' \oplus \sigma_e \oplus F'' : F' \oplus \sigma_e T \oplus F'' \longrightarrow F'\oplus T \oplus F'',$$
        where $e$ is any edge of $T$. The \textit{degeneracies} of $F$ are the morphisms which can be written as a composition of elementary degeneracies.
        \item The \textit{elementary faces} of $F$ are the maps of the form
        $$ F' \oplus \partial_{\bullet} \oplus F'' : F'\oplus \partial_{\bullet} T \oplus F'' \longrightarrow F' \oplus T \oplus F'',$$
        where $\partial_{\bullet} T$ is either an elementary inner or leaf face of $T$, or the root face of $T$; we also allow for the map $\partial_{\bullet} T \to T$ to be $\emptyset \to \eta$. The \textit{faces} of $F$ are the morphisms which can be written as a composition of elementary faces.
    \end{enumerate}
\end{definition}

\begin{proposition}
    Any morphism $F_1 \to F_2$ in $\For$ can be written uniquely up to isomorphism as a composition of the form
    $$ F_1 \xlongrightarrow{\sigma} F_1' \xlongrightarrow{\cong} F_2' \xlongrightarrow{\partial} F_2,$$
    where $\sigma$ is a degeneracy and $\partial$ is a face map.
\end{proposition}
\begin{proof}
    See \citep[Lem. 3.1.3]{HeutsHinnichMoerdijk} for a proof.
\end{proof}

Taking presheaves on $\For$ leads to the category of \textit{forest sets}, which we will write as $\fSets$. Within this setting, we can also extend the forest operation $\oplus$ to a bifunctor
$$- \oplus - : \fSets \times \fSets \longrightarrow \fSets$$
by left Kan extending along each variable. In addition to this, taking left and right Kan extensions along $u$ leads to a diagram of adjoint functors the category of dendroidal sets $\dSets$ to the category of forest sets $\fSets$:
$$
\begin{tikzcd}
\dSets \arrow[rr, "u_!", bend left] \arrow[rr, "u_*", bend right] &  & \fSets \arrow[ll, "u^*"'].
\end{tikzcd}
$$

In general, it is difficult to give a practical and simple expression for the left Kan extension $u_! X$ for $X \in \dSets$, since, for a fixed tree $T$, the overcategory $\For_{/T}$ can be quite complicated\footnote{We note however that if $T= \Delta^k$ is a linear tree, than $T$ only admits maps from other linear trees. Therefore $u_!$ is easy to describe when restricted to $\sSets$.}. However, one can easily check using the adjointness of the pair $(u^\ast, u_\ast)$ that $u_\ast X$ evaluates on a forest $\bigoplus_{i=1}^n T_i$ to $\prod_{i=1}^n X_{T_i}$, and sends a morphism $(f, \{ \alpha_i\}) \colon \bigoplus_{i \in I} T_i \to \bigoplus_{j \in J} S_j$ to the composition of
$$ \prod_{j\in J} X_{S_j} \xrightarrow{\mathrm{project \: onto \:  im}(f)} \prod_{j \in \mathrm{im}(f)} X_{S_j} \xrightarrow{\prod_{i \in I} (\alpha_i)^*} \prod_{i \in I} X_{T_i}.$$

As one would expect, the relevant notions for the homotopy theory of dendroidal sets that we have already discussed in the previous section have an analogue for forest sets. Since these definitions are quite similar to the dendroidal case, we will quickly go through them, and the interested reader can read about the details in Section 3 of \cite{HeutsHinnichMoerdijk}.

\begin{itemize}[label=$\diamond$]
    \item Any forest $F \in \For$ admits a \textit{boundary} $\partial F \in \fSets$, given as the colimit 
    $$ \partial F = \colim_{F' \to F} F'$$
    taken over all face maps $F' \to F$ in $\For$. The forest boundary satisfies the equations 
    $$\partial(F \oplus G) = \partial F \oplus G \cup_{\partial F \oplus \partial G} F \oplus \partial G \hspace{1.5em } \mathrm{and} \hspace{1.5em} \partial( u^*F ) = u^*(\partial F)$$
    for any forests $F, G$. 

    \item The horn inclusions have immediate generalizations, and they can either be \textit{inner, leaf} or \textit{root horns}. The only difference compared to $\dSets$ is that now any non-empty forest $F$ admits a root horn $\Lambda^r F$ for each choice of root edge $r$.
    
    \item A monomorphism of forest sets $f\colon X \to Y$ is said to be a \textit{normal monomorphism} if, given any forest $F \in \For$, the group of automorphisms $\mathsf{Aut}_{\For}(F)$ acts freely on the set $Y_F -\mathrm{im}(f)$. In the special case when $X = \emptyset$, we say that $Y$ is a \textit{normal forest set}.

    Normal forest sets admit a skeletal filtration according to \cite[Prop. 3.3.5]{HeutsHinnichMoerdijk}, which in particular shows that the class of forest normal monomorphisms is the saturated class generated by the boundary inclusions $\partial F \to F$.

    \item The definition of being \textit{forest inner/leaf/root anodyne} can be done in exactly the same way as we did for dendroidal sets using saturated classes. A similar remark also applies for defining the respective classes of fibrations. 
    
    \item The Boardman--Vogt tensor product of dendroidal sets can be extended to forests sets by setting
    $$ F \otimes G = \left( \bigoplus_{i=1}^n T_i \right) \otimes \left( \bigoplus_{j=1}^m S_j \right) = \bigoplus_{i=1}^n \bigoplus_{j=1}^m u_! \left( T_i \otimes S_j\right) $$
    for $F, G \in \For$, and stating that it commutes with colimits on both variables. 
\end{itemize}

\subsection{The tensor product of forest sets and forest root anodynes}

\leavevmode

In this section we will focus on proving some technical results about the theory of the forest sets that we will frequent in the subsequent sections. We will mostly focus on studying the class of forest root anodynes and on how it interacts with the Boardman--Vogt tensor product. The reader should feel free to skip this section, except possibly for the definition of a forest spine in \cref{definition_forest_spine}, and refer back to it whenever some result is used.

We begin our analysis by presenting some of the already-known interactions of the Boardman--Vogt tensor product with the different classes of anodyne morphisms we have presented. Before stating it, we define a certain special class of forest root anodynes.

\begin{definition}
    The class of \textit{forest unitary root anodynes} is the smallest saturated class of forest morphisms which contains all inner horn inclusions and all root horn inclusions $\Lambda^r F \to F$, where $r$ is a root of $F$ which has only one incoming edge.
\end{definition}

\begin{remark}
    One should be warned about the usage of the expression \textit{root anodyne} in \cite{HeutsHinnichMoerdijk} for two reasons: they not only exclude inner anodynes from the generators of the saturated class, but they only include \textit{unitary} root horn inclusions in their generators. We decided to stick to a broader class of root anodynes since this makes the definition of forest right fibration cleaner.
\end{remark}

\begin{proposition}
    Let $f\colon M \to N$ be a monomorphism of simplicial sets and $g\colon X \to Y$ a normal monomorphism of forest sets, and consider its pushout-product
    \begin{equation}
     f \hat{\otimes} g : M \otimes Y \cup_{M \otimes X} N \otimes X \longrightarrow N \otimes Y
    \label{pushout_product}
    \end{equation}
    with respect to the Boardman--Vogt tensor product of forest sets. Then the following hold:
    \begin{enumerate}[label=(\alph*)]
        \item The puhsout-product $ f \hat{\otimes} g$ is a normal monomorphism of forest sets.
        \item If $g$ is an inner anodyne, then $ f \hat{\otimes} g$ is also an inner anodyne.
        \item If $g$ is a unitary root anodyne, then $ f \hat{\otimes} g$ is also a unitary root anodyne.
        \item If $f$ is a right anodyne, then $ f \hat{\otimes} g$ is a unitary root anodyne.
    \end{enumerate}
    \label{closure_BV_pushout_product}
\end{proposition}

\begin{proof}
    All these statements are proven at different stages of \cite{HeutsHinnichMoerdijk}: case (a) is Proposition 3.4.1 (i), case (b) is the saturated version of Proposition 3.6.2, and (c) and (d) are unmarked versions of some of the variations mentioned in Proposition 4.2.7.
\end{proof}

\begin{remark}
    There is also a version of case (c) of the previous result where one allows more generally for $g$ to be a forest root anodyne, which is \cref{saturated_rightanodynes}. This is not covered in \cite{HeutsHinnichMoerdijk}, and the proof of it will be given later in this section.
\end{remark}

In the next few results, we will try find certain criteria for showing that certain classes of forest normal monomorphisms contain the class of forest root anodynes.

\begin{definition}
    Let $F \in \For$ and let $R \subseteq R(F)$ be a non-empty subset of roots. We write $\Lambda^R F \subseteq F$ for the \textit{generalized root horn}, given as the union of all elementary faces of $F$, except for the root faces indexed by the elements of $R$. 
    \label{generalized_horn}
\end{definition}

\begin{lemma}
   Let $F \in \For$ and consider a non-empty subset $R \subseteq R(F)$. Then the inclusion $\Lambda^R F \to F$ is a forest root anodyne. 
   \label{generalized_horn_root}
\end{lemma}

\begin{proof}
    If $R$ is a singleton, then $\Lambda^R F$ is the usual root horn inclusion, which is root anodyne by definition. If $R = R' \cup \{ r \}$ for some strictly smaller set $R' \subseteq R$, then there exists a pushout diagram in $\fSets$
\[\begin{tikzcd}[cramped]
	{\Lambda^{R' }(\partial_r F)} & {\Lambda^{R}F} \\
	{\partial_r F} & {\Lambda^{R'}F}
	\arrow[from=1-1, to=1-2]
	\arrow[from=1-1, to=2-1]
	\arrow[from=1-2, to=2-2]
	\arrow[from=2-1, to=2-2]
	\arrow["\lrcorner"{anchor=center, pos=0, rotate=180, scale=1.5}, draw=none, from=2-2, to=1-1]
\end{tikzcd}\]
and, since the left map is root anodyne by induction on the size of $R'$, we conclude that $\Lambda^R F \to \Lambda^{R'} F$ must also be root anodyne. Consequently, so is the composite $\Lambda^R F \to \Lambda^{R'} F \to F$, as we wanted to show.
\end{proof}

Next we define two versions of the spine of a forest, which will be useful when we talk about forest spaces.

\begin{definition}
    Let $F$ be a forest. We define the \textit{forest spine} of $F$ to be the forest set $\fSp(F) \subseteq F$ constructed inductively as follows:
    \begin{enumerate}[label=(\roman*)]
        \item If $F$ is a either a corolla or $\eta$, then $\fSp(F) = F$.
        \item If $G$ and $H$ are forests, then $\fSp(G \oplus H) \cong \fSp(G) \oplus \fSp(H)$.
        \item If $T$ is a tree written as a grafting $C_k \circ (T_1, \ldots, T_k)$, then
        $$ \fSp\left( T \right) = C_k \cup_{k \cdot \eta} \left( \fSp(T_1) \oplus \cdots \oplus \fSp(T_k) \right).$$
    \end{enumerate}
    \label{definition_forest_spine}
    Alternatively, $\fSp(F)$ can also be defined as the colimit over all subforests $G \subseteq F$ where each component has at most one vertex. There is a canonical inclusion $\fSp(F) \to F$ of the forest spine, for any forest $F$.
\end{definition}

We will use the definition above of spine in order to introduce the concept of a Segal forest space. However, in future chapters it will also be convenient to have a slight variation of this construction, which we explain now.

\begin{definition}
    Let $F$ be a forest. We define the \textit{wide spine} of $F$ to be the forest set $\wfSpine(F) \subseteq F$ constructed inductively as follows:
    \begin{enumerate}[label=(\roman*)]
        \item If $F$ is a either a corolla or $\eta$, then $\wfSpine(F) = F$.
        \item If $G$ and $H$ are forests, then $\wfSpine(G \oplus H) \cong \wfSpine(G) \oplus \wfSpine(H)$.
        \item If $T$ is a tree written as a grafting $C_k \circ (T_1, \ldots, T_k)$, then
        $$ \wfSpine\left( T \right) = C_k \cup_{\eta^{\amalg k}} \left( \wfSpine(T_1) \amalg \cdots \amalg \wfSpine(T_k) \right).$$
    \end{enumerate}
    There are canonical inclusions $\wfSpine(F) \to \fSp(F) \to F$ of the spines, for any forest $F$.
\end{definition}

\begin{example}
    In order to illustrate the difference between $\fSp(F)$ and $\wfSpine(F)$, consider the tree represented below
    \begin{center}
        \begin{tikzpicture}[scale=0.75]
            \draw (0,1)--(0,0)--(0,-1);
            \filldraw[black] (0,0) circle (1.5pt);
            \draw (1,1)--(0,0)--(-1,1);
            \filldraw[black] (1,1) circle (1.5pt);
            \filldraw[black] (-1,1) circle (1.5pt);
            \draw (-1.5,2)--(-1,1)--(-0.5,2);
            \draw (1,2)--(1,1);
        \end{tikzpicture}
    \end{center}
The forest spine and the wide spine of this tree are given by the forest sets
$$ C_3 \cup_{ 3 \cdot \eta} \left( C_2 \oplus \eta \oplus C_1 \right) \hspace{2em} \mathrm{and} \hspace{2em} C_3 \cup_{ \eta^{\amalg 3}} \left( C_2 \amalg \eta \amalg C_1 \right)$$
respectively. In general, we note that if $F$ is written as a concatenation of trees $\bigoplus_{i \in I} T_i$, then we have the identities
$$ \fSp(F) = \bigoplus_{i \in I} \fSp(T_i) \hspace{2em} \mathrm{and} \hspace{2em} \wfSpine(F) = \bigoplus_{i \in I} u_! \mathsf{Sp}(T_i).$$
since $\wfSpine(T) = u_! \mathsf{Sp}(T)$ for any tree $T$.
\label{different_spines}
\end{example}

The forest spine interacts well with the class of inner anodynes (as it already happens for their dendroidal and simplicial counterparts), and in a lot of situations they tend to be more workable in practice than inner horn inclusions.

For the next few results, we define a \textit{basic forest} to be a forest $F$ such that its components are either corollas $C_k$ for $k \geq 0$, or edges $\eta$. We also note that the inclusion of the roots $\rho(F) \to F$ always factors through the forest spine of $F$, leading to a map $\rho(F) \to \fSp(F)$.

\begin{lemma}
    Let $F$ be a forest. Then the inclusion $\rho(F) \to \fSp(F)$ can be written as a composition of pushouts of maps of the form $\rho(G) \to G$, where $G$ is a basic forest.
    \label{basic_to_F}
\end{lemma}

\begin{proof}
    We proceed inductively on the number $n(F)=\lvert V(F) \rvert - \lvert R(F) \rvert$, which is always non-negative if we make the convention that components of the form $\eta$ contribute with 0 to $n(F)$. 
    
    The case when $n(F)=0$ happens exactly when $F$ is a basic forest, for which our statement holds since $\fSp(F)=F$ in this scenario. If $F$ is a forest of the form $T_1 \oplus \cdots \oplus T_m$ with each component a tree, then writing $T_i$ as a grafting gives us the formula for the forest spine
    $$ \fSp\left( T_i \right) = C_{k(i)} \cup_{k(i) \cdot \eta} \left( \fSp(T_i^1) \oplus \cdots \oplus \fSp(T_i^{k(i)}) \right) $$
    where $k(i) \geq 0$. Writing $F'$ for the forest $\bigoplus_{i=1}^m \bigoplus_{j=1}^{k(i)} T_i^j$ obtained by taking all the elementary root faces of $F$, we note that the induction hypothesis applies to this forest, since $F'$ has less vertices but more components than $F$ (and thus $n(F') < n(F)$). In light of the square in the diagram below
\[\begin{tikzcd}
	& {\rho(F')} & {\fSp(F')} \\
	{\rho(F)} & {\bigoplus_{i=1}^m C_{k(i)}} & {\fSp(F)}
	\arrow[from=1-2, to=1-3]
	\arrow[from=1-2, to=2-2]
	\arrow[from=1-3, to=2-3]
	\arrow["{\mathrm{root}}", from=2-1, to=2-2]
	\arrow[from=2-2, to=2-3]
        \arrow["\lrcorner"{anchor=center, pos=0, rotate=180, scale=1.5}, draw=none, from=2-3, to=1-2]
\end{tikzcd}\]
    being a pushout square (this can be checked using the axioms for the forest spine), we conclude that the composition of the bottom row is of the desired form, as we wanted to show.
\end{proof}

The next result is fundamental in reducing arguments concerning forest right anodynes to only checking something concerning basic forests and spine inclusions. Before proving it, we recall that, given classes of morphisms $\mathcal{M} \subseteq \mathcal{N}$ in a category $\mathcal{C}$, we will say that $\mathcal{M}$ satisfies the \textit{right cancellation property} within $\mathcal{N}$ if, given morphisms $A \xrightarrow{i} B \xrightarrow{j} C$ satisfying $i, j, ji \in \mathcal{N}$ and $ji, i \in \mathcal{M}$, then $j$ is also in $\mathcal{M}$.

\begin{proposition}
    Let $\mathcal{A}$ be a saturated class of normal monomorphisms in $\fSets$ which satisfies the right cancellation property. Suppose additionally that $\mathcal{A}$ contains one of the following sets of morphisms:
    \begin{enumerate}[label=(\alph*)]
        \item The root inclusions $\rho(F) \to F$ for any forest $F$.
        \item The root inclusions $\rho(G) \to G$ for any basic forest $G$, together with all spine inclusions $\fSp(F) \to F$ for any forest $F$.
    \end{enumerate}
    Then $\mathcal{A}$ contains all forest root anodynes.
    \label{closure_forest_root anodynes}
\end{proposition}

\begin{proof}
    We begin by showing that under the conditions of the first paragraph, $\mathcal{A}$ contains (a) if and only if it contains (b). For that purpose, consider the string
    $$ \rho(F) \xlongrightarrow{i} \fSp(F) \xlongrightarrow{j} F.$$

    If $\mathcal{A}$ satisfies (a), then both $i$ and $ji$ are in $\mathcal{A}$, the former due to Lemma \ref{basic_to_F}; by right cancellation, we get that $j \in \mathcal{A}$. If $\mathcal{A}$ satisfied (b) instead, then $j \in \mathcal{A}$, and another application of Lemma \ref{basic_to_F} shows that $i \in \mathcal{A}$; therefore $ji \in \mathcal{A}$.

    Due to $\mathcal{A}$ being saturated, it suffices that we show that the root horn inclusions of any forest $F$ are in $\mathcal{A}$. In order to do so, we prove the following stronger statement: for any forest $F$ and $R$ any set of roots of $F$, all the maps in the string
    $$ \rho(F) \xlongrightarrow{j_1} V_1 \xlongrightarrow{j_2} V_2 \xlongrightarrow{j_3} \Lambda^R F \xlongrightarrow{j_4} F$$
    belong to $\mathcal{A}$, where:
    \begin{itemize}[label=$\diamond$]
        \item The map $j_1$ attaches a subset of inner faces of $F$.
        \item The map $j_2$ attaches some of the leaf faces of $F$.
        \item The map $j_3$ attaches the remaining faces of $F$.
        \item The map $j_4$ is the usual horn inclusion into $F$.
    \end{itemize}

    This implies the required by taking $R$ to be a singleton. We will show this step by step by working inductively along the number $n(F)=\lvert V(F) \rvert - \lvert F \rvert$ we already defined during the proof of Lemma \ref{basic_to_F}. The case when $n(F)=0$ holds by hypothesis, since this corresponds to basic forests.

    \underline{\textit{The maps of type $j_1$ are in $\mathcal{A}$}}: if $V_1$ is obtained by adding a single inner face $\partial_e F$, then $j_1$ is just the map $\rho(\partial_e F) \to \partial_e F$ since contracting edges doesn't change the root set of $F$, and thus it is in $\mathcal{A}$ by induction. If $V_1$ is obtained from $V_1'$ by adding a single inner face $\partial_e F$, then we canb encode in a pushout diagram
\[\begin{tikzcd}
	{V_1' \cap \partial_e F} & {V_1'} \\
	{\partial_e F} & {V_1.}
	\arrow[from=1-1, to=1-2]
	\arrow[from=1-1, to=2-1]
	\arrow[from=1-2, to=2-2]
	\arrow[from=2-1, to=2-2]
        \arrow["\lrcorner"{anchor=center, pos=0, rotate=180, scale=1.5}, draw=none, from=2-2, to=1-1]
\end{tikzcd}\]

However, the forest set $V_1' \cap \partial_e F$ is given by attaching some set of inner faces of $\partial_e F$ to the root forest $\rho(\partial_e F)$, which implies that the map $V_1'\cap \partial_e F \to \partial_e F$ is of type $j_4 j_3 j_2$. Since $n(\partial_e F) < n(F)$, this map is in $\mathcal{A}$ by induction, and therefore so is $V_1'\to V_1$ by $\mathcal{A}$ being saturated. Consequently, so is the composition of $\rho(F) \to V_1' \to V_1$.

\vspace{0.5em}

    \underline{\textit{The maps of type $j_2$ are in $\mathcal{A}$}}: the argument for this type of morphisms is similar to the previous case, where we now attach each leaf face one at a time, so we omit it.

\vspace{0.5em}

    \underline{\textit{The maps of type $j_3$ are in $\mathcal{A}$}}: to show that $j_3 \in \mathcal{A}$, we proceed by induction on the cardinality of $R(F) - R$.

    For the base case, consider the diagram

\[\begin{tikzcd}
	{\rho(F)} & {V_1} & {V_2} & {\Lambda^R F,} \\
	&& {V_2'}
	\arrow["{j_1}", from=1-1, to=1-2]
	\arrow["{j_2}", from=1-2, to=1-3]
	\arrow["{i_2}"', from=1-2, to=2-3]
	\arrow["{j_3}", from=1-3, to=1-4]
	\arrow["{i_3}"', from=2-3, to=1-4]
\end{tikzcd}\]
where $j_1$ adds all inner faces of $F$ in $V_2$, and $j_2$ adds all leaf faces of $F$ in $V_2$, both of which are in $\mathcal{A}$ by the previous two cases. By the right cancellation property, we only need to show that $j_3 j_2 \in \mathcal{A}$ in order to prove that $j_3 \in \mathcal{A}$.

With that purpose in mind, we decompose $j_3 j_2$ via $i_2$, which adds \textit{all} the remaining inner faces of $F$ not in $V_1$, followed by $i_3$, which attaches \textit{all} the leaf faces of $F$. Due to we being in the base case $R=R(F)$, we note that $i_3$ is just the identity and so is in $\mathcal{A}$. Moreover, both $j_1$ and $i_2 j_1$ are in $\mathcal{A}$ since they are both of type $j_1$; by the right cancellation property, $i_2 \in \mathcal{A}$. Putting these two considerations together, we deduce that $j_3 j_2 = i_3 i_2 \in \mathcal{A}$, and thus by another application of the right cancellation property we have that $j_3 \in \mathcal{A}$.

Suppose now that $R$ is obtained from $R'$ by removing a root $r$, so that we have $\Lambda^{R'} F \to F$ in $\mathcal{A}$ by induction. Then we deduce from the pushout diagram below
\[\begin{tikzcd}
	{\Lambda^{R}(\partial_r F)} & {\Lambda^{R'} F} \\
	{\partial_r F} & {\Lambda^{R} F,}
	\arrow[from=1-1, to=1-2]
	\arrow[from=1-1, to=2-1]
	\arrow[from=1-2, to=2-2]
	\arrow[from=2-1, to=2-2]
        \arrow["\lrcorner"{anchor=center, pos=0, rotate=180, scale=1.5}, draw=none, from=2-2, to=1-1]
\end{tikzcd}\]
that $\Lambda^{R'} F \to \Lambda^R F$ is in $\mathcal{A}$, because $n(\partial_r F) < n(F)$ and $\mathcal{A}$ is saturated. Finally, another application of the right cancellation property to $\Lambda^{R'} F \to \Lambda^R F \to F$ shows that $\Lambda^R F$.

\vspace{0.5em}

    \underline{\textit{The maps of type $j_4$ are in $\mathcal{A}$}}: this follows directly from the right cancellation property applied to $ \rho(F) \xlongrightarrow{j_3 j_2 j_1} \Lambda^R F \xlongrightarrow{j_4} F.$

\end{proof}

Our last main result we will need is an extension of \cref{closure_BV_pushout_product} for the case of a pushout-product of a simplicial monomorphism and a forest root anodyne.

\begin{proposition}
    Let $n \geq 1$ and $F$ a forest. For any root $r \in R(F)$, the pushout-product
    $$ \partial \Delta^n \otimes F \cup_{\partial \Delta^n \otimes \Lambda^r F} \Delta^n \otimes \Lambda^r F \longrightarrow \Delta^n \otimes F$$
    is a forest root anodyne.
    \label{root_anodyne_bv_tensor}
\end{proposition}

Using the saturatedness of the class of root anodynes and the fact that the Boardman--Vogt tensor product commutes with colimits, the following corollary follows immediately.

\begin{corollary}
    Let $f\colon M \to N$ be a monomorphism of simplicial sets and $g\colon X \to Y$ a forest root anodyne. Then its pushout-product
    \begin{equation}
     f \hat{\otimes} g : M \otimes Y \cup_{M \otimes X} N \otimes X \longrightarrow N \otimes Y
    \end{equation}
    is a forest root anodyne.
    \label{saturated_rightanodynes}
\end{corollary}

Given $X, Y \in \fSets$, we will write $\Hom(X, Y) \in \sSets$ for the simplicial set $\Hom(X, Y)_k = \fSets( \Delta^k \otimes X, Y)$ for each $k \geq 0$.

\begin{proposition}
    The class of forest root anodynes satisfies the right cancellation property within the class of forest normal monomorphisms.
    \label{forest_right_cancellation}
\end{proposition}

\begin{proof}
    The proof follows the one in \cite{StevensonInnerFibrations} and \cite{BarataRightCancelation} for the case of simplicial and dendroidal inner anodynes respectively; however, these simplify somewhat for the case of root anodynes, so we will present this simpler version. We will use throughout the following property of the Boardman--Vogt tensor product: if $X \in \fSets$ and $ A \in \sSets$, then there exists a projection map $\pi_X \colon X \otimes A \to X$ onto the forest component.

    We want to consider a lifting problem in $\fSets$ of the bottom square in
\begin{equation}
\begin{tikzcd}
A \arrow[d, "i"'] \arrow[rd, "fi"] &                  \\
B \arrow[d, "j"'] \arrow[r, "f"]   & X \arrow[d, "p"] \\
C \arrow[r, "g"]                   & Y            
\end{tikzcd},
\label{desired_diagram}
\end{equation}
    with $p$ a forest right fibration, and the maps $i, j, ji$  are normal monomorphisms with $i$ and $ji$ being forest root anodynes. We first consider a diagonal lift $s \colon C \to X$ for the outermost square: this almost defines a lift for our desired square, except that $sj=f$ might not hold.
    
    Using the closure properties in \cref{closure_BV_pushout_product}, one can check that, since $i$ is a root anodyne and $p$ is a right fibration, the canonical map
    \begin{equation}
    \Hom(B,X) \longrightarrow \Hom(A,X) \times_{\Hom(A,Y)} \Hom(B,Y),
    \label{pullback_homs}
    \end{equation}
    is a trivial Kan fibration. The morphisms $f, sj\colon B \to X$ lie on the fiber of \eqref{pullback_homs} over the pair $(fi, gj)$, and therefore the dashed lift $H \colon J \otimes B \to X$ in the diagram below
    \[\begin{tikzcd}
\partial J \otimes B \cup_{\partial J \otimes A} J \otimes A \arrow[d] \arrow[r] & X \arrow[d, "p"] \\
 J \otimes B \arrow[r] \arrow[ru, "H", dashed]                  & Y               
\end{tikzcd}\]
    exists, where the top map is $(f \amalg sj, fi\pi_A)$ and the bottom one is $gj\pi_B$. This follows by adjunction and using that \eqref{pullback_homs} is a trivial Kan fibration. 
    
    Since the simplicial inclusion $0 \to J$ is a right anodyne (and therefore also a forest root anodyne) its pushout-product with $j \colon B \to C$ is a (unitary) root anodyne. Therefore, the diagram below
    $$
\begin{tikzcd}
0 \otimes C \cup_{0 \otimes B} J \otimes B \arrow[rr, "{(s,H)}"] \arrow[d] &                  & X \arrow[d, "p"] \\
J \otimes C \arrow[r, "\pi_C"]                                             & C \arrow[r, "g"] & Y               
\end{tikzcd}
$$
also admits a diagonal lift $K \colon J \otimes C \to X$, since $p$ is a right fibration. One can now check that $K(1,-)\colon C \to X$ gives a lift for the bottom square of \eqref{desired_diagram}, as we wanted to show.
\end{proof}

\begin{lemma}
    Let $n \geq 1$ and $G$ a basic forest. Then the pushout-product
    $$ \partial \Delta^n \otimes G \cup_{\partial \Delta^n \otimes \rho(G)} \Delta^n \otimes \rho(G) \longrightarrow \Delta^n \otimes G$$
    is a forest root anodyne.
    \label{basic_forest_case}
\end{lemma}

\begin{proof}
    Our proof will follow the one for Proposition 4.2.7 in \cite{HeutsHinnichMoerdijk} and we will just indicate the main steps of this proof and leave the details for this reference. Firstly, we claim that it suffices to consider the case when $G$ is a tree, and therefore $\rho(G) \to G$ becomes the inclusion of the root edge of a certain corolla $C_k$. This reduction is spelled out at the end of the proof of \citep[Prop. 3.6.2]{HeutsHinnichMoerdijk}, and essentially uses that root anodynes are closed under puhsouts, together with the fact that $\otimes$ and $\oplus$ commute with colimits on each variable.    

    We analyse the required map by attaching at each stage the shuffles of $\Delta^n \otimes C_k$, which are the representables forests appearing in a certain colimit description of this forest set. For our specific case, these are very easy to describe: for instance, for $\Delta^3 \otimes C_2$ these are given as follows:

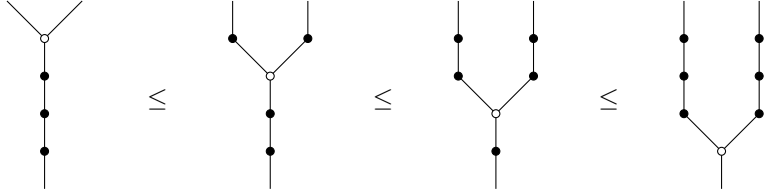
\begin{center}
    \begin{tikzpicture}
    \draw (-4,1.5) -- (-4,-0.5);
    \draw ({-4-sin(30)}, 2) -- (-4,1.5) -- ({-4+sin(30)},2);
    \filldraw[black] (-4,0) circle (1.5pt);
    \filldraw[black] (-4,0.5) circle (1.5pt);
    \filldraw[black] (-4,1) circle (1.5pt);
    \filldraw[color=black, fill=white] (-4,1.5) circle (1.5pt);

    \draw (-1,1) -- (-1,-0.5);
    \draw ({-1-sin(30)}, 1.5) -- (-1,1) -- ({-1+sin(30)},1.5);
    \draw ({-1-sin(30)}, 2) -- ({-1-sin(30)}, 1.5);
    \draw ({-1+sin(30)},2) -- ({-1+sin(30)},1.5);
    \filldraw[black] (-1,0) circle (1.5pt);
    \filldraw[black] (-1,0.5) circle (1.5pt);
    \filldraw[black] ({-1-sin(30)},1.5) circle (1.5pt);
    \filldraw[black] ({-1+sin(30)},1.5) circle (1.5pt);
    \filldraw[color=black, fill=white] (-1,1) circle (1.5pt);

    \draw (2,0.5) -- (2,-0.5);
    \draw ({2-sin(30)}, 1) -- (2,0.5) -- ({2+sin(30)},1);
    \draw ({2-sin(30)}, 2) -- ({2-sin(30)}, 1);
    \draw ({2+sin(30)},2) -- ({2+sin(30)},1);
    \filldraw[black] (2,0) circle (1.5pt);
    \filldraw[black] ({2-sin(30)},1) circle (1.5pt);
    \filldraw[black] ({2+sin(30)},1) circle (1.5pt);
    \filldraw[black] ({2-sin(30)},1.5) circle (1.5pt);
    \filldraw[black] ({2+sin(30)},1.5) circle (1.5pt);
    \filldraw[color=black, fill=white] (2,0.5) circle (1.5pt);

    \draw (5,0) -- (5,-0.5);
    \draw ({5-sin(30)}, 0.5) -- (5,0) -- ({5+sin(30)},0.5);
    \draw ({5-sin(30)}, 2) -- ({5-sin(30)}, 0.5);
    \draw ({5+sin(30)}, 2) -- ({5+sin(30)}, 0.5);
    \filldraw[color=black, fill=white] (5,0) circle (1.5pt);
    \filldraw[black] ({5-sin(30)},0.5) circle (1.5pt);
    \filldraw[black] ({5+sin(30)},0.5) circle (1.5pt);
    \filldraw[black] ({5-sin(30)},1) circle (1.5pt);
    \filldraw[black] ({5+sin(30)},1) circle (1.5pt);
    \filldraw[black] ({5-sin(30)},1.5) circle (1.5pt);
    \filldraw[black] ({5+sin(30)},1.5) circle (1.5pt);
    
   \node [label={[xshift=-2.5cm, yshift=0.3cm]$\leq$}] {};
   \node [label={[xshift=0.5cm, yshift=0.3cm]$\leq$}] {};
   \node [label={[xshift=3.5cm, yshift=0.3cm]$\leq$}] {};
\end{tikzpicture}
\end{center}

    The general pattern of the shuffles is easy to deduce from this one: there are $n+1$ shuffles in $\Delta^n \otimes C_k$, and they are given by "propagating" a copy of $C_k$ along the linear tree $\Delta^n$. We can order the shuffles as it appears in the picture above: the minimal element is the tree with the corolla a leaf vertex, and the maximal element is the one with the corolla as the root vertex. 

    We attach these shuffles one by one using this ordering, which leads to a filtration
    $$ A_0 \subseteq A_1 \subseteq \cdots \subseteq A_n \subseteq A_{n+1} = \Delta^n \otimes C_k,$$
    where $A_0 = \partial \Delta^n \otimes G \cup_{\partial \Delta^n \otimes \eta} \Delta^n \otimes \eta$. In \cite{HeutsHinnichMoerdijk} it is proven that all the stages $A_i \subseteq A_{i+1}$ are inner anodynes as long as $i \neq n$, and their proof also work in our context. The only place where their proof doesn't apply in our scenario is for $A_n \subseteq A_{n+1}$, where we attach the shuffle $S$. In any case, we shall now explain that it is a forest root anodyne, which will then finish the proof.

    In order to gain a better grasp of the leaf and inner faces of $S$, we define first a finite filtration 
    $$A_n = A_n^0 \subseteq A_n^1 \subseteq \cdots \subseteq A_{n+1}$$
    by attaching all \textit{prunings} of $S$. These are the trees obtained from $S$ by successively removing leaf vertices of $S$; these are ordered under containment, with the maximal element being $S$.

    We now further filter each inclusion $A_n^i \subseteq A_n^{i+1}$ given by attaching a pruning $P$ as follows: writing $I(P)$ for the set of inner edges of $P$, we will look at subsets $J \subseteq I(P)$ and at the tree $P^J$ given by contracting all edges \textit{not} in $J$. Extending the containment order of the power set of $I(P)$ in some way, this leads to another finite filtration
    $$ A_n^i = A_n^{i,0} \subseteq A_n^{i,1} \subseteq \cdots \subseteq  A_n^{i+1}$$
    where at the stage corresponding to $J$ we attach $P^J$. We claim that $A_{n}^{i,j} \subseteq A_n^{i,j+1}$ is forest root anodyne.

    Suppose that the inclusion in question corresponds to attaching $P^J$ for some set $J$. Of course, if $P^J$ factors through the previous stage $A_n^{i,j}$, then we have nothing to show. This in particular happens when $J$ doesn't contain any of the incoming edges of the root vertex of $S$: indeed, if that were the case, then to form $P^J$ from $P$ we would need to contract all the incoming edges of the root vertex of $S$. However, the Boardman--Vogt interchange relation would then imply that $P^J$ would have to factor through $A_n$. 

    Supposing then that $J$ contains at lest one incoming edge of the root vertex of $S$, we make the following observations:
    \begin{itemize}
        \item Any leaf face of $P^J$ must have appeared in a smaller pruning, and therefore such a face is contained in $A_n^i$.
        \item Any inner face of $P^J$ must be contained in $A_n^{i,j}$ since we have ordered according to the size of $J$.
        \item The root face of $P^J$ is not in $A_n^{i,j}$: since it contains at least one incoming edge of $S$, it can't be in a previous shuffle. It can't also be contained in a pruning smaller than $P$, and it can't factor through an earlier $P^{J'}$ with $J'$ smaller than $J$, by definition of these forests.
    \end{itemize}

    We conclude from this that $A_{n}^{i,j} \subseteq A_n^{i,j+1}$ is the pushout of the root horn inclusion $\Lambda^r(P^J) \to P^J$, which finishes the proof.
\end{proof}

\begin{proof}[Proof of \cref{root_anodyne_bv_tensor}]
    Consider the class $\mathcal{M}_n$ of forest normal monomorphisms $j \colon A \to B$ such that the pushout-product
    \begin{equation}
    \partial \Delta^n \otimes B \cup_{\partial \Delta^n \otimes A} \Delta^n \otimes A \longrightarrow \Delta^n \otimes B
    \label{eqpp}
    \end{equation}
    is a forest root anodyne. It's a routine check that $\mathcal{M}_n$ is saturated, and it also satisfies the right cancellation property by \cref{forest_right_cancellation}.
    
    \cref{basic_forest_case} ensures that $\mathcal{M}_n$ contains all root inclusions $\rho(G) \to G$ whenever $G$ is a basic forest. Moreover, it also contains all forest spine inclusions $\mathsf{fSp}(F) \to F$: these are examples of inner anodynes by \citep[Prop. 3.6.6]{HeutsHinnichMoerdijk}, and \cref{closure_BV_pushout_product} (b) yields that its pushout-product with $\partial \Delta^n \to \Delta^n$ is again inner anodyne, and thus in $\mathcal{M}_n$. \cref{closure_forest_root anodynes} guarantees that any saturated class satisfying the right cancellation property which contains these two types of morphisms must in particular contain all the forest root anodynes, which is exactly what we wanted to show.
\end{proof}

\subsection{The contravariant model structure on $\fSets_{/V}$}

\leavevmode

Turning now to the homotopy theory of $\fSets$, Heuts--Hinich--Moerdijk establish in \cite{HeutsHinnichMoerdijk} the existence of a model structure on $\fSets$, which is closely related to the Cisinski--Moerdijk model structure $\dSets_{\mathsf{opd}}$ on dendroidal sets constructed in \cite{CisinskiMoerdijk}, as well as the Joyal model structure $\sSets_{\mathsf{cat}}$ on simplicial sets, in the sense that there exist Quillen pairs
\[\begin{tikzcd}
	{\sSets_{\mathsf{cat}}} & {\dSets_{\mathsf{opd}}} & {\fSets_{\mathsf{opd}}}
	\arrow["{\iota_!}", shift left, from=1-1, to=1-2]
	\arrow["{\iota^\ast}", shift left, from=1-2, to=1-1]
	\arrow["{u^*}", shift left, from=1-2, to=1-3]
	\arrow["{u_*}", shift left, from=1-3, to=1-2]
\end{tikzcd}\]
relating all of these model categories. In fact, the operadic model structure on $\fSets$ that is defined in \cite{HeutsHinnichMoerdijk} makes $(u_\ast, u^\ast)$ into a Quillen equivalence; the definition we will work with here is a modified version of their model structure which doesn't make this pair into a Quillen equivalence anymore, see \cref{not_quillen_equivalence} for more on this.

We begin by first explaining what are the different classes of maps which define the model category $\fSets_{\mathsf{opd}}$. Recall the definition of the simplicial set $\Hom(X,Y)$ for any forest sets $X$ and $Y$ presented before \cref{forest_right_cancellation}.

\begin{definition}
    We say $X \in \fSets$ is \textit{operadically local} if the following conditions hold:
    \begin{enumerate}[label=(\alph*)]
        \item For every normal forest morphism $A \to B$, the pullback map
        $$ \Hom(B, X) \longrightarrow \Hom(A,X)$$
        is a categorical fibration of simplicial sets.

        \item For every forest $F$ and any choice of inner edge $e \in E(F)$, the map
        $$ \Hom(F,X) \longrightarrow \Hom(\Lambda^e F, X)$$
        is a categorical trivial fibration of simplicial sets.
    \end{enumerate}

    In particular, the first condition implies that $\Hom(A,X)$ is an $\infty$-category, whenever $A$ is a normal forest set and $X$ is operadically local.

   Before the next definition, let us say that a  \textit{normalization} of a forest map $f \colon X \to Y$ is a forest map $f'\colon X' \to Y'$ which fits into a commutative diagram
    $$
\begin{tikzcd}
X' \arrow[r, "f'"] \arrow[d] & Y' \arrow[d] \\
X \arrow[r, "f"]             & Y           
\end{tikzcd}
    $$
    such that $X'$ and $Y'$ are normal forest sets, and the vertical maps have the right lifting property with respect to all normal monomorphisms.
\end{definition}

\begin{definition}
    A forest map $X \to Y$ is an \textit{operadic weak equivalence} if it admits a normalization $X'\to Y'$ such that, for every $E$ operadically local forest set, the map
    \begin{equation}
    \Hom(Y', E) \longrightarrow \Hom(X', E)
    \label{operadic_w_e}
    \end{equation}
    is a categorical equivalence between $\infty$-categories.
\end{definition}

By Lemma 3.7.11 in \cite{HeutsHinnichMoerdijk}, the choice of normalization is irrelevant: that is, if \eqref{operadic_w_e} is an equivalence for \textit{some} choice of normalization of $X \to Y$, then it is an equivalence for \textit{all} choices of normalization.

\begin{definition}
    We say a forest map $A \to B$ is an \textit{operadically anodyne map} if it belongs to the saturated class generated by the maps
    $$\Lambda^e F \to F \hspace{2em} \mathrm{and} \hspace{2em} 0 \otimes F \cup J \otimes \partial F \to J \otimes F$$
    for any forest $F$ and any inner edge $e$. 
    \label{def_operadic_anodynes}
\end{definition}

The next result is a slightly altered version of Theorem 3.7.9 and its subsequent corollary in \cite{HeutsHinnichMoerdijk}.
\begin{theorem}
    There exists a cofibrantly generated and left proper model structure on $\fSets$ satisfying the following properties:
    \begin{enumerate}[label=(\alph*)]
        \item The cofibrations are the normal monomorphisms of forest sets.
        \item The weak equivalences are the operadic weak equivalences. In particular, the weak equivalences between normal forest sets are the maps $A \to B$ such that, for any operadically local forest set $X$, the map
        $$ \Hom(B,X) \longrightarrow \Hom(A,X)$$
        is a categorical equivalence of $\infty$-categories.
        \item The fibrant objects are the operadically local forest sets.
        \item The fibrations between fibrant objects are exactly the maps which have the right lifting property with respect to the operadically anodyne maps.
    \end{enumerate}
    
    We write $\fSets_{\mathsf{opd}}$ for the corresponding model category. Moreover, the adjoint pairs
    \[\begin{tikzcd}
	{\fSets_{\mathsf{opd}}} & {\dSets_{\mathsf{opd}}}
	\arrow["{u^*}", shift left, from=1-1, to=1-2]
	\arrow["{u_*}", shift left, from=1-2, to=1-1]
\end{tikzcd} \hspace{1.5em} and \hspace{1.5em} \begin{tikzcd}
	{\sSets_{\mathsf{cat}}} & {\fSets_{\mathsf{opd}}}
	\arrow["{(u\iota)_!}", shift left, from=1-1, to=1-2]
	\arrow["{(u \iota)^*}", shift left, from=1-2, to=1-1]
\end{tikzcd}\]
    define Quillen adjunctions.
    \label{model_structure_operadic_forest}
\end{theorem}

\begin{proof}
    As we mentioned above, this is essentially what is proved in \cite{HeutsHinnichMoerdijk}, except that in their work there is a further localization making the derived counit into a weak equivalence when defining the operadically local forest sets. 

   One can check that the proof in \cite{HeutsHinnichMoerdijk} for the existence of the operadic model structure goes through when one makes this small change. Alternatively, the machinery of $\mathcal{A}$-model structures which we will explain next also gives a way of constructing this model structure, see \cref{example_operadic} for an elaboration on this.

    The only statement which is not proved is the one concerning the Quillen adjunction with $\sSets_{\mathsf{cat}}$. Using the colimit formula for the left Kan extension and using the fact that linear trees are the only forests admitting maps to $\eta$, one easily checks that the value of the left adjoint $(u \iota)_!$ on a simplicial set $X$ is
    $$(u \iota)_! (X)_F = 
    \begin{cases}
    X_k &\mathrm{if} \: F = Delta^k \: \text{for some} \: k,\\
    \emptyset  &\mathrm{otherwise.} 
    \end{cases}
    $$

    It follows from this explicit description that $(u \iota)_!$ preserves cofibrations. If $f \colon X \to Y$ is a categorical trivial cofibration of simplicial sets, then for any $E$ operadically local forest set, we know that the map
    $$ \tau\left( \Hom(Y, (u \iota)^* E) \right) \longrightarrow  \tau\left( \Hom(X, (u \iota)^* E) \right)$$
    is an equivalence of categories, since $(u \iota)^* E$ is an $\infty$-category and $f$ is a categorical equivalence. Here $\tau \colon \sSets \to \mathsf{Cat}$ is the left adjoint to the usual nerve functor for categories. By adjunction, this condition implies that 
     $$ \tau\left( \Hom((u \iota)_! Y, E) \right) \longrightarrow  \tau\left( \Hom((u \iota)_! X, E) \right)$$
     is also an equivalence of categories, where $\tau$ is the left adjoint to the usual categorical nerve. Since $(u \iota)_! f$ is already a normalization of itself, this shows that 
     $$ \Hom((u \iota)_! Y, E) \longrightarrow \Hom((u \iota)_! X, E) $$
     is a categorical equivalence, and therefore $(u \iota)_! f$ is a trivial cofibration in $\fSets_{\mathsf{opd}}$, as we wanted to show.
\end{proof}

\begin{remark}
    As we mentioned at the start of the proof, we have opted out from localizing with respect to the counit map. Informally speaking, this condition makes forest sets into objects that behave in much the same as dendroidal sets, since the former are exactly the presheaves that treat $\oplus$ as a coproduct, see \citep[Rem. 3.2.3]{HeutsHinnichMoerdijk} for a more sheaf-theoretic version of this statement. 
    
    This further localization is not ideal for our intent of studying right modules, since we will see further on that the forest set associated with an operadic right module is not often in the image of $u_\ast$. For instance, the modules coming from embedding calculus give rise to true forest sets, essentially because the map 
    $$ \mathsf{Emb} \left( \underline{k} \times  D^d, M^d \right) \longrightarrow \prod_{i=1}^k \mathsf{Emb}(D^d, M^d)$$
    that keeps track of each individual embedded disk is far from being a homeomorphism (not even a weak homotopy equivalence).
    \label{not_quillen_equivalence}
\end{remark}

\begin{remark}
It is natural to ask whether the adjoint pair $(u_!, u^*)$ also defines a Quillen adjunction, but this is not true: the main problem lies in the fact $u_!$ doesn't preserve normal monomorphisms, and the previous argument fails since such an explicit description for $u_!$ fails to exist. For an explicit counterexample to $u_!$ being a left Quillen functor see \citep[Rem. 3.3.4]{HeutsHinnichMoerdijk}.
    \label{operadic_model_remark}
\end{remark}

In order to prepare ourselves for the construction of the contravariant model structure, we remark that there exists an alternative method for constructing the operadic model structure $\fSets_{\mathsf{opd}}$, using the \textit{$\mathcal{A}$-model structures} developed in Chapter 9 of \cite{HeutsMoerdijkDendroidal}. This is essentially a general framework that, given a certain set of forest monomorphisms $\mathcal{A}$, produces a model structure on $\fSets$ where one localizes with respect to $\mathcal{A}$. In fact, this different strategy will be the one we will make use of for defining the contravariant model structure, and therefore we shall briefly go over it now.

\begin{notation}
Throughout this discussion, we fix a forest set $V$ and a set $\mathcal{A}$ of forest normal monomorphisms in $\fSets_{/V}$ between normal forest sets.   
\label{notation_A_model}
\end{notation}

\begin{warning}
    If $V$ is a dendroidal set (or afterwards a dendroidal space), we will simplify the notation and always write $\fSets_{/V}$ for the slice over the forest set $u_\ast V$.
\end{warning}

We will begin by defining the main classes for this model structure on the slice category $\fSets_{/V}$. Given any two objects $p \colon X \to V$ and $q \colon Y \to V$, we define the simplicial set $\Hom_V(X,Y)$ as the fiber of 
$$ \Hom(X,Y) \xlongrightarrow{ q \circ -} \Hom(X,V)$$
over $p \in \Hom(X,V)_0$.

\begin{definition}
    Let $V$ and $\mathcal{A}$ be as in \cref{notation_A_model}.
    \begin{itemize}
        \item We say that an object $X \to V$ is \textit{$\mathcal{A}$-local} if, for any forest normal monomorphism $A \to B$ over $V$, the induced map
        $$ \Hom_V(B, X) \longrightarrow \Hom_V(A,X)$$
        is a categorical fibration in $\sSets$, which is moreover trivial if $A \to B$ belongs to $\mathcal{A}$. In particular, taking $A = \emptyset$, one sees that $\Hom_V(B,X)$ is an $\infty$-category whenever $B$ is normal and $X$ is $\mathcal{A}$-local.
        \item We say that a forest map $f \colon A \to B$ over $V$ is an \textit{$\mathcal{A}$-weak equivalence} if, given any normalization $A' \to B'$ over $V$ of $f$, the map
        $$ \Hom_V(B', X) \longrightarrow \Hom_V(A',X)$$
        is a categorical equivalence between $\infty$-categories, whenever $X \to V$ is $\mathcal{A}$-local.
    \end{itemize}
\end{definition}

\begin{definition}
    Let $V$ and $\mathcal{A}$ be as in \cref{notation_A_model}. We define the class of \textit{$\mathcal{A}$-anodynes} to be the saturated class generated by the following types of morphisms:
    \begin{enumerate}[label=(\roman*)]
    \item For $A \to B$ a forest normal monomorphism over $V$ and $\Lambda^i [n] \to \Delta^n$ any inner horn inclusion, the pushout-product
    $$ \Lambda^i [n] \otimes B \cup \Delta^n \otimes A \longrightarrow \Delta^n \otimes B.$$
    \item For $A \to B$ a forest normal monomorphism over $V$, the pushout-product
    $$ 0 \otimes B \cup J \otimes A \longrightarrow J \otimes B.$$
    \item For $A \to B$ a forest normal monomorphism over $V$ in $\mathcal{A}$ and $n \geq 0$, the pushout-product
    $$ \partial \Delta^n \otimes B \cup \Delta^n \otimes A \longrightarrow \Delta^n \otimes B.$$
    \end{enumerate}
    \label{definition_A_anodyne}
\end{definition}

Before stating the main result, we will need to be a bit more restrictive on the kind of sets $\mathcal{A}$ we will consider. With that in mind, we will say that $\mathcal{A}$ is \textit{admissible} if, for any $i \colon A \to B$ in $\mathcal{A}$ and $n \geq 0$, the pushout-product
$$ \partial \Delta^n \otimes B \cup_{\partial \Delta^n \otimes A} \Delta^n \otimes A \longrightarrow \Delta^n \otimes B$$
is in $\mathcal{A}$.

\begin{theorem}
    Let $\mathcal{A}$ and $V$ be as in \cref{notation_A_model}. If $\mathcal{A}$ is admissible, then there exists a left proper, cofibrantly generated model structure on the slice category $\fSets_{/V}$ such that:
        \begin{enumerate}[label=(\alph*)]
        \item The cofibrations are the normal monomorphisms over $V$.
        \item The weak equivalences are exactly the $\mathcal{A}$-weak equivalences.
        \item The fibrant objects are the $\mathcal{A}$-local objects.
        \item The fibrations between fibrant objects are exactly the morphisms having the right lifting property with respect to all $\mathcal{A}$-anodynes.
    \end{enumerate}

    We will call this model structure the $\mathcal{A}$-model structure on $\fSets_{/V}$.
\end{theorem}

\begin{proof}
    The proof of the dendroidal version of this result is contained in the proofs of Theorems 9.99 and 9.55 in \cite{HeutsMoerdijkDendroidal}. One can check that the forest variants of the more technical lemmas used in proving these results are also true, and most of these are actually proven in \cite{HeutsHinnichMoerdijk}. Therefore, the arguments present in these sources can still be used to show the existence of this model structure on $\fSets_{/V}.$
\end{proof}

\begin{example}
    Take $V$ to be the terminal forest set so that $\fSets_{/V} \cong \fSets$, and consider the set $\mathcal{A}$ of forest inner anodynes $A \to B$ with $A$ and $B$ finite normal forest sets. Then it follows from \cref{closure_BV_pushout_product} (b) that $\mathcal{A}$ is admissible, and the associated $\mathcal{A}$-model structure is $\fSets_{\mathsf{opd}}$: indeed, the cofibrations are the normal monomorphisms for both model categories, and, by comparing definitions, we see that the operadically local objects are exactly the $\mathcal{A}$-local objects.
    \label{example_operadic}
\end{example}

\begin{theorem}
    Let $V$ be a forest set. The slice category $\fSets_{/V}$ admits a left proper, cofibrantly generated model structure such that the following conditions hold:
    \begin{enumerate}[label=(\alph*)]
        \item The cofibrations are the normal monomorphisms over $V$.
        \item The weak equivalences between normal forest sets are the maps $A \to B$ over $V$ such that, for any right fibration $X \to V$, the map
        $$ \Hom_V(B,X) \longrightarrow \Hom_V(A,X)$$
        is a categorical equivalence of $\infty$-categories.
        \item The fibrant objects are the forest right fibrations over $V$.
        \item The fibrations between fibrant objects are the forest right fibrations.
    \end{enumerate}

    We will call this model structure the contravariant model structure on $\fSets_{/V}$, and write $\left( \fSets_{/V} \right)_{\ctv}$  for the associated model category.
    \label{ctv_character}
\end{theorem}

\begin{proof}
Letting $\mathcal{R}$ be the saturated class generated by the horn inclusions $\Lambda^a(F) \to F$ over $V$, where $a$ is either an inner edge or a root of $F$, we define $\mathcal{A}$ to be the set of maps in $\mathcal{R}$ which have both source and target finite forest sets. Then $\mathcal{A}$ is an admissible set due to \cref{root_anodyne_bv_tensor}, so we can talk about the $\mathcal{A}$-model structure on $\fSets_{/V}$.

All the properties we have stated are clear from the characterization of the $\mathcal{A}$-model structure except for (c) and (d). We will check (c) since the verification of (d) is exactly the same.

An object being fibrant in the $\mathcal{A}$-model structure is the same as saying that it has the right lifting property with respect to the maps in \cref{definition_A_anodyne}. The maps of type (i) and (ii) are always root anodyne: for instance, this follows from \cref{closure_BV_pushout_product} (d) in both cases (for (i) one can be a bit more precise and see that it is actually inner anodyne, but we won't need this). As for the maps of type (iii), these are forest root anodynes by \cref{root_anodyne_bv_tensor}. It is clear now from this discussion that the right fibrations over $V$ are examples of $\mathcal{A}$-local objects. The reverse also holds: it suffices to consider the maps in (iii) with $n=0$ and $A \to B$ any inner/root horn inclusion over $V$.
\end{proof}

The following corollary yields a class of forest maps which are examples of contravariant trivial cofibrations, namely the forest right anodynes.

\begin{corollary}
    Let $f \colon X \to Y$ be a forest map over a forest set $V$. If $f$ is a forest root anodyne, then $f$ is a trivial cofibration in $\left( \fSets_{/V}\right)_{\ctv}$.
    \label{root_anodyne_ctv_we}
\end{corollary}

\begin{proof}
    By \citep[Lem. 8.42]{HeutsMoerdijkDendroidal}, it suffices to check that $f$ has the left lifting property with respect to fibrations between fibrant objects. By \cref{ctv_character}, these are exactly the right fibrations.
\end{proof}

Although we won't make use of it in the rest of the text, one can prove a certain converse to \cref{root_anodyne_ctv_we}, in the sense that any contravariant trivial cofibration with contravariantly fibrant target is also a forest root anodyne.

\begin{lemma}
    Let $V$ be a forest set and $f \colon X \to Y$ a morphism in $\fSets_{/V}$. If the reference map $Y \to V$ is a right fibration and $f$ is a contravariant trivial cofibration, then $f$ is a forest root anodyne.
    \label{root_anodynes_help_lemma}
\end{lemma}

\begin{proof}
    The proof is a typical instance of the retract argument in model categories: we begin by factoring $f$ as $X \xrightarrow{i} W \xrightarrow{p} Y$, where $i$ is a forest root anodyne and $p$ is a forest right fibration. Consider the diagram
    \begin{equation}
    \begin{tikzcd}
X \arrow[d, "f"'] \arrow[r, "i"] & W \arrow[d, "p"] \\
Y \arrow[r, equal] & Y.              
\end{tikzcd}
    \label{lift1}    
    \end{equation}

    The composite $W \xrightarrow{p} Y \to V$ exhibits $W$ as a forest set over $V$, and in fact it is also a right fibration by our initial hypothesis. Thus, it follows from \cref{ctv_character} that $p$ is also a contravariant fibration over $V$, and therefore a lift $s\colon Y \to W$ in \eqref{lift1} exists. This witnesses $f$ as a retract of the root anodyne $i$, which concludes the proof.
\end{proof}

The next result identifies the contravariant model structure as a certain left Bousfield localization of the relative operadic model category. Its proof mostly rests on the following lemma.
\begin{lemma}
    Let $p \colon X \to V$ be a map of forest sets. Then $p$ is a forest right fibration if and only if $p$ is an operadic fibration and the map
    \begin{equation}
     \mathsf{Map}_V( F, X) \longrightarrow \mathsf{Map}_V(\rho(F),X)
     \label{localitiy}
    \end{equation}
    is a weak homotopy weak equivalence of simplicial sets, for every root inclusion $\rho(F) \to F$ over $V$. Here $\Map_V(X,Y)$ denotes the (derived) mapping space with respect to the operadic model structure on $\fSets_{/V}$.
    \label{character_fibrant_bosufield}
\end{lemma}

\begin{proof}
Suppose first that $p \colon X \to V$ is an operadic fibration satisfying the mapping space condition, and consider the class $\mathcal{M}$ of normal monomorphisms $A \to B$ over $V$ such that the map obtained after applying $\mathsf{Map}_V(-,X)$ is a weak homotopy equivalence. This is a saturated class which satisfies the right cancellation property and by hypothesis contains all root inclusions; therefore, the relative version
of \cref{closure_forest_root anodynes} shows that $\mathcal{M}$ contains all root horn inclusions over $V$, that is,
$$ \mathsf{Map}_V(F, X) \xlongrightarrow{\simeq} \mathsf{Map}_V(\Lambda^r F,X)$$
is a weak homotopy equivalence for any root horn inclusion over $V$. As $X \to V$ is an operadic fibration, this map is modelled by a trivial Kan fibration
$$ \left( \fSets_{/V} \right)\left( \mathcal{X}^\bullet \otimes F, X \right) \longrightarrow  \left( \fSets_{/V} \right)\left( \mathcal{X}^\bullet \otimes \Lambda^r F, X \right),$$
 where $ \mathcal{X}^\bullet$ is a cosimplicial resolution of the unit with $\mathcal{X}^0 = \eta$, which in particular implies the desired lifting property of $p$ with respect to $\Lambda^r F \to F$.

On the other hand, suppose $X \to V$ is a right fibration; in particular this means that $p$ is an operadic fibration, so it suffices to shows that it satisfies the condition on the mapping spaces. One can find a cosimplicial resolution $\mathcal{X}^\bullet$ of the terminal forest set for the operadic model structure, where each stage $\mathcal{X}^k$ is a simplicial set: just choose any cosimplicial resolution in $\sSets_{\mathsf{opd}}$ first and then apply $(u \iota)_!$, using that this functor is left Quillen by \cref{model_structure_operadic_forest}. Consequently, we can model \eqref{localitiy} by the map
\begin{equation}
 \fSets_{/V}( \mathcal{X}^{\bullet} \otimes F, X) \longrightarrow \fSets_{/V}( \mathcal{X}^\bullet \otimes \rho(F), X).
 \label{locality_2}
\end{equation}

Associated with this, we consider the following (equivalent) lifting problems
\[\begin{tikzcd}[cramped]
	{\partial \Delta^n} & {\fSets_V(\mathcal{X}^\bullet \otimes F, X)} & {\partial \mathcal{X}^n \otimes F \cup \mathcal{X}^n \otimes \rho(F)} & X \\
	{\Delta^n} & {\fSets_V(\mathcal{X}^\bullet \otimes \rho(F), X)} & {\mathcal{X}^n \otimes F} & V
	\arrow[from=1-1, to=1-2]
	\arrow[from=1-1, to=2-1]
	\arrow[from=1-2, to=2-2]
	\arrow[from=1-3, to=1-4]
	\arrow[from=1-3, to=2-3]
	\arrow[from=1-4, to=2-4]
	\arrow[dashed, from=2-1, to=1-2]
	\arrow[from=2-1, to=2-2]
	\arrow[dashed, from=2-3, to=1-4]
	\arrow[from=2-3, to=2-4]
\end{tikzcd}\]
where $\partial \mathcal{X}^n$ is the colimit of $\mathcal{X}^k$ along all proper monomorphisms $[k] \to [n]$. The left vertical map in the right diagram is a forest root anodyne by \cref{closure_forest_root anodynes}, and therefore the dashed arrow exists. By looking at the equivalent left square, this shows that \eqref{locality_2} is a trivial Kan fibration, and consequently \eqref{localitiy} is a weak homotopy equivalence.    
\end{proof}

\begin{proposition}
    Let $V$ be a forest set. Then the contravariant model structure on $\fSets_{/V}$ is a left Bousfield localization of the relative operadic model structure at the root inclusions $\rho(F) \to F$ over $V$.
    \label{ctv_as_bousfield}
\end{proposition}
\begin{proof}
Clearly both model structures share the normal monomorphisms over $V$ as the class of cofibrations. \cref{character_fibrant_bosufield} is exactly the statement that the model categories also share the same fibrant objects, which shows that they must have the same underlying model structure.
\end{proof}

For our final results on this section, we want to study how exactly the contravariant model structure depends on the object we are slicing over, in particular up to operadic weak equivalence of forest sets. 

\begin{lemma}
    Let $\varphi \colon V \to W$ be a forest map. Then the adjunction
    \[\begin{tikzcd}
	{\left( \fSets_{/V} \right)_\ctv} & {\left( \fSets_{/W}\right)_{\ctv}}
	\arrow["{\varphi_!}", shift left, from=1-1, to=1-2]
	\arrow["{\varphi^*}", shift left, from=1-2, to=1-1]
\end{tikzcd}\]
defines a Quillen pair.
\label{lemma_quillen_pair}
\end{lemma}

\begin{proof}
    By Lemma \citep[Lem. 8.42]{HeutsMoerdijkDendroidal}, it suffices to check that $\varphi_!$ preserves cofibrations, and $\varphi^*$ sends fibrations between fibrant objects to fibrations.

    For the first one, this is immediate since the left adjoint is just composing the reference maps with $\varphi$. As for the right adjoint $\varphi^*$, this functor is given by pulling back along $\varphi$; by the characterization in \cref{ctv_character} of fibrant objects and fibrations between such objects, we conclude that $\varphi^*$ preserves both of these classes, since right fibrations are stable under pullback.
\end{proof}

\begin{proposition}
    Let $\varphi \colon V \to W$ be a forest map.
    \begin{enumerate}[label=(\alph*)]
        \item If $\varphi$ is a forest right fibration, then the contravariant model structure $\left( \fSets_{/ V} \right)_\ctv$ coincides with the relative model structure on $\left( \fSets_{/ W} \right)_{/ \varphi}$ given by equipping $\fSets_{/W}$ with the contravariant model structure.
        \item If $\varphi$ is an operadic weak equivalence between operadically local objects, then the Quillen adjunction
        \[\begin{tikzcd}
	{\left( \fSets_{/V} \right)_\ctv} & {\left( \fSets_{/W}\right)_{\ctv}}
	\arrow["{\varphi_!}", shift left, from=1-1, to=1-2]
	\arrow["{\varphi^*}", shift left, from=1-2, to=1-1]
\end{tikzcd}\]
is a Quillen equivalence.
    \end{enumerate}
    \label{contravariance_operad_invariance}
\end{proposition}

\begin{proof}
    For (a), it suffices to check that the fibrant objects on both categories coincide, since it is clear that the class of cofibrations is the class of normal monomorphisms in both situations. Unwinding the definitions, we see that a fibrant object in $\left( \fSets_{/ W} \right)_{/ \varphi}$ is the data of forest maps $p \colon X \to W$ and $f \colon X \to V$ related via the equation $p = \varphi f$. This in particular exhibits $f$ as a morphism in $\fSets_{/W}$, and the fibrancy condition translates into $f$ being a contravariant fibration in $\left( \fSets_{/W} \right)_{\ctv}$. 
    
    By Theorem \ref{ctv_character}, $\varphi$ is a contravariantly fibrant object in $\fSets_{/W}$, which additionally implies that the same also holds for $p$, since $f$ is a contravariant fibration in $\fSets_{/W}$. Consequently, $f$ defines a fibration between fibrant objects in $\left( \fSets_{/W} \right)_{\ctv}$, and another application of \cref{ctv_character} it follows that $f$ is a right fibration of forest sets; since these are the fibrant objects in $(\fSets_{/V})_{\ctv}$, we are done.

    For (b), the fact that $(\varphi_!, \varphi^*)$ is a Quillen adjunction is the content of Lemma \ref{lemma_quillen_pair}, so we only have to check that under the additional hypothesis it is also a Quillen equivalence.

    We first claim that we can assume that $\varphi$ is a trivial fibration in $\left( \fSets_{/W} \right)_{\ctv}$. Indeed, we note first that $\varphi$ being a contravariantly trivial fibration is equivalent to $\varphi$ being an operadic trivial fibration over $V$, as both of the model categories in question share exactly the same class of cofibrations; therefore it suffices to show that we can assume that $\varphi$ is an operadic trivial fibration. Since $\varphi$ is by hypothesis an operadic weak equivalence between operadically local objects, we can apply Brown's lemma to find a factorization of $\varphi$ as
    $$ X \xlongrightarrow{i} Z \xlongrightarrow{t} Y$$
    with $t$ an operadic trivial fibration, together with another map $p \colon Z \to X$ satisfying $pi=\mathrm{id}_X$. From this last equation we see that if we first prove (b) for $t$, then this also implies the same statement but for $\varphi$. As $t$ is an operadic trivial fibration, this shows the desired reduction.
    
    Under this assumption, we can use Example 8.47 (i) of \cite{HeutsHinnichMoerdijk}, to conclude that $\varphi$ defines a Quillen equivalence
    \[\begin{tikzcd}
	{\left( \fSets_{/W} \right)_{/ \varphi}} & {\left( \fSets_{/W}\right)_{/ \mathrm{id}_W},}
	\arrow["{\varphi_!}", shift left, from=1-1, to=1-2]
	\arrow["{\varphi^*}", shift left, from=1-2, to=1-1]
\end{tikzcd}\]
where both categories have the model structure explained in the statement of part (a). Noticing now that the right hand side is isomorphic to $(\fSets_{/W})_{\ctv}$ and that part (a) identifies the left hand side with $\left( \fSets_{/V} \right)_{\ctv}$ finishes the proof.
\end{proof}

An easy corollary of this result is that the contravariant model structure over an $\infty$-operad $V$ is unique under operadic weak equivalences in $\dSets$, up to Quillen equivalence. 

\begin{corollary}
    Suppose $\varphi \colon V \to W$ is a weak equivalence between $\infty$-operads in $\dSets_{\mathsf{opd}}$. Then there is a Quillen equivalence 
    \[\begin{tikzcd}
	{\left( \fSets_{/V} \right)_\ctv} & {\left( \fSets_{/W}\right)_{\ctv}.}
	\arrow["{(u_\ast \varphi)_!}", shift left, from=1-1, to=1-2]
	\arrow["{(u_\ast \varphi)^*}", shift left, from=1-2, to=1-1]
\end{tikzcd}\]
\end{corollary}

\begin{proof}
    This is immediate from Proposition \ref{contravariance_operad_invariance}, the fact that $u_\ast$ is a right Quillen functor for the operadic model structure and that $\infty$-operads are exactly the fibrant objects in $\dSets_{\mathsf{opd}}$.
\end{proof}

\subsection{An application to right modules over an operad}

\leavevmode

In this section we will prove a technical result relating the contravariant model structure on $\fSets_{/V}$ and its relation to operadic right modules, which will become more relevant later on when we consider forest \textit{spaces}. This will involve the construction of a functor analogous to the nerve construction but for the category of right $\mathcal{P}$-modules, which is of central importance for this work and especially in the next section.

\begin{notation}
    The following notation will be useful from now on: given a forest $F$, recall that we write $\rho(F)$ for the forest of roots of $F$. Given an operad $\mathcal{P}$, a right $\mathcal{P}$-module $M$ and $ k \geq 1$, we define
    $$ M(k) = \coprod_{(c_1, \ldots, c_k)} M(c_1, \ldots, c_k),$$
    where the coproduct is indexed over all distinct $k$-tuples of colours of $\mathcal{P}$. In particular, this coincides with the already-existing notation for uncoloured operads.

    We extend this notation to forests of edges by setting $M(k \cdot \eta) = M(k)$, which consequently defines $M( \rho(F))$ for any $F \in \For$. Finally, we note that there is a map natural in $F$
    $$ M(\rho(F)) \longrightarrow N\mathcal{P}_{\rho(F)}$$
    given as follows: any $m \in M(\rho(F))$ (that is, an element $m \in M(c_1, \ldots, c_k)$ in the module for a certain tuple $(c_1, \ldots, c_k)$ of colours) is sent to the respective tuple of colours $(c_1, \ldots, c_k)$. With this notation, we note, for any forest $F$, there is a morphism
    $$ \cdot_F : N\mathcal{P}_F \times_{N\mathcal{P}_{\rho(F)}} M(\rho(F)) \longrightarrow M(\lambda(F)),$$
    where $\lambda(F) \in \For$ is the forest of all leaves of $F$: an element of $N\mathcal{P}(F)$ can be thought of as the forest $F$ decorated with operations at each vertex and the edges labelled with the input and oupt edges accordinlgy. Then, we act with these on $M(\rho(F))$, starting from the root vertices until we end up at the leaves of $F$.
    \label{notation_forests_modules}
\end{notation}

For the next definition, recall that a coloured operad $\mathcal{P}$ is \textit{closed} if, for any colour $c$, the space of constants $\mathcal{P}(-;c)$ is a singleton. 

\begin{proposition}
    Let $\mathcal{P}$ be a closed operad. Then there is a functor
    $$ N_{\mathcal{P}} : \RMod{\mathcal{P}} \longrightarrow \fSets_{/ N\mathcal{P}}$$
    which, for $M$ a right $\mathcal{P}$-module, is given via the pullback
\[\begin{tikzcd}[cramped]
	{\left( N_{\mathcal{P}}M \right)_F} & {M(\rho(F))} \\
	{N\mathcal{P}_F} & {N\mathcal{P}_{\rho(F)}}
	\arrow[from=1-1, to=1-2]
	\arrow[from=1-1, to=2-1]
	\arrow["\lrcorner"{anchor=center, pos=0, scale=1.5}, draw=none, from=1-1, to=2-2]
	\arrow[from=1-2, to=2-2]
	\arrow[from=2-1, to=2-2]
\end{tikzcd}\]
    when evaluated at a forest $F$. In this diagram, the right map is the one mentioned in \cref{notation_forests_modules}, and the bottom one is induced by the root inclusion $\rho(F) \to F$.
\end{proposition}

\begin{proof}
    We will explain how the generating morphisms mentioned in \cref{arrows_For} act on $N_{\mathcal{P}} M$. We will also omit the functor $u_\ast$ from the notation from now on.

    Suppose $T \to S$  is a forest map between trees which is either an isomorphism, an elementary degeneracy or face map of trees which is \textit{not} the root face, then the action is
    $$
\begin{tikzcd}
M(\rho(S)) \arrow[d, equal]\arrow[r] & N\mathcal{P}_{\rho(S)}   \arrow[d, equal] & N\mathcal{P}_{S} \arrow[d, "N\mathcal{P}^*(T \to S)"] \arrow[l] \\
M(\rho(T))                 \arrow[r]               & N\mathcal{P}_{\rho(T)}                               & N\mathcal{P}_{T}         \arrow[l]                           
\end{tikzcd}
    $$
    where, since $\rho(S) = \rho(T) = \eta$, it makes sense to set the left and middle maps to be the respective identities.

     For a summand inclusion $F \to F \oplus G$, we note first that there are natural maps
     \begin{equation}
     N\mathcal{P}_{G} \times_{N\mathcal{P}_{\rho(G)}} M(\rho(F) \oplus \rho(G)) \longrightarrow M(\rho(F) \oplus \lambda(G)) \longrightarrow M(\rho(F)).
    \label{rewriting_action_Gbar}
    \end{equation}
     The first of these maps is $\cdot_G$ as we have defined\footnote{Here there is an extra $\rho(F)$ term, but this is irrelevant.} before the proposition; the second one comes from the observation that there is always a map $f^ * \colon M(\ell \cdot \eta) \to M(k \cdot \eta)$ for every injective function $f \colon \underline{k} \to \underline{\ell}$. Indeed, the domain is the copropduct of the sets $M(c_1, \ldots, c_\ell)$, and we always have a map
     $$ M(c_1, \ldots, c_\ell) \longrightarrow M(c_{f(1)}, \ldots, c_{f(\ell)})$$
     by acting with the constants of $\mathcal{P}$ with output colour $c_i$ with $i \not\in \mathrm{im}(f)$. With this in mind, $N_\mathcal{P} M$ evaluated on $F \oplus G$ is given by
    \begin{equation}
    \left( N_{\mathcal{P}}M\right)_{F \oplus G} = N\mathcal{P}_F \times_{N\mathcal{P}_{\rho(F)}} \times \left[ N\mathcal{P}_{G} \times_{N\mathcal{P}_{\rho(G)}} M(\rho(F) \oplus \rho(G)) \right].
    \label{rewriting_NP_expression}
    \end{equation}
     Acting with \eqref{rewriting_action_Gbar} on the right hand side of this expression yields the desired map:
    $$
    (N_{\mathcal{P}}M)_{F \oplus G}  \longrightarrow N\mathcal{P}_F \times_{N\mathcal{P}_{\rho(F)}} \times M(\rho(F)) = (N_{\mathcal{P}} M)_F.
    $$

    Finally, let $T$ be a tree obtained as a grafting $C_k \circ (T_1, \ldots, T_k)$, so that the root face map can be written as $T_1 \oplus \cdots \oplus T_k \rightarrow T$. By the Segal condition for the presheaf $N \mathcal{P}$, the set $(N_{\mathcal{P}} M)_T$ is isomorphic to
    $$ \left[ \prod_{i=1}^k N\mathcal{P}_{T_i} \times_{N\mathcal{P}_{k \cdot \eta}} N\mathcal{P}_{C_k} \right] \times_{N\mathcal{P}_{\rho(T)}} M(\rho(T)) =  \left( \prod_{i=1}^k N\mathcal{P}_{T_i} \right) \times_{N\mathcal{P}_{k \cdot \eta}} \left[ N\mathcal{P}_{C_k} \times_{N\mathcal{P}_{\rho(C_k)}} M(\rho(C_k)) \right],$$
    where in the last step we have replaced $T$ by $C_k$ via the inclusion of the root vertex of the former. Using the map $\cdot_{C_k}$ on the right hand side gives rise to the required structure map
    $$ (N_\mathcal{P} M)_T \longrightarrow \left( \prod_{i=1}^n N\mathcal{P}_{T_i} \right) \times_{N\mathcal{P}_{k \cdot \eta}} M(k \cdot \eta) = (N_{\mathcal{P}} M)_{T_1 \oplus \cdots \oplus T_k}.$$
    
\end{proof}

\begin{remark}
    The reason for asking for $\mathcal{P}$ to be closed is so that we can define how the summand inclusions act on $N_{\mathcal{P}} M$: in essence, this map corresponds to $M$ defining a colax monoidal presheaf with respect to the tensor product in $\Prop(\mathcal{P})$. Recall that we had already observed in \cref{colax_monoidal} that right modules over a closed operad always define a colax monoidal presheaves on $\Prop(\mathcal{P})$.
    
    For our purposes and examples we are mostly interested in, the condition for being closed is always satisfied, so this added constraint is not an important issue. However, if one wanted to define a similar functor for non-closed operads, then one can try to work with an altered version of $\For$ where one doesn't allow for summand inclusions (in the language of \cref{forest_definition}, this amounts to saying that if $f \colon I \to J$ is an injective map, then it is also surjective). The main consequence of this is that then there exists no map in $\fSets$
    $$ \coprod_{i=1}^n F_i \longrightarrow \bigoplus_{i=1}^n F_i$$
    for $F_i$ representable presheaves, which is troublesome for developing the theory of forest sets, especially when wants to compare with the dendroidal framework. With this in mind, we preferred to instead restrict the range of the operads we are working with.
\end{remark}

Let us spell out some basic properties of this nerve functor that we will need further on.

\begin{lemma}
    Let $\mathcal{P}$ a closed operad. Then the nerve functor $N_{\mathcal{P}} \colon \RMod{\mathcal{P}} \to \fSets_{/ N \mathcal{P}}$ is the right adjoint in an adjunction
    \[\begin{tikzcd}
	{\RMod{\mathcal{P}}} & {\fSets_{/N\mathcal{P}}}. & {}
	\arrow["{N_{\mathcal{P}}}"', shift right, from=1-1, to=1-2]
	\arrow["{\tau_{\mathcal{P}}}"', shift right, from=1-2, to=1-1]
\end{tikzcd}\]
    Additionally, the functor $N_{\mathcal{P}}$ commutes with coproducts in $\RMod{\mathcal{P}}$.
    \label{corpoduct_nerve}
\end{lemma}
\begin{proof}
    It suffices to specify what $\tau_{\mathcal{P}}$ does on representables $\alpha\colon F \to N\mathcal{P}$, because the condition that it commutes with colimits will then determine the rest of the functor. Given such a presheaf $\alpha$, we set 
    $$\tau_{\mathcal{P}}(\alpha) = \Free_{\mathcal{P}} \left( \rho(\alpha) \colon \rho(F) \to N\mathcal{P} \right)$$ 
    to be the free right $\mathcal{P}$-module on the set of colours determined by the restriction of $\alpha$ to the roots of $\alpha$, as discussed in \cref{examples_right_modules}. The fact that this is indeed the left adjoint follows essentially from the observation that $N_{\mathcal{P}} M$ is determined on a forest $F$ by its restriction to $\rho(F)$, if we consider it as a forest set over $N\mathcal{P}$.

    The fact that $N_{\mathcal{P}}$ commutes with coproducts is easily verified, using that in $\Sets$ pulling back along a morphism commutes with coproducts.
\end{proof}

In the rest of this section we will want to have a better understanding of the counit of the adjunction $(N_{\mathcal{P}}, \tau_{\mathcal{P}})$. Observe that the unit of this adjunction will always be an isomorphism, as one can easily verify that $N_{\mathcal{P}}$ is fully faithful, so we won't pursue its study any further.

Informally speaking, the counit of this adjunction will associate to a forest map $X \to N\mathcal{P}$ the free right $\mathcal{P}$-module generated by such data (more correctly, the right fibration associated with it), and therefore one should not expect it to be an isomorphism in general. However, we will now see that for the representables $F \to N\mathcal{P}$ the counit is a contravariant weak equivalence, which for our purposes is good enough.

Before stating and proving the next result, let us briefly describe the operadic right $\mathcal{P}$-module $\tau_{\mathcal{P}}(\alpha)$ for $\alpha\colon F \to N{\mathcal{P}}$ a representable, as this will be crucial in a moment. By definition, this is the free right $\mathcal{P}$-module on the colours determined by the map $\rho(\alpha) \colon \rho(F) \to N\mathcal{P}$ as we discussed in \cref{examples_right_modules} using representable presheaves. In terms of modules, we have the following:
\begin{itemize}[label=$\diamond$]
    \item If $r_1, \ldots, r_k$ are the roots of $F$, then there exists a \textit{generator} element of the module
    $$m_{\mathrm{gen}} \in \tau_{\mathcal{P}}(\alpha) (\alpha(r_1), \ldots, \alpha(r_k)),$$ 
    together with  $\sigma^* m_{\mathrm{gen}} \in \tau_{\mathcal{P}}(\alpha) (\alpha(r_{\sigma(1)}), \ldots, \alpha(r_{\sigma(k)}))$ for each $\sigma \in \Sigma_k$, which comes from the object in $\mathbf{Env}(\mathcal{P})$ as determined by $\rho(\alpha)$.
    \item The remaining elements of the module are freely generated by letting $\mathcal{P}$ act on the permutations of the generator $m_{\mathrm{gen}}$ under the structure maps of the right action. By the axioms of a right $\mathcal{P}$-module, such an element will be determined by some tuple of operations $(p_1, \ldots, p_k)$ in $\mathcal{P}$, with $p_i$ having output $\alpha(r_i)$. In fact, by the equivariance of operadic right modules with respect to the action of $\Sigma_n$, we can take $m_{\mathrm{gen}}$ instead of a permutation of it.
\end{itemize}

\begin{proposition}
    Let $\mathcal{P}$ be a closed $\Sigma$-free operad and $\alpha \colon k \cdot \eta \to N\mathcal{P}$ a forest map for some integer $k \geq 1$. Then the counit map
    $$ k \cdot \eta \longrightarrow N_{\mathcal{P}} \tau_{\mathcal{P}}(\alpha)$$
    is a forest root anodyne. In particular, it is a contravariant trivial cofibration in the overcategory $\fSets_{/N\mathcal{P}}$.
    \label{counit_roots}
\end{proposition}

\begin{proof}
    The fact that it is a contravariant trivial cofibration is a direct consequence of \cref{root_anodyne_ctv_we}, so we will focus on showing that the required map is root anodyne. Our proof is an adaptation of the case for dendroidal left anodynes presented in \citep[Prop. 9.67]{HeutsMoerdijkDendroidal}.

    Fix a forest $F \in \For$ and consider an element of $N_{\mathcal{P}} \tau_{\mathcal{P}}(\alpha)_F$, which, by the pullback description of the nerve of an operadic right module, is the data of a pair $(a,m)$ given by:
    \begin{itemize}[label=$\diamond$]
        \item An element $a \in N\mathcal{P}_F$, which in particular evaluates to the tuple of colours $(a_1, \ldots, a_n)$ when restricted to the forest of roots $\rho(F)$ (these colours are elements in the image of $\rho(\alpha)$ by the construction of $\tau_{\mathcal{P}}$).
        \item An element $m \in \tau_{\mathcal{P}}(\alpha)(a_1, \ldots, a_n)$. By the description of $\tau_{\mathcal{P}}$ we presented above, this is equivalent to giving the generating element $m_{\mathrm{gen}} \in \tau_{\mathcal{P}}(\alpha)(\alpha(r_1), \ldots, \alpha(r_k))$ of the module, as well as arity $n(i)$ operations
        $$ p_i \in \mathcal{P}(a^1_i, \ldots, a^{n(i)}_i; \alpha(r_i))$$
        for each $1 \leq i \leq k$, such that $m = m_{\mathrm{gen}} \cdot (p_1, \ldots, p_k)$.
    \end{itemize}

    We claim that from this information we can extract a different pair $(a', m_{\mathrm{gen}}) \in N_{\mathcal{P}} \tau_{\mathcal{P}}(\alpha)_{F'}$ in the following way: we can see each operation $p_i$ as an element of the presheaf $N\mathcal{P}$ evaluated at the corolla $C_{n(i)}$. Together with $a \in N\mathcal{P}_F$ and the Segal condition, these gather into $a' \in N\mathcal{P}_{F'}$, where $F'$ is obtained from $F$ by attaching the corollas $C_{n(i)}$ and therefore has roots $r_1, \ldots, r_k$. As $m_{\mathrm{gen}} \in \tau_{\mathcal{P}}(\alpha(r_1), \ldots, \alpha(r_k))$, so that $(a', m_{\mathrm{gen}})$ indeed lies in $N_{\mathcal{P}} \tau_{\mathcal{P}}(\alpha)_{F'}$.

    Let us define an \textit{admissible pair} to be any element of the presheaf $N_{\mathcal{P}} \tau_{\mathcal{P}}(\alpha)$ of the form $(a,m_{\mathrm{gen}})$. Then we can rephrase the previous paragraph as saying that any pair $(a,m) \in N_{\mathcal{P}} \tau_{\mathcal{P}} (\alpha)_F$ is in the image of an admissible pair $(a',m_{\mathrm{gen}}) \in N_{\mathcal{P}} \tau_{\mathcal{P}}(\alpha)_{F'}$ via a face map: more explicitly, the forest face map $F \to F'$ is given by applying the root faces associated to all the roots of $F'$. This observation allows us to filter the counit map as
    $$ k \cdot \eta = R_0 \subseteq R_1 \subseteq R_2 \subseteq \cdots \subseteq R_\ell \subseteq \cdots \subseteq N_{\mathcal{P}} \tau_{\mathcal{P}} (\alpha)$$
    where $R_\ell$ is the forest set generated by the non-degenerate admissible pairs indexed by forests with at most $\ell$ vertices. The forest set $N\mathcal{P}$ is a normal forest set since $\mathcal{P}$ is $\Sigma$-free, and we can use this to write each stage of the filtration via a pushout diagram in $\fSets$
    \begin{equation}
\begin{tikzcd}[cramped]
	{\underset{a}{\coprod} \: R_{\ell-1} \cap F} & {R_{\ell-1}} \\
	{\underset{a}{\coprod} \: F} & {R_\ell}.
	\arrow[from=1-1, to=1-2]
	\arrow[from=1-1, to=2-1]
	\arrow[from=1-2, to=2-2]
	\arrow[from=2-1, to=2-2]
	\arrow["\lrcorner"{anchor=center, pos=0, rotate=180, scale=1.5}, draw=none, from=2-2, to=1-1]
\end{tikzcd}
    \label{pushout_NF}
    \end{equation}
   
    Here the coproducts are indexed over all non-degenerate admissible pairs $(a,m_{\mathrm{gen}}) \in N_{\mathcal{P}} \tau_{\mathcal{P}}(\alpha)_F$ with $F$ having exactly $\ell$ vertices. On the one hand, the action of the elementary inner faces and leaf faces of $F$ will appropriately change $a$ but fix $m_{\mathrm{gen}}$, so these types of morphisms will preserve admissible pairs. On the other hand, any root face of a component of $F$ will originate a non-admissible pair, since the root faces act via the right module action, according to the definition of the nerve functor $N_{\mathcal{P}}$, and we assumed that the elements we are considering are non-degenerate.
    
     Putting these things together, we can identify each normal monomorphism $R_{n-1} \cap F \to F$ with the generalized root horn inclusion $\Lambda^R F \to F$ defined in \cref{generalized_horn}, where $R = R(F)$ is the set of all roots of $F$. By \cref{generalized_horn_root}, the left map in \eqref{pushout_NF} is therefore a root anodyne, and consequently the same holds for $R_{\ell-1} \to R_\ell$, which finishes the proof.
\end{proof}

The following corollary gives an extension of the previous result to the scenario when the source of $\alpha$ is any representable presheaf $F$.

\begin{corollary}
    Let $\mathcal{P}$ be a closed $\Sigma$-free operad and $\alpha \colon F \to N\mathcal{P}$ a forest map for $F \in \For$. Then the counit map
    $$ F \longrightarrow N_{\mathcal{P}} \tau_{\mathcal{P}}(\alpha)$$
    is a forest root anodyne. In particular, it is a contravariant trivial cofibration in the overcategory $\fSets_{/N\mathcal{P}}$.
    \label{counit_forest}
\end{corollary}

\begin{proof}
    Consider the commutative diagram below:
    $$
\begin{tikzcd}
\rho(F) \arrow[d] \arrow[r] & N_{\mathcal{P}} \tau_{\mathcal{P}}(\rho(\alpha)) \arrow[d] \\
F \arrow[r]                 & N_{\mathcal{P}} \tau_{\mathcal{P}}(\alpha).                
\end{tikzcd}
    $$

    By \cref{counit_roots}, the top horizontal map is a forest root anodyne. Moreover, the same also holds for both vertical maps: for the left one this is clear by \ref{closure_forest_root anodynes}, and for the right one it is actually an isomorphism, because the description of $\tau_{\mathcal{P}}$ we gave above clearly only depends on the value of $\alpha$ on $\rho(F)$. Consequently, by the right cancellation property for forest root anodynes, the same holds for the bottom map, as we wanted to show.
\end{proof}

\section{Operadic right modules via forest spaces}

We now turn to goal of contextualizing the theory of simplicial operadic right modules within the formalism of forest \textit{spaces}. With that in mind, we will dedicate Section 4.1 to reviewing some classical model category theory of the category of forest spaces $\fSpaces$, such as the existence of a projective and Reedy model structure, together with their complete Segal variants. After that, we will dedicate Section 4.2 to introducing the contravariant model structure within the formalism of forest spaces and proving some of its basic properties. In particular, in Section 4.4. we will investigate the more difficult question of how the contravariant model structure depends on the forest space we are slicing over. Finally, in Section 4.5 we prove \cref{eqv_ctv_RMod}, which establishes, for a certain class of well-behaved simplicial operads $\mathcal{P}$, a Quillen equivalence between the model category of right $\mathcal{P}$-modules, and the projective contravariant model structure on $\fSpaces_{/N\mathcal{P}}$. We conclude this work by extracting some consequences pertaining to this equivalence.

In Section 4.3 we also briefly go over the Day convolution product on $\fSpaces$ coming from the concatenation of forests  and how it interacts with the different model structures considered in the preceding sections, relating it to the tensor product we mentioned at the end of Section 2.

\subsection{The category of forest spaces and its homotopy theory}

\leavevmode

The purpose of this section is to recall some aspects developed in \cite{HeutsHinnichMoerdijk} about the category of \textit{forest spaces}, that is, the category of simplicial presheaves $\Fun(\For^{\mathsf{op}}, \sSets)$, which we will denote by $\fSpaces$. From now on, we will make very little use of the category of \textit{dendroidal spaces} $\mathsf{dSpaces}$, but we note that most of the discussion in this section also applies to this presheaf category with little change.

\begin{notation}
    In order to not overwhelm the notation, we will write $X(F) \in \sSets$ for the evaluation of a forest space $X$ at the forest $F$.
\end{notation}

Besides the projective model structure that $\fSpaces$ inherits from being a category of simplicial presheaves, it also inherits a Reedy model structure due to $\For$ being a generalized Reedy category, as defined in \cite{BergerMoerdijkReedy} by Berger--Moerdijk: we let the degree function be the number of edges, so that the degree lowering/raising morphisms are the maps which are injective/surjective at the level of the edges.

\begin{proposition}
    The identity functor induces a Quillen equivalence 
\[\begin{tikzcd}
	{\fSpaces_{\mathsf{P}}} & {\fSpaces_{\mathsf{R}}}
	\arrow["{\mathrm{id}_!}", shift left, from=1-1, to=1-2]
	\arrow["{\mathrm{id}^*}", shift left, from=1-2, to=1-1]
\end{tikzcd}\]
between the projective model structure $\fSpaces_{\mathsf{P}}$ and the Reedy model structure $\fSpaces_{\mathsf{R}}$.
\label{projective_reedy}
\end{proposition}

\begin{proof}
   By Corollary 10.12 of \cite{HeutsMoerdijkDendroidal}, it holds that, for any generalized Reedy category, the Reedy (trivial) fibrations are also projective fibrations, which shows that the adjunction above is a Quillen adjunction. It is also an equivalence since the classes of weak equivalences in both model categories are the same.
\end{proof}

An important consequence for us of the Quillen equivalence above is that the (derived) mapping spaces for the Reedy model structure and for the projective one coincide, so we don't need extra care in distinguishing these in what follows. We will write $\Map(X,Y)$ for such a mapping space, whenever $X$ and $Y$ are forest spaces.

It is easy to characterize the Reedy cofibrations in $\fSpaces$. Before doing so, we can associate to a given  forest space $X$ and simplicial set $K$ a new forest space $X \boxtimes K$, given by 
$$ (X \boxtimes K)(F) = X(F) \times K,$$
for $F \in \For$. The forest action on $X \boxtimes K$ is induced by the one on $X$. 

\begin{proposition}
    A forest map $f \colon X \to Y$ of forest spaces is a Reedy cofibration if and only if the map $X_k \to Y_k$ is a normal monomorphism of forest sets for each $k \geq 0$.

    In particular, normal forest sets and maps of normal monomorphisms of forest sets are Reedy cofibrant and Reedy cofibrations in $\fSpaces$, respectively.
    \label{characterization_Reedy_cofibration}
\end{proposition}

\begin{proof}
    The Reedy cofibrations in $\fSpaces$ are generated by the maps
    \begin{equation}
    F \boxtimes \partial \Delta^n \cup_{\partial F \boxtimes \partial \Delta^n} \partial F \boxtimes \Delta^n \longrightarrow F \boxtimes \Delta^n
    \label{generators_Reedy_cofibrations}
    \end{equation}
    for every forest $F$ and $n \geq 1$. Since forests are normal, it follows that the maps of the form \eqref{generators_Reedy_cofibrations} are levelwise normal monomorphisms. In fact, the class of such maps is easily verified to be saturated, so in particular must contain all Reedy cofibrations.

    In the other direction, suppose $f$ is a levelwise normal monomorphism. We write $L_F X \in \sSets$ for the latching object of $X$ evaluated at $F$, which can be described as the space of degenerate $F$-forests of $X$. By the definition of Reedy cofibrations, we need to verify that the latching map
    \begin{equation}
     L_F Y \cup_{L_F X} X(F) \longrightarrow Y(F)
    \label{latching_map}
     \end{equation}
    is a projective cofibration in the functor category $\Fun\left(\mathsf{Aut}_{\For}(F)^{\mathsf{op}}, \sSets\right)$. This involves checking two conditions, namely that \eqref{latching_map} is a monomorphism and the right $\mathsf{Aut}_{\For}(F)$-action is free on the complement of its image.

    We will first focus on showing that \eqref{latching_map} is a monomorphism, for which we need to check that, given $a \in X(F)_k$ such that $f(a) \in Y(F)_k$ is a degenerate element, $a$ is also degenerate. By hypothesis, let $\sigma \colon F \to G$ be a forest degeneracy such that $\sigma^{\ast} b = f(a)$ for some $b \in Y(G)_k$. From this data we can form the diagram
\[\begin{tikzcd}[cramped]
	F & {X_k} \\
	G & P \\
	&& {Y_k}
	\arrow["a", from=1-1, to=1-2]
	\arrow["\sigma"', from=1-1, to=2-1]
	\arrow["\psi", from=1-2, to=2-2]
	\arrow["f", bend left, from=1-2, to=3-3]
	\arrow[from=2-1, to=2-2]
	\arrow["b"', bend right, from=2-1, to=3-3]
	\arrow["\lrcorner"{anchor=center, pos=0, rotate=180, scale=1.5}, draw=none, from=2-2, to=1-1]
	\arrow["\varphi"', dashed, from=2-2, to=3-3]
\end{tikzcd}\]
where $P$ is defined as a pushout. As $\sigma$ is a degeneracy, it admits a section and therefore the same holds for $\psi$, which implies the surjectivity of this map. Moreover, as $f$ is a monomorphism, the same holds for $X_k \to P$, and thus $\psi$ is actually an isomorphism. The inverse of $\psi$ produces a decomposition of $a \colon F \to X_k$ of the form
$$ F \xlongrightarrow{\sigma} G \longrightarrow X_k,$$
which shows that $a$ is a degenerate element.

For the freeness condition, we observe first that we have an inclusion$X(F) \subseteq L_F Y \amalg_{L_F X} X(F)$ of subspaces of $Y(F)$, because the morphism $X(F) \to L_F Y \amalg_{L_F X} X(F)$ is a pushout of the monomorphism $L_F X \to L_F Y$, and thus is also a monomorphism. Therefore, after taking complements we get
$$ Y(F) - L_F Y \amalg_{L_F X} X(F) \subseteq Y(F) - X(F),$$
which of course is a $\mathsf{Aut}_{\For}(F)$-equivariant inclusion. As the rightmost set admits a free action, the same holds for the left hand side, which shows the required.
\end{proof}

The other important model structures on $\fSpaces$ we will want to consider are analogous to Rezk's complete Segal spaces model structure on $\mathsf{sSpaces}$ \cite{RezkSegalComplete}, and the version of Cisinski--Moerdijk for dendroidal spaces \cite{CisinskiMoerdijkSegal}. For the next definition, we recall that in \cref{definition_forest_spine} we defined the spine of a forest.

\begin{definition}
    Let $X$ be a forest space. We say that $X$ is satisfies the \textit{Segal condition} if, given any forest $F$, the map induced by the inclusion of the forest spine
    $$ \Map(F, X) \longrightarrow \Map(\fSp(F), X)$$
    is a weak homotopy equivalence of simplicial sets. We will write $\fSpaces_{\mathsf{PS}}$ and $\fSpaces_{\mathsf{RS}}$ for the left Bousfield localization of the projective and Reedy model structures for the Segal condition, respectively. We say $X$ is a \textit{Segal forest space} if it is fibrant in $\fSpaces_{\mathsf{RS}}$.
\end{definition}

We can also define what it means for a forest space to be complete.

\begin{definition}
    Let $X$ be a forest space. We say that $X$ is \textit{complete} if the inclusion\footnote{Here $J$ is seen as a forest set via $(u \iota)_!$, see Section 3 for the definition of this functor.} $\eta \to J$ of one of the objects of $J$ induces a weak homotopy equivalence
    $$ \Map(J, X) \longrightarrow \Map(\eta, X)$$
    of mapping spaces. We will write $\fSpaces_{\mathsf{PSC}}$ and $\fSpaces_{\mathsf{RSC}}$ for the left Bousfield localizations of $\fSpaces_{\mathsf{PS}}$ and $\fSpaces_{\mathsf{RS}}$ at the completeness condition, respectively. We will say that $X$ is a \textit{complete Segal forest space} if it is a fibrant object in $\fSpaces_{\mathsf{RSC}}$.
\end{definition}

\begin{remark}
We observe that the completeness condition only depends on the underlying simplicial space of the forest space $X$ in question: indeed, sincee the map we localize for is a simplicial map, it immediately follows that $X$ is a complete forest space if and only if its underlying simplicial space is complete in the sense of Rezk. 
\end{remark}

One can relate all these different model structures via \cref{projective_reedy}, which allows us to conclude that the following diagram of adjunctions, consisting only of the identity functors,

\[\begin{tikzcd}[cramped]
	{\fSpaces_{\mathsf{P}}} & {\fSpaces_{\mathsf{PS}}} & {\fSpaces_{\mathsf{PSC}}} \\
	{\fSpaces_{\mathsf{R}}} & {\fSpaces_{\mathsf{RS}}} & {\fSpaces_{\mathsf{RSC}}}
	\arrow[shift left, from=1-1, to=1-2]
	\arrow["\simeq"', shift right, from=1-1, to=2-1]
	\arrow[shift left, from=1-2, to=1-1]
	\arrow[shift left, from=1-2, to=1-3]
	\arrow["\simeq"', shift right, from=1-2, to=2-2]
	\arrow[shift left, from=1-3, to=1-2]
	\arrow["\simeq"', shift right, from=1-3, to=2-3]
	\arrow[shift right, from=2-1, to=1-1]
	\arrow[shift left, from=2-1, to=2-2]
	\arrow[shift right, from=2-2, to=1-2]
	\arrow[shift left, from=2-2, to=2-1]
	\arrow[shift left, from=2-2, to=2-3]
	\arrow[shift right, from=2-3, to=1-3]
	\arrow[shift left, from=2-3, to=2-2]
\end{tikzcd}\]

\vspace{0.1em}

\noindent is made up of Quillen adjunctions, and the vertical functors define Quillen equivalences. The following result gives us a comparison between the operadic model structure $\fSets_{\mathrm{opd}}$ and the model structure on complete Segal forest spaces $\fSpaces_{\mathsf{RSC}}$.
\begin{theorem}
    There is a Quillen equivalence
    \[\begin{tikzcd}
	{\fSets_{\mathsf{opd}}} & {\fSpaces_{\mathsf{RSC}}},
	\arrow["{\mathrm{con}}", shift left, from=1-1, to=1-2]
	\arrow["{\mathrm{ev}_0}", shift left, from=1-2, to=1-1]
\end{tikzcd}\]
where $\mathrm{con}(X)=X$ sees a forest set as a discrete forest space, and $\mathrm{ev}_0(X)_F = X(F)_0$.
\label{quillen_eqv_segal_complete}
\end{theorem}

\begin{proof}
    This is the content of Theorem 3.9.4 of \cite{HeutsHinnichMoerdijk}, with the small caveat already made in \cref{not_quillen_equivalence}. This extra localization made by Heuts--Hinnich--Moerdijk translates into condition $(f\beta')$ that appears in their proof. However, their argument at the start of Section 3.9 still holds without this extra restriction, and will correspond to the aforementioned Quillen equivalence.
\end{proof}

Since our interest for forest spaces mostly stems from the presheaves which actually come from dendroidal spaces via the functor $u_\ast$, we would like to be assured that the different notions of fibrancy in $\fSpaces$ that we mentioned coincide with the ones that are already known for $\mathsf{dSpaces}$. This is guaranteed by the next lemma.

\begin{lemma}
    Let $\ast$ stand for $\mathsf{P}, \mathsf{R}, \mathsf{PS}, \mathsf{RS}, \mathsf{PSC}$ or $\mathsf{RSC}$. Then the adjoint pair
    \begin{equation}
    \begin{tikzcd}
	{\mathsf{fSpaces}_{\ast}} & {\mathsf{dSpaces}_{\ast}}
	\arrow["{u_*}", shift left, from=1-2, to=1-1]
	\arrow["{u^\ast}", shift left, from=1-1, to=1-2]
\end{tikzcd}
\label{quillen_adjunction}
\end{equation}
defines a Quillen adjunction.
\label{everything_is_the_same}
\end{lemma}

\begin{proof}
    It is clear from the explicit formula for $u_\ast$ that it preserves projective (trivial) fibrations. As for the Reedy model structure, the condition on \cref{characterization_Reedy_cofibration} characterizing forest Reedy cofibrations also characterizes dendroidal Reedy cofibrations, by replacing all forests with trees. From this observation, it is clear that $u^*$ preserves Reedy (trivial) cofibrations. This deals with the cases when $\ast$ is $\mathsf{P}$ or $\mathsf{R}$.

    These Quillen adjunctions induce Quillen pairs, where on the left hand side of \eqref{quillen_adjunction} we localize for the Segal condition, and on the right hand side we localize for 
    $$ u^* \fSp(F) \longrightarrow u^*(F).$$
    However, this is the same as localizing for the Segal condition: this is a consequence of $u^*$ commuting with coproducts and $\fSp(-)$ commuting with $\oplus$. Therefore the Quillen adjunction we want also exists for when $\ast$ is $\mathsf{PS}$ and $\mathsf{RS}$. Finally, the case for the completeness condition is immediate since we are localizing for the same simplicial map $\eta \to J$.

\end{proof}

\subsection{First properties of the contravariant model structure on $\fSpaces_{/V}$}

\leavevmode

In this section we will construct a version of the contravariant model structure for forest spaces, which will be of central importance when comparing the theory of right fibrations of forest spaces to the theory of right modules over an operad. Having this in mind, let us begin by explaining what a right fibration of forest spaces should be.

\begin{definition}
    Let $p \colon X \to V$ be a map of forest spaces. We say that $p$ is a \textit{right fibration} if, for any forest $F$ and morphism $\alpha \colon F \to V$, the induced map between mapping spaces
    \begin{equation}
    \alpha^\ast : \Map_V(F,X) \longrightarrow \Map_V(\rho(F),X)
    \label{localitiy_contravariant}
    \end{equation}
    is a weak homotopy equivalence.
    \label{definition_right_fib}
\end{definition}

The mapping spaces above should be taken with respect to either the relative projective or Reedy model structures on $\fSpaces_{/V}$; again by \cref{projective_reedy}, the choice of model structure doesn't impact the homotopy type of $\Map_V(X,Y)$.

We would like to have a more practical description of right fibrations of forest spaces for when we try to relate the contravariant model structure to the theory of operadic right modules via  $N_{\mathcal{P}}$. The next easy lemma will be useful for computing mapping spaces. 

\begin{lemma}
    Let $V$ be a forest space, together with forest maps $p \colon X \to V$ and $q \colon Y \to V$. Then the following statements hold:
    \begin{enumerate}[label=(\alph*)]
        \item If $q$ is a Reedy fibration and $X$ is a Reedy cofibrant object, then 
        $$\Map_V(X,Y)_k = \left( \fSpaces_{/V} \right) (X \boxtimes \Delta^k, Y)$$
        for every $k \geq 0$.
        \item If in (a) we also assume that $V$ is Reedy fibrant, then $\Map_V(X,Y)$ is the fiber of
        $$ \Map(X,Y) \longrightarrow \Map(X,V)$$
        over $p \in \Map(X,V)_0$.
    \end{enumerate}

    The same conclusions hold if our assumptions are taken with respect to the projective model structure instead of the Reedy model structure.

    \label{mapping_spaces_reedy}
\end{lemma}

\begin{proof}
    We will focus on the Reedy statements, as the projective case is entirely analogous.

    For part (a), the cosimplicial object $\Delta^{\bullet}$ constitutes a cosimplicial resolution of the terminal forest space: in fact, this holds more generally for any category of simplicial presheaves, see Example 11.15 in \cite{HeutsMoerdijkDendroidal}. Consequently, if $X$ is Reedy cofibrant then $X \boxtimes \Delta^{\bullet}$ is a cosimplicial resolution of $X$. Therefore, the formula in (a) determines a model for the mapping space $\Map_V(X,Y)$ whenever $q \colon Y \to V$ is a Reedy fibration.

    In particular, this shows that one can compute $\Map_V(X,Y)$ as the fiber of the simplicial map
    \begin{equation}
    \fSpaces(X \boxtimes \Delta^{\bullet}, Y) \xlongrightarrow{q \circ -}\fSpaces(X \boxtimes \Delta^{\bullet}, V)
    \label{q_ast}
    \end{equation}
    over the composition of $X \boxtimes \Delta^{\bullet} \to X \xrightarrow{p} V$, where $X \boxtimes \Delta^{\bullet} \to X$ is the map exhibiting the source as a cosimpicial resolution of $X$. If we now assume that $V$ is Reedy fibrant, then so is $Y$ since $q$ is a Reedy fibration. Consequently, the source and target simplicial sets in \eqref{q_ast} serve as models for the respective mapping spaces, which shows (b).
\end{proof}

We can now give a different characterization of right fibrations which avoids mentioning mapping spaces and is easier to verify in practice.

\begin{proposition}
    Let $p \colon X \to V$ be a map of forest spaces. Then $p$ is a right fibration if and only if, for every $F \in \For$, the diagram
  \begin{equation}
\begin{tikzcd}
X(F) \arrow[d] \arrow[r] & X(\rho(F)) \arrow[d] \\
V(F) \arrow[r]           & V(\rho(F))    
\end{tikzcd}
\label{right_fibration_hopullback}
\end{equation}
    is a homotopy pullback square.
    \label{characterization_ctv_fibrant}
\end{proposition}

\begin{proof}
    Since we are interested in the homotopy type of $X$ evaluated at a forest $F$ and projective weak equivalences are levelwise weak homotopy equivalences, we can replace $p$ by a projective fibration over $V$, which we assume from now on to be the case.

     Since $\rho(F)$ and $F$ are projectively cofibrant as they are representable presheaves, we are in position of applying \cref{mapping_spaces_reedy} so that the map in \eqref{localitiy_contravariant} is obtained as the map between the vertical fibers of
     $$
\begin{tikzcd}
{\fSpaces(F \boxtimes \Delta^{\bullet}, X)} \arrow[d] \arrow[r] & {\fSpaces(\rho(F) \boxtimes \Delta^{\bullet}, X)} \arrow[d] \\
{\fSpaces(F \boxtimes \Delta^{\bullet}, V)} \arrow[r]     & {\fSpaces(\rho(F) \boxtimes \Delta^{\bullet}, V)},   
\end{tikzcd}
     $$
     which in turn coincides with the square in the statement. To check that the pullback condition is equivalent to \eqref{localitiy_contravariant} being a weak homotopy equivalence, it suffices to check that the vertical fibers coincide with the homotopy fibers. This follows from $p$ being a projective fibration, so that the vertical maps are Kan fibrations of simplicial sets.
\end{proof}

If $V$ is a dendroidal Segal space, then $V$ as a forest space will be local for the forest spine inclusions $\fSp(F) \to F$, but also for the wide spine $\wfSpine(F) \to F$. The next result, which we will only need for Chapter 4, says that this locality condition can be lifted along right fibrations $X \to V$ from the target to the source.

\begin{proposition}
    Suppose $V$ is a dedroidal Segal space and $X \to V$ is a Reedy right fibration over $V$, with $X$ a forest space. 
    \begin{enumerate}[label=(\alph*)]
        \item  Let $F$ be a forest together with a choice of leave $\ell \in L(F)$. Given $k \geq 0$, the diagram 
\[\begin{tikzcd}[cramped]
	{X(F \circ_{\ell} C_k)} & {X(C_k)} \\
	{X(F)} & {X(\eta)}
	\arrow[from=1-1, to=1-2]
	\arrow[from=1-1, to=2-1]
	\arrow[from=1-2, to=2-2]
	\arrow[from=2-1, to=2-2]
\end{tikzcd}\]
is homotopy cartesian. Here $F \circ_\ell C_k \in \For$ is the forest obtained from $F$ by grafting a copy of the corolla $C_k$ along $\ell$.
    \item The forest space $X$ is local with respect to the wide spine, that is, for any forest $F$ the map
    $$ X(F) \longrightarrow X(\wfSpine(F))$$
     is a weak homotopy equivalence. 
     \item Suppose $V$ is a reduced dendroidal space. If $F = \bigoplus_{i=1}^k F_i$ such that $F_i$ is a linear tree, then there is a weak homotopy equivalence
     $$ X(F) \xlongrightarrow{\simeq} X(k \cdot \eta)$$
     induced by the inclusion of the roots of $F$.
    \end{enumerate}
    \label{right_fibration_Segal}
\end{proposition}

\begin{proof}
    For (a), consider the diagram
\[\begin{tikzcd}[cramped]
	{X(F \circ_\ell C_k)} & {X(F)} & {X(\rho(F))} \\
	{V(F \circ_\ell C_k)} & {V(F)} & {V(\rho(F))}
	\arrow[from=1-1, to=1-2]
	\arrow[from=1-1, to=2-1]
	\arrow[from=1-2, to=1-3]
	\arrow[from=1-2, to=2-2]
	\arrow[from=1-3, to=2-3]
	\arrow[from=2-1, to=2-2]
	\arrow[from=2-2, to=2-3]
\end{tikzcd}\]
Then all the squares above are homotopy pullbacks, since the whole rectangle and the rightmost square are homotopy cartesian because of the right fibration condition. We can complete the leftmost square into a cube
\[\begin{tikzcd}[scale=0.7, cramped]
	& {X(F\circ_\ell C_k)} && {X(F)} \\
	{X(C_k)} && {X(\eta)} \\
	& {V(F \circ_{\ell} C_k)} && {V(F)} \\
	{V(C_k)} && {V(\eta)}
	\arrow[from=1-2, to=1-4]
	\arrow[from=1-2, to=2-1]
	\arrow[from=1-2, to=3-2]
	\arrow[from=1-4, to=2-3]
	\arrow[from=1-4, to=3-4]
	\arrow[from=2-1, to=2-3]
	\arrow[from=2-1, to=4-1]
	\arrow[from=2-3, to=4-3]
	\arrow[from=3-2, to=3-4]
	\arrow[from=3-2, to=4-1]
	\arrow[from=3-4, to=4-3]
	\arrow[from=4-1, to=4-3]
\end{tikzcd}\]
which has the bottom and front faces also homotopy cartesian: the former due to $V$ coming from a dendroidal Segal space, and the latter by the right fibration condition. As we have already seen that the back face is a homotopy pullback, we conclude that the same holds for the top face, as we wanted to show.

Part (b) is a consequence (in fact, it is equivalent) to part (a). Indeed, one proceeds inductively on the number of vertices of a forest $F$. The base case is $F = n \cdot \eta$ for which there is nothing to prove since $\wfSpine(F) = F$ in this situation. For the induction step, we use that any forest $F \neq n \cdot \eta$ must be of the form $F = F' \circ_\ell C_k$ for some corolla $C_k$, and then we can apply (a).

    For (c), $V$ being reduced implies that $V(F) \simeq V(\rho(F))$, because $V$ sends $\oplus$ to products and $V(L) \simeq V(\eta)$ for linear trees. The root inclusion for $F$ yields a diagram 
\[\begin{tikzcd}[cramped]
	{X(F) } & {X(k \cdot \eta)} \\
	{V(F)} & {V(k \cdot \eta)}
	\arrow[from=1-1, to=1-2]
	\arrow[from=1-1, to=2-1]
	\arrow[from=1-2, to=2-2]
	\arrow[from=2-1, to=2-2]
\end{tikzcd}\]
which is homotopy cartesian by the right fibration condition. The top map is a weak homotopy equivalence by putting together the observations we made.
\end{proof}

\begin{remark}
    Assuming that $V$ is a dendroidal Segal space and not just a forest Segal space is essential in the proposition above. The reason is that a general forest Segal space $V$ will not be local with respect to the wide spine, but it is if $V$ comes from a dendroidal Segal space.
\end{remark}

The next alternative characterization of Reedy fibrations which are also right fibrations will be quite useful in later parts of the text.
\begin{lemma}
    Let $p\colon X \to V$ a Reedy fibration of forest spaces. Then the following statements are equivalent:
    \begin{enumerate}[label=(\alph*)]
        \item The map $X \to V$ is a right fibration of forest spaces.
        \item For any forest $F$, the map
        $$ X(F) \longrightarrow X(\rho(F)) \times_{V(\rho(F))} V(F)$$
        is a trivial Kan fibration.
        \item For any forest $F$ and any inner or root horn inclusion $\Lambda^a F \to F$, the map
        $$ X(F) \longrightarrow X(\Lambda^a F) \times_{V(\Lambda^a F)} V(F)$$
        is a trivial Kan fibration.
    \end{enumerate}
    \label{def_Reedy_right_fibration}
\end{lemma}

\begin{proof}
    The Reedy fibrancy condition on $p$ yields that the maps in (b) and (c) are Kan fibrations by Proposition 10.10 of \cite{HeutsMoerdijkDendroidal}. The equivalence of (a) and (b) is then immediate from  \cref{characterization_ctv_fibrant}, using that the Kan--Quillen model structure is right proper and that the right map in \eqref{right_fibration_hopullback} is a Kan fibration (so the strict pullback coincides with the homotopy pullback). The equivalence of (b) and (c) can be checked using the characterization of the class of forest root anodynes we presented in \cref{closure_forest_root anodynes}.
\end{proof}

For the next example, we note that the dendroidal nerve construction can be extended to the simplicial context, yielding in that scenario a functor
$$ N : \sOp \longrightarrow \mathsf{dSpaces}$$
given by $(N\mathcal{P})_k = N(\mathcal{P}_k)$, where here we see a simplicial operad as a presheaf $\mathcal{P}_{\bullet} \colon \Delta^{\mathsf{op}} \to \Op$ with constant set of colours. Similarly, the nerve functor for operadic right modules we constructed in the previous sections also admits a simplicial version
\[\begin{tikzcd}
	{\RMod{\mathcal{P}}} & {\fSpaces_{/N\mathcal{P}}},
	\arrow["{N_{\mathcal{P}}}"', shift right, from=1-1, to=1-2]
	\arrow["{\tau_{\mathcal{P}}}"', shift right, from=1-2, to=1-1]
\end{tikzcd}\]
with a similar restriction that $\mathcal{P}$ needs to a closed simplicial operad. The definition for $N_{\mathcal{P}}M$ is entirely similar, given as the strict pullback
\begin{equation}
\begin{tikzcd}[cramped]
	{\left( N_{\mathcal{P}}M \right)(F)} & {M(\rho(F))} \\
	{N\mathcal{P}(F)} & {N\mathcal{P}(\rho(F)).}
	\arrow[from=1-1, to=1-2]
	\arrow[from=1-1, to=2-1]
	\arrow["\lrcorner"{anchor=center, pos=0, scale=1.5}, draw=none, from=1-1, to=2-2]
	\arrow[from=1-2, to=2-2]
	\arrow[from=2-1, to=2-2]
\end{tikzcd}
\label{ho_pullback_question}
\end{equation}
\begin{example}
    As the diagram \eqref{ho_pullback_question} defining $N_{\mathcal{P}}$ indicates, the nerve functor for operadic right modules provides plentiful examples of right fibrations, \textit{except} that the strict pullback on question might not coincide with the homotopy version. We present here two situations where these are the same.
    \begin{enumerate}[label=(\alph*)]
        \item Suppose $\mathcal{P}$ is a closed simplicial operad with all spaces of operations $\mathcal{P}(c_1, \ldots, c_n; c_0)$ being Kan complexes. Then we claim that, for any simplicial right $\mathcal{P}$-module $M$, the forest map $N_{\mathcal{P}}M \to N\mathcal{P}$ is a right fibrations of forest spaces. 
    
    Indeed, by the definition of the nerve functor and \cref{characterization_ctv_fibrant}, it is sufficient to check that the square \eqref{ho_pullback_question} is a homotopy pullback square. As $N\mathcal{P}(\rho(F))$ is discrete and from our assumptions on the spaces of operations, we conclude that the map $N\mathcal{P}(F) \to N\mathcal{P}(\rho(F))$ is a Kan fibration; since the Kan--Quillen model structure is right proper, this implies the intended.

    \item Suppose now that $\mathcal{P}$ is only a closed simplicial operad and that $M$ is a projectively fibrant right $\mathcal{P}$-module. Then the same conclusion as in (a) holds. To see this, note that $M(\rho(F))$ is a Kan complex for any choice of forest $F$, and therefore the right map in \eqref{ho_pullback_question} is a Kan fibration, as the base is discrete, and the rest of argument goes as in (a). It is also clear from this discussion that, under this hypothesis, $N_{\mathcal{P}}M \to N\mathcal{P}$ is also a projective fibration.
    \label{examples_fibrations_right}
    \end{enumerate}

    Let $\mathcal{P}$ be a closed simplicial operad with all spaces of operations $\mathcal{P}(c_1, \ldots, c_n; c_0)$ being Kan complexes. Then we claim that, for any simplicial right $\mathcal{P}$-module $M$, the forest map $N_{\mathcal{P}}M \to N\mathcal{P}$ is a right fibrations of forest spaces. 
    
    Indeed, by the definition of the nerve functor and \cref{characterization_ctv_fibrant}, it is sufficient to check that the square \eqref{ho_pullback_question} is a homotopy pullback square. As $N\mathcal{P}(\rho(F))$ is discrete and from our assumptions on the spaces of operations, we conclude that the map $N\mathcal{P}(F) \to N\mathcal{P}(\rho(F))$ is a Kan fibration; since the Kan--Quillen model structure is right proper, this implies the intended.
    \label{operadic_nerve_right_fibrations}
\end{example}

Our strategy for defining the contravariant model structure for the slice categories $\fSpaces_{/V}$ will be different from how we defined it for the discrete case in previous sections. Since we would like this model structure to reflect in some way the theory of right fibrations and taking inspiration from \cref{definition_right_fib}, we make the following definition.

\begin{definition}
    Let $V$ be a forest space and let $\ast$ stand for either $\mathsf{P}$ or $\mathsf{R}$. We define the \textit{contravariant model structure} on $\fSpaces_{/V}$ to be the left Bousfield localization of the relative model category $ \left( \fSpaces _{/V} \right)_\ast$ corresponding to $\ast$, with respect to the set of all root inclusions $\rho(F) \to F$ over $V$, for $F \in \For$. We will write
    $$\left( \fSpaces_{/V} \right)_{\ast, \ctv}$$
    for the contravariant model category. If $\ast = \mathsf{P}$ we call this the \textit{projective contravariant model category}, and if $\ast = \mathsf{R}$ we will call it the \textit{Reedy contravariant model structure}.
\end{definition}

One could naturally ask about defining versions of the contravariant model structure where we take as initial input the versions of the projective and Reedy model structures which have bee localized for the additional Segal and completeness conditions. These will not play a very important role in what follows, but they will be useful for some technical arguments. The next lemma will be sufficient for our purposes.

\begin{lemma}
    Let $V$ be a forest space and let $\mathcal{M}$ be $\fSpaces_{/V}$ equipped with the Reedy contravariant model structure.
    \begin{enumerate}[label=(\alph*)]
        \item If $V$ is a Segal forest space, then $\mathcal{M}$ is a left Bousfield localization of the relative model structure
        $$ \left( \fSpaces_{/V} \right)_{\mathsf{RS}}$$
    at the contravariant condition.
    \item If $V$ is a complete Segal forest space, then $\mathcal{M}$ is a left Bousfield localization of the relative mode structure
        $$ \left( \fSpaces_{/V} \right)_{\mathsf{RSC}}$$
    at the contravariant condition.
    \end{enumerate}
    \label{localizing_contravariant}
\end{lemma}

\begin{proof}
    Since all model structures in question have the same class of cofibrations, it suffices to check that the fibrant objects are the same. This is the verification which is done in \citep[Thm. 13.6]{HeutsMoerdijkDendroidal}, which is simple to adapt to our context by respectively changing "trees" and "left" with "forests" and "right" through their argument.
\end{proof}

\begin{notation}
    Following the definition of the projective and Reedy contravariant model categories, we will also write 
    $$ \left( \fSpaces_{/V} \right)_{\mathsf{RS}, \ctv} \hspace{1em} \text{and} \hspace{1em} \left( \fSpaces_{/V} \right)_{\mathsf{RSC}, \ctv}$$
    for the contravariant variants mentioned in \cref{localizing_contravariant}.
\end{notation}

Our next result builds a version of the Quillen equivalence of \cref{quillen_eqv_segal_complete} comparing this newly-defined contravariant model structure with the one we already studied for $\fSets$.

\begin{theorem}
    Let $V$ be a complete Segal forest space. Then there is a Quillen equivalence
    \[\begin{tikzcd}
	{\left( \fSets_{/ \mathsf{ev}_0(V)} \right)_{\ctv} } & {\left( \fSpaces_{/V} \right)_{\mathsf{R},\: \ctv} }.
	\arrow["{\mathrm{con}}", shift left, from=1-1, to=1-2]
	\arrow["{\mathrm{ev}_0}", shift left, from=1-2, to=1-1]
\end{tikzcd}\]
\label{dsicrete_space}
\end{theorem}

\begin{proof}

    By case (k) of \citep[Ex. 8.47]{HeutsMoerdijkDendroidal}, the Quillen equivalence of \cref{quillen_eqv_segal_complete} carries over to the relative setting
    \[\begin{tikzcd}
	{\left( \fSets_{/ \mathsf{ev}_0(V)} \right)_{\mathsf{opd}} } & {\left( \fSpaces_{/V} \right)_{\mathsf{RSC}} }
	\arrow["{\mathrm{con}}", shift left, from=1-1, to=1-2]
	\arrow["{\mathrm{ev}_0}", shift left, from=1-2, to=1-1]
\end{tikzcd}\]
since $V$ is by definition a fibrant object of the model category on the right. We can localise on the left for the root inclusions $\rho(F) \to F$ over $\mathsf{ev}_0(V)$, which by \cref{character_fibrant_bosufield} leads to the Quillen equivalence
\[\begin{tikzcd}
	{\left( \fSets_{/ \mathsf{ev}_0(V)} \right)_{\ctv} } & {\left( \fSpaces_{/V} \right)_{\mathsf{RSC},\: \mathsf{con}(\ctv)} }.
	\arrow["{\mathrm{con}}", shift left, from=1-1, to=1-2]
	\arrow["{\mathrm{ev}_0}", shift left, from=1-2, to=1-1]
\end{tikzcd}\]

For the right hand side, localizing for $\mathsf{con}(\ctv)$ is the same as localizing for the contravariant condition, as any root inclusion $\rho(F) \to F$ over $V$ automatically factors through $\mathsf{ev}_0(V)$. Finally, we can now replace the $\mathsf{RSC}$ condition for the Reedy condition $\mathsf{R}$, due to \cref{localizing_contravariant}(b).
\end{proof}

We end this section with a comparison result between the projective contravariant model structure and the Reedy contravariant model structure.

\begin{proposition}
    Let $V$ be a projectively fibrant forest space, and let $\varphi\colon V \to V'$ be a Reedy fibrant replacement of $V$. Then the adjunction induced by this fibrant replacement
    \[\begin{tikzcd}
	{\left( \fSpaces_{/ V} \right)_{\mathsf{P}, \ctv} } & {\left( \fSpaces_{/V'} \right)_{\mathsf{R},\: \ctv} }
	\arrow["{\varphi_!}", shift left, from=1-1, to=1-2]
	\arrow["{\varphi^*}", shift left, from=1-2, to=1-1]
\end{tikzcd}\]
defines a Quillen equivalence.
\label{right_fibr_fibrant}
\end{proposition}

\begin{proof}
    We first claim that such a Quillen equivalence exists before localizing with respect to the contravariant model structure. To see this, notice that this adjoint pair can be decomposed as the composition of the following adjunctions:
\[\begin{tikzcd}[cramped]
	{(\fSpaces_{/V})_{\mathsf{P}}} & {(\fSpaces_{/V'})_{\mathsf{P}}} & {(\fSpaces_{/V'})_{\mathsf{R}}}.
	\arrow["{\varphi_!}", shift left, from=1-1, to=1-2]
	\arrow["{\varphi^*}", shift left, from=1-2, to=1-1]
	\arrow["{\mathrm{id}_!}", shift left, from=1-2, to=1-3]
	\arrow["{\mathrm{id}^*}", shift left, from=1-3, to=1-2]
\end{tikzcd}\]

The leftmost adjunction is a Quillen equivalence, since $\varphi$ a projective equivalence between projectively fibrant objects and by \citep[Ex. 8.47(i)]{HeutsHinnichMoerdijk}. The other pair also defines a Quillen equivalence, by a relative version of \cref{projective_reedy}.

Consequently, after localizing for the contravariant model structure we get a Quillen equivalence
\[\begin{tikzcd}
	{\left( \fSets_{/ V} \right)_{\mathsf{P}, \ctv} } & {\left( \fSpaces_{/V'} \right)_{\mathsf{R},\: \varphi_!(\ctv)} },
	\arrow["{\varphi_!}", shift left, from=1-1, to=1-2]
	\arrow["{\varphi^*}", shift left, from=1-2, to=1-1]
\end{tikzcd}\]
where on the right hand side we are localizing with respect to all root inclusions $\rho(F) \to F$ over $V'$ that factor along $\varphi$. Therefore, we will have shown the intended once we verify that the rightmost model category corresponds to the Reedy contravariant model structure on $\fSpaces_{/V'}$, and to do so it is enough to check that the fibrant objects of both model categories are the same, since the classes of cofibrations clearly coincide.

On the one hand, if we localize $\left( \fSpaces_{/V'} \right)_{\mathsf{R}}$ for the contravariant condition, then \cref{right_fibr_fibrant} identifies its fibrant objects as the Reedy fibrations $X \to V'$ such that, for the diagram below,
\begin{equation}
\begin{tikzcd}
X(F) \arrow[d] \arrow[r] & X(\rho(F)) \arrow[d] \\
V'(F) \arrow[r]           & V'(\rho(F))    
\end{tikzcd}
\label{help_right_fibration}
\end{equation}
the map between the vertical homotopy fibers is a weak equivalence, for any choice of basepoint in $V'(F)$. On the other hand, if we instead localise for $\varphi_!(\ctv)$, then the fibrant objects are the Reedy fibrations such that the map between the vertical homotopy fibers of \eqref{help_right_fibration} is a weak equivalence, \textit{but} only for the basepoints in the image of $f\colon V(F) \to V'(F)$. However, $f$ being a projective equivalence implies that $\pi_0\left( V(F) \right) \to \pi_0 \left( V'(F) \right)$ is an isomorphism, so this condition for the vertical fibers actually holds for all basepoints of $V(F)$. Hence Reedy $\ctv$-fibrant objects coincide with Reedy $\varphi_!(\ctv)$-fibrant objects in this case.
\end{proof}

\subsection{The concatenation product of forest spaces}

\leavevmode

In this short section we will introduce a tensor product on the category of forest spaces called the \textit{concatenation product}. In particular, we will see that all the model categories we presented in the sections above can be enhanced to monoidal model categories with respect to this monoidal structure. Finally, we will further explore how the concatenation product relates to the product $\boxplus$ we had already defined in Section 2 for operadic right modules.

We recall that the category $\For$ admits a bifunctor
$$ - \oplus - : \For \times \For \longrightarrow \For$$
which corresponds to the concatenation of forests. In fact, this operation endows $\For$ with a symmetric monoidal structure, with unit\footnote{This is the only reason why we included the empty forest as a forest, so that this monoidal structure has a unit.} given by the empty forest $\emptyset$. As a consequence, the category of forest spaces admits a closed symmetric monoidal structure $(\fSpaces, \oplus, \emptyset)$ via Day convolution, and the tensor product is exactly the bifunctor
$$ - \oplus - : \fSpaces \times \fSpaces \longrightarrow \fSpaces$$
obtained by the formula
\begin{equation}
 (F \boxtimes K) \oplus (G \boxtimes L) = (F \oplus G) \boxtimes (K \times L)
 \label{boxes_interaction}
 \end{equation}
for forests $F, G$ and simplicial sets $K, L$, and extended to $\fSpaces$ via left Kan extension. We will call this the \textit{concatenation product}. It will also be useful to recall that the coend formula for the Day convolution allows us to write the simplicial set $(X \oplus Y)(F)$ as the quotient of the coproduct
\begin{equation}
\coprod_{F_1, F_2 \in \For} \For(F, G_1 \oplus G_2) \times X(G_1) \times Y(G_2)
\label{Day_coproduct_formula}
\end{equation}
by the following relation: if $(\varphi, x', y') \in \For(F, G_1 \oplus G_2) \times X(G'_1) \times Y(G'_2)$ and $\alpha_i \colon G_i \to G'_i$, then $((\alpha_1 \oplus \alpha_2) \varphi, x', y')$ is identified with $(\varphi, \alpha_1^* x', \alpha_2^* y')$. As we already did for right modules, this tensor product can be further simplified; the proof of it is exactly the same as the one in \cref{coend_simplify_modules}, so we omit it.

\begin{proposition}
   If $X$ and $Y$ are forest spaces, then, for every forest $F$, there is a natural isomorphism
        $$ (X \oplus Y)(F) \cong \coprod_{F_1, F_2} X(F_1) \times Y(F_2),$$
        where the coproduct is taken over subforests $F_1, F_2$ of $F$  inducing an isomorphism $F_1 \oplus F_2 \xrightarrow{\cong} F$.
    \label{coend_fibrations}
\end{proposition}

Since we are ultimately interested in slice categories of $\fSpaces$, the following result will be important. We will say that a forest space $V$ is \textit{reduced} if $V(\emptyset)=\Delta^0$.

\begin{lemma}
    Let $V$ be a dendroidal space such that $u_\ast V$ is reduced. Then the following hold:
    \begin{enumerate}[label=(\alph*)]
        \item The forest space $u_\ast V$ is a monoid object in the symmetric monoidal category $(\fSpaces, \oplus, \emptyset)$.
        \item The slice category 
    $$\fSpaces_{/V}$$
    carries a closed symmetric monoidal structure induced by the concatenation product $\oplus$.  
    \end{enumerate}
\end{lemma}

\begin{proof}
    The second statement follows from the first one, as the reader can easily verify that slicing a closed symmetric monoidal category over a monoid object leads to a closed symmetric monoidal category. The product and internal hom are defined in the remark following the proof.

    For the first part, we need to construct the multiplication map $\mu \colon u_\ast V \oplus u_\ast V \to u_\ast V$. By the adjunction $(u^*, u_\ast)$ together with the property $u^*(X \oplus Y) \cong u^*(X) \amalg u^*(Y)$, the existence of $\mu$ is equivalent to constructing a morphism of the form
    $$ u^* u_* V \amalg u^* u_* V \longrightarrow V,$$
    which we define to the codiagonal of the counit $u^* u_* V \to V$. The unit of the monoid structure $\eta\colon \emptyset \to u_\ast V$ is the unique map of this form, since $u_\ast V$ is assumed to be reduced.
\end{proof}

\begin{remark}
    Still using the notation in the proof, the concatenation product of $\alpha_i \colon X_i \to u_\ast V$ for $i=1,2$ is
    $$ X_1 \oplus X_2 \xlongrightarrow{\alpha_1 \oplus \alpha_2} u_\ast V \oplus u_\ast V \xlongrightarrow{\mu} u_\ast V,$$
    and if $[-,-]$ denotes the internal hom in $\fSpaces$, then its sliced version $[X_1, X_2]_V \in \fSpaces_{/V}$ is the pullback
\[\begin{tikzcd}[cramped]
	{[X_1, X_2]_V} && {[X_1,X_2]} \\
	u_\ast V & {[u_\ast V, u_\ast V]} & {[X_1,u_\ast V]}
	\arrow[from=1-1, to=1-3]
	\arrow[from=1-1, to=2-1]
	\arrow["\lrcorner"{anchor=center, pos=0, scale=1.5}, draw=none, from=1-1, to=2-2]
	\arrow["{(\alpha_2)_\ast}", from=1-3, to=2-3]
	\arrow["{\hat{\mu}}", from=2-1, to=2-2]
	\arrow["{(\alpha_1)^*}", from=2-2, to=2-3]
\end{tikzcd}\]
where $\hat{\mu}$ is the forest morphism adjoint to the monoid multiplication map $\mu$. 

\end{remark}

\begin{remark}
    As we have mentioned before, the empty forest should be thought of as an extra element we have added to $\For$ to ensure the existence of a unit for the concatenation product. This has the consequence in the previous proof of us having to restrict our attention to reduced forest spaces, but this is pretty harmless. Indeed, the full subcategory of $\fSpaces$ on all reduced forest spaces is reflective, and the left adjoint is given via the pushout diagram
\[\begin{tikzcd}[cramped]
	\emptyset & X \\
	{\Delta^0} & {\mathsf{red}(X).}
	\arrow[from=1-1, to=1-2]
	\arrow[from=1-1, to=2-1]
	\arrow[from=1-2, to=2-2]
	\arrow[from=2-1, to=2-2]
	\arrow["\lrcorner"{anchor=center, pos=0, rotate=180, scale=1.5}, draw=none, from=2-2, to=1-1]
\end{tikzcd}\]

This only changes the space $X(\emptyset)$ and leaves $X(F)$ unaltered for $F$ non-empty. Since the evaluation of $X$ at the empty forest is irrelevant information for us, this shows that the restriction to reduced forest spaces is not relevant.
\label{empty_forest}
\end{remark}

As an easy application of the formula \cref{coend_fibrations} we prove the following lemma, which is the forest version of part (b) of \cref{coend_simplify_modules}.

\begin{lemma}
    Let $V$ be a dendrodial space such that $u_\ast V$ is reduced. Then, if $X \to V$ and $Y \to V$ are projective right fibrations over $V$, then their concatenation product $X \oplus Y \to V$ over $V$ is a projective right fibration.
    \label{concantenation_fibrations}
\end{lemma}

\begin{proof}
    The fact that it is still a projective fibration follows immediately from \cref{coend_fibrations}. The right fibration property is also clear: indeed, for each $F$, we need to verify the condition of the diagram in \cref{characterization_ctv_fibrant}, which corresponds to the outer rectangle in
    \[\begin{tikzcd}[cramped]
	{\coprod_{F_1, F_2} X(F_1)\times Y(F_2)} & {\coprod_{F_1, F_2} V(F_1)\times V(F_2)} & {V(F)} \\
	{\coprod_{F_1, F_2} X(\rho(F_1)) \times Y(\rho(F_2))} & {\coprod_{F_1, F_2} V(\rho(F_1)) \times V(\rho(F_2))} & {V(\rho(F))}
	\arrow[from=1-1, to=1-2]
	\arrow[from=1-1, to=2-1]
	\arrow[from=1-2, to=1-3]
	\arrow[from=1-2, to=2-2]
	\arrow[from=1-3, to=2-3]
	\arrow[from=2-1, to=2-2]
	\arrow[from=2-2, to=2-3]
\end{tikzcd}\]
being a homotopy pullback. The left square is a homotopy pullback since $X \to V$ and $Y \to V$ are right fibrations, and the right square is certainly s homotopy pullback also as the top and horizontal maps are coproducts of isomorphisms due to $V$ be a dendroidal space. 
\end{proof}

We now turn to the question of how the concatenation product intertwines with the model structures we have considered up to now.

\begin{proposition}
    Let $V$ be a dendroidal space such that $u_\ast V$ is reduced. Then the following four model categories
    $$ \fSpaces_\mathsf{P} \hspace{2em} \fSpaces_\mathsf{R} \hspace{2em} \left(\fSpaces_{/V} \right)_{\mathsf{P}} \hspace{2em} \left(\fSpaces_{/V} \right)_{\mathsf{R}}$$
    have the structure of a monoidal model category with respect to the concatenation product.
    \label{monoidal_without_ctv}
\end{proposition}

\begin{proof}
    We will only focus on proving the absolute versions of the statements, as the relative versions have an entirely analogous proof. Secondly, since the unit $\emptyset$ is both projective and Reedy cofibrant, we only need to check that pushout-product axiom with respect to $\oplus$ is satisfied; as all the model structures are cofibrantly generated and $- \oplus -$ commutes with colimits on both variables, it suffices to consider the generating (trivial) cofibrations.

    For the projective model structure, the generating cofibrations and trivial cofibrations are
    $$ F \boxtimes \Lambda^i [n] \longrightarrow F \boxtimes \Delta^n \hspace{1em} \mathrm{and} \hspace{1em}  F \boxtimes \partial \Delta^n \longrightarrow F \boxtimes \Delta^n$$
    respectively, and the verification of the pushout-product axiom goes through as we did in \cref{monoidality_projective}, using equation \eqref{boxes_interaction}.

     As for the Reedy model structure, we will need to consider pushout-products of maps of the form
    $$\psi_{F, A \to B} : F \boxtimes A \cup_{\partial F \boxtimes A} \partial F \boxtimes B \longrightarrow F \boxtimes B $$
    where $A \to B$ is some simplicial map and $F$ is a forest. After some manipulation and using that $\boxplus$ coincides with $\oplus$ on representables, one easily concludes that the pushout-product of $\psi_{F, A \to B}$ with $\psi_{G, A' \to B'}$ is
    $$(F \oplus G) \boxtimes K \cup_{\partial(F \oplus G) \boxtimes K} \partial(F \oplus G) \boxtimes (B \times B') \longrightarrow (F \oplus G) \boxtimes (B \times B'),$$
    where $K = (A \times B') \cup_{A \times A'} (B \times A')$. The necessary verification needed for the Reedy model structure is now easily done using this formula, with the simplicial maps being either horn inclusions or boundary inclusions of simplices.    
\end{proof}

\begin{proposition}
    Let $V$ be a dendroidal space such that $u_\ast V$ is reduced. Then the following two model categories
    $$  \left(\fSpaces_{/V} \right)_{\mathsf{P}, \ctv} \hspace{2em} \left(\fSpaces_{/V} \right)_{\mathsf{R}, \ctv}$$
    have the structure of a monoidal model category with respect to the concatenation product.
    \label{monoidal_with_ctv}
\end{proposition}

\begin{proof}
    The considerations we made at the start of the previous proof also apply here. In this case, we only need to verify that the pushout-product of a cofibration with a trivial cofibration is a trivial cofibration, as the class of cofibrations is preserved under left Bousfield localization.
   
   The generating projective contravariant trivial cofibrations are generated by the generating projective trivial cofibrations and the map
   \begin{equation}
    \rho(G) \boxtimes \Delta^m \cup_{\rho(G) \boxtimes \partial \Delta^m} G \boxtimes \partial \Delta^m \longrightarrow G \boxtimes \Delta^m
    \label{generator}
   \end{equation}
   for every forest $G$ and $m \geq 0$. After some computations, we see that the pushout-product of such a map with $F \boxtimes \partial \Delta^n \to F \boxtimes \Delta^n$ is
   $$ (G \oplus F) \boxtimes \partial(\Delta^m \times \Delta^n) \cup_{(\rho(G) \oplus F) \boxtimes \partial( \Delta^m \times \Delta^n)} (\rho(G) \oplus F) \boxtimes (\Delta^m \times \Delta^n) \longrightarrow (G \oplus F) \boxtimes (\Delta^m \times \Delta^n).$$

    We already know that this map is a cofibration, so it suffices to check that it is a weak equivalence in the projective contravariant model structure. For that, it will be sufficient to verify that the forest map $\rho(G) \oplus F \to G \oplus F$ induces a weak homotopy equivalence
    $$ \Map_V(G \oplus F, X) \longrightarrow \Map_V(\rho(G) \oplus F, X)$$
    between mapping spaces, where $X \to V$ is a projective right fibration. Here the domain is projectively cofibrant and the source is a projective fibration, so the argument in the proof of \cref{characterization_ctv_fibrant} shows that we can equivalently show that the left hand side square in
\[\begin{tikzcd}[cramped]
	{X(G \oplus F)} & {X(\rho(G) \oplus F)} & {X(\rho(G \oplus F))} \\
	{V(G \oplus F)} & {V(\rho(G) \oplus F)} & {V(\rho(G \oplus F))}
	\arrow[from=1-1, to=1-2]
	\arrow[from=1-1, to=2-1]
	\arrow[from=1-2, to=1-3]
	\arrow[from=1-2, to=2-2]
	\arrow[from=1-3, to=2-3]
	\arrow[from=2-1, to=2-2]
	\arrow[from=2-2, to=2-3]
\end{tikzcd}\]
    is a homotopy pullback square, where the right hand square comes from applying $\rho(-)$ to the $F$ component. Note that the outer rectangle and right hand square are homotopy pullbacks, by the characterization of right fibrations in \cref{characterization_ctv_fibrant}. It follows that then so is the left hand square, as we wanted to show.

    We can proceed similarly for the Reedy case, and we get in that case that map we want to analyse in the pushout-product, with respect to $\boxtimes$, of the forest map
    \begin{equation}
     G \oplus \partial F \cup_{\rho(G) \oplus \partial F} \rho(G) \oplus F \longrightarrow G \oplus F
    \label{need_help}
     \end{equation}
    and the simplicial boundary inclusion $\partial( \Delta^m \times \Delta^n) \to \Delta^m \times \Delta^n$. It suffices to check that the map \eqref{need_help} is a trivial cofibration in the Reedy contravariant model structure over $V$. This is an immediate consequence of the square below
\[\begin{tikzcd}[cramped]
	{\rho(G)} & { G \oplus \partial F \cup_{\rho(G) \oplus \partial F} \rho(G) \oplus F} \\
	G & {G \oplus F}
	\arrow[from=1-1, to=1-2]
	\arrow[from=1-1, to=2-1]
	\arrow[from=1-2, to=2-2]
	\arrow[from=2-1, to=2-2]
	\arrow["\lrcorner"{anchor=center, pos=0, rotate=180, scale=1.5}, draw=none, from=2-2, to=1-1]
\end{tikzcd}\]
    being a pushout square, together with $\rho(G) \to G$ being a trivial cofibration in the Reedy contravariant model structure.
    \end{proof}

Finally, we give a result relating the concatenation product of forest spaces with the tensor product of simplicial operadic modules.

\begin{proposition}
    Let $\mathcal{P}$ be a closed simplicial operad. Then the nerve functor for operadic right modules defines a strong monoidal functor
    $$ N_{\mathcal{P}} : \left( \RMod{\mathcal{P}}, \boxplus, \mathcal{Y}_{\emptyset}\right) \longrightarrow \left( \fSpaces_{/ N\mathcal{P}}, \oplus, \emptyset \right).$$
    In particular, the left adjoint $\tau_\mathcal{P} \colon \fSpaces_{/N\mathcal{P}} \to \RMod{\mathcal{P}}$ defines an oplax monoidal functor with respect to the monoidal structures above.
    \label{NP_strong}
\end{proposition}

\begin{proof}

We need to produce a natural isomorphism
$$ N_\mathcal{P}(M) \oplus N_\mathcal{P}(L) \longrightarrow N_{\mathcal{P}}(M \boxplus L)$$
for every right $\mathcal{P}$-modules $M$ and $L$. This is given by the following chain of isomorphisms:
\begin{align*}
    N_\mathcal{P}(M \boxplus L)(F) &= (M \boxplus L)(\rho(F)) \times_{N\mathcal{P}(\rho(F))} N\mathcal{P}(F) \\
    &\cong \coprod_{F_1, F_2} M(\rho(F_1)) \times L(\rho(F_2)) \times_{N\mathcal{P}(\rho(F_1)) \times N\mathcal{P}(\rho(F_2))} N\mathcal{P}(F_1) \times N\mathcal{P}(F_2) \\
    &\cong \coprod_{F_1, F_2} M(\rho(F_1)) \times_{N\mathcal{P}(\rho(F_1)} N\mathcal{P}(F_1) \times L(\rho(F_2)) \times_{N\mathcal{P}(\rho(F_2))}  N\mathcal{P}(F_2) \\
    &\cong \coprod_{F_1, F_2} N_{\mathcal{P}}M(F_1) \times N_{\mathcal{P}}L(F_2) \\
    &\cong (N_\mathcal{P}(M) \oplus N_{\mathcal{P}}L)(F)
\end{align*}
Here the coproduct is over $F_1$ and $F_2$ as in \cref{coend_fibrations}, and we have used the isomorphism $N\mathcal{P}(F) \xrightarrow{\cong} N\mathcal{P}(F_1) \times N\mathcal{P}(F_2)$ in the second step. As for the unit, there is also an obvious isomorphism $\emptyset \xrightarrow{\cong} N_\mathcal{P}(\emptyset)$.
\end{proof}

\subsection{Dependence of the contravariant model structure on the base forest space}

\leavevmode

    In the previous sections we provided some first basic results on the behaviour of the contravariant model structure on $\fSpaces_{/V}$ and how it relates to the version we had previously considered for forest sets. A natural question that arises is about how this new model category depends on the choice of the object we are slicing over; in this regard, we have already given an answer to this same question for $\fSets_{/V}$ in \cref{contravariance_operad_invariance}, where we showed that we produce Quillen equivalent model categories whenever we slice over operadically equivalent operadically fibrant forest sets. 
    
    This same question for the overcategories of $\fSpaces$ seems to be technically more intricate: the main technical hurdle we will have to tackle when analysing these problems will be in showing that certain forest maps $f \colon X \to Y$ over a fixed $V$ are in fact contravariant weak equivalences in $\fSpaces_{/V}$. The idea we will employ to overcome our difficulties is that of \textit{simplicial diagrams of model categories}, which is explained in \citep[Sec. 13.2]{HeutsMoerdijkDendroidal} and which we recall in Appendix A. The main insight one should retain from this approach is the following: the presheaf category $\fSpaces$ can also be seen as the presheaf category $\Fun(\Simplicial^{\mathsf{op}}, \fSets)$. After taking the slice by $V$, this identification sends the category $\fSpaces_{/V}$ to a simplicial diagram\footnote{We have omitted the degeneracy maps from the diagram.}
\[\begin{tikzcd}[cramped]
	{\fSets_{/V_0}} & {\fSets_{/V_1}} & {\fSets_{/V_2}} & \cdots
	\arrow[shift right, from=1-2, to=1-1]
	\arrow[shift left, from=1-2, to=1-1]
	\arrow[shift right=2, from=1-3, to=1-2]
	\arrow[shift left=2, from=1-3, to=1-2]
	\arrow[from=1-3, to=1-2]
	\arrow[shift left=3, from=1-4, to=1-3]
	\arrow[shift left, from=1-4, to=1-3]
	\arrow[shift right=3, from=1-4, to=1-3]
	\arrow[shift right, from=1-4, to=1-3]
\end{tikzcd}\]
of overcategories constructed out of $\fSets$. In particular, any map $f \colon X \to Y$ over $V$ will induce a collection of morphisms $f_i \colon X_i \to Y_i$ in $\fSets_{/V_i}$ for each $i \geq 0$. If it were the case that each $f_i$ was a contravariant weak equivalence in $\fSets_{/V_i}$, then it is natural to ask if $f$ is a contravariant weak equivalence in $\fSpaces_{/V}$. We will show in \cref{bootstrap} that this is indeed the case for any choice of $V$. The main upshot of such a result is that we usually have more criteria for checking that each $f_i$ is a contravariant equivalence: for instance, this happens if $f_i$ is a forest root anodyne, by \cref{root_anodyne_ctv_we}.

We will thus begin by showing that we can indeed use our knowledge of the discrete version of discrete contravariant model structure to bootstrap information regarding its simplicial version. For the next proposition, we will freely use the language of Appendix A.

\begin{proposition}
    Let $V$ be a forest space, and consider the the simplicial diagram of model categories $\mathcal{M}_{\bullet} = \{\left( \fSets_{/V_i} \right)_{\ctv} : i \geq 0 \}$. 
    
    Then the categorical identification $\Tot(\mathcal{M}) \cong \fSpaces_{/V}$ of \cref{example_2} can be enhanced to a Quillen equivalence of model categories
    $$ \left( \fSpaces_{/V} \right)_{\mathsf{R}, \ctv} \simeq \Tot\left( \mathcal{M} \right)_{\mathsf{w}}.$$
    \label{identification_simplicial_diagram_ctv}
\end{proposition}

\begin{proof}
The proof is similar to the one in \citep[Thm. 13.23]{HeutsHinnichMoerdijk} for the covariant model structure on dendroidal spaces. We will use throughout the equivalence of categories 
 \[\begin{tikzcd}
	{\Tot(\mathcal{M})} & {\fSpaces_{/V} }
	\arrow["{L}", shift left, from=1-1, to=1-2]
	\arrow["{R}", shift left, from=1-2, to=1-1]
\end{tikzcd}\]
that we discussed in \cref{example_2}. Using an argument similar to the one in \cref{characterization_Reedy_cofibration}, one can check that the cofibrations in $\Tot(\mathcal{M})_{\mathsf{w}}$ (that is, the Reedy cofibrations) are exactly the same as the Reedy cofibrations in $\fSpaces_{/V}$. Therefore, it suffices to consider whether the equivalence $(L,R)$ identifies the fibrant objects of both categories.

By \cref{def_Reedy_right_fibration}, one can characterize the fibrant objects in $\left( \fSpaces_{/V} \right)_{\mathsf{R}, \ctv}$ as the forest maps $X \to V$ having the right lifting property with respect to the following classes of morphisms:
\begin{enumerate}[label=(\alph*)]
    \item For each $n \geq 0$ and $0 \leq k \leq n$, the map
    $$ F \boxtimes \Lambda^k[n] \cup_{\partial F \boxtimes \Lambda^k [n]} \partial F \boxtimes \Delta^n \longrightarrow F \boxtimes \Delta^n$$
    for each forest $F$. This accounts for the map being a Reedy fibration.
    \item For each $n \geq 0$, the map
    $$ F \boxtimes \partial \Delta^n \cup_{\Lambda^a(F) \boxtimes \partial \Delta^n} \Lambda^a(F) \boxtimes \Delta^n \longrightarrow F \boxtimes \Delta^n$$
    for each forest $F$, and inner or root horn inclusion $\Lambda^a(F) \to F$. This accounts for the map being a right fibration.
\end{enumerate}

From this point on the proof of the theorem consists in checking how having the right lifting property with respect to (a) and (b) translates into within the category $\Tot(\mathcal{M})_{\mathsf{w}}$ and that it indeed corresponds to the notion of fibrancy we wanted. The proof follows by replacing left fibrations with right fibrations in \citep[Thm. 13.23]{HeutsHinnichMoerdijk}, and using the characterization in \cref{ctv_character} of the contravariantly fibrant objects in $\fSets \to V_i$ as the right fibrations over $V_i$.
\end{proof}

\begin{corollary}
    Let $V$ be a forest space, and $f \colon X \to Y$ a forest map over $V$ such that each map $f_i \colon X_i \to Y_i$ is a contravariant weak equivalence in $\left( \fSets_{/V_i} \right)_{\ctv}$. Then $f$ is a contravariant weak equivalence in $\left( \fSpaces_{/V} \right)_{\mathsf{R}, \ctv}$.
    \label{bootstrap}
\end{corollary}

\begin{proof}
    The condition on the morphisms $f_i$ implies that, under the identification of \cref{identification_simplicial_diagram_ctv} and using the language in its proof, $R(f)$ will correspond to a weak equivalence in $\Tot(\mathcal{M})_{\mathsf{R}}$ by the characterization in \cref{Reedy_simplicial_diagrams}. As we indicated in the sketch of the proof of \cref{RSC_Totalization}, $\Tot(\mathcal{M})_{\mathsf{w}}$ is a left Bousfield localization of $\Tot(\mathcal{M})_{\mathsf{R}}$, and therefore $R(f)$ remains a weak equivalence in the localized setting $\Tot\left( \mathcal{M} \right)_{\mathsf{w}}$. Together with \cref{identification_simplicial_diagram_ctv}, this finishes the proof.
\end{proof}

Having done all the technical preliminaries, we will focus for the rest of this section on giving some possible answers to the question of the homotopy dependence of the contravariant model structure on the forest space we are slicing over.

\begin{proposition}
    Let $\varphi\colon V \to W$ be a weak equivalence between fibrant objects in $\fSpaces_{\mathsf{RS}}$, that is, $\varphi$ is a levelwise weak equivalence between Segal forest spaces. Then $\varphi$ induces a Quillen equivalence
    \[\begin{tikzcd}
	{\left( \fSpaces_{/ V} \right)_{\mathsf{R}, \ctv} } & {\left( \fSpaces_{/W} \right)_{\mathsf{R},\: \ctv} }
	\arrow["{\varphi_!}", shift left, from=1-1, to=1-2]
	\arrow["{\varphi^*}", shift left, from=1-2, to=1-1]
\end{tikzcd}\]
    between the Reedy contravariant model structures.
    \label{ctv_Segal_eqv}
\end{proposition}

\begin{proof}
    Firstly, one can assume that $\varphi$ is a trivial fibration in $\fSpaces_{\mathsf{RS}}$ by Brown's lemma. From this extra fibrancy condition, it follows that the Quillen adjunction
    \[\begin{tikzcd}
	{\left( \fSpaces_{/ V} \right)_{\mathsf{RS}} } & {\left( \fSpaces_{/W} \right)_{\mathsf{RS}} }
	\arrow["{\varphi_!}", shift left, from=1-1, to=1-2]
	\arrow["{\varphi^*}", shift left, from=1-2, to=1-1]
\end{tikzcd}\]
is actually a Quillen equivalence between the relative model structures. Therefore, in a similar fashion to what happened in the proof of \cref{right_fibr_fibrant}, we see that this also defines a Quillen equivalence after localizing for the contravariant condition, as long as we ensure that on the right hand side of the adjunction this is the same as localizing with respect to $\varphi_!(\ctv)$. This can be argued by noting that a dashed lift in
\[\begin{tikzcd}[cramped]
	& V \\
	F & W
	\arrow["\varphi", from=1-2, to=2-2]
	\arrow[dashed, from=2-1, to=1-2]
	\arrow[from=2-1, to=2-2]
\end{tikzcd}\]
always exists, since $F$ is cofibrant for any of the model structures considered, and $\varphi$ is a trivial fibration.

To eliminate the Segal condition that appears throughout, one now appeals to part (a) of \cref{localizing_contravariant}, which provides the identification of
$$ \left( \fSpaces_{/V}\right)_{\mathsf{R}, \ctv} \hspace{1em} \text{and} \hspace{1em} \left( \fSpaces_{/W}\right)_{\mathsf{R}, \ctv}$$
with the $\mathsf{RS}$ variants of the contravariant model structure we have considered, since $V$ and $W$ are assumed to be Segal forest spaces.
\end{proof}

For a wide variety of purposes the homotopy invariance result we have just proven is enough, as it for instance captures the homotopy invariance for forest Segal spaces (thus in particular for $\infty$-operads) connected by levelwise weak equivalences. However, a situation which is not covered by \cref{ctv_Segal_eqv} is that of completion of forest spaces.

Recall that a \textit{completion} of a forest Segal space $X$ is a map (or zig-zag of such)
$$ \varphi : X \longrightarrow \overline{X}$$
such that $\varphi$ is a complete weak equivalence and $\overline{X}$ is a complete Segal forest space. We will now show how we can improve on \cref{ctv_Segal_eqv}
by using our knowledge of simplicial diagrams of model categories, and show that the contravariant model structure is invariant under completions.
\begin{proposition}
    Let $\varphi \colon V \to W$ be a forest map between Segal forest spaces, such that $\varphi_i \colon V_i \to W_i$ is an operadic equivalence for each simplicial degree $i \geq 0$. Then there is a Quillen equivalence
    \[\begin{tikzcd}
	{\left( \fSpaces_{/ V} \right)_{\mathsf{R}, \ctv} } & {\left( \fSpaces_{/W} \right)_{\mathsf{R},\: \ctv} }.
	\arrow["{\varphi_!}", shift left, from=1-1, to=1-2]
	\arrow["{\varphi^*}", shift left, from=1-2, to=1-1]
\end{tikzcd}\]
\label{bootstrap_ctv_eqv}
\end{proposition}

\begin{proof}
    Each $\varphi_i$ is an operadic equivalence between operadically fibrant objects. By \cref{contravariance_operad_invariance} (b), we get Quillen equivalences
    \[\begin{tikzcd}
	{\left( \fSets_{/V_i} \right)_\ctv} & {\left( \fSets_{/W_i}\right)_{\ctv}}
	\arrow["{\varphi_!}", shift left, from=1-1, to=1-2]
	\arrow["{\varphi^*}", shift left, from=1-2, to=1-1]
\end{tikzcd}\]
The desired conclusion is now a combination of \cref{identification_simplicial_diagram_ctv} together with \cref{comparison_totalizations}.
\end{proof}

\begin{corollary}
    Let $\varphi \colon V \to W$ be a weak equivalence in $\fSpaces_{\mathsf{RSC}}$ between Segal forest spaces. Then there is a Quillen equivalence
    \[\begin{tikzcd}
	{\left( \fSpaces_{/ V} \right)_{\mathsf{R}, \ctv} } & {\left( \fSpaces_{/W} \right)_{\mathsf{R},\: \ctv} }.
	\arrow["{\varphi_!}", shift left, from=1-1, to=1-2]
	\arrow["{\varphi^*}", shift left, from=1-2, to=1-1]
\end{tikzcd}\]
\label{complete_weak_equivalence}
\end{corollary}

\begin{proof}
    If $V$ and $W$ were complete, then this means that they would be fibrant in $\fSpaces_{\mathsf{RSC}}$, and therefore $\varphi$ is a projective weak equivalences by general theory of left Bousfield localizations. Hence we are in the situation of \cref{ctv_Segal_eqv} and the result follows.
    
    If $V$ and $W$ are not complete, then we can consider completions of $V$ and $W$. Appealing now to \cref{bootstrap_ctv_eqv}, it will suffice that we show that there is a completion $V \to \overline{V}$ such that $V_i \to \overline{V}_i$ is an operadic weak equivalence for each $i \geq 0$ -- the fact that this implies the result for all completions follows from another application of \cref{ctv_Segal_eqv}. Such a completion can always be constructed: this follows from an adaption to the forest setting of the proof of \citep[Cor 12.21]{HeutsMoerdijkDendroidal}, which is immediate.
\end{proof}

\begin{corollary}
    Let $V$ be a Segal forest space and suppose $\varphi \colon V \to \overline{V}$ is a completion. Then there is a Quillen equivalence
    \[\begin{tikzcd}
	{\left( \fSpaces_{/ V} \right)_{\mathsf{R}, \ctv} } & {\left( \fSpaces_{/\overline{V}} \right)_{\mathsf{R},\: \ctv} }.
	\arrow["{\varphi_!}", shift left, from=1-1, to=1-2]
	\arrow["{\varphi^*}", shift left, from=1-2, to=1-1]
\end{tikzcd}\]
between the respective contravariant model structures.
\label{completion_invariance}
\end{corollary}

\begin{remark}
    All these results have projective versions by \cref{right_fibr_fibrant}.
\end{remark}

\subsection{Operadic right modules as right fibrations in $\fSpaces$}

\leavevmode

We recall that in \cref{operadic_nerve_right_fibrations} we saw that simplicial operadic right modules give us examples of right fibrations under certain restrictions on the operad $\mathcal{P}$; we will also see soon that one can get a similar result without restricting the operad, but instead asking for the module to be projectively fibrant. In any case, from these observations it is natural that we try to have a more meaningful understanding of this relation, which we will do in this section. We will focus on using our knowledge on the contravariant model structure in order to analyse the adjunction
\[\begin{tikzcd}
	{\RMod{\mathcal{P}}} & {\fSpaces_{/N\mathcal{P}}.} & {}
	\arrow["{N_{\mathcal{P}}}"', shift right, from=1-1, to=1-2]
	\arrow["{\tau_{\mathcal{P}}}"', shift right, from=1-2, to=1-1]
\end{tikzcd}\]
coming from the nerve of an operadic right module. 

\begin{theorem}
    Let $\mathcal{P}$ be a closed simplicial $\Sigma$-free operad. Then the adjoint pair defined by the nerve functor $N_{\mathcal{P}}$ can be upgraded to a Quillen equivalence
    \[\begin{tikzcd}
	{\RMod{\mathcal{P}} } & {\left( \fSpaces_{/N\mathcal{P}} \right)_{\mathsf{P}, \ctv}} & {}
	\arrow["{N_{\mathcal{P}}}"', shift right, from=1-1, to=1-2]
	\arrow["{\tau_{\mathcal{P}}}"', shift right, from=1-2, to=1-1]
\end{tikzcd}\]
    between the projective model structure on $\RMod{\mathcal{P}}$, and the projective contravariant model structure on $\fSpaces_{/N\mathcal{P}}$.
    \label{eqv_ctv_RMod}
\end{theorem}

We will start the proof of \cref{eqv_ctv_RMod} by showing that the nerve adjunction defines a Quillen pair, even without the $\Sigma$-free condition.

\begin{lemma}
    Let $\mathcal{P}$ be a closed simplicial operad. Then the adjoint pair
    \[\begin{tikzcd}
	{\RMod{\mathcal{P}}} & {\left( \fSpaces_{/N\mathcal{P}} \right)_{\mathsf{P}, \ctv}} & {}
	\arrow["{N_{\mathcal{P}}}"', shift right, from=1-1, to=1-2]
	\arrow["{\tau_{\mathcal{P}}}"', shift right, from=1-2, to=1-1]
\end{tikzcd}\]
is a Quillen adjunction. Moreover, $N_{\mathcal{P}}$ also detects weak equivalences.

\label{Modules_Quillen_Adjunction}

\end{lemma}

\begin{proof}
    It is immediate that
    \[\begin{tikzcd}
	{\RMod{\mathcal{P}}} & {\left( \fSpaces_{/N\mathcal{P}} \right)_{\mathsf{P}}} & {}
	\arrow["{N_{\mathcal{P}}}"', shift right, from=1-1, to=1-2]
	\arrow["{\tau_{\mathcal{P}}}"', shift right, from=1-2, to=1-1]
\end{tikzcd}\]
    is a Quillen pair, as $N_{\mathcal{P}}$ will preserve weak equivalence and fibrations, since these are defined levelwise for both model categories. To check that this Quillen adjunction still holds after localizing for the contravariant condition, it suffices to check that $N_{\mathcal{P}}$ sends fibrant modules to right fibrations over $N\mathcal{P}$ which are also projective fibrations. This is exactly the content of \cref{examples_fibrations_right}(b).
    
    For the statement concerning the detection of weak equivalences, note that it is certainly true before localizing for the contravariant condition, since the projective weak equivalences are defined levelwise. In order to deal with the localization condition, we make the following two observations:
    \begin{enumerate}[label=(\alph*)]
        \item  If $\varphi\colon M \to L$ is a map between \textit{fibrant} right $\mathcal{P}$-modules, then $N_{\mathcal{P}}(\varphi)$ is a contravariant equivalence if and only if it is a projective equivalence: indeed, $N_{\mathcal{P}}(\varphi)$ is a map between contravariantly fibrant objects, and, for any left Bousfield localization, a map between locally fibrant objects is a weak equivalence if and only if it is a local weak equivalence.
        \item If $M$ and $L$ are not necessarily fibrant modules, consider the diagram below where the horziontal maps are fibrant replacements:
\[\begin{tikzcd}[cramped]
	M & {M'} \\
	L & {L'}
	\arrow["\simeq", tail, from=1-1, to=1-2]
	\arrow["\varphi"', from=1-1, to=2-1]
	\arrow["{\varphi'}", from=1-2, to=2-2]
	\arrow["\simeq", tail, from=2-1, to=2-2]
\end{tikzcd}\]

As the top and bottom maps are projective weak equivalences, they are sent to projective weak equivalences via $N_\mathcal{P}$, which are in particular contravariant weak equivalences. Therefore $\varphi$ is a weak equivalence if and only if $\varphi'$ is a weak equivalence, and $N_{\mathcal{P}}(\varphi)$ is a contravariant weak equivalence if and only if $N_\mathcal{P}(\varphi')$ is a contravariant weak equivalence.
    \end{enumerate}

Chaining these properties together we get
\begin{align*}
    \varphi \: \text{is a weak equivalence} &\overset{(b)}{\iff} \varphi' \: \text{is a weak equivalence} \\
    &\iff  N_{\mathcal{P}}(\varphi') \: \text{is a projective weak equivalence} \\
    &\overset{(a)}{\iff} N_{\mathcal{P}}(\varphi') \: \text{is a contravariant weak equivalence} \\
    &\overset{(b)}{\iff} N_{\mathcal{P}}(\varphi) \: \text{is a contravariant weak equivalence} \\
\end{align*}
\end{proof}

\cref{Modules_Quillen_Adjunction} gets us rather close to finishing the proof of \cref{eqv_ctv_RMod}: the fact that the right Quillen functor $N_{\mathcal{P}}$ preserves and detects weak equivalences makes us only have to check that, for any forest morphism $X \to N\mathcal{P}$ with $X$ a normal forest space, the counit map
$$ X \longrightarrow N_{\mathcal{P}} \tau_{\mathcal{P}}X$$
is a contravariant weak equivalence of forest spaces over $N\mathcal{P}$ (the condition for the unit map then follows from the triangle relations for the adjunction). The next lemma about the projective model structure on $\fSpaces$ allows us to reduce this problem to a question concerning representable presheaves.

\begin{lemma}
    Let $X$ be a forest space which is projectively cofibrant. Then $X$ is the retract of forest space $X'$ such that, for each $i \geq 0$, the forest set $X'_i$ is a coproduct of forests.
    \label{coprod_representables}
\end{lemma}

\begin{proof}
    This follows from \citep[Ex. 13.36]{HeutsMoerdijkDendroidal}.
\end{proof}

\begin{proof}[Proof of \cref{eqv_ctv_RMod}:]
  As we have already discussed, \cref{Modules_Quillen_Adjunction} implies that we only have to check that 
  $$ X \longrightarrow  N_{\mathcal{P}}\tau_{\mathcal{P}}X$$
  is a contravariant weak equivalence, whenever $X \to N\mathcal{P}$ has projectively cofibrant source. By \cref{coprod_representables}, we can assume that $X$ additionally satifies the condition that each $X_i$ is a coproduct of forests. Moreover, from \cref{bootstrap} it suffices to only check that the forest map
  $$ X_i \longrightarrow  N_{\mathcal{P}_i}\tau_{\mathcal{P}_i} X_i$$
  is a contravariant weak equivalence over $N\mathcal{P}_i$, for each value of $i$. This holds by \cref{counit_forest}, together with the fact that $N_{\mathcal{P}_i}$ preserves coproducts, as shown in \cref{corpoduct_nerve}.
\end{proof}

The first corollary we want to discuss is how this equivalence is not only one of model categories, but actually one between \textit{monoidal} model categories. For us, a Quillen equivalence of monoidal model categories corresponds to the notion of "weak monoidal Quillen equivalence" as discussed by Schwede--Shipley \citep[Def. 3.6]{SchwedeShiplyMonoidalEqv}.

\begin{corollary}
    Let $\mathcal{P}$ be a unital $\Sigma$-free simplicial operad. Then the adjunction
\[\begin{tikzcd}
	{\RMod{\mathcal{P}} } & {\left( \fSpaces_{/N\mathcal{P}} \right)_{\mathsf{P}, \ctv}} & {}
	\arrow["{N_{\mathcal{P}}}"', shift right, from=1-1, to=1-2]
	\arrow["{\tau_{\mathcal{P}}}"', shift right, from=1-2, to=1-1]
\end{tikzcd}\]
    defines a Quillen equivalence of monoidal model categories, where both categories are equipped with the respective concatenation product.
    \label{monoidal_eqv}
\end{corollary}

\begin{proof}
    We have just shown that under these circumstances the adjunction indeed defines a Quillen equivalence of model categories, so we just need to lift this to an equivalence further respecting the monoidal structure.

     We showed in \cref{NP_strong} that $N_{\mathcal{P}}$ is a strong monoidal functor, which in particular means that its left adjoint $\tau_\mathcal{P}$ is oplax monoidal. The last thing we need to check is that the structure map
     $$ \tau_\mathcal{P}(X \oplus Y) \longrightarrow \tau_\mathcal{P} X \boxplus \tau_\mathcal{P} Y$$
     is a weak equivalence of right $\mathcal{P}$-modules, whenever $X, Y \in \fSpaces_{N \mathcal{P}}$ are cofibrant. We can rewrite this map as described in the following commutative diagram  below
\[\begin{tikzcd}
	{\tau_\mathcal{P}(X \oplus Y)} && {\tau_\mathcal{P}X \boxplus \tau_\mathcal{P}Y} \\
	{\tau_\mathcal{P}(N_\mathcal{P} \tau_\mathcal{P}X \oplus N_\mathcal{P} \tau_\mathcal{P}Y)} && {\tau_\mathcal{P}N_\mathcal{P} (\tau_\mathcal{P}X \boxplus \tau_\mathcal{P}Y)}
	\arrow[from=1-1, to=1-3]
	\arrow[from=1-1, to=2-1]
	\arrow["\cong", from=2-1, to=2-3]
	\arrow[from=2-3, to=1-3]
\end{tikzcd}\]
     where the left and right morphisms are respectively constructed using the unit and counit of the adjunction $(\tau_\mathcal{P}, N_\mathcal{P})$, and the bottom map is an isomorphism since $N_\mathcal{P}$ is strong monoidal. As the counit is also always an isomorphism, it is clear that the right map is an isomorphism. The fact that the leftmost map is a weak equivalence is a consequence of $- \oplus -$ being a Quillen bifunctor (this is part of the proof of \cref{monoidal_without_ctv} and \cref{monoidal_with_ctv}), the unit of the nerve adjunction being a contravariant equivalence on cofibrant objects (this is part of the proof of \cref{eqv_ctv_RMod}), and $\tau_\mathcal{P}$ being left Quillen.
\end{proof}

The following important corollary allows us to compute mapping spaces of right modules in terms of right fibrations of forest spaces.

\begin{corollary}
    Let $\mathcal{P}$ be a closed simplicial $\Sigma$-free operad, and suppose $M$ and $L$ are simplicial right $\mathcal{P}$-modules which are projectively fibrant. Then the simplicial map between mapping spaces
    $$ \Map_{\mathcal{P}}(M,L) \longrightarrow \Map_{N \mathcal{P}}( N_\mathcal{P} M, N_{\mathcal{P}} L)$$
    is a weak homotopy equivalence. Here the left hand side is computed using the projective model structure on $\RMod{\mathcal{P}}$, and the right hand side using the projective contravariant model structure on $\fSpaces_{\mathcal{P}}$.
    \label{compute_map_spaces}
\end{corollary}

\begin{proof}
    This follows immediately from \cref{eqv_ctv_RMod}.
\end{proof}

\begin{remark}
    By \cref{right_fibr_fibrant}, we can also consider the rightmost mapping space in the Reedy contravariant model structure in the statement above, if we first replace $N\mathcal{P}$ by an equivalent Reedy fibrant object. 
    
    This result has a very important consequence for the computation of mapping spaces of right $\mathcal{P}$-modules: whereas the projective model structure on $\RMod{\mathcal{P}}$ has very few cofibrant objects, the cofibrant objects in the Reedy contravariant model structure are just the maps $X \to N\mathcal{P}$ with $X$ normal, which in our context is automatic, as $N\mathcal{P}$ is normal due to $\mathcal{P}$ being a $\Sigma$-free operad. Therefore, one should expect that computing the mapping spaces in the context of $\fSpaces_{/N\mathcal{P}}$ is a much more tractable problem, as the resolutions we have to take of the source and target will be smaller and easier to handle.
\end{remark}

Another easy consequence of the results just proven is that the homotopy theory of the category of right $\mathcal{P}$-modules only depends on the simplicial operad up to Dwyer-Kan equivalence: this is the notion of equivalence of simplicial operads explained in \citep[Section 14.3]{HeutsHinnichMoerdijk}, and essentially is the operadic version of fully faithful and essentially surjective maps of operads.

\begin{corollary}
    Let $\Phi \colon \mathcal{P} \to \mathcal{Q}$ be an operad map between closed simplicial $\Sigma$-free operads. If $\Phi$ is a Dwyer-Kan equivalence, then there is a Quillen equivalence
    \[\begin{tikzcd}
	{\RMod{\mathcal{P}}} & {\RMod{\mathcal{Q}}}
	\arrow["{\Phi^*}"', shift right, from=1-1, to=1-2]
	\arrow["{\Phi_!}"', shift right, from=1-2, to=1-1]
\end{tikzcd}\]
between the respective projective model structures.
\label{DK_eqv}
\end{corollary}

\begin{proof}
    Firstly, it is clear that $\Phi^*$ preserves (trivial) fibrations since these are given levelwise, so that the pair in question defines a Quillen adjunction. We can expand the given adjunction to the following square of Quillen pairs
\[\begin{tikzcd}[cramped]
	{\RMod{\mathcal{P}}} & {\RMod{\mathcal{Q}}} \\
	{\left( \fSpaces_{/N\mathcal{P}}\right)_{\mathsf{P}, \ctv}} & {\left( \fSpaces_{/N\mathcal{Q}}\right)_{\mathsf{P}, \ctv}.}
	\arrow["{\Phi_!}", shift left, from=1-1, to=1-2]
	\arrow["{N_{\mathcal{P}}}", from=1-1, to=2-1]
	\arrow["{\Phi^*}", shift left, from=1-2, to=1-1]
	\arrow["{N_\mathcal{Q}}"', shift right, from=1-2, to=2-2]
	\arrow["{\tau_{\mathcal{P}}}", shift left=2, from=2-1, to=1-1]
	\arrow["{(N\Phi)_!}", shift left=2, from=2-1, to=2-2]
	\arrow["{\tau_\mathcal{Q}}"', shift right, from=2-2, to=1-2]
	\arrow["{(N\Phi)^*}", from=2-2, to=2-1]
\end{tikzcd}\]

By \citep[Lem. 14.23]{HeutsHinnichMoerdijk}, $N\Phi$ is a complete weak equivalence between complete Segal dendroidal spaces, and therefore the same holds for $u_\ast N \Phi$. Using \cref{eqv_ctv_RMod} together with \cref{complete_weak_equivalence} shows that the left, right and bottom adjunctions are Quillen equivalences, which implies the desired conclusion.
\end{proof}

We end this section with a strictification result for operadic right modules, based on a similar result for dendroidal left fibrations and simplicial algebras present in \citep[Thm. 14.44]{HeutsMoerdijkDendroidal}. We note that its proof will involve some results about simplicial operads which we have not discussed here, but the interested reader should look at the reference for the results we are leaving out. We will spell out the main steps of the proof in any case.

For the statement of the result, we note that there exists a version of the homotopy coherent nerve $w^*$ for simplicial operads, which takes part in an adjunction
\[\begin{tikzcd}
	{\mathsf{sOp}} & {\mathsf{dSpaces}.}
	\arrow["{w^*}"', shift right, from=1-1, to=1-2]
	\arrow["{w_!}"', shift right, from=1-2, to=1-1]
\end{tikzcd}\]
For more details, see Section 14.6 of \cite{HeutsMoerdijkDendroidal}.
\begin{theorem}
    Let $V$ be a normal closed $\infty$-operad. Then there exists a natural zig-zag of Quillen equivalences
\[\begin{tikzcd}[cramped]
	{\left( \fSets_{/V}\right)_{\ctv}} & \cdots & {\mathsf{RMod}_{w!(V)}.}
	\arrow[shift right, from=1-1, to=1-2]
	\arrow[shift right, from=1-2, to=1-1]
	\arrow[shift right, from=1-2, to=1-3]
	\arrow[shift right, from=1-3, to=1-2]
\end{tikzcd}\]
\label{strictification}
\end{theorem}

\begin{proof}
    The first step of the chain of Quillen equivalences is of the form
    \[\begin{tikzcd}
	{\left( \fSets_{/V} \right)_{\ctv}} & {\left( \fSpaces_{/X} \right)_{\mathsf{R}, \ctv},}
 \arrow[shift right, from=1-1, to=1-2]
	\arrow[shift right, from=1-2, to=1-1]
\end{tikzcd}\]
which we can achieve by choosing a complete Segal forest space $X$ such that $X_0= u_\ast V$: this can be done by first using the complete Segal dendroidal space $X'$ constructed in Example 12.24 of \cite{HeutsMoerdijkDendroidal} and setting $X = u_\ast X'$. The fact that this is a Quillen equivalence is then a consequence of \cref{dsicrete_space}.

The next two equivalences take the form
\begin{equation}
\begin{tikzcd}[cramped]
	{\left( \fSpaces_{/X}\right)_{\mathsf{R}, \ctv}} & {\left( \fSpaces_{/Y}\right)_{\mathsf{P},\ctv}} & {\left( \fSpaces_{/N\tau(Y')}\right)_{\mathsf{P},\ctv}.}
	\arrow[shift right, from=1-1, to=1-2]
	\arrow[shift right, from=1-2, to=1-1]
	\arrow[shift right, from=1-2, to=1-3]
	\arrow[shift right, from=1-3, to=1-2]
\end{tikzcd}
\label{Qeqv}
\end{equation}
for a certain dendroidal space $Y'$ and forest space $Y$ which we now explain. Firstly, construct a Segal dendroidal space $Y'$ such that the following two conditions are satisfied:
\begin{enumerate}[label=(\alph*)]
    \item The unit map $Y' \to N\tau(Y')$ is a complete weak equivalence of dendroidal spaces. This can be achieved by setting $Y'$ to be a cofibrant replacement of $X'$ in the sparse model structure of \citep[Thm. 14.17]{HeutsMoerdijkDendroidal} and applying \citep[Lem. 14.24]{HeutsMoerdijkDendroidal}. We can also ensure that $Y'$ is Segal by \citep[Lem. 14.16]{HeutsMoerdijkDendroidal}.
    \item There exists a complete weak equivalence $X'\to Y'$. As the sparse model structure as the complete weak equivalences as its weak equivalences, this follows from part (a).
    \item There is a (zig-zag) of Quillen equivalences between the simplicial operads $\tau (Y')$ and $w_!(V)$. This is \citep[Cor. 14.30]{HeutsMoerdijkDendroidal}.
\end{enumerate}

Having chosen this, set $Y = u_\ast Y'$. Then the leftmost adjunction in \eqref{Qeqv} is an equivalence by \cref{complete_weak_equivalence} as $X$ and $Y$ are Segal forest spaces and by property (b); as everything is Reedy fibrant, we are also free to change to the projective model structure. The second equivalence is then another application of \cref{complete_weak_equivalence}, together with property (a).

Finally, we claim there are Quillen equivalences
\begin{equation}
\begin{tikzcd}[cramped]
	{\left( \fSpaces_{/N\tau(Y)}\right)_{\mathsf{P},\ctv}} & {\mathsf{RMod}_{\tau(Y)}} & {\mathsf{RMod}_{w_!(V)}.}
	\arrow[shift right, from=1-1, to=1-2]
	\arrow[shift right, from=1-2, to=1-1]
	\arrow[shift right, from=1-2, to=1-3]
	\arrow[shift right, from=1-3, to=1-2]
\end{tikzcd}
\label{Qeqv2}
\end{equation}

In order to show this, we require that $Y'$ satisfies the following additional property:
\begin{enumerate}
    \item[(c)] There exists a (zig-zag) of Dwyer--Kan equivalences between the simplicial operads $\tau (Y')$ and $w_!(V)$. This is guaranteed by\citep[Cor. 14.30]{HeutsMoerdijkDendroidal}.
\end{enumerate}

The first adjunction is an equivalence by \cref{eqv_ctv_RMod}, and the last one by property (c) and \cref{DK_eqv}.
\end{proof}

\section{A dendroidal approach to Goodwillie--Weiss manifold calculus}

In this final chapter of the article we aim at showing how our formalism of right fibrations of forest spaces can be effectively used to study the tower for smooth embeddings studied by Goodwillie and Weiss. We construct this tower using only a certain filtration on the category of forests and show that it is indeed equivalent to the usual embedding tower. The rest of this part of the article is then dedicated to showing how one can compute the layers of the tower using a latching-matching description reminiscent of the one in \cite{GopplWeiss}. We finish with some connectivity estimates obtained via this description.

\subsection{The operadic tower for right modules}

Recall that we have associated to every simplicial operad $\mathcal{P}$ its symmetric monoidal envelope $\mathsf{Env}(\mathcal{P})$, which is a symmetric monoidal category encoding the operadic composition internal to the operad $\mathcal{P}$ via its own internal composition. Having this category at hand, we further constructed the category $\RMod{\mathcal{P}}$ of right $\mathcal{P}$-modules as the category of simplicial presheaves on $\mathsf{Env}(\mathcal{P})$. 

The aim of this section is to study a filtration on the category of right $\mathcal{P}$-modules
$$ \RMod{\mathcal{P}}^{\smallleq 0} \subseteq \RMod{\mathcal{P}}^{\smallleq 1} \subseteq \cdots \subseteq \RMod{\mathcal{P}}^{\smallleq k} \subseteq \cdots \subseteq \RMod{\mathcal{P}}$$
which, when restricted to the context of the modules $\mathbb{E}_M$ from embedding theory, recovers the Goodwillie--Weiss tower. We will build this filtration by constructing it first at the level of the envelope category $\mathsf{Env}(\mathcal{P})$.

Let $\mathcal{P}$ be a simplicial operad. Notice that it admits a unique map to the commutative operad $\mathsf{Com}$, as the latter defines the terminal operad in $\mathsf{sOp}$. Recalling from \cref{examples_right_modules} (b) that $\mathsf{Env}(\mathsf{Com})$ is the category of finite sets $\Fin$, we see that there exists a natural functor $\mathsf{Env}(\mathcal{P}) \to \Fin$, which ignores the information of $\mathsf{Env}(\mathcal{P})$ coming from the operad $\mathcal{P}$. We will filter $\mathsf{Env}(\mathcal{P})$ via the cardinality filtration on $\Fin$.

\begin{definition}
    Let $k \geq 0$ and write $\Fin^{\smallleq k} \subseteq \Fin$ for the full subcategory on finite sets of cardinality at most $k$. We set $\mathsf{Env}(\mathcal{P})^{\smallleq k} \subseteq \mathsf{Env}(\mathcal{P})$ to be the full subcategory determined by the pullback diagram
\[\begin{tikzcd}[cramped]
	{\mathsf{Env}(\mathcal{P})^{\smallleq k}} & {\,\mathsf{Env}(\mathcal{P})} \\
	{\Fin^{\smallleq k}} & {\Fin}
	\arrow[from=1-1, to=1-2]
	\arrow[from=1-1, to=2-1]
	\arrow[from=1-2, to=2-2]
	\arrow[from=2-1, to=2-2]
    \arrow["\lrcorner"{anchor=center, pos=0, scale=1.5}, draw=none, from=1-1, to=2-2]
\end{tikzcd}\]
in $\mathsf{sCat}$. We define the category of \textit{$k$-truncated right $\mathcal{P}$-modules} $\RMod{\mathcal{P}}^{\smallleq k}$ to be the category of simplicial presheaves on $\mathsf{Env}(\mathcal{P})^{\smallleq k}$. 
\label{filt_def}
\end{definition}

\begin{remark}
    We can give a more explicit description of the category $\mathsf{Env}(\mathcal{P})^{\smallleq k}$. Recall that a morphism $(I, \alpha) \to (J, \beta)$ is the data of a function between finite sets $f \colon I \to J$, together with operations
$$ p_j \in \mathcal{P}\left( f^{-1}(\alpha(j)); \beta(j) \right)$$
for each $j \in J$. The truncation we have imposed corresponds to the cardinality restriction $\lvert I \rvert, \lvert J \rvert \leq k$, which in particular also limits the allowed arities of $p_j$, since the \textit{sum} of the arities of all the $p_j$'s must be at most $k$ also.
\end{remark}

\begin{remark}
The category $\mathsf{Env}(\mathcal{P})^{\smallleq 0}$ is clearly just the trivial category; therefore from now on we will ignore this stage of the filtration and always assume that it starts at $\mathsf{Env}(\mathcal{P})^{\smallleq 1}$. 

We can offer a different perspective on this via a comment similar to \cref{empty_forest}: given $M \in \RMod{\mathcal{P}}$, we can always change $M(0)$ to be $\Delta^0$ -- here $M(0)$ is the space coming from evaluation on the empty set -- and this still defines a right module. We assume from now on that this is always the case.
\label{filtration_zero}
\end{remark}

\begin{remark}
We note that the first stage of the filtration $\mathsf{Env}(\mathcal{P}^{\smallleq 1})$ is given by the underlying category $\mathsf{Cat}(\mathcal{P})$ of the given operad. Note that this is different from $\mathsf{Env}(\mathsf{Cat}(\mathcal{P}))$. 

For $k \geq 2$, a similar succinct expression for $\mathsf{Env}(\mathcal{P})^{\smallleq k}$ seems difficult. However, we can say that in general $\mathsf{Env}\left(\mathcal{P}\right)^{\smallleq k}$ is different from $\mathsf{Env}\left( \mathcal{P}^{\smallleq k} \right)$, where $\mathcal{P}^{\smallleq k}$ denotes the truncation of $\mathcal{P}$ at the operations of arity at most $k$. Indeed, we have the following descriptions:
\begin{itemize}
    \item Presheaves on $\mathsf{Env}\left( \mathcal{P}\right)^{\smallleq k}$ will correspond to considering the truncated sequence $\{ M(\ell) : 0 \leq \ell \leq k\}$, and its action by the symmetric sequence $\{ \mathcal{P}(\ell) : \ell \geq 0 \}$. In particular, we only allow for partial composition maps
    \begin{equation}
     M(\ell) \times \mathcal{P}(n) \xlongrightarrow{\cdot_i} M(\ell+n-1)
     \label{action_eq}
    \end{equation}
    satisfying $\ell+n-1 \leq k$.
    \item Presheaves on $\mathsf{Env}\left( \mathcal{P}^{\smallleq k} \right)$ will correspond to considering the truncated sequence $\{ M(\ell) : 0 \leq \ell \leq k\}$, and its action by the \textit{truncated}  sequence $\{ \mathcal{P}(\ell) : 0 \leq \ell \leq k \}$. We still have the same action maps \eqref{action_eq} as above, but the difference is that now we are also filtering the operad $\mathcal{P}$.
\end{itemize}
The main point from these remarks is that, if we had considered $\mathsf{Env}\left( \mathcal{P}^{\smallleq k} \right)$, we would have an additional filtration on $\mathcal{P}$ which is unnecessary for the study of right $\mathcal{P}$-modules.
\label{truncated_operads}
\end{remark}

The truncated variants of the symmetric monoidal envelope assemble themselves into a filtration
$$ \mathsf{Cat}(\mathcal{P}) = \mathsf{Env}(\mathcal{P})^{\smallleq 1} \subseteq \mathsf{Env}(\mathcal{P})^{\smallleq 2} \subseteq \cdots \subseteq \mathsf{Env}(\mathcal{P})^{\smallleq k} \subseteq \cdots \subseteq \mathsf{Env}(\mathcal{P})$$
with corresponding filtered colimit recovering $\mathsf{Env}(\mathcal{P})$. At the level of operadic modules, we have a cofiltration given by the restriction functors
\[\begin{tikzcd}[cramped]
	{\RMod{\mathcal{P}}} & \cdots & {\RMod{\mathcal{P}}^{\smallleq k+1}} & {\RMod{\mathcal{P}}^{\smallleq k}} & \cdots & {\RMod{\mathcal{P}}^{\smallleq 1}} 
	\arrow[from=1-1, to=1-2]
	\arrow["{J^*_{k+1}}", from=1-2, to=1-3]
	\arrow["{J^*_k}", from=1-3, to=1-4]
	\arrow["{J^*_{k-1}}", from=1-4, to=1-5]
	\arrow["{J^*_1}", from=1-5, to=1-6]
\end{tikzcd}\]
induced by the inclusions $J_k \colon \mathsf{Env}(\mathcal{P})^{\smallleq k} \to \mathsf{Env}(\mathcal{P})^{\smallleq k+1}$. Finally, we also get functors $U_k^* \colon \RMod{\mathcal{P}} \to \RMod{\mathcal{P}}^{\smallleq k}$ coming from the inclusions $U_k \colon \mathsf{Env}(\mathcal{P})^{\smallleq k} \to \mathsf{Env}(\mathcal{P})$.

As we have previously considered, we will always endow $\RMod{\mathcal{P}}$ with the projective model structure, and therefore we will put the same model structure on $\RMod{\mathcal{P}}^{\smallleq k}$, as these are still categories of presheaves. Moreover, we will also consider the adjoints to $U_k^*$
$$
\begin{tikzcd}
\RMod{\mathcal{P}}^{\smallleq k} \arrow[rr, "(U_k)_\ast", bend right] \arrow[rr, "(U_k)_!", bend left] &  & \RMod{\mathcal{P}} \arrow[ll, "U_k^*"']
\end{tikzcd}
$$
constructed via left and right Kan extension. From the fact that $(U_k)_!$ will preserve representables we can conclude that $\left( (U_k)_!, U_k^* \right)$ is a Quillen adjunction, but it is not clear whether $\left( U_k^*, (U_k)_\ast \right)$ is also a Quillen pair because it is not clear that $U_k^*$ will preserve projective cofibrations. In any case, since the truncation functor $U_k^*$ clearly preserves all levelwise weak equivalences, it is already a derived functor by \citep[Prop. 8.4.8]{HirschhornLocalizations}. Similar considerations also hold for $J_k^*$ of course. 

In the next result, we will show that $(U_k)_*$ can be easily described as a limit over a cubical diagram.

\begin{notation}
Before the next few results, let us fix some notation:
\begin{itemize}
    \item For a right $\mathcal{P}$-module $M$, we will write $M^{\smallleq k}  \in \RMod{\mathcal{P}}^{\smallleq k}$ for its $k$-truncation $U_k^*M$.
    \item If $k \geq 1$ and $A$ is a finite set, we will write $\mathsf{P}(A; k)$ for the poset
        $$ \mathsf{P}(A;k) = \left\{ I \subseteq A \: : \:  \lvert I \rvert \leq k\right\}$$
    of subsets of $A$ of cardinality at most $k$, ordered using the set containment.
\end{itemize}
 \label{notation_labels_modules}
\end{notation}

Suppose now that $\mathcal{P}$ is a unital simplicial operad, and fix $(I, \alpha) \in \mathsf{Env}(\mathcal{P})$, which we also see as a representable presheaf. Then, given $k \geq 1$, there exists a functor 
$$ \Psi : \mathsf{P}(I; k) \longrightarrow \mathsf{Env}(\mathcal{P})^{\smallleq k}_{/ (I, \alpha)}$$
 sending $I'\subseteq I$ to the map $(I', \alpha') \to (I, \alpha)$ indexed by the inclusion of sets $I' \subseteq I$ together with the operadic identities for the colours of $I'$, and the constants for the colours in the complement of $I'$ (this is why we require $\mathcal{P}$ to be unital). We will make use of this functor in the next lemma.
 
\begin{lemma}
    Let $k \geq 1$ and $\mathcal{P}$ a unital simplicial operad. Then the following statements hold:
    \begin{enumerate}[label=(\roman*)]
        \item For any representable $(I, \alpha) \in \mathsf{Env}(\mathcal{P})$, the $k$-truncated right module $(I, \alpha)^{\smallleq k} \in \RMod{\mathcal{P}}^{\smallleq k}$ is given as a colimit of representables
        $$(I, \alpha)^{\smallleq k} = \colim_{I' \in \mathsf{P}(I; k)} \: (I', \alpha'),$$
        where $\alpha'$ is the restriction of $\alpha$ to $I'$. 
        \item If $M \in \RMod{\mathcal{P}}^{\smallleq k}$, then the right $\mathcal{P}$-module $(U_k)_\ast M$ is given by
        $$ \left( (U_k)_\ast M \right)(I, \alpha) = \lim_{I' \in \mathsf{P}(I; k)^{\mathsf{op}}} \: M(I', \alpha')$$
        when evaluated at any $(I, \alpha) \in \mathsf{Env}(\mathcal{P})$. 
    \end{enumerate}
    \label{formula_Uk}
\end{lemma}

\begin{proof}
    Statement (ii) is clear once one shows (i), since $(U_k)_\ast$ is right adjoint to $(U_k)^*$. For (i), we can write $(I, \alpha)^{\smallleq k}$ as a colimit of representables, leading to the formula.
    $$(I, \alpha)^{\smallleq k} = \colim_{(J, \beta)} \: (J, \beta),$$
    where the indexing category for the colimit is $\mathsf{Env}(\mathcal{P})^{\smallleq k}_{/ (I, \alpha)}$. Consider the functor
    $$ \Psi : \mathsf{P}(I; k) \longrightarrow \mathsf{Env}(\mathcal{P})^{\smallleq k}_{/ (I, \alpha)}$$
     that we introduced before the lemma. We claim that $\Psi$ is right cofinal, that is, the slice categories $\Psi_{(J, \beta) \to (I, \alpha)/}$ are non-empty and connected, which shows (i).

   In order to check that it is non-empty, we need to verify that the left arrow in the diagram
\[\begin{tikzcd}[cramped]
	{(J, \beta)} & {(I',\alpha')} \\
	{(I, \alpha)}
	\arrow["{(g, q_\ast)}", dashed, from=1-1, to=1-2]
	\arrow["{(f, p_\ast)}"', from=1-1, to=2-1]
	\arrow["{\Psi(I')}", dashed, from=1-2, to=2-1]
\end{tikzcd}\]
   with $\lvert J \rvert \leq k$ admits the dashed factorization. Letting $J \xrightarrow{g} I' \xrightarrow{\iota} I$ be the surjective-injective factorization of $f$, we note that $I' \in \mathsf{P}(I; k)$. Setting $q_\ast$ to be the operations of $p_\ast$ which are not constants, we can define a morphism $(g, q_\ast) \colon (J, \beta) \to (I', \alpha')$ which factors $(f, p_\ast)$ via a map in the image of $\Psi$.
   
   To check that the category is connected, it suffices to show that the solid diagram over $(I, \alpha)$ below
\[\begin{tikzcd}[cramped]
	(J, \beta) & {(I', \alpha')} \\
	{(I'', \alpha'')} & {(I' \cap I'', \alpha' \cap \alpha'' )}
	\arrow["{(f, p_\ast)}", from=1-1, to=1-2]
	\arrow["{(g, q_\ast)}"', from=1-1, to=2-1]
	\arrow["{(h, r_\ast)}", dashed, from=1-1, to=2-2]
	\arrow[dashed, from=2-2, to=1-2]
	\arrow[dashed, from=2-2, to=2-1]
\end{tikzcd}\]
can be completed via the dashed arrows, with the diagram with source $(I' \cap I'', \alpha' \cap \alpha'' )$ being in the image of $\Psi$. The existence of the left and top morphisms imply that $p_i = q_i$ for $i \in I' \cap I''$ and otherwise all other operations are constants. Therefore the factorization $(h, r_\ast) \colon (J, \beta) \to (I' \cap I'', \alpha' \cap \alpha'' )$ exists, as we wanted to show.
\end{proof}

For the sake of convenience, we gather below all the previous homotopical considerations about the operadic filtration in the following proposition.

\begin{proposition}
    Suppose $k \geq 1$ and let $\mathcal{P}$ be a simplicial operad.  Consider the following diagram of adjunctions
\[\begin{tikzcd}[cramped]
	{\RMod{\mathcal{P}}} & {\RMod{\mathcal{P}}^{\smallleq k}} & {\RMod{\mathcal{P}}^{\smallleq k+1}}
	\arrow["{U_k^*}", from=1-1, to=1-2]
	\arrow["{(U_k)_!}"', bend right=60, from=1-2, to=1-1]
	\arrow["{(U_k)_\ast}", bend left=40, from=1-2, to=1-1]
	\arrow["{(J_k)_!}", bend left=55, from=1-2, to=1-3]
	\arrow["{(J_k)_*}"', bend right=40, from=1-2, to=1-3]
	\arrow["{J_k^*}"', from=1-3, to=1-2]
\end{tikzcd}\]
    induced from $U_k$ and $J_k$. Then:
    \begin{enumerate}[label=(\roman*)]
        \item The pairs $((U_k)_!, U_k^*)$ and $((J_k)_!, J_k^*)$ define Quillen adjunctions.
        \item The restrictions functors $U_k^*$ and $J_k^*$ preserve projective weak equivalences.
    \end{enumerate}

    \label{Quillen_adj_projective_modules}
\end{proposition}

\begin{remark}
    The restriction functors $U_k^*$ and $J_k^*$ will not preserve projective cofibrations, which is a general phenomenon of the restriction functor of presheaf categories. This is also not true if we restrict ourselves to unital operads.
\end{remark}

Let us now discuss mapping spaces. Since the restriction functors preserve projective weak equivalences, these are already both left and right derived. Therefore, given right $\mathcal{P}$-modules $M, L \in \RMod{\mathcal{P}}$, we get morphisms between the respective derived mapping spaces
$$ \Map_{\mathcal{P}}(M, L) \xlongrightarrow{U_{k+1}^*} \Map_{\mathcal{P}}^{\smallleq k+1}(M^{\smallleq k+1}, L^{\smallleq k+1}) \xlongrightarrow{J_k^*} \Map^{\smallleq k}_{\mathcal{P}}(M^{\smallleq k}, L^{\smallleq k})$$
with composite given by $U_k^*$. Here we have written $\Map^{\smallleq k}_{\mathcal{P}}(-, -)$ for the derived mapping space computed in $\RMod{\mathcal{P}}^{\smallleq k}$, and from now one we will write $M$ instead of $M^{\smallleq k}$ when considering mapping spaces. 

Therefore these restrictions give rise to a tower of mapping spaces, and the convergence properties of this tower are clarified in the next result.

\begin{proposition}
    Suppose $\mathcal{P}$ is a simplicial operad, and let $M, L$ be operadic right $\mathcal{P}$-modules. Then the $k$-truncated categories of right modules $\RMod{\mathcal{P}}^{\smallleq k}$ induce a tower of derived mapping spaces
\[\begin{tikzcd}[cramped]
	&& \vdots \\
	&& {\Map_{\mathcal{P}}^{\smallleq 3}(M, L)} \\
	&& {\Map_{\mathcal{P}}^{\smallleq 2}(M, L)} \\
	{\Map_{\mathcal{P}}(M,L)} && {\Map_{\mathcal{P}}^{\smallleq 1}(M, L)}
	\arrow[from=1-3, to=2-3]
	\arrow[from=2-3, to=3-3]
	\arrow[from=3-3, to=4-3]
	\arrow[bend left=10, from=4-1, to=2-3]
	\arrow[from=4-1, to=3-3]
	\arrow[from=4-1, to=4-3]
\end{tikzcd}\]
which converges to  $\Map_{\mathcal{P}}(M,L)$, that is, the map to the homotopy limit of the tower
$$ \Map_{\mathcal{P}}(M,L) \longrightarrow \holim_k \Map_{\mathcal{P}}^{\smallleq k}(M, L)$$
is a weak homotopy equivalence.
\label{operadic_tower}
\end{proposition}

\begin{proof}
    The only aspect which we haven't mentioned yet concerns the convergence of the tower. As the filtered colimit of the full subcategories $\mathsf{Env}(\mathcal{P})^{\smallleq k} \subseteq \mathsf{Env}(\mathcal{P})$ is exactly $\mathsf{Env}(\mathcal{P})$, the convergence is ensured by Lemma 3.3.1 of \cite{GopplWeiss}.
\end{proof}

In the rest of this chapter we want to develop a careful study of this tower for any simplicial unital $\Sigma$-free operad $\mathcal{P}$ and any right modules $M, N$ over it. Recall that this general setting has the Goodwillie--Weiss manifold calculus as a particular instance, and we will address this specific case in in the last section of this article.

\begin{remark}
    Although we have said that manifold calculus is an example of the scenario that we are considering, we observe that the convergence of our tower is not equivalent to the convergence of the Goodwillie--Weiss tower, which in fact doesn't converge in general. The reason is that the term $\mathsf{Map}_{\mathbb{E}_d}(\mathbb{E}_M, \mathbb{E}_N)$ is the limit of the embedding tower, and not the embedding space $\mathsf{Emb}(M, N)$ as it appears in manifold calculus.
\end{remark}

The standard strategy for analysing a tower of spaces in homotopy theory usually goes through three steps, namely:
\begin{itemize}
    \item \textit{A description of the first stage of the tower}: for us, this would mean having a handle on
    $\Map_{\mathcal{P}}^{\smallleq 1}(M, L)$ for any right $\mathcal{P}$-modules $M$ and $L$.
    \item \textit{A description of the layers of the tower}: that is, we need to compute the homotopy fibers of the connecting maps
    $$ \Map_{\mathcal{P}}^{\smallleq k}(M, L) \longrightarrow \Map_{\mathcal{P}}^{\smallleq k-1}(M, L)$$
    for all $k \geq 2$.
    \item \textit{A discussion of the convergence of the tower}: in our situation this is part of the content of \cref{operadic_tower}.
\end{itemize}

The problem concerning the layers of the tower will be the main focus of the remaining sections since the flexibility of the category of forest spaces and the different model structures on it will be crucial. However, we can already address the problem concerning the first stage.

\begin{proposition}
    Let $\mathcal{P}$ be a simplicial operad. Then there is a Quillen equivalence
    $$\RMod{\mathcal{P}}^{\smallleq 1} \simeq \mathsf{Psh}(\mathsf{Cat}(\mathcal{P}))$$
    between the category of 1-truncated $\mathcal{P}$-modules, and the category of simplicial presheaves on $\mathsf{Cat}(\mathcal{P})$, both equipped with the projective model structure.
    \label{first_stage}
\end{proposition}

\begin{proof}
    This is clear from the identification $\mathsf{Env}(\mathcal{P})^{\smallleq 1} \simeq \mathsf{Cat}(\mathcal{P})$ that we observed in \cref{truncated_operads}.
\end{proof}

In general $\mathsf{Cat}(\mathcal{P})$ is quite explicit and manageable: for instance, for uncoloured operads it will be a category with a single object, so it is often a group $G$. Supposing that we are in the latter, then we are left with computing a derived mapping space in the category of $G$-spaces, which is feasible. We give below evidence of this in the context of embedding calculus. 

\begin{example}
    Let $d, k \geq 1$, and consider manifolds $M^d$ and $N^{d+k}$. We note that $\mathbb{E}_M$ and $\mathbb{E}_N$ can be seen as modules over the $\mathbb{E}^\mathsf{fr}_d$-operad: for the latter, we consider the pullback of $\mathbb{E}_N \in \mathsf{RMod}_{\mathbb{E}^\mathsf{fr}_{d+k}}$ along the inclusion $\mathbb{E}^\mathsf{fr}_d \to \mathbb{E}^\mathsf{fr}_{d+k}$ coming from the topological inclusion $D^d \to D^{d+k}$. We will still refer to this $\mathbb{E}^\mathsf{fr}_d$-module as $\mathbb{E}_N$.

    Firstly, we note that the presheaf $\mathbb{E}_N^{\smallleq 1}$ is determined by the space $\mathbb{E}_N(1) \simeq \mathrm{Emb}(D^d, N)$, together with an action via precomposition by $\mathrm{Emb}(D^d, D^d)$. It is well-known that there is a weak homotopy equivalence
    $$ O(d) \xlongrightarrow{\simeq} \mathrm{Emb}(D^d, D^d)$$
    via the inclusion of the orthogonal group. This allows us to replace the action by $\mathrm{Emb}(D^d, D^d)$ by the $O(d)$-action, and therefore we have the chain of equivalences below:
    \begin{align*}
        \mathsf{Map}^{\smallleq 1}_{\mathbb{E}^\mathsf{fr}_d}(\mathbb{E}_M, \mathbb{E}_N) &\simeq \mathsf{Map}^{O(d)}\left( \mathrm{Emb}(D^d, M) , \mathrm{Emb}(D^d, N)\right) \\
        &\simeq \mathsf{Map}^{O(d)}\left( \mathrm{Fr}_d(M), \mathrm{Fr}_d(N) \right). 
    \end{align*}
    Here we write $\mathsf{Map}^{O(d)}(-,-)$ for the mapping space for $O(d)$-spaces. In the last step we have used that the derivative at the center of a disk is an $O(d)$-equivariant weak homotopy equivalence
    $$ \mathrm{Emb}(D^d, N^{d+k}) \xlongrightarrow{\simeq} \mathrm{Fr}_d(N)$$
    to the total space of the principal $O(d)$-bundle of $d$-frames in $N$. The evaluation of the differential of an immersion induces a map between the respective framed bundles, thus leading to a map
    \begin{equation}
    \mathrm{Imm}(M, N) \longrightarrow \mathsf{Map}^{O(d)}\left( \mathrm{Fr}_d(M), \mathrm{Fr}_d(N) \right).
    \label{smale_hirsch}
    \end{equation}
    The Smale--Hirsch theorem \citep[Thm. 1.3]{WeissImmersionTheory} states that \eqref{smale_hirsch} is a weak homotopy equivalence whenever $M$ is compact with $\dim M < \dim N$, or $M$ is open with no dimension restrictions.

    \label{T1_example}
\end{example}

\subsection{The forest tower for right fibrations}

In the previous section we defined a filtration on the symmetric monoidal envelope $\mathsf{Env}(\mathcal{P})$ of a simplicial operad $\mathcal{P}$, and saw how this consequently leads to a tower of spaces approximating the derived mapping space $\mathsf{Map}_{\mathcal{P}}(M,L)$ for any right $\mathcal{P}$-modules. Although we were able to prove convergence of said tower in \cref{operadic_tower}, as well as give a convenient description of the first stage in \cref{first_stage}, we still haven't addressed the question of computing the layers. 

We claim that the computation of the layers of the operadic tower is not as easy to carry out if one adopts the framework of right $\mathcal{P}$-modules as presheaves on $\mathsf{Env}(\mathcal{P})$. Such a computation essentially entails studying how these presheaves extend along the inclusion
\begin{equation}
\mathsf{Env}(\mathcal{P})^{\smallleq k-1} \longrightarrow \mathsf{Env}(\mathcal{P})^{\smallleq k}
\label{associated_graded_env}
\end{equation}
and proving that, up to homotopy, only a small amount of information about $\mathsf{Env}(\mathcal{P})^{\smallleq k}$ needs to be captured by the layer in question. The main obstruction in following up with this strategy is one that we had already encountered before, namely that the $\mathsf{Env}(\mathcal{P})$ mixes up the combinatorics of the forest category $\For$ with the operations of $\mathcal{P}$.

With these remarks in mind, in this section we aim to produce a filtration $\FilFor{k}$ for the forest category $\For$ which parallels the operadic one, and compare the towers obtained via these procedures. In Section 4.3 we will come back to the question about the layers, where we will see that the Segal condition will reduce the analysis of the inclusion $\For^{\smallleq k-1} \subseteq \For^{\smallleq k}$ to that of understanding the forest of edges $k \cdot \eta$. For the moment, we will begin with introducing the filtration we will consider.

\begin{definition}
    Let $F$ be a forest. We say that $F$ lies in \textit{filtration $k$} if it has exactly $k$ tree components, that is, $F$ is of the form $\bigoplus_{i=1}^k F_i$ where each $F_i$ is a tree. Given $k \geq 0$, we set the \textit{category of $k$-forests} $\FilFor{k} \subseteq \For$ to be the full subcategory spanned by the forests of filtration at most $k$.
\end{definition}

\begin{remark}
    The only forest with filtration 0 is the empty forest, and therefore $\For^{\smallleq 0}$ is the trivial category. The observations we presented in \cref{empty_forest} and \cref{filtration_zero} also apply here, and therefore we will from now on assume that all forests have positive filtration.
\end{remark}

\begin{example}
    As the forests with exactly one component are the trees, we conclude that $\For^{\smallleq 1}$ coincides with the dendroidal category $\Tree$.
\end{example}

\begin{remark}
Before continuing, let us briefly comment on another reasonable filtration one could study on $\For$. This filtration originates from the work of Göppl and Weiss in \cite{GopplWeiss} on computing derived mapping spaces $\Map(\mathcal{P}, \mathcal{Q})$ between simplicial operads $\mathcal{P}, \mathcal{Q}$. They approach this question via dendroidal spaces by using a filtration on $\Tree$
$$ \Tree^{\smallleq 0} \subseteq \Tree^{\smallleq 1} \subseteq \cdots \subseteq \Tree^{\smallleq k} \subseteq \cdots \subseteq \Tree,$$
where $T \in \Tree^{\smallleq k}$ if the vertices of $T$ have at most $k$ incoming edges. From the operadic point of view, this will mimic the arity filtration of an operad. We can extend this to a filtration $\mathsf{fil}^{\mathsf{GW}} \colon \For \to \mathbb{Z}$ via the rule 
$$\mathsf{fil}^{\mathsf{GW}}(F \oplus G) = \mathsf{max} \left( \mathsf{fil}^{\mathsf{GW}}(F), \mathsf{fil}^{\mathsf{GW}}(G) \right)$$
for any $F, G \in \For$. In other words, we have that $\mathsf{fil}^{\mathsf{GW}}(F) \leq k$ if and only if the vertices of the forest $F$ have at most $k$ incoming edges.

Notice that the Göppl--Weiss filtration is actually closer to the filtration on $\mathsf{Env}(\mathcal{P})$ that we considered in \cref{filt_def}. Indeed, we have previously seen that the following identifications of the first stage
$$ \mathsf{Env}(\mathcal{P})^{\smallleq 1} = \mathsf{Cat}(\mathcal{P}) \hspace{2em} \mathrm{and} \hspace{2em} \For^{\smallleq 1} = \Tree$$
which capture different kinds of information about the operad. On the other hand, it it clear that $\mathsf{fil}^{\mathsf{GW}}(F) \leq 1$ if $F$ is a linear forest with possibly some stumps, which is the type of data captured by $\mathsf{Cat}(\mathcal{P}) = \mathsf{Env}(\mathcal{P})^{\smallleq 1}$. One would be naturally led to conclude that we have chosen the wrong filtration for $\For$, but we will argue later on in \cref{table_filtrations} that our choice is much more efficient, especially in what concerns the layers of the tower.

\end{remark}

The filtration of the forest category we have constructed yields inclusions of full subcategories
$$ \Tree = \FilFor{1} \subseteq \FilFor{2} \subseteq \cdots \subseteq \FilFor{k} \subseteq \cdots \subseteq \For$$
with filtered colimit given by $\For$ once again. We will write $j_k \colon \FilFor{k} \to \FilFor{k+1}$ and $u_k \colon \FilFor{k} \to \For$ for the obvious inclusions. Note that $u_1$ is exactly what we have called $u$ in Chapter 3 and Chapter 3. At the level of presheaves, we will write $\FilfSpaces{k}$ for the category of simplicial presheaves on $\FilFor{k}$ and, if $V \in \fSpaces$, we will write $\FilfSpaces{k}_{/V}$ for the category of presheaves on $\FilFor{k}$ lying over $V^{\smallleq k} \coloneqq u_k^*V$. 

In what concerns the restriction along $u_k$ and $j_k$, as well as the  Kan extensions along them, the
 considerations we made in Section 4.1 for the projective model structure still holds. An analogue of \cref{Quillen_adj_projective_modules} is provided by the next lemma.

\begin{notation}
 Let $F$ be a forest and $k \geq 1$.
 \begin{itemize}
     \item We write $\lvert F \rvert$ for the number of components of $F$, that is, $\lvert F \rvert$ is the filtration of $F$.
     \item We adopt the notation $\mathsf{P}(F;k)$ for the poset
     $$ \mathsf{P}(F;k) = \mathsf{P}(\underline{\lvert F \rvert};k) = \left\{ I \subseteq \underline{\lvert F \rvert} : I \neq \underline{\lvert F \rvert}\: \mathrm{and} \:  \lvert I \rvert \leq k \right\}$$
     ordered under set inclusion. Given $I \in \mathsf{P}(F; k)$, we use the notation $F_I \in \For$ for the subforest of $F$ formed by the components indexed by $I$.
 \end{itemize}
\end{notation}

\begin{lemma}
    Let $k \geq 1$ and $F$ any forest. Then the following statements hold:
    \begin{enumerate}[label=(\roman*)]
        \item The $k$-truncated presheaf $F^{\smallleq k}$ is given by the colimit formula 
        $$F^{\smallleq k} = \colim_{I \in \mathsf{P}(F; k)} F_I.$$
        \item If $X \in \FilfSpaces{k}$, then the forest space $(u_k)_\ast X$ is given by
        $$ \left((u_k)_\ast X \right) (F) = \lim_{I \in \mathsf{P}(F; k)^{\mathsf{op}}} X(F_I)$$
        when evaluated at any forest $F$.
    \end{enumerate}
    \label{formula_truncated_forest}
\end{lemma}

\begin{proof}
    The proof of the lemma is analogous \cref{formula_Uk}, and in this case it suffices to show that for the functor
$$ \Phi : \mathsf{P}(F; k) \longrightarrow \FilFor{k}_{/F}$$
sending a set $I$ to the forest $F_I \subseteq F$, the respective slice categories $\Phi_{(G \to F)/}$ are non-empty and connected. The proof is quite similar, so we only point out the crucial steps one needs to check:
\begin{itemize}
    \item  For non-emptiness, given the left arrow in the diagram
\[\begin{tikzcd}[cramped]
	{G} & {F_I} \\
	{F}
	\arrow[dashed, from=1-1, to=1-2]
	\arrow[from=1-1, to=2-1]
	\arrow["{\Phi(I)}", dashed, from=1-2, to=2-1]
\end{tikzcd}\]
   with $\lvert G \rvert \leq k$, the dashed factorization exists.
   \item For connectedness, given the solid diagram over $F$ below
\[\begin{tikzcd}[cramped]
	G & {F_I} \\
	{F_J} & {F_{I \cap J}}
	\arrow["{f}", from=1-1, to=1-2]
	\arrow["{g}"', from=1-1, to=2-1]
	\arrow["{h}", dashed, from=1-1, to=2-2]
	\arrow[dashed, from=2-2, to=1-2]
	\arrow[dashed, from=2-2, to=2-1]
\end{tikzcd}\]
with $\lvert G \rvert \leq k$, we can complete it as shown by the dashed arrows, with the arrows with source $F_{I \cap J}$ being in the image of $\Phi$.
\end{itemize}

The proof of these results comes down to the simple observation that any morphism of forests $G \to F$ with  $G \in \FilFor{k}$ has to factor through a subforest of $F$ with at most $k$ components via the obvious inclusion.
\end{proof}

Similarly to what we did for \cref{Quillen_adj_projective_modules}, we will collect in the proposition below the main points of the discussion up to now.

\begin{proposition}
    Suppose $k \geq 1$ and let $V$ be a forest space.  Consider the following diagram of adjunctions
    \[\begin{tikzcd}[cramped]
	{\left( \fSpaces_{/V} \right)_{\mathsf{P}}} & {\left( \FilfSpaces{k}_{/V} \right)_{\mathsf{P}}} & {\left( \FilfSpaces{k+1}_{/V} \right)_{\mathsf{P}}}
	\arrow["{u_k^*}", from=1-1, to=1-2]
	\arrow["{(u_k)_!}"', bend right=40, from=1-2, to=1-1]
	\arrow["{(u_k)_\ast}", bend left=40, from=1-2, to=1-1]
	\arrow["{(j_k)_!}", bend left=40, from=1-2, to=1-3]
	\arrow["{(j_k)_*}"', bend right=40, from=1-2, to=1-3]
	\arrow["{j_k^*}"', from=1-3, to=1-2]
\end{tikzcd}\]
    induced from $u_k$ and $j_k$. Then:
    \begin{enumerate}[label=(\roman*)]
        \item The pairs $((u_k)_!, u_k^*)$ and $((j_k)_!, j_k^*)$ define Quillen adjunctions.
        \item The restriction functors $u_k^*$ and $j_k^*$ preserve projective weak equivalences.
    \end{enumerate}
    \label{projective_forests}
\end{proposition}

We will also need to understand how these adjunctions behave with respect to the contravariant model structure. To begin with, we will endow $\fSpaces^{\smallleq k}_{/ V}$ with the \textit{projective contravariant model structure}, which is the left Bousfield localization of the projective model structure at the root inclusions 
$$\rho(F) \longrightarrow F$$
but now we restrict to forests $F$ that have at most $k$ components. The next result states that our considerations for the projective model structure also hold after this localization. We note that the statement that the restriction functors preserves weak equivalences is less immediate in the contravariant setting.

\begin{proposition}
    Suppose $k \geq 1$ and let $V$ be a forest space.  Consider the following diagram of adjunctions
    \[\begin{tikzcd}[cramped]
	{\left( \fSpaces_{/V} \right)_{\mathsf{P}, \ctv}} & {\left( \FilfSpaces{k}_{/V} \right)_{\mathsf{P}, \ctv}} & {\left( \FilfSpaces{k+1}_{/V} \right)_{\mathsf{P}, \ctv}}
	\arrow["{u_k^*}", from=1-1, to=1-2]
	\arrow["{(u_k)_!}"', bend right=40, from=1-2, to=1-1]
	\arrow["{(u_k)_\ast}", bend left=40, from=1-2, to=1-1]
	\arrow["{(j_k)_!}", bend left=40, from=1-2, to=1-3]
	\arrow["{(j_k)_*}"', bend right=40, from=1-2, to=1-3]
	\arrow["{j_k^*}"', from=1-3, to=1-2]
\end{tikzcd}\]
induced from $u_k$ and $j_k$. Then:
    \begin{enumerate}[label=(\roman*)]
        \item The pairs $((u_k)_!, u_k^*)$ and $((j_k)_!, j_k^*)$ define Quillen adjunctions.
        \item The restriction functors $u_k^*$ and $j_k^*$ preserve contravaraint weak equivalences.
    \end{enumerate}
\label{ctv_forests}
\end{proposition}

\begin{proof}
    We will only consider the pair $((u_k)_!, u_k^*)$. We already know that the adjunction is Quillen for the projective model structure by \cref{projective_forests}, and so we get a further Quillen pair
\[\begin{tikzcd}[cramped]
	{\left( \FilfSpaces{k}_{/V} \right)_{\mathsf{P}}} & {\left( \fSpaces_{/V} \right)_{\mathsf{P}, \ctv}} 
	\arrow["{(u_k)_!}", shift left, from=1-1, to=1-2]
	\arrow["{u_k^*}", shift left, from=1-2, to=1-1]
\end{tikzcd}\]
after we localize for the contravariant condition on the right hand side. To show that this adjunction also holds after we similarly localize on the left, it suffices to check that $(u_k)_!(\rho(F) \to F)$ is a contravariant weak equivalence in $\fSpaces_{/V}$ for all $F \in \For^{\smallleq k}_{/V}$. This is clear as $(u_k)_!$ preserves representable presheaves, and $\rho(F) \to F$ is a contravariant weak equivalence for any forest $F$. This concludes the proof of (i).

For (ii), we need to show that a contravariant weak equivalence $\varphi \colon A \to B$ in $\fSpaces_{/V}$ gives rise to a contravariant weak equivalence $\varphi^{\smallleq k} \colon A^{\smallleq k} \to B^{\smallleq k}$ after restricting to filtration degree at most $k$. We remark first that we can assume that $\varphi$ is a projective cofibration: this follows from the factorization of $\varphi$ into a trivial cofibration followed by a trivial fibration, and from part (i) we know that $u_k^*$ preserves contravariant trivial fibrations.

Consider the class $\mathcal{M}_k$ of forest maps over $V$ defined below
$$\mathcal{M}_k = \left\{ \varphi \colon A \to B : \varphi^{\smallleq k} \; \text{is a trivial cofibration in } \left( \FilfSpaces{k}_{/V} \right)_{\mathsf{R},\ctv} \right\}.$$
Here the model structure at the end of the definition of $\mathcal{M}_k$ comes from first considering the Reedy model structure on $\FilfSpaces{k}_{/V}$ that uses the number of edge of a forest as a degree function, and then further localizing for the conravariant condition for forests with at most $k$ components. It is clear that $\mathcal{M}_k$ is a saturated class, so as long as we show that all of the generating trivial cofibrations of $\left( \fSpaces_{/V} \right)_{\mathsf{P}, \ctv}$ are in $\mathcal{M}_k$, then by the previous paragraph we will have proven statement (ii).

Since $u_k^*(F \boxtimes \Lambda^i [n] \to F \boxtimes \Delta^n) = F^{\smallleq k} \boxtimes \Lambda^i [n] \to F^{\smallleq k} \boxtimes \Delta^n$, it is clear that the generating morphism $F \boxtimes \Lambda^i [n] \to F \boxtimes \Delta^n$ is in $\mathcal{M}_k$. As for the generating trivial cofibrations coming from the Bousfield localization, it suffices to check that $\rho^{\smallleq k}(F) \to F^{\smallleq k}$ is a Reedy contravariant trivial cofibration for any forest $F \in \For$.

Appealing to \cref{formula_truncated_forest}, we know that 
$$F^{\smallleq k} = \colim_{I \in \mathsf{P}(F; k)} F_I$$
and we claim that this colimit coincides with its homotopy colimit taken with respect to the projective/Reedy model structure. Indeed, it suffices to prove that diagrams above are Reedy cofibrant, and the latching maps in this case are of the form
$$ \colim_{J \subseteq I, \: J \neq I} F_J \longrightarrow F_I$$
These are easily checked to be normal monomorphisms, so they define Reedy cofibrations, thus showing the Reedy cofibrancy of the diagram. A similar argument also works for $\rho(F)^{\smallleq k}$.

From this, we have the following diagram
\[\begin{tikzcd}
	{\Map^{\smallleq k}_{V}(F^{\smallleq k}, X)} & {\Map^{\smallleq k}_{V}(\hocolim_{I} F_I, X)} & {\holim_I \Map^{\smallleq k}_{V}(F_I, X)} \\
	{\Map^{\smallleq k}_{V}(\rho(F)^{\smallleq k}, X)} & {\Map^{\smallleq k}_{V}(\hocolim_{I} \rho(F_I), X)} & {\holim_I \Map^{\smallleq k}_{V}(\rho(F_I), X)}
	\arrow["\simeq", from=1-1, to=1-2]
	\arrow[from=1-1, to=2-1]
	\arrow["\simeq", from=1-2, to=1-3]
	\arrow[from=1-2, to=2-2]
	\arrow[from=1-3, to=2-3]
	\arrow["\simeq", from=2-1, to=2-2]
	\arrow["\simeq", from=2-2, to=2-3]
\end{tikzcd}\]
for any $X \in \FilfSpaces{k}_{/V}$. If $X \to V$ is a fibrant object in $\left( \FilfSpaces{k}_{/V} \right)_{\mathsf{P}, \ctv}$, the rightmost map is a weak homotopy equivalence, and thus so are the remaining vertical maps. In particular, the left map being a weak equivalence shows that $\rho(F)^{\smallleq k} \to F^{\smallleq k}$ is a contravariant weak equivalences. It is also clear that it is a Reedy cofibration, so we see that it is in the class $\mathcal{M}_k$.
\end{proof}

\begin{remark}
    As we had already noted for the operadic case, the adjunctions $(u_k^*, (u_k)_*)$ and $(j_k^*, (j_k)_*)$ will not be Quillen for both the projective and contravariant model structures by essentially the same reason as that previous scenario.
\end{remark}

With these results at hand, we can construct the desired tower for forests spaces, as well as show that it indeed converges. In its proof we will use the following criteria that appears in \cite{GopplWeiss}, which describes the category of presheaves on a category in terms of certain closure properties. 

\begin{lemma}
    Let $\mathcal{C}$ be a category, and let $\mathcal{A}$ be a class of objects of the category of simplicial presheaves $\mathsf{Psh}(\mathcal{C})$. Suppose that $\mathcal{A}$ satisfies the following properties:
    \begin{enumerate}[label=(\roman*)]
        \item $\mathcal{A}$ is invariant under projective weak equivalences.
        \item $\mathcal{A}$ contains the representable presheaves on $\mathcal{C}$.
        \item $\mathcal{A}$ is closed under arbitrary disjoint unions and homotopy pushouts.
    \end{enumerate}
    Then $\mathcal{A}$ coincides with the presheaf category $\mathsf{Psh}(\mathcal{C})$.
    \label{closure_classes}
\end{lemma}

\begin{proof}[Proof sketch:]
    This is the content of \citep[Prop. 3.1.5]{GopplWeiss}, so we defer the exact details to it and explain the main ideas.

    By property (i), it will suffice to check that the projectively cofibrant presheaves $X$ are in $\mathcal{A}$. These can be written as a filtered colimit of a diagram of the form
    $$ X_0 \longrightarrow X_1 \longrightarrow X_2\longrightarrow \cdots \longrightarrow X$$
    with each morphism $X_{i-1} \to X_i$ being a pushout of a disjoint union of maps of the form
    $$ c \boxtimes \partial \Delta^i \longrightarrow c \boxtimes \Delta^i$$
    for $c \in \mathcal{C}$. One argues inductively that $X_i$ is in $\mathcal{A}$ for all $i \geq 0$: the base case is just a disjoint union of representable functors, and therefore is in $\mathcal{A}$ by (iii).

    For the inductive step, from property (iii) it will suffice to show that $c \boxtimes \partial \Delta^i$ and $c \boxtimes \Delta^i$ are in $\mathcal{A}$, since the pushout diagram describing $X_{i-1} \to X_i$ as a pushout is also a homotopy pushout. For the former, this holds by the induction hypothesis, for the latter this is due to property (ii). Finally, for the final step one argues that the filtered colimit describing $X$ can also be written as a homotopy pushout of objects in $\mathcal{A}$, and another application of (iii) finishes the proof.
\end{proof}

\begin{proposition}
    Suppose $V$ is a forest space and let $X , Y \in \fSpaces_{/V}$. Then the $k$-truncated categories of forest spaces $\FilfSpaces{k}_{/V}$ induce a tower of derived mapping spaces computed with projective contravariant model structure
\[\begin{tikzcd}[cramped]
	&& \vdots \\
	&& {\Map_{V}^{\smallleq 3}(X, Y)} \\
	&& {\Map_{V}^{\smallleq 2}(X, Y)} \\
	{\Map_{V}(X,Y)} && {\Map_{V}^{\smallleq 1}(X, Y)}
	\arrow[from=1-3, to=2-3]
	\arrow[from=2-3, to=3-3]
	\arrow[from=3-3, to=4-3]
	\arrow[bend left=10, from=4-1, to=2-3]
	\arrow[from=4-1, to=3-3]
	\arrow[from=4-1, to=4-3]
\end{tikzcd}\]
which converges to  $\Map_{V}(X,Y)$, that is, the map to the homotopy limit of the tower
$$ \Map_{V}(X,Y) \longrightarrow \holim_k \Map_{V}^{\smallleq k}(X, Y)$$
is a weak homotopy equivalence.
\label{forest_tower}  
\end{proposition}

\begin{proof}
    Fix $Y \in \fSpaces_{/V}$ and consider the class $\mathcal{A}$ of objects $X \in \fSpaces_{/V}$ such that the map to the limit of the tower
    $$\Map_{V}(X,Y) \xlongrightarrow{\simeq} \holim_k \Map_V^{\smallleq k}(X,Y)$$
    is a weak homotopy equivalence. As the derived mapping spaces and homotopy limits are invariant under weak equivalences, we will assume from now on that $Y$ is fibrant, that is, a projective right fibration over $V$.

    By \cref{closure_classes}, we will be done once we show that the class $\mathcal{A}$ satisfies the following properties:
    \begin{itemize}
        \item \textit{$\mathcal{A}$ is invariant under projective weak equivalences}: the statement is certainly invariant contravariant weak equivalences, and projective weak equivalences are a special instance of these.
        \item \textit{Representable presheaves are in $\mathcal{A}$}: for the case of a representable presheaf $X = F$ the morphism to the limit becomes
        $$Y(F) \longrightarrow \holim_k \Map_V^{\smallleq k}(F, Y),$$
        and since we assumed $Y$ to be fibrant, the term inside the homotopy limit is equivalent to the homotopy fiber of $Y(F^{\smallleq k}) \to V(F^{\smallleq k})$. Since $F^{\smallleq k} \cong F$ for all $k$ larger than the filtration degree of $F$, we see that the limit in question eventually stabilizes to $Y(F)$.
        \item \textit{$\mathcal{A}$ is closed under arbitrary disjoint unions and homotopy pushouts}: this follows from the fact that the derived mapping space sends coproducts on the first variable to products, and homotopy pushouts to homotopy pullbacks.
    \end{itemize}
\end{proof}

\begin{remark}
    The work of Göppl concerns the convergence of a tower of derived mapping spaces computed with the projective model structure
    $$ \Map(X,Y) \longrightarrow \holim_k \Map^{\smallleq k}(X,Y)$$
    associated to some filtration of $\mathcal{C}$ by full subcategories $\mathcal{C}_k$ with filtered colimit $\mathcal{C}$. What we have shown above is that this argument is still valid if we replace the projective model structure with any left Bousfield localization of it.
\end{remark}

We conclude this section by connecting the operadic and the forest towers. Suppose that we are given a $k$-truncated module $M \in \RMod{\mathcal{P}}^{\smallleq k}$, and recall the nerve functor for operadic right modules $N_{\mathcal{P}} \colon \RMod{\mathcal{P}} \to \fSpaces_{/ N_d \mathcal{P}}$ that we constructed whenever $\mathcal{P}$ is a unital simplicial operad. Since $M$ is $k$-truncated, it will not define a presheaf on $\For$, but we claim that it does define one on $\For^{\smallleq k}$ via the formula $$ N_{\mathcal{P}}M(F) = M(\rho(F)) \times_{(N_d \mathcal{P})(\rho(F))} (N_d \mathcal{P})(F)$$
that we already had before. Indeed, this follows from the fact that $\lvert \rho(F) \rvert = \lvert F \rvert$ for any forest $F$, and the action of a root face $\partial_{r} G \to G$ such that $\lvert \partial_r G \rvert, \lvert G \rvert \leq k$ corresponds to one of the right actions from the module structure on $M$ which is already detected when restricted to $\mathsf{Env}(\mathcal{P})^{\smallleq k}$. 

To summarize the last paragraph, the nerve construction for right modules can be restricted along the truncation functors, leading to the dashed arrow
\[\begin{tikzcd}[cramped]
	{\RMod{\mathcal{P}}} & {\RMod{\mathcal{P}}^{\smallleq k}} \\
	{\fSpaces_{/ N_d \mathcal{P}}} & {\FilfSpaces{k}_{/N_d \mathcal{P}}}
	\arrow["{U_k^*}", from=1-1, to=1-2]
	\arrow["{N_{\mathcal{P}}}"', from=1-1, to=2-1]
	\arrow["{N_{\mathcal{P}}}"', dashed, from=1-2, to=2-2]
	\arrow["{u_k^*}", from=2-1, to=2-2]
\end{tikzcd}\]
making the square above commute. This functor also admits a left adjoint 
$$\tau_{\mathcal{P}} \colon \FilfSpaces{k}_{/N_d \mathcal{P}} \to \RMod{\mathcal{P}}^{\smallleq k}$$
which sends a $k$-truncated map $X \to (N_d \mathcal{P})^{\smallleq k}$ to the $k$-truncated $\mathcal{P}$ generated by such data. The next result establishes the compatibility of the nerve with the operadic and forest filtrations.

\begin{proposition}
    Let $k \geq 1$ and $\mathcal{P}$ a unital simplicial operad. Then:
    \begin{enumerate}[label=(\roman*)]
        \item The equalities $N_{\mathcal{P}} U_k^* = u_k^* N_{\mathcal{P}}$ and $N_{\mathcal{P}} (U_k)_\ast = (u_k)_* N_{\mathcal{P}}$ hold.
        \item The adjunction associated to the nerve functor
        \[\begin{tikzcd}[cramped]
	{\RMod{\mathcal{P}}^{\smallleq k}} & {\left( \FilfSpaces{k}_{/N_d \mathcal{P}} \right)_{\mathsf{P}, \ctv}}
	\arrow["{N_{\mathcal{P}}}"', shift right, from=1-1, to=1-2]
	\arrow["{\tau_{\mathcal{P}}}"', shift right, from=1-2, to=1-1]
\end{tikzcd}\]
defines a Quillen pair, and $N_{\mathcal{P}}$ detects contravariant weak equivalences. If $\mathcal{P}$ is also a $\Sigma$-free operad, then this adjunction is a Quillen equivalence.
    \end{enumerate}
\label{compatible_Quillen}
\end{proposition} 

\begin{proof}
    It will be convenient throughout the proof to keep track of the truncation, so we will write $N_{\mathcal{P}}^{k}$ for the $k$-truncated nerve. In what concerns (i), the first equality is exactly the commutativity of the diagram before the proposition. For the other one, this follows from a direct computation using the expressions in \cref{formula_Uk} and \cref{formula_truncated_forest}. As for (ii), $(\tau^k_{\mathcal{P}}, N^k_{\mathcal{P}})$ being a Quillen pair and $N^k_{\mathcal{P}}$ detecting contravariant weak equivalences can be proved by the same argument as we used in the proof of \cref{Modules_Quillen_Adjunction}.
    
Finally, and assuming from now on that $\mathcal{P}$ is a $\Sigma$-free operad, by the detection property we have just proved we can reduce checking that $N_{\mathcal{P}}$ defines a Quillen equivalence to verifying that the unit of the adjunction
    $$X \longrightarrow N_{\mathcal{P}}^k \tau_{\mathcal{P}}^k X$$
    is a contravariant weak equivalence, whenever $X$ is cofibrant in $\left( \FilfSpaces{k}_{/N_d \mathcal{P}} \right)_{\mathsf{P}, \ctv}$. 
    
    In order to do so, we can fit this arrow into the diagram below
\begin{equation}
\begin{tikzcd}[cramped]
	X & {N_{\mathcal{P}}^k \tau_{\mathcal{P}}^k X} \\
	{u_k^*(u_k)_! X} & { N_{\mathcal{P}}^k \tau_{\mathcal{P}}^k u_k^*(u_k)_!X}
	\arrow[from=1-1, to=1-2]
	\arrow["\cong"', from=1-1, to=2-1]
	\arrow["\cong"', from=1-2, to=2-2]
	\arrow[from=2-1, to=2-2]
\end{tikzcd}
\label{square_units}
\end{equation}
    with the vertical maps being isomorphisms due to the fully faithfullness of $u_k$. We can rewrite the bottom corner presheaf as $u_k^* N_{\mathcal{P}} \tau_{\mathcal{P}} (u_k)_! X$ using (i). Moreover, the bottom horizontal map is the result of applying $u_k^*$ to the counit map
    $$\varepsilon_{(u_k)_!(X)} : (u_k)_!(X) \to N_{\mathcal{P}} \tau_{\mathcal{P}} (u_k)_! (X).$$ 
    Since $(u_k)_!$ is left Quillen by \cref{ctv_forests}, we know that $\varepsilon_{(u_k)_!X}$ is a contravariant weak equivalence by the proof of \cref{eqv_ctv_RMod}. This finishes the proof of (ii), since $u_k^*$ preserves contravariant weak equivalences by \cref{ctv_forests}.
\end{proof}

A consequence of the compatibility between the nerve for right modules and the filtrations we have been analysing is that there exists a map of derived mapping spaces
$$ \Map_{\mathcal{P}}(M,L) \longrightarrow \Map_{N_d \mathcal{P}}^{\smallleq k} \left( \mathbb{R}N_{\mathcal{P}}M, \mathbb{R}N_{\mathcal{P}}L \right)$$
for every $k \geq 1$, which leads to the following theorem.
\begin{theorem}
    Suppose $\mathcal{P}$ is a unital simplicial operad, and let $M, L$ be operadic right $\mathcal{P}$-modules. Then the nerve of operadic right $\mathcal{P}$-modules assembles into a map between the operadic and forest towers
\[\begin{tikzcd}[cramped]
	{\Map_{\mathcal{P}}(M,L)} & {\Map_{N_d\mathcal{P}}\left( \mathbb{R}N_{\mathcal{P}} M, \mathbb{R}N_{\mathcal{P}}L \right)} \\
	{\left\{ \Map^{\smallleq k}_{\mathcal{P}}(M,L) \right\}_{k \geq 1}} & {\left\{ \Map^{\smallleq k}_{N_d\mathcal{P}}\left( \mathbb{R}N_{\mathcal{P}} M, \mathbb{R}N_{\mathcal{P}} L \right) \right\}_{k \geq 1}}
	\arrow["{N_{\mathcal{P}}}", from=1-1, to=1-2]
	\arrow[from=1-1, to=2-1]
	\arrow[from=1-2, to=2-2]
	\arrow["{N_{\mathcal{P}}}", from=2-1, to=2-2]
\end{tikzcd}\]
If $\mathcal{P}$ is additionally $\Sigma$-free, then $N_{\mathcal{P}}$ induces a Quillen equivalence between each stage of the filtration.
\label{main_thm}
\end{theorem}

\begin{proof}
    This is immediate from part (ii) of \cref{compatible_Quillen}.
\end{proof}

Restricting our attention to the objects with filtration one gives another interpretation of the first stage of the tower, complementing the result in \cref{first_stage}. From the point of view of dendroidal homotopy theory, it also clarifies the role the contravariant model structure plays in dendroidal spaces, instead of forest spaces.

\begin{corollary}
    Let $\mathcal{P}$ be a unital simplicial operad. Then the nerve for operadic right modules defines a Quillen adjunction
    \[\begin{tikzcd}[cramped]
	{\mathsf{Psh}(\mathsf{Cat}(\mathcal{P}))_{\mathsf{P}}} & {\left( \dSpaces_{/N_d \mathcal{P}} \right)_{\mathsf{P}, \ctv}.}
	\arrow["{N_{\mathcal{P}}}"', shift right, from=1-1, to=1-2]
	\arrow["{\tau_{\mathcal{P}}}"', shift right, from=1-2, to=1-1]
\end{tikzcd}\]
    If $\mathcal{P}$ is additionally $\Sigma$-free, then it can be upgraded to a Quillen equivalence.
    \label{first_stage_forests}
\end{corollary}

\subsection{The layers of the forest tower}

In the previous section we constructed the truncation tower for forest spaces over a fixed forest space $V$, and we argued that this filtration is moreover compatible with the operadic filtration for simplicial right $\mathcal{P}$-modules, ultimately yielding an equivalence of towers at least when $\mathcal{P}$ is unital and $\Sigma$-free. In addition to this, \cref{first_stage} and \cref{first_stage_forests} give descriptions of the first stage of both towers.

Following the strategy explained at the end of the first section, we now wish to study the layers of the forest tower, that is, for $X, Y \in \fSpaces_{/V}$ we are interested in the space
$$ \mathsf{Fib}_k(X,Y) \coloneqq \mathrm{hofib} \left( \Map_V^{\smallleq k}(X,Y) \longrightarrow \Map_{V}^{\smallleq k-1}(X,Y) \right)$$
for $k \geq 2$. The goal of this section will be to simplify our mapping spaces by replacing the target of the inclusion $\FilFor{k-1} \to \FilFor{k}$ by a smaller category.

\begin{notation}
    From now on we will need to work with the category of \textit{closed forests} as we will need our forests to have stumps for some of the proofs in Section 4.3.1, for instance \cref{convex_maps}. This means that we will replace the forest category $\For$ with the full subcategory on all forests $F$ where all edges are the outcoming edge of some vertex. In particular, by \textit{leaves} we will always mean the outcoming edges of the stumps. However, in order not to overwhelm the notation below, we have adopted to still write $\For$ and $\fSpaces$ for this closed context.

    We point out that everything we have done until now works without change: the main difference is that we don't have leaf faces anymore, but we still have the inner and root faces available. As we are only interested in right fibrations, this change is immaterial.
\end{notation}

\subsubsection{Reduction to $k \cdot \eta$}

We begin by introducing a category of forests which lies between $\For^{\smallleq k-1}$ and $\FilFor{k}$.

\begin{definition}
    Let $k \geq 2$. We define the category of \textit{basic k-forests} to be the full subcategory $\mathsf{Bas}^{\smallleq k}$ of $\FilFor{k}$ spanned by the the objects of $\FilFor{k-1}$, together with the forest of edges $k \cdot \eta$. This category sits between the two consecutive stages of the forest filtration
    $$ \FilFor{k-1} \longrightarrow \mathsf{Bas}^{\smallleq k} \xlongrightarrow{\mathsf{bas}} \FilFor{k}$$
    via the inclusions above.
\end{definition}

The first main step in our desired reduction is that restriction of forest spaces from $\FilFor{k}$ to $\mathsf{Bas}^{\smallleq k}$ will not change the homotopy type of the derived mapping space. The desired simplification of the layers of the tower will directly follow.

Recall from now on that a dendroidal space $V$ is \textit{reduced} whenever $V(C_1) \simeq V(\eta)$. 

\begin{proposition}
    Let $V$ be a reduced dendroidal Segal space, and let $X, Y \in \fSpaces_{/V}$. For $k \geq 2$, the map induced by the restriction to the subcategory of basic $k$-forests $\mathsf{Bas}^{\smallleq k}$
$$\Map_V^{\smallleq k}(X,Y) \xlongrightarrow{\mathsf{bas}^*} \Map_V^{\smallleq k}(\mathsf{bas}^*X,\mathsf{bas}^*Y) $$
    is a weak homotopy equivalence.
    \label{reduction_basic_forests}
\end{proposition}

\begin{corollary}
    Let $V$ be a dendroidal Segal space, and let $X, Y \in \fSpaces_{/V}$. The $k^{th}$-layer of the forest tower for $\Map_V(X,Y)$ is given by 
    $$ \mathsf{Fib}_k(X,Y) \simeq \mathrm{hofib} \left( \Map_V^{\smallleq k}(\mathsf{bas}^* X,\mathsf{bas}^*Y) \longrightarrow \Map_{V}^{\smallleq k-1}(X,Y) \right)$$
    for $k \geq 2$.
    \label{first_layer_desc}
\end{corollary}

\begin{proof}
    Since $\mathsf{Bas}^{\smallleq k}$ includes all forests of filtration $ \leq k-1$, we have the following commutative diagram
\[\begin{tikzcd}[cramped]
	{\Map_V^{\smallleq k}(X,Y)} & {\Map_V^{\smallleq k}(\mathsf{bas}^*X,\mathsf{bas}^*Y)} \\
	{\Map_V^{\smallleq k-1}(X,Y)} & {\Map_V^{\smallleq k-1}(X,Y)}
	\arrow["\simeq", from=1-1, to=1-2]
	\arrow[from=1-1, to=2-1]
	\arrow[from=1-2, to=2-2]
	\arrow["{\simeq }", from=2-1, to=2-2]
\end{tikzcd}\]
which is a homotopy cartesian square by \cref{reduction_basic_forests}. The corollary follows by taking vertical homotopy fibers.
\end{proof}

In order to prove \cref{reduction_basic_forests}, we filter $\mathsf{Bas}^{\smallleq k} \to \FilFor{k}$ by an intermediate full subcategory
$$ \mathsf{Bas}^{\smallleq k} \xlongrightarrow{\alpha} \mathsf{Ext}^{\smallleq k} \xlongrightarrow{\beta} \FilFor{k}$$
which we construct in \cref{extended_forests_definition}, and show that restriction along the two inclusions above also induces weak homotopy equivalences on derived mapping spaces. The proofs of these results will involve some manipulation and simplification of homotopy limits and homotopy colimts, and for that reason we now recall the definition of homotopy left/right cofinal functors.

\begin{definition}
    Let $\varphi \colon \mathcal{C} \to \mathcal{D}$ be a functor.
    \begin{itemize}
        \item We say that $\varphi$ is a \textit{homotopy left cofinal functor} if the slice categories $\varphi_{/d}$ are contractible for every object $d \in \mathcal{D}$.
        \item We say that $\varphi$ is a \textit{homotopy right cofinal functor} if the slice categories $\varphi_{d/}$ are contractible for every object $d \in \mathcal{D}$.
    \end{itemize}
\end{definition}

The importance of homotopy left/right cofinal functors is that homotopy limits and colimits are invariant upon restricting along them. More precisely, we have the following for simplifying homotopy limits and colimits of contravariant functors.

\begin{proposition}
    Let $\varphi \colon \mathcal{C} \to \mathcal{D}$ be a functor, and consider a presheaf $X \colon \mathcal{D}^{\mathsf{op}} \to \mathcal{M}$ with values in some  model category $\mathcal{M}$.
    \begin{enumerate}[label=(\roman*)]
        \item If $\varphi$ is a homotopy left cofinal functor, then the map below
        $$ \hocolim_{\mathcal{C}^\op} \: \varphi^* X \longrightarrow \hocolim_{\mathcal{D}^\op} \: X$$
        is a weak equivalence in $\mathcal{M}$.
        \item If $\varphi$ is a homotopy right cofinal functor, then the map below
        $$ \holim_{\mathcal{C}^\op} \: \varphi^* X \longrightarrow \holim_{\mathcal{D}^\op} \: X$$
        is a weak equivalence in $\mathcal{M}$.
        
    \end{enumerate}
    \label{limit_simplifying}
\end{proposition}

\begin{proof}
    See \citep[Thm. 19.6.3]{HirschhornLocalizations}.
\end{proof}

We will also need to have a practical way of checking when a functor is homotopy left/right cofinal. The well-known criteria below are probably the easiest and most useful manner for verifying this, and we will make use of it quite often.

\begin{lemma}
    Let $\varphi \colon \mathcal{C} \to \mathcal{D}$ be a functor.
    \begin{enumerate}[label=(\roman*)]
        \item If $\varphi$ is fully faithful and admits a right adjoint, then $\varphi$ is a homotopy left cofinal functor.
        \item If $\varphi$ is fully faithful and admits a left adjoint, then $\varphi$ is a homotopy right cofinal functor.
    \end{enumerate}
    \label{adjoint_functor_simplifying}
\end{lemma}

\begin{proof}
    If $\lambda \colon \mathcal{D} \to \mathcal{C}$ is a left adjoint to $\varphi$, then the counit map $d \to \varphi \lambda (d)$ is an isomorphism due to $\varphi$ being fully faithful. Consequently, this morphism defines a terminal object in $\varphi_{d/}$, which shows that this category is contractible. A dual argument shows (ii).
\end{proof}

Let us now focus on analysing the inclusion $\beta \colon \mathsf{Ext}^{\smallleq k} \to \FilFor{k}$, for which we will first need to set up some definitions and prove an auxiliary lemma.

\begin{definition}
 We say that a morphism of forests is a \textit{convex map} if it satisfies the following conditions:
\begin{itemize}
    \item If $f \colon T \to S$ is a dendroidal morphism, then it is convex whenever it  maps the root edge of $T$ to the root edge of $S$, and  given any collection of ordered edges $e_1 < e_{12} < e_2$ of the tree $S$, the edge $e_{12}$ is in the image of $f$ whenever $e_1$ and $e_2$ are also in the image of $f$.
    \item If $f \colon F \to G$ is a forest morphism, then it is convex when the dendroidal map defined by each component of $F$ is convex.
\end{itemize}
\end{definition}

As a quick example to keep in mind, root faces and inner faces don't define convex maps.

\begin{remark}
    Our definition of convex maps is close to the one in  \citep[Rem. 2.9]{BonventrePereira}, as well as to the definition of subtrees present in \citep[Def. 3.1.17]{GopplWeiss}. Our definition diverges in that we imposed the root condition, because we need the wide forest spine to be a natural construction with respect to convex maps, and it is not for the root face of a forest (one can verify this for the tree in \cref{different_spines} for instance).
\end{remark}

\begin{definition}
    Let $k \geq 2$ and fix $F \in \FilFor{k}$. 
    \begin{itemize}
        \item Define the category of \textit{extended k-forests} $\mathsf{Ext}^{\smallleq k} \subseteq \FilFor{k}$ to be the full subcategory spanned by the objects of $\FilFor{k-1}$, the forest of edges $k \cdot \eta$, and its degeneracies (these are exactly the forests of linear trees, so we will call them the \textit{linear $k$-forests} from now on).
        \item Define the category $\mathsf{Conv}_F^{\smallleq k} \subseteq  \mathsf{Ext}^{\smallleq k}_{/ F}$ to be the full subcategory spanned by the maps of forests $G \to F$ which are convex. Note that a morphism in this category is automatically convex also.
    \end{itemize}
    \label{extended_forests_definition}
\end{definition}

As pointed out in \citep[Rem. 3.5]{BonventrePereiraEquivariant}, the interest in convex maps stems from the fact that the spine construction is functorial for these type of maps, which makes it possible to apply the Segal condition in certain situations. This is also our motivation for considering them, but we will need to be careful about what kind of spine we choose.

In order to see this, recall that in \cref{different_spines} we saw that for the given tree $T$ the forest and wide spines are
$$ C_3 \cup_{ 3 \cdot \eta} \left( C_2 \oplus \eta \oplus C_1 \right) \hspace{2em} \mathrm{and} \hspace{2em} C_3 \cup_{ \eta^{\amalg 3}} \left( C_2 \amalg \eta \amalg C_1 \right)$$
respectively. Even though $T$ is a tree -- and therefore has filtration degree 1 -- its forest spine is a pushout of forests of filtration degree higher than 1, namely $3 \cdot \eta$ appears in this colimit. However, the same kind of problem doesn't happen for the wide spine of $T$: for instance, $3 \cdot \eta$ is now replaced by $\eta^{\amalg 3}$, which is a coproduct of presheaves which lie in filtration 1.

In summary, the following always holds: given $F$ of filtration $k$, its wide spine $\wfSpine(F)$ is a colimit of representables in $\For^{\smallleq k}$, whereas the forest spine $\fSp(F)$ is not compatible with the filtration we have been considering. 

The next lemma establishes funtoriality of the wide spine with respect to convex maps.

\begin{lemma}
    Let $k \geq 2$ and fix $F \in \For^{\smallleq k}$.
    \begin{enumerate}[label=(\roman*)]
        \item The inclusion of categories $\mathsf{Conv}_F^{\smallleq k} \to \mathsf{Ext}^{\smallleq k}_{/ F}$ is fully faithful and admits a left adjoint. In particular, it is a homotopy right cofinal functor.
        \item Suppose $\varphi \colon F_1 \to F_2$ is a convex map of forests. Then it induces a morphism of forest sets
        $$ \wfSpine(\varphi) : \wfSpine(F_1) \longrightarrow \wfSpine(F_2)$$
        between the wide spines.
    \end{enumerate}

    \label{convex_maps}
\end{lemma}

\begin{proof}
    Throughout the proof we will use the following notation. Let $T$ be a tree, and $e_1 < e_2$ distinct edges in $T$ such that there is no $e \in E(T)$ such that $e_1 < e < e_2$. Given a linear tree $\ell$, we write 
    $$ T/e_1 \circ_{\ell} e_2/T$$
    for the tree obtained from $T$ by replacing the vertex incident to $e_1$ and $e_2$ with the linear path $\ell$.

    Let $L \colon  \mathsf{Ext}^{\smallleq k}_{/ F} \to \mathsf{Conv}_F^{\smallleq k}$ be the desired left adjoint, which we now wish to describe. We will do this by producing, for an object $f \colon G \to F$ in the source category, a natural decomposition
\[\begin{tikzcd}[cramped]
	G & {G'} \\
	F
	\arrow[from=1-1, to=1-2]
	\arrow["f"', from=1-1, to=2-1]
	\arrow["f'", from=1-2, to=2-1]
\end{tikzcd}\]
    with $G' \to F$ a convex map. We then set $L(G \to F)$ to be $G' \to F$, and it is clear that this defines a left adjoint, for instance by verifying the triangle identities. We will write $G = \bigoplus_{i \in I} G_i$ for the decomposition of $G$ into trees, with $r_i$ the root of $G_i$.
    
    We begin by explaining how one can make $f'$ respect the root condition when $F$ is a tree with root edge $r_F$. Firstly, we note that if $f(r_i) = r_F$ for some $i$, then the independence condition forces $G$ to be a tree. Therefore, we may assume from now on that $f(r_i) \neq r_F$ for all $i \in I$. In this scenario, let $T$ be the subtree of $F$ with leaves $f(r_i)$ for all $i \in I$ and root $r_F$. We define $G'$ to be the tree obtained from $G$ by grafting $T$ along the set of roots $R(G)$. There is an inclusion of forests $G \to G'$, and it is easy to see that $f$ extends to a map on $G'$ which now satisfies the root condition. Below we illustrate this construction.
    \[\begin{tikzpicture}[scale=0.45]
        \draw (0,-2.5)--(0,-2);
        \draw (-1,-1)--(0,-2)--(1,-1)--(-1,-1);
        
        \node[scale=0.6] at (1.2,-2) {$\oplus$};
        \node[scale=1] at (2.2,-2) {$\cdots$};
        \node[scale=0.6] at (3.2,-2) {$\oplus$};
        \node[scale=0.6] at (0,-1.4) {$G_1$};
        \node[scale=0.6] at (4.4,-1.4) {$G_n$};

        \draw (4.4,-2.5)--(4.4,-2);
        \draw (3.4,-1)--(4.4,-2)--(5.4,-1)--(3.4,-1);

        \draw (9,-0.5)--(9,0);
        \draw (8,1)--(9,0)--(10,1)--(8,1);
        
        \node[scale=0.6] at (9,0.6) {$G_1$};
        \node[scale=0.6] at (12.4,0.6) {$G_n$};
        \node[scale=1] at (10.7,0) {$\cdots$};
        \draw[->] (5.5,-2)--(8,-2);

        \draw (12.4,-0.5)--(12.4,0);
        \draw (11.4,1)--(12.4,0)--(13.4,1)--(11.4,1);
        \draw (10.7,-3)--(10.7,-3.5);
        \draw (10.7,-3)--(13.4,-0.5)--(8,-0.5)--(10.7,-3);

        \node at (10.7,-1.5) {$T$};

        \draw (0.5,-6.5)--(0.5,-6);
        \draw (-0.5,-5)--(0.5,-6)--(1.5,-5)--(-0.5,-5);
        \draw (4.4,-6.5)--(4.4,-6);
        \draw (3.4,-5)--(4.4,-6)--(5.4,-5)--(3.4,-5);
        \draw (2.45,-9)--(2.45,-9.5);
        \draw (2.45,-9)--(5.4,-6.5)--(-0.5,-6.5)--(2.45,-9);

        \node at (2.45,-7.5) {$T$};

        \node[scale=0.5] at (0.5,-5.3) {$f(G_1)$};
        \node[scale=0.5] at (4.4,-5.3) {$f(G_n)$};
        \node[scale=1] at (2.45,-6) {$\cdots$};
        \draw[->] (2.45,-3)--(2.45,-4.5);
        \node[scale=0.7] at (2.15,-3.75) {$f$};
        \draw[->] (10,-4)--(6,-6.5);
    \end{tikzpicture}\]
 If $F$ is not a tree, then we do the construction just described to each component of $F$, and the map $f' \colon G' \to F$ thus constructed sends root edges to root edges.
    
    Suppose now that we are given two distinct comparable edges $g_1 < g_2$ in $G$ for which there exists no other edge $g$ such that $g_1 < g < g_2$. Let us say that the pair $(g_1, g_2)$ has \textit{convex image} if, for any edge $f(g_1) < e < f(g_2)$ in $F$, we have that $e$ is in the image of $f$. Our objective is to naturally alter $G$ over $F$ in such a way that all pairs of comparable edges $(g_1, g_2)$ have convex image. 
    
    In order to do this, suppose that $(g_1, g_2)$ does not have convex image, and let $\ell_{12}$ be the path in $F$ from $f(g_1)$ to $f(g_2)$. The remaining edges of $\ell_{12}$ cannot be in the image of $f$, as otherwise this would imply that there exists at least one edge between $g_1$ and $g_2$ in $G$. There is\footnote{From now on we really need to assume that we are working with closed trees, as otherwise the map $\lambda \colon G \to G'$ we will define might not exist. This is relevant whenever $e_1$ and $e_2$ are incident to a vertex $v$ which also has others edges incident to it.} a dendroidal map $\lambda \colon C_1 \to \ell_{12}$ mapping the extreme edges of $C_1$ to the extreme edges of $\ell_{12}$. Thinking of $C_1$ as the vertex incident to $e_1$ and $e_2$ in $G$, the map $\lambda$ induces a diagram on graftings
\[\begin{tikzcd}[cramped]
	{ G= G/g_2 \circ_{C_1} g_1/G} & {G'= G/g_2 \circ_{\ell_{12}} g_1/G} \\
	{ F = F/f(g_2) \circ_{\ell_{12}} f(g_1)/F }
	\arrow["\lambda", from=1-1, to=1-2]
	\arrow["f"', from=1-1, to=2-1]
	\arrow["f'", from=1-2, to=2-1]
\end{tikzcd}\]
where now for the diagonal map $G' \to F$ we have reduced by one the number of pairs with convex image. Since there are only finitely many edges in $G$, this procedure must eventually halt; noticing moreover that this process always leads to root-preserving maps (if $f$ already satisfies this condition), we define $L(G \to F)$ to be the end result, thus showing (i).

For the functoriality with respect to convex maps of forests, we note that this class of maps is generated under composition by:
\begin{enumerate}[label=(\alph*)]
    \item The maps of forests $\mathrm{id}_H \oplus f \oplus \mathrm{id}_{H'} \colon H \oplus T \oplus H' \to H \oplus S \oplus H'$, where $f \colon T \to S$ is a convex dendroidal map.
    \item The component inclusion of forests.
\end{enumerate}

For the case (a), we can assume that the components $H$ and $H'$ don't exist; then the naturality follows from the observation that $f$ must be generated under composition by leaf faces, degeneracies and root faces when they exist in the category $\Tree$, and all of these are examples of convex maps. As for (b), this is clear since $\wfSpine(-)$ commutes with the concatenation of forests. 
\end{proof}

\begin{proposition}
     Suppose $V$ is a dendroidal Segal space, and let $X, Y \in \fSpaces_{/V}$. For $k \geq 2$, the map induced by restriction along the inclusion $\beta \colon \mathsf{Ext}^{\smallleq k} \to \FilFor{k}$
     $$ \Map_V^{\smallleq k}(X, Y) \longrightarrow \Map_V^{\smallleq k}(\beta^* X, \beta^* Y)$$
     is a weak homotopy equivalence.
     \label{red_ext}
\end{proposition}

\begin{proof}
Since the statement is invariant under contravariant weak equivalences, we can always replace $X \to V$ and $Y \to V$ by Reedy right fibrations. It suffices to show that the derived unit map $X \to \mathbb{R}\beta_\ast \beta^* X$ is a contravariant weak equivalence. We claim that the map we get after evaluating the unit map on a forest $F$
$$ X(F) \longrightarrow \Tilde{X}(F) \coloneqq (\mathbb{R}\beta_\ast \beta^* X) (F) = \holim_{\left( \mathsf{Ext}^{\smallleq k}_{/F} \right)^\op} \: X(G)$$
defines a weak homotopy equivalence over $V(F)$, which will immediately yield the desired result. Moreover, since $\mathsf{Ext}^{\smallleq k}$ contains all the corollas and all the forests of edges $\ell \cdot \eta$ for $ \ell \leq k$, we only need to show that $\tilde{X}$ satisfies the wide Segal condition. By \cref{convex_maps} (i) and \cref{limit_simplifying} (ii), there is a weak homotopy equivalence
$$ \holim_{\left( \mathsf{Ext}^{\smallleq k}_{/F} \right)^\op} \: X(G) \xlongrightarrow{\simeq} \holim_{\left(  \mathsf{Conv}^{\smallleq k}_{F} \right)^\op} \: X(G),$$
and by \cref{convex_maps} (ii) the restriction of $\tilde{X}$ along the wide spine inclusion $\wfSpine(F) \to F$ coincides with the map below
\begin{equation}
 \holim_{\left(  \mathsf{Conv}^{\smallleq k}_{F} \right)^\op} \: X(G) \xlongrightarrow{\simeq} \holim_{\left(  \mathsf{Conv}^{\smallleq k}_{F} \right)^\op} \: X(\wfSpine(G)).
\label{Segal_convex}
 \end{equation}
As $X$ satisfies the wide Segal condition by \cref{right_fibration_Segal}, each term of the homotopy limit above is a weak homotopy equivalence, and consequently so is \eqref{Segal_convex}, as we wanted to show.
\end{proof}

We now aim to show how can simplify our presheaves by removing all extended versions of $k \cdot \eta$ from $\mathsf{Ext}^{\smallleq k}$. For that purpose, we define, for each linear $k$-forest $F$, full subcategories
$$ \mathcal{C}_F^{\mathrm{r}} \subseteq \mathcal{C}_F \subseteq \mathsf{Bas}^{\smallleq k}_{/F}$$
as follows:
\begin{itemize}
    \item The subcategory $\mathcal{C}_F$ is spanned by all maps $G \to F$ with $G$ a linear, such that root edges are sent to root edges.
    \item The subcategory $\mathcal{C}_F^{\mathrm{r}}$ is spanned by all maps $\ell \cdot \eta \to F$ that send all root edges to root edges. 
\end{itemize}

Notice that if $G \to F$ with $G \in \mathsf{Bas}^{\smallleq k}$ and $F$ a linear $k$-forest, then $G$ is automatically linear also.

\begin{lemma}
    Let $F$ be a linear $k$-forest. Then the following hold:
    \begin{enumerate}[label=(\roman*)]
        \item The inclusion $\mathcal{C}_F \to \mathsf{Bas}^{\smallleq k}_{/F}$ is homotopy right cofinal.
        \item The inclusion $\varphi \colon \mathcal{C}_F^{\mathrm{r}} \to \mathcal{C}_F$ induces a weak homotopy equivalence
        $$ \holim_{\left( \mathcal{C}_F^{\mathrm{r}} \right)^\op} \varphi^* X \xlongrightarrow{\simeq} \holim_{\left( \mathcal{C}_F \right)^\op} X $$
        for any presheaf $X$ satisfying $X(G) \simeq X(\rho(G))$ for $G$ a linear forest.
    \end{enumerate}
    \label{simplify_extension}
\end{lemma}

\begin{proof}
    The proof of (i) is similar to the one for \cref{convex_maps}. In this case one notices that any object $f \colon G \to F$ in $\mathsf{Bas}^{\smallleq k}_{/F}$ factors as
    \[\begin{tikzcd}[cramped]
	G & {G'} \\
	F
	\arrow[from=1-1, to=1-2]
	\arrow["f"', from=1-1, to=2-1]
	\arrow["f'", from=1-2, to=2-1]
\end{tikzcd}\]
with $f'$ in $\mathcal{C}_F$. This can be accomplished as follows: consider a component $f_i \colon G_i \to F_j$ of $f$, and assume that the image of the root edge $f_i(r_i)$ is not the root edge of $F_j$, call it $r_F$. In this situation we can consider the path $\ell_i$ between these edges. Set $G_i'$ to the be the tree obtained by grafting $\ell_i$ to the root of $G_i$. Proceeding like this for each component leads to the desired factorization.

For (ii) we appeal to \citep[Prop. A.4] {DuggerCombinatorial} to simplify the homotopy limit as we stated. A very similar argument is also present in the proof of Lemma 3.1.20 of \cite{GopplWeiss}. We observe that there is an adjunction
\[\begin{tikzcd}[cramped]
	{\mathcal{C}_F^{\mathrm{r}}} & {\mathcal{C}_F}
	\arrow["\varphi", shift left, from=1-1, to=1-2]
	\arrow["{\rho(-)}", shift left, from=1-2, to=1-1]
\end{tikzcd}\]
with the right adjoint sending an object $G \to F$ to its restriction to the root edges $\rho(G) \to F$. The unit of this adjoint pair is an isomorphism, and the counit when evaluated on $X$ is the root inclusion map
$$ X(G) \longrightarrow X(\rho(G)),$$
which is a weak homotopy equivalence for $G$ linear. The conclusion now follows by applying Dugger's result that we mentioned above.
\end{proof}

\begin{proposition}
     Suppose $V$ is a reduced dendroidal Segal space, and let $X, Y \in \fSpaces_{/V}$. For $k \geq 2$, the restriction along the inclusion $\alpha$
     $$ \mathsf{Bas}^{\smallleq k} \xlongrightarrow{\alpha} \mathsf{Ext}^{\smallleq k} \xlongrightarrow{\beta} \FilFor{k}$$
     induces a morphism between derived mapping spaces
     $$ \Map_V^{\smallleq k}(\beta^* X, \beta^* Y) \longrightarrow \Map_V^{\smallleq k}(\alpha^* \beta^* X, \alpha^* \beta^* Y)$$
     which is a weak homotopy equivalence.
     \label{red_extcor}
\end{proposition}

\begin{proof}
    Following exactly the same line of reasoning as in the proof of \cref{red_ext}, we can first replace $X \to V$ by a Reedy right fibration and show that we have a weak homotopy equivalence
    \begin{equation}
    X(F) \longrightarrow \holim_{\left( \mathsf{Bas}^{\smallleq k}_{/F} \right)^\op} \: X(G)
    \label{expss_bas}
    \end{equation}
    induced by the restriction along $\beta\alpha$, for all choices of $F \in \mathsf{Ext}^{\smallleq k}$. The only cases where one needs to check something is when $F \not\in \mathsf{Bas}^{\smallleq k}$, that is, when $F$ is a linear $k$-forest.

    In light of \cref{red_ext}, we can replace the category $\mathsf{Bas}^{\smallleq k}_{/F}$ in \eqref{expss_bas} by $\mathcal{C}_F^{\mathrm{r}}$, since $X$ is a presheaf satisfying the conditions of \cref{red_ext} (ii) because of \cref{right_fibration_Segal} (c). As the inclusion of the roots $\rho(F) \to F$ defines a terminal object in $\mathcal{C}_F^{\mathrm{r}}$, equation \eqref{expss_bas} further simplifies to 
    $$ X(F) \longrightarrow X(\rho(F)),$$
    which is a weak homotopy equivalence again by \cref{right_fibration_Segal} (c).
\end{proof}

\begin{proof}[Proof of \cref{reduction_basic_forests}]
    Considering the full subcategories we have analysed, we can rewrite the map given by restricting along $\mathsf{bas} = \beta \alpha$ as the composition of
    $$ \Map_V^{\smallleq k}(X,Y) \xrightarrow{\beta^*} \Map_V^{\smallleq k}(\beta^* X, \beta^* Y) \xrightarrow{\alpha^*} \Map_V^{\smallleq k}(\alpha^* \beta^* X, \alpha \beta^* Y),$$
    and these maps are all weak homotopy equivalences by the combination of \cref{red_ext} and \cref{red_extcor}, as we wanted to show.
\end{proof}

\begin{remark}
    Before continuing, let us return to the question of the other filtration $\mathsf{fil}^{\mathsf{GW}}(-)$ on $\For$, arising from the dendroidal filtration discussed in \cite{GopplWeiss}. We claim that the arguments above lead to a more complicated analysis of $\mathsf{For}^{\smallleq k-1} \subseteq \mathsf{For}^{\smallleq k}$.

    Indeed, we can still use the Segal condition to simplify the layers. However, one needs to be careful that the $k$-corolla $C_k$ satisfies $\mathsf{fil}^{\mathsf{GW}}(C_k) = k$, and hence we also need to add this tree when computing the layers. The following more heuristic remarks might be illuminating to understand this discrepancy:
    \begin{itemize}
        \item The component filtration arises from the filtration of a right $\mathcal{P}$-module $M$ by the cardinality of the symmetric sequence it defines \textit{but} we don't filter the operad $\mathcal{P}$ in any way.
        \item The Göppl--Weiss filtration corresponds to filtering both the sequences defined by the module \textit{and} the operad via cardinality, at the same time.
    \end{itemize}

    The main point is that filtering the operad itself is not needed in the context of operadic right modules and would lead to much more complicated layers, even though the first stage is simpler. In our situation, having a clearer description of the layers is more preferable, and that's the reason for opting for the component filtration in place of the Göppl--Weiss filtration.
    
    We summarize this discussion in the table below.
    \begin{center}
    \setlength{\tabcolsep}{1em}
\begin{tabular}{ l l l } 
\textsc{Filtration}              & \textsc{First Stage}              & \textsc{Layer}       \\
\hline\hline
\small Component filtration    & \small $\Tree$ (Hard)           & \small Adds $k \cdot \eta$ (Easy)               \\
\hline
\small Göppl--Weiss filtration & \small $\mathsf{\Delta}$ (Easy) & \small Adds $k \cdot \eta$ and $C_k$ (Hard)    \\
\hline
\end{tabular}
\end{center}

\label{table_filtrations}

\end{remark}

\subsubsection{A latching-matching decomposition of the layers}

In the previous section we reduced the computation of the layers of the forest tower to computing the homotopy fiber of the map
$$ \Map_V^{\smallleq k}(\mathsf{bas}^* X, \mathsf{bas}^*Y) \longrightarrow \Map_V^{\smallleq k-1}(X,Y),$$
where on the right hand side we see $X$ and $Y$ as presheaves on $\FilFor{k-1}$, and on the left hand side as presheaves on $\mathsf{Bas}^{\smallleq k}$, which is the full subcategory of $\FilFor{k}$ on $\FilFor{k-1}$ and the forest $k \cdot \eta$. From this simplification one would expect that the behaviour of the layer is concentrated in the object $k \cdot \eta$ and its automorphisms $\Sigma_k$ given by the symmetric group on $k$ letters. 

Our first objective will be to make this intuition precise via the following proposition. The exact definitions that we will need will be made further on in the text.

\begin{notation}
    For $k \geq 1$ and $V \in \fSpaces$, we will abbreviate the space $V(k \cdot \eta)$ to $V(k)$. This space defines a presheaf on the the symmetric group $\Sigma_k$.
\end{notation}

\begin{proposition}
    Let $V$ be a reduced dendroidal Segal space, and suppose $X, Y \in \fSpaces_{/V}$. Then there are latching-matching functors
   $$ \partial_k X : [2] \longrightarrow \mathsf{Psh}_{V(k)}(\Sigma_k) \hspace{2em} \mathrm{and} \hspace{2em} \partial_k' X : [1] \longrightarrow \mathsf{Psh}_{V(k)}(\Sigma_k)$$
   such that the following commutative diagram
\begin{equation}
\begin{tikzcd}[cramped]
	{\Map_V^{\smallleq k}(\mathsf{bas}^*X, \mathsf{bas}^*Y) } & {\Map_{V(k)}^{[2] \times \Sigma_k}(\partial_k X, \partial_k Y)} \\
	{\Map_V^{\smallleq k-1}(X, Y)} & {\Map_{V(k)}^{[1] \times \Sigma_k}(\partial_k' X, \partial_k' Y)}
	\arrow["\partial_k", from=1-1, to=1-2]
	\arrow[from=1-1, to=2-1]
	\arrow[from=1-2, to=2-2]
	\arrow["{\partial'_k}", from=2-1, to=2-2]
\end{tikzcd}
\label{latc_match_square}
\end{equation}
is homotopy cartesian. Here the right hand side is computed using the projective model structure, and the right vertical map is induced by $[1] \to [2]$ sending $0$ to $0$, and $1$ to $2$.
\label{latch_match_decomp}
\end{proposition}

We begin by properly defining the operators $\partial_k$ and $\partial'_k$. These are functors which, given the truncated presheaf $X^{\smallleq k-1}$, associate to it the space $X(k)$, but also its latching and matching objects, which respectively form the initial and terminal way of extending $X$ to $k \cdot \eta$ by only using the information up to filtration $k-1$.

\begin{definition}
    Suppose $V$ is a dendroidal space and let $X \in \fSpaces_{/V}$. For $k \geq 2$ we make the following definitions:
    \begin{enumerate}[label=(\roman*)]
        \item The \textit{latching object} of $X$ is the evaluation at $k \cdot \eta$ of the homotopical left Kan extension of $X^{\smallleq k-1}$ along the inclusion $\FilFor{k-1} \to \mathsf{Bas}^{\smallleq k}$. More explicitly, it is given by
        $$ \mathsf{Latch}_k(X)  = \hocolim_{\left( \FilFor{k-1}_{k \cdot \eta/} \right)^\op} X(G)$$
        where $G \in  \FilFor{k-1}_{k \cdot \eta /}$. There is a canonical map $\mathsf{Latch}_k(X) \to X(k)$.
        \item The \textit{matching object} of $X$ is the evaluation at $k \cdot \eta$ of the homotopical right Kan extension of $X^{\smallleq k-1}$ along the inclusion $\FilFor{k-1} \to \mathsf{Bas}^{\smallleq k}$. More explicitly, it is given by
        $$ \mathsf{Match}_k(X)  = \holim_{\left( \FilFor{k-1}_{/ k \cdot \eta} \right)^\op} X(G)$$
        where $G \in  \FilFor{k-1}_{/ k \cdot \eta}$. There is a canonical map $X(k) \to \mathsf{Match}_k(X)$.
        \item We define the functors
        $$ \partial_k X : [2] \longrightarrow \mathsf{Psh}_{V(k)}^{[2]}(\Sigma_k) \hspace{2em} \mathrm{and} \hspace{2em} \partial'_k X : [1] \longrightarrow \mathsf{Psh}_{V(k)}(\Sigma_k)$$
        to be respectively given by the latching-matching diagrams 
        $$\mathsf{Latch}_k(X) \to X(k) \to \mathsf{Match}_k(X) \hspace{2em} \mathrm{and} \hspace{2em} \mathsf{Latch}_k(X) \to \mathsf{Match}_k(X).$$
    \end{enumerate}
\end{definition}

\begin{remark}
    One thing we didn't explain in the definition is how $\mathsf{Latch}_k(X)$ and $\mathsf{Match}_k(X)$ define spaces over $V(k)$. For the latching object this is easy, as one just take the composite $\mathsf{Latch}_k (X) \to X(k) \to V(k)$. For the matching object, we note that the matching map $V(k) \simeq \mathsf{Match}_k(V)$ is a weak equivalence if $V$ is a dendroidal space. This follows easily from the formula for the matching object in \cref{matching_latching_simplification}, for instance.
\end{remark}

Now that we have defined all the terms in \cref{latch_match_decomp}, we are ready to prove this result.

\begin{proof}[Proof of \cref{latch_match_decomp}:]
    We use the same strategy as in the proof of \cref{forest_tower}. Here we should be careful about our notation: we are not actually working with $X, Y$ when computing these mapping spaces, but only considering their restrictions to $\mathsf{Bas}^{\smallleq k}$.
    
    Fixing $Y \in \mathsf{Psh}_V(\mathsf{Bas}^{\smallleq k})$ throughout, let $\mathcal{A}$ be the class of $X \in \mathsf{Psh}_V(\mathsf{Bas}^{\smallleq k})$ for which the diagram 
    \begin{equation}
    \begin{tikzcd}[cramped]
	{\Map_V^{\smallleq k}(X, Y) } & {\Map_{V(k)}^{[2] \times \Sigma_k}(\partial_k X, \partial_k Y)} \\
	{\Map_V^{\smallleq k-1}(X, Y)} & {\Map_{V(k)}^{[1] \times \Sigma_k}(\partial_k' X, \partial_k' Y)}
	\arrow["\partial_k", from=1-1, to=1-2]
	\arrow[from=1-1, to=2-1]
	\arrow[from=1-2, to=2-2]
	\arrow["{\partial'_k}", from=2-1, to=2-2]
\end{tikzcd}
\label{VID}
\end{equation}
    is homotopy cartesian. By \cref{closure_classes} we just need to show that $\mathcal{A}$ satisfies the following properties:
    \begin{itemize}
        \item \textit{$\mathcal{A}$ is invariant under projective weak equivalences}: everything is defined in terms of derived limits, colimits and mapping spaces, so this property is clear.
        \item \textit{Representable presheaves are in $\mathcal{A}$}: this involves more work, so we defer its verification to \cref{representables}.
        \item \textit{$\mathcal{A}$ is closed under arbitrary disjoint unions and homotopy pushouts}: this is not difficult to check explicitly by verifying on homotopy fibers, and using that the derived mapping spaces sends homotopy pushouts and arbitrary coproducts to homotopy pullbacks and products respectively. This argument is completely written down in Lemmas 4.2.1 and 4.2.2 of \cite{GopplWeiss}.
    \end{itemize}    
\end{proof}

For handling the representable case, we will first prove two auxiliary lemmas. The first one is quite straightforward, and it provides a convenient description of the homotopy fiber of the rightmost map in \eqref{latc_match_square} in terms of the total homotopy fiber of a square whose vertices are easier to understand.

\begin{remark}
    From now on we will write the mapping space for $\Sigma_k$-equivariant maps over $V(k)$ as $\mathsf{Map}_{V}^{\Sigma_k}(-,-)$ instead of $\mathsf{Map}_{V(k)}^{\Sigma_k}(-,-)$.
\end{remark}

\begin{lemma}
    Let $V$ be a dendroidal Segal space, and suppose $X, Y \in \fSpaces_{/V}$. For $k \geq 2$, the homotopy fiber of the map
    $$\Map_V^{[2] \times \Sigma_k}(\partial_k X,\partial_k Y) \longrightarrow \Map_V^{[1] \times \Sigma_k}(\partial'_k X, \partial'_k Y)$$
    is weakly homotopy equivalent to the total homotopy fiber of the square
\begin{equation}
\begin{tikzcd}[cramped]
	{\Map^{\Sigma_k}_V(X(k), Y(k))} & {\Map^{\Sigma_k}_V(\mathsf{Latch}_k(X), Y(k))} \\
	{\Map^{\Sigma_k}_V(X(k), \mathsf{Match}_k(Y))} & {\Map^{\Sigma_k}_V(\mathsf{Latch}_k(X), \mathsf{Match}_k(Y))}
	\arrow[from=1-1, to=1-2]
	\arrow[from=1-1, to=2-1]
	\arrow[from=1-2, to=2-2]
	\arrow[from=2-1, to=2-2]
\end{tikzcd}
\label{tot_fiber_ho}
\end{equation}
obtained using the respective latching and matching maps.
\label{useful_square}
\end{lemma}

\begin{proof}
    This can be easily verified by, for instance, computing the vertical homotopy fibers of the diagram.
\end{proof}

\begin{remark}
    Putting together \cref{latch_match_decomp} and \cref{useful_square} together, we can form the following interpretation of the fiber $\mathsf{Fib}_k(X,Y)$ for $k \geq 2$. The information of filtration $k-1$ determines the solid $\Sigma_k$-equivariant diagram over $V(k)$ shown below 
\[\begin{tikzcd}[cramped]
	{\mathsf{Latch}_k(X)} & {X(k)} & {\mathsf{Match}_k(X)} \\
	{\mathsf{Latch}_k(Y)} & {Y(k)} & {\mathsf{Match}_k(Y).}
	\arrow[from=1-1, to=1-2]
	\arrow[from=1-1, to=2-1]
	\arrow[from=1-2, to=1-3]
	\arrow[dashed, from=1-2, to=2-2]
	\arrow[from=1-3, to=2-3]
	\arrow[from=2-1, to=2-2]
	\arrow[from=2-2, to=2-3]
\end{tikzcd}\]
    The space $\mathsf{Fib}_k(X,Y)$ will then be the space of $\Sigma_k$-equivariant maps $X(k) \to Y(k)$ making all the squares above commute up to homotopy.
\end{remark}

The second lemma computes the latching and matching objects for the special instance where the presheaf in question is representable.

\begin{lemma}
    Suppose $k \geq 2$ and let $F \in \mathsf{Bas}^{\smallleq k}$. Then the following hold:
    \begin{enumerate}[label=(\roman*)]
        \item If $F = k \cdot \eta$, then $\mathsf{Latch}_k(F) = \emptyset$.
        \item If $F \neq k \cdot \eta$, then the latching map $\mathsf{Latch}_k(F) \to F(k)$ is a weak homotopy equivalence.
        \item If $F \neq k \cdot \eta$ admitting no morphisms of the form $k \cdot \eta \to F$, then $\mathsf{Latch}_k(F) = F(k) = \emptyset$.
    \end{enumerate}
    \label{latch_match_representables}
\end{lemma}

\begin{proof}
    For (i), one notices that the latching object of $X$ is defined as a homotopy colimit of objects of the form $\For(G, k \cdot \eta)$, indexed by morphisms $k \cdot \eta \to G$ with $G \in \FilFor{k-1}$. If this colimit was non-empty, then in particular there would be forest morphisms $k \cdot \eta \to G \to k \cdot \eta$, which is easily checked to be impossible due to $G$ having filtration strictly lower than $k$.
    
    For statement (ii), we simplify the homotopy colimit defining the latching object. Firstly, notice that we have equivalences
    $$ \mathsf{Latch}_k(F) = \hocolim_{\left( \FilFor{k-1}_{k \cdot \eta/} \right)^\op} \For^{\smallleq k-1}(G, F) \simeq \hocolim_{G \in \For^{\smallleq k-1}} (k \cdot \eta \to G \to F),$$
    and this last homotopy colimit is exactly computing the classifying space for the opposite of the category $\FilFor{k-1}_{k \cdot \eta/ /F}$. In order to simplify it, we can consider the functor
    $$ \mathcal{C} \longrightarrow \FilFor{k-1}_{k \cdot \eta/ /F}$$
    associated to the inclusion of the full subcategory $\mathcal{C} \subseteq \FilFor{k-1}_{k \cdot \eta/ /F}$ on objects of the form $k \cdot \eta = k \cdot \eta \to F$. Notice that $\mathcal{C}$ is equivalent to the set $F(k)$, seen here as a discrete category. Moreover the inclusion $\mathcal{C} \subseteq \FilFor{k-1}_{k \cdot \eta/ /F}$ admits a right adjoint which sends $k \cdot \eta \xrightarrow{g} G \xrightarrow{f} F$ to $k \cdot \eta = k \cdot \eta \xrightarrow{fg} F$. Putting these remarks together and on account of \cref{adjoint_functor_simplifying}, we have the weak homotopy equivalences
    $$ \mathsf{Latch}_k(F) \simeq \hocolim_{G \in \For^{\smallleq k-1}} (k \cdot \eta \to G \to F) \simeq \hocolim_{\mathcal{C}^\op} \ast \simeq F(k), $$
    as we wanted to show.

    Finally, for (iii), the lack of maps of the form $k \cdot \eta \to F$ immediately gives us that $F(k)$ is empty, and the same also holds for the latching space $\mathsf{Latch}_k(F)$ since it maps to $F(k)$.
\end{proof}

\begin{lemma}
    Using the notation of the proof of \cref{latch_match_decomp}, the class $\mathcal{A}$ contains the representable presheaves on $\mathsf{Bas}^{\smallleq k}$.
    \label{representables}
\end{lemma}

\begin{proof}
    In this proof we will write $F$ instead of $X$ since we are dealing with representables. We divide the proof into three cases.
    \begin{enumerate}[label=\protect\circled{\arabic*}]
        \item \textit{$F = k \cdot \eta$}: in this scenario the left vertical map of \eqref{VID}, which is
    $$ \Map_V^{\smallleq k}(F, \mathsf{bas}^*Y) \longrightarrow \Map_V^{\smallleq k-1}(F,Y),$$
    will coincide with the matching map $Y(k) \to \mathsf{Match}_k(Y)$, so we just need to show that the same also holds for the right hand map of the same diagram, which we will analyse using \eqref{tot_fiber_ho}. 
    
    By \cref{latch_match_representables} (i) we know that $\mathsf{Latch}_k(F) = \emptyset$ . Using the description in \cref{useful_square}, we conclude that the homotopy fiber of the right hand side map of the square of \cref{latch_match_decomp} is also described as the total homotopy fiber of
\[\begin{tikzcd}[cramped]
	{\Map^{\Sigma_k}_V(F(k), Y(k))} & {\Map^{\Sigma_k}_V(\emptyset, Y(k))} \\
	{\Map^{\Sigma_k}_V(F(k), \mathsf{Match}_k(Y))} & {\Map^{\Sigma_k}_V(\emptyset, \mathsf{Match}_k(Y)).}
	\arrow[from=1-1, to=1-2]
	\arrow[from=1-1, to=2-1]
	\arrow["\simeq", from=1-2, to=2-2]
	\arrow[from=2-1, to=2-2]
\end{tikzcd}\]
As $F(k \cdot \eta) \cong \Sigma_k$, we can replace $F(k \cdot \eta)$ with $\ast$ and remove the $\Sigma_k$-equivariant condition on the left vertical map above. Consequently, the total homotopy fiber of this square will be the matching map $Y(k) \to \mathsf{Match}_k(Y)$). This shows that \eqref{latc_match_square} is homotopy cartesian, concluding the verification of the first case.
    \item \textit{$F \neq k \cdot \eta$ and there is some map $k \cdot \eta \to F$}: the left map from \eqref{latc_match_square} becomes a weak homotopy equivalence since $F \in \FilFor{k-1}$, so we will be done once we check that the rightmost map also has trivial homotopy fiber.

    By \cref{latch_match_representables} (ii), we have a weak homotopy equivalence $\mathsf{Latch}_k(F) \simeq F(k)$. Consequently, the horizontal maps of \eqref{tot_fiber_ho} are weak homotopy equivalences, and therefore this square is homotopy cartesian. The content of \cref{useful_square} now yields the desired conclusion. 

    \item \textit{$F \neq k \cdot \eta$ and there is no map $k \cdot \eta \to F$}: here one can apply \cref{latch_match_representables} to explicitly describe the square \eqref{tot_fiber_ho}, and the result immediately follows.
    \end{enumerate}
    \end{proof}

    Chaining together \cref{first_layer_desc} and \cref{latch_match_decomp} yields the following corollary, which is our final description of the layers
    \begin{corollary}
        Let $V$ be a reduced dendroidal Segal space, and suppose $X, Y \in \fSpaces_{/V}$. For $k \geq 2$, the layer $\mathsf{Fib}_k(X,Y)$ can be described as the total homotopy fiber of
        $$ 
        \begin{tikzcd}[cramped]
	{\Map^{\Sigma_k}_V(X(k), Y(k))} & {\Map^{\Sigma_k}_V(\mathsf{Latch}_k(X), Y(k))} \\
	{\Map^{\Sigma_k}_V(X(k), \mathsf{Match}_k(Y))} & {\Map^{\Sigma_k}_V(\mathsf{Latch}_k(X), \mathsf{Match}_k(Y)),}
	\arrow[from=1-1, to=1-2]
	\arrow[from=1-1, to=2-1]
	\arrow[from=1-2, to=2-2]
	\arrow[from=2-1, to=2-2]
\end{tikzcd}
        $$
        where the vertical and horizontal arrows are induced by the latching map $\mathsf{Latch}_k(X) \to X(k)$, and the matching map $Y(k) \to \mathsf{Match}_k(Y)$. 
        \label{better_layer_description}
    \end{corollary}

    \begin{remark}
        An analysis of operadic right modules and its associated tower is also present in the literature in the recent work of Krannich and Kupers \cite{KrannichKupersRightModules}. Their article deals with the similar problems to the ones we address in this article, and in particular a description of the layers like the one in \cref{better_layer_description} is also present in their work. For instance, a similar latching-matching decomposition for the layers in described in \citep[Sec.\ 4.4]{KrannichKupersRightModules}.
        
        The way this article differs from their work is in the methods we use. Whereas Krannich--Kupers work in a model-independent fashion, we instead opt to use the dendroidal and forest formalisms. We believe that our approach has the advantage of being quite explicit and some of the technical problems are sometimes easier to solve in the dendroidal setting. An example of this is the reduction problems we discuss at the start of Section 4.4, which in the model-independent strategy relies on the more intricate machinery of assembly for generalized $\infty$-operads.
    \end{remark}

    We will use the corollary above in the next section to give connectivity estimates for the layers of the forest tower in the case of embedding calculus. Such estimates are easy to provide once one has some understanding over the maps defining the square, which in turn comes down to getting a grasp of the latching and matching maps
    $$ \mathsf{Latch}_k(X) \longrightarrow X(k) \hspace{2em} \mathrm{and} \hspace{2em} Y(k) \longrightarrow \mathsf{Match}_k(Y).$$
    We will end this section by showing that one can already simplify the categories $\FilFor{k-1}_{k \cdot \eta / }$ and $\FilFor{k-1}_{/ k \cdot \eta }$ in order to simplify the latching and matching spaces, respectively.

    \begin{notation}
        Let $k \geq 2$. We define the category of \textit{k-factorizations} 
        $$\mathsf{Fact}_k \subseteq \FilFor{k-1}_{k \cdot \eta /}$$ 
        to be the full subcategory spanned by the objects $f \colon k \cdot \eta \to G$ with $G$ having no unary vertices and having exactly $k$ leaves. 
    \end{notation}

        The category $\mathsf{Fact}_k$ is supposed to parameterize all ways of factoring the arity $k$ component of a module $M$ via its operadic right action. 
        
        Indeed, recall that the module structure is detected at the level of forests by the root faces, as we saw in the definition of the nerve functor $N_{\mathcal{P}}$. Given $f \colon k \cdot \eta \to G$ in $\mathsf{Fact}_k$, we can repeatedly apply the root faces to $G$ -- which from the module viewpoint corresponds to repeatedly applying the right actions -- until we reach a forest of edges. This will be exactly the forest of leaves of $G$ and will of the form $k \cdot \eta$ due to the definition of $\mathsf{Fact}_k$, and is the forest parameterizing the arity $k$ component of a module.
    
    \begin{lemma}
        Let $V$ be a reduced forest space and $X \in \fSpaces_{/V}$. For $k \geq 2$, the following descriptions of the latching and matching objects hold:
        \begin{enumerate}[label=(\roman*)]
            \item The matching object $\mathsf{Match}_k(X)$ can be written as a cubical diagram
            $$ \mathsf{Match}_k(X) \xlongrightarrow{\simeq} \holim_{I \subsetneq \underline{k}} X(I \cdot \eta)$$
            over the opposite of the poset of proper non-empty subsets of $\underline{k}$.
            \item The latching object $\mathsf{Latch}_k(X)$ can be written as a homotopy colimit
            $$ \mathsf{Latch}_k(X) \xlongrightarrow{\simeq} \hocolim_{\mathsf{Fact}_k^\op} X(G)$$
            over the category of $k$-factorizations $\mathsf{Fact}_k$. 
        \end{enumerate}
        \label{matching_latching_simplification}
    \end{lemma}
    \begin{proof}
        For the case of the matching object, we note that the inclusion of the full subcategory of $\FilFor{k-1}_{/k \cdot \eta}$ spanned by the forests $\ell \cdot \eta$ with $\ell < k$ is homotopy right cofinal. This is a consequence of the existence of the factorization below for the vertical map
\[\begin{tikzcd}[cramped]
	G & {\ell \cdot \eta} \\
	{k \cdot\eta}
	\arrow["f", from=1-1, to=1-2]
	\arrow[from=1-1, to=2-1]
	\arrow[from=1-2, to=2-1]
\end{tikzcd}\]
This is in turn due to the image of $G$ under such $f$ factoring through a strict subforest of $k \cdot \eta$, because $G$ has filtration smaller than $k$.

        For the latching object, we follow the same strategy as in the proof of \cref{simplify_extension}, which consists in considering an extra full subcategory
        $$ \mathsf{Fact}_k \subseteq \mathcal{C}_{k \cdot \eta} \subseteq \FilFor{k-1}_{k \cdot \eta /}$$
        where in $\mathcal{C}_{k \cdot \eta}$ we allow for the objects $k \cdot \eta \to G$ with $G$ with $k$ leaves but possibly with unary vertices. We first show that the $\mathcal{C}_{k \cdot \eta} \subseteq \FilFor{k-1}_{k \cdot \eta/}$ admits a right adjoint. This follows by from the decomposition of map $f \colon k \cdot \eta \to G$ as below
\[\begin{tikzcd}[cramped]
	& {k \cdot \eta} \\
	{G'} & G
	\arrow["f'"', from=1-2, to=2-1]
	\arrow["f", from=1-2, to=2-2]
	\arrow[from=2-1, to=2-2]
\end{tikzcd}\]
    with $f' \colon k \cdot \eta \to G'$ as required. This can be constructed by contracting all inner edges of $G$ without changing the image of $f$, until it satisfies the required property (this again uses that we are working with closed forests). For the reduction from $\mathcal{C}_{k \cdot \eta}$ to $\mathsf{Fact}_k$, we use the same argument and criteria as we used at the end of the proof of \cref{simplify_extension}.
    \end{proof}

\begin{remark}
    In the next section we will apply our results to the specific case of embedding calculus, but we would like to offer some remarks on the more general setting we have been working on. 
    
    A basic question about the layers of a tower concerns the computation of the connectivity of its layers, in our case the spaces $\mathsf{Fib}_k(X,Y)$ for each $k \geq 2$. In order to understand how this can be computed, we recall that for $A$ is a $d$-dimensional CW complex and and $B$ is a $m$-connected space, the mapping space $\mathsf{Map}(A,X)$ is $(m-d)$-connected. We can apply this to our situation via \cref{better_layer_description}, which implies that
    $$ \mathrm{con}\left( \mathsf{Fib}_k(X,Y) \right) \geq  \mathrm{con} \left( Y(k) \to \mathsf{Match}_k(Y) \right) - \mathrm{hodim}_{\mathsf{Latch}_k(X)}\left( X(k) \right) -1$$
    where $\mathrm{hodim}_{\mathsf{Latch}_k(Y)}(Y(k))$ is the relative homotopical dimension of the pair of spaces $(Y(k), \mathsf{Latch}_k(Y))$.

    The connectivity of $Y(k) \longrightarrow \mathsf{Match}_k(Y)$ is quite easy to compute in general. Indeed, we have just shown that $\mathsf{Match}_k(Y)$ is constructed from a certain cubical diagram built out of the set $\underline{k}$, and $X(k)$ is the evaluation of $X$ at the initial object of this cube. Thus, we can compute the connectivity of $Y(k) \to \mathsf{Match}_k(Y)$ via the cubical version of Blakers--Massey \citep[Thm. 2.4]{GoodwillieII}, which is quite optimal.

    The question for the latching object is more intricate. To see this, we will first consider the case of operads instead of modules, as this is well understood in the literature. We can define $\mathsf{Latch}_k(\mathcal{P})$ for any operad $\mathcal{P}$, which will be the space that captures all the factorizations of a arity $k$ operation by operations of lower arity via the composition of $\mathcal{P}$. The space we are interested in
    \begin{equation}
     \mathrm{hocof} \left( \mathsf{Latch}_k(\mathcal{P}) \longrightarrow \mathcal{P}(k)\right)
     \label{cofiber}
    \end{equation}
    has already appeared in the literature on operads and, up to a suspension, coincides with the space of indecomposable operations of $\mathcal{P}(k)$. By \citep[Sec. 4]{HeutsMoerdijkPartition} and \citep[Rem. 4.3]{HeutsKoszul}, this is also the space of arity $k$ cooperations in the cooperad $B \mathcal{P}$ coming from the bar construction of $\mathcal{P}$. The results in \cite{HeutsKoszul} can then be used to give estimates for the homology of \eqref{cofiber}, which is what one would want.

    A similar theory and description for operadic right modules seems possible -- see \cite{MalinTaggartKoszul} for work in this direction -- and would lead to a completely algebraic approach to understanding \eqref{cofiber}. However, in this article we will instead analyse the cofiber in the context of embedding calculus via a geometric approach using the Fulton--MacPherson compactification of smooth manifolds.
    \label{general_connectivity}
    \end{remark}

\subsection{Embedding calculus via the forest tower}

In the last section of this article we will apply the description of the layers given in \cref{better_layer_description} to the particular case of manifold calculus in order to obtain connectivity estimates for these spaces. We will also draw some conclusions about the connectivity of the space of embeddings $\mathrm{Emb}(M,N)$ from these computations, by appealing to the convergence results of Goodwillie--Klein.

\begin{notation}
For the remainder of this section we make the following conventions:
\begin{itemize}
    \item From now on we fix smooth manifolds $M^d$ and $N^{d+n}$, with $d\geq 1$ and $n \geq 0$.
    \item As we have already done in \cref{T1_example}, we will need to understand $\mathbb{E}_N$ as an $\mathbb{E}_d^{\mathsf{fr}}$-module, and we will still denote this module by $\mathbb{E}_N$.
    \item We are always working within the model of forest spaces from this points on, and so we will drop any references to the nerve functors to make the notation clearer.
    \item To make our notation coincide with that of manifold calculus, we will define
    $$ T_k \mathrm{Emb}(M,N) \coloneqq \Map^{\smallleq k}_{\mathbb{E}^{\mathsf{fr}}_d}\left( \mathbb{E}_M, \mathbb{E}_N \right) \hspace{1em} \mathrm{and} \hspace{1em} T_\infty \mathrm{Emb}(M,N) \coloneqq \Map_{\mathbb{E}^{\mathsf{fr}}_d}\left( \mathbb{E}_M, \mathbb{E}_N \right)$$
    for each $k \geq 1$.
\end{itemize}  
\end{notation}

We need to first address an important technical question before applying the results we have proved. We have assumed throughout Section 4.3 that we are working over a \textit{reduced} dendroidal Segal space, and the reduced condition was essential for replacing the inclusion $\FilFor{k-1} \subseteq \FilFor{k}$ to $\mathsf{Bas}^{\smallleq k} \subseteq \FilFor{k}$. However, the operad $\mathbb{E}^{\mathsf{fr}}_d$ does not define a reduced dendroidal Segal space and is not even complete: indeed, the underlying category of this operad is the monoid $\mathrm{Emb}(D^d, D^d)$, which is weakly equivalent to the orthogonal group $O(d)$.

Recalling from \cref{completion_invariance} that the contravariant model structure is invariant under completion of the base forest space, our strategy for solving this issue will be to construct a completion 
$$ \mathbb{E}_d^{\mathsf{fr}} \longrightarrow \mathbb{E}_{d, hO}^{\mathsf{fr}}$$
of the framed $\mathbb{E}_d$-operad, which will in particular be a reduced dendroidal Segal space. Moreover, at the level of modules a similar construction also holds, which will in this case produce a map
$$ \mathbb{E}_M \longrightarrow \mathbb{E}_{M, hO}$$
playing a similar role to the completion above.

\begin{remark}
    This issue with the reduced condition also appears when working with Lurie's model for $\infty$-operads. In that context, this is solved in Sections 2.3.2 and 2.3.4 of \cite{LurieHA} by the machinery of assembly for generalized $\infty$-operads, and the exact case which is of interest to us is explained in detail in \citep[Prop. 2.2]{HorelKrannichKupers}. Our approach using the dendroidal formalism also appears in the recent preprint \cite{BoavidaCiriciHorelEquivariant}.
\end{remark}

\begin{notation}
    Let $F$ be a forest and $d \geq 1$. We define $O(F)$ to be the group $O(d)^{\times E(F)}$, where $E(F)$ is the set of edges of $F$. Note that any forest morphism $F \to G$ will induce a map $O(F) \to O(G)$ via the function on the set of edges $E(F) \to E(G)$.
\end{notation}

We begin by observing that, for any tree $T$, the space $\mathbb{E}_d^{\mathsf{fr}}(T)$ admits an $O(T)$-action, constructed as follows:
 \begin{itemize}
        \item Let $T$ be a $k$-corolla and $(g_0, \ldots, g_k) \in O(T)$, with $g_0$ indexed by the root edge. Then the action of $(g_0, \ldots g_k)$ on $f \in \mathbb{E}_d^{\mathsf{fr}}(T)$ has component $1 \leq i \leq k$ given by $g_0 f_i g_i^{-1}$.
        \item For a general tree $T$, we can use the Segal property to decompose $\mathbb{E}_d^{\mathsf{fr}}(T)$ as a product of corollas, and we use the previous case to construct the $O(T)$-action.
        \item The face and degeneracy maps for $\mathbb{E}_d^{\mathsf{fr}}$ are constructed using the group multiplication and the identity of $O(d)$.
    \end{itemize}

More generally, for the right $\mathbb{E}_d^{\mathsf{fr}}$-module $\mathbb{E}_M$ we can endow $\mathbb{E}_M(F)$ with an $O(F)$-action for any forest $F$. This is done in the same fashion as above.

\begin{definition}
Let $d \geq 1$. We define the dendroidal space $\mathbb{E}_{d, hO}^{\mathsf{fr}}$ via the formula
$$ \mathbb{E}_{d,hO}^{\mathsf{fr}}(T) \coloneqq \mathbb{E}_d^{\mathsf{fr}}(T)_{hO(T)}$$
 for any tree $T$. For a general smooth manifold $M^d$, we similarly construct a forest space $\mathbb{E}_{M, hO}$ as
 $$ \mathbb{E}_{M, hO}(F) \coloneqq \mathbb{E}_M(F)_{hO(F)}$$
 for any forest $F$.
\end{definition}

It follows that there exist dendroidal and forest maps $\psi_d \colon \mathbb{E}_d^{\mathsf{fr}} \to \mathbb{E}_{d, hO}^{\mathsf{fr}}$ and $\psi_M \colon \mathbb{E}_M \to \mathbb{E}_{M, hO}$. It will be important to observe for the next proof that these morphisms fit into fiber sequences of spaces
\begin{equation}
\mathbb{E}_d^{\mathsf{fr}}(T) \to \mathbb{E}_{d, hO}^{\mathsf{fr}}(T) \to BO(T) \hspace{1.5em} \mathrm{and} \hspace{1.5em} \mathbb{E}_M(F) \to \mathbb{E}_{M, hO}(F) \to BO(F)
\label{fibre_sequence}
\end{equation}
for any tree $T$ and forest $F$. The first statement of the next proposition is part of the content of \citep[Lem. 2.2]{BoavidaCiriciHorelEquivariant}.

\begin{proposition}
    Let $d \geq 1$ and consider a smooth manifold $M^d$.
    Then:
    \begin{enumerate}[label=(\roman*)]
        \item  The dendroidal space $\mathbb{E}^{\mathsf{fr}}_{d, hO}$ is a complete dendroidal Segal space, and the map $\psi_d \colon \mathbb{E}^{\mathsf{fr}}_{d} \to \mathbb{E}_{d, hO}^{\mathsf{fr}}$ is a complete weak equivalence.
        \item There is a homotopy pullback square of forest spaces
\[\begin{tikzcd}[cramped]
	{\mathbb{E}_{M}} & {\mathbb{E}_{M, hO}} \\
	{\mathbb{E}_{d}^{\mathsf{fr}}} & {\mathbb{E}_{d, hO}^{\mathsf{fr}}}
	\arrow["{\psi_M}", from=1-1, to=1-2]
	\arrow[from=1-1, to=2-1]
	\arrow[from=1-2, to=2-2]
	\arrow["{\psi_d}", from=2-1, to=2-2]
\end{tikzcd}\]
and the right vertical map $\mathbb{E}_{M, hO} \to \mathbb{E}_{d, hO}^{\mathsf{fr}}$ defines a right fibration.
    \end{enumerate}   
    \label{completing_E_d}
\end{proposition}

\begin{proof}
We will begin by showing that $\mathbb{E}_{d, hO}^{\mathsf{fr}}$ defines a complete dendroidal Segal space. Given any tree $T$, we have a commutative diagram
\[\begin{tikzcd}[cramped]
	{\mathbb{E}_d^{\mathsf{fr}}(T)} & {\mathbb{E}_{d, hO}^{\mathsf{fr}}(T)} & {BO(T)} \\
	{\mathbb{E}_d^{\mathsf{fr}}(\mathsf{Sp}(T))} & {\mathbb{E}_{d, hO}^{\mathsf{fr}}(\mathsf{Sp}(T))} & {\prod_{v \in V(T)} BO(C_{\lvert v \rvert})}
	\arrow[from=1-1, to=1-2]
	\arrow["\simeq"', from=1-1, to=2-1]
	\arrow[from=1-2, to=1-3]
	\arrow[from=1-2, to=2-2]
	\arrow["\simeq"', from=1-3, to=2-3]
	\arrow[from=2-1, to=2-2]
	\arrow[from=2-2, to=2-3]
\end{tikzcd}\]
with the horizontal lines being fiber sequences. Then the leftmost and rightmost maps are weak homotopy equivalences, since $\mathbb{E}_d^{\mathsf{fr}}$ satisfies the Segal condition and by definition of $O(T)$ respectively. Therefore the middle map is also such an equivalence, which shows the Segal condition for $\mathbb{E}_{d, hO}^{\mathsf{fr}}$. As for completeness, the underlying simplicial space of $\mathbb{E}_{d, hO}^{\mathsf{fr}}$ is the constant simplicial space on $BO(d)$, and it is clear that this is complete.

We now prove that $\psi_d \colon \mathbb{E}^{\mathsf{fr}}_{d} \to \mathbb{E}_{d, hO}^{\mathsf{fr}}$ defines a complete weak equivalence, for which by \citep[Cor. 12.34]{HeutsMoerdijkDendroidal} it will suffice to check that $\psi_d$ is essentially surjective and fully faithful. The former is immediate, and the latter holds since the square 
\[\begin{tikzcd}[cramped]
	{\mathbb{E}_d^{\mathsf{fr}}(C_k)} & {\mathbb{E}_{d, hO}^{\mathsf{fr}}(C_k)} \\
	{\prod_{i=0}^k \mathbb{E}_d^{\mathsf{fr}}(\eta)} & {\prod_{i=0}^k \mathbb{E}_{d, hO}^{\mathsf{fr}}(\eta)}
	\arrow[from=1-1, to=1-2]
	\arrow[from=1-1, to=2-1]
	\arrow[from=1-2, to=2-2]
	\arrow[from=2-1, to=2-2]
\end{tikzcd}\]
is homotopy cartesian for all $k \geq 0$. This is due to the left and right bottom vertices being contractible and $BO(C_k)$ respectively, together with the fiber sequence for the homotopy orbits.

For (ii), one can verify that the square is homotopy cartesian by taking the homotopy fibers of the horizontal maps. Finally, consider for any forest $F$ the commutative cube
\[\begin{tikzcd}[cramped]
	{\mathbb{E}_M(F)} && {\mathbb{E}_{M,hO}(F)} \\
	& {\mathbb{E}_M(\rho(F))} && {\mathbb{E}_{M,hO}(\rho(F))} \\
	{\mathbb{E}^{\mathsf{fr}}_d(F)} && {\mathbb{E}^{\mathsf{fr}}_{d, hO}(F)} \\
	& {\mathbb{E}^{\mathsf{fr}}_d(\rho(F))} && {\mathbb{E}^{\mathsf{fr}}_{d, hO}(\rho(F))}
	\arrow[from=1-1, to=1-3]
	\arrow[from=1-1, to=2-2]
	\arrow[from=1-1, to=3-1]
	\arrow[from=1-3, to=2-4]
	\arrow[from=1-3, to=3-3]
	\arrow[from=2-2, to=2-4, crossing over]
	\arrow[from=2-4, to=4-4]
	\arrow[from=3-1, to=3-3]
	\arrow[from=3-1, to=4-2]
	\arrow[from=3-3, to=4-4]
	\arrow[from=4-2, to=4-4]
	\arrow[from=2-2, to=4-2, crossing over]
\end{tikzcd}\]
The left and front faces are homotopy cartesian, respectively because $\mathbb{E}_M \to \mathbb{E}_d$ is a right fibration and by taking the horizontal fibers of the front face. By a similar reason, the back face is also a homotopy pullback square, and the same holds for the right face by taking the vertical fibers of the cube and using that $\mathbb{E}_d(F) \to \mathbb{E}_{d,hO}(F)$ is $\pi_0$-surjective. By the characterization in \cref{characterization_ctv_fibrant}, this implies that $\mathbb{E}_{M, hO} \to \mathbb{E}_{d, hO}$ is a right fibration.
\end{proof}

\begin{corollary}
    For $1 \leq k \leq \infty$, there is a weak homotopy equivalence 
    $$ T_k \mathrm{Emb}(M,N) \simeq \Map^{\smallleq k}_{\mathbb{E}_d^{\mathsf{fr}}}(\mathbb{E}_M, \mathbb{E}_N) \simeq \Map^{\smallleq k}_{\mathbb{E}_{d,hO}^{\mathsf{fr}}}\left( \mathbb{E}_{M, hO}, \mathbb{E}_{N, hO} \right)$$
    induced by the pullback along $\psi_d \colon \mathbb{E}_d^{\mathsf{fr}} \to \mathbb{E}_{d, hO}^{\mathsf{fr}}$.
    \label{completion_mapping_spaces}
\end{corollary}

\begin{proof}
    By \cref{completing_E_d} (ii) we know that $\psi^* \mathbb{E}_{M, hO}^{\mathsf{fr}} \simeq \mathbb{E}_M^{\mathsf{fr}}$ and $\psi^* \mathbb{E}_{N, hO}^{\mathsf{fr}} \simeq \mathbb{E}_N^{\mathsf{fr}}$, and this also holds for the $k$-truncated versions of these presheaves. By \cref{completion_invariance} and \cref{completing_E_d} (i), the pullback along $\psi^*$ yields weak homotopy equivalences
    $$ \Map^{\smallleq k}_{\mathbb{E}_d^{\mathsf{fr}}}(\mathbb{E}_M, \mathbb{E}_N) \simeq \Map^{\smallleq k}_{\mathbb{E}_d^{\mathsf{fr}}}(\psi^* \mathbb{E}_{M, hO}, \psi^* \mathbb{E}_{N, hO}) \simeq \Map^{\smallleq k}_{\mathbb{E}_{d, hO}^{\mathsf{fr}}}( \mathbb{E}_{M, hO}, \mathbb{E}_{N, hO}),$$
    where we took the non-derived version of $\psi^*$ since all the presheaves are already fibrant by \cref{completing_E_d} (ii).
\end{proof}

In order to simplify the notation from now on, we will write $\tilde{\mathbb{E}}_d$ and $\tilde{\mathbb{E}}_M$ instead of $\mathbb{E}_{d, hO}^{\mathsf{fr}}$ and $\mathbb{E}_{M, hO}^{\mathsf{fr}}$ respectively.

The benefit of \cref{completion_mapping_spaces} is twofold. On the one hand and as we have mentioned before, $\tilde{\mathbb{E}}_d$ is a reduced dendroidal Segal space, and therefore we are in a position for applying the results of the previous sections. On the other hand, taking homotopy orbits for the action by the orthogonal groups greatly simplifies $\mathbb{E}_M$, as it eliminates the immersion component.

\begin{lemma}
    Let $d \geq 1$ and consider a smooth manifold $M^d$. Then there is a weak homotopy equivalence
    $$ \tilde{\mathbb{E}}_M(k) \simeq \mathsf{Conf}_k(M)$$
    for every $k \geq 1$.
    \label{embeddings_simplifies_to_configurations}
\end{lemma}

\begin{proof}
    Taking the derivative at the center of each disk yields a weak homotopy equivalence
    \begin{equation}
     \mathrm{Emb}( \underline{k} \times D^d, M) \xlongrightarrow{\simeq} \mathsf{Conf}_k(M) \times_{M^{\times k}} \mathrm{Imm}(\underline{k} \times D^d, M)
    \label{eq1}
     \end{equation}
    for each $k \geq 1$. Moreover, there is a fiber sequence
    \begin{equation}
    O(d)^{\times k} \longrightarrow \mathrm{Imm}( \underline{k} \times D^d, M) \longrightarrow M^{\times k}
    \label{eq2}
    \end{equation}
    coming from the action of $O(d)^{\times k}$ on the total space via precomposition. By taking the homotopy orbits with respect to the $O(d)^{\times k}$-action in the right hand side of \eqref{eq1} gives us the configuration space $\mathsf{Conf}_k(M)$ because of the fibration of\eqref{eq2}. This finishes the proof.
\end{proof}

The matching and latching objects of $\tilde{\mathbb{E}}_M$ are now easy to describe.

\begin{lemma}
    Consider a smooth manifold $M^d$. For each $k \geq 2$, the matching map $\tilde{\mathbb{E}}_M(k) \to \mathsf{Match}_k(\tilde{\mathbb{E}}_M)$ is weakly equivalent to
    $$ \mathsf{Conf}_k(M) \longrightarrow \holim_{I \subsetneq \underline{k}} \mathsf{Conf}_I(M).$$
In particular, the connectivity of the matching map $\tilde{\mathbb{E}}_M(k) \to \mathsf{Match}_k(\tilde{\mathbb{E}}_M)$ is given by $k(\dim M -2 )-(\dim M -3)$. 
    \label{cube_identification}
\end{lemma}

\begin{proof}
    The identification of the matching map with the cubical limit above is a consequence of \cref{embeddings_simplifies_to_configurations}, together with the general description of the 
    matching object in \cref{matching_latching_simplification} (i). The connectivity estimate is a classic result proved using excision and the Fadell--Neuwirth fibrations for configuration spaces, see for instance \citep[Ex. 6.2.9]{MunsonVolic}.
\end{proof}

Before giving a topological model for the latching map, recall from Section 4.3 that $\mathsf{Latch}_k(\tilde{\mathbb{E}}_M)$ was defined as a certain homotopy colimit indexed over over the category of $k$-factorizations $\mathsf{Fact}_k$. This is the full subcategory of $\FilFor{k-1}_{k \cdot \eta /}$ on the maps $k \cdot \eta \to G$ with $G$ having exactly $k$ leaves and no unary vertices.

We note that the constraint on the number of leaves leads to a accompanying constraint on the kinds of morphisms
\[\begin{tikzcd}[cramped]
	& {k \cdot \eta} \\
	G && {G'}
	\arrow[from=1-2, to=2-1]
	\arrow[from=1-2, to=2-3]
	\arrow["f", from=2-1, to=2-3]
\end{tikzcd}\]
we allow in $\mathsf{Fact}_k$. For instance, $f$ cannot be a leaf face, but it is certainly allowed to be any inner face or also any root face. In fact, the essential observation for understanding $\mathsf{Latch}_k(\tilde{\mathbb{E}}_M)$ is that $\mathsf{Fact}_k^{\mathsf{op}}$ can be described\footnote{Here $\mathsf{fBin}_k$ is the category of binary forests introduced in \cref{examples_right_modules} \ref{example_FM_manifolds} for explaining the stratification of the boundary of the Fulton--MacPherson compactification of a smooth manifold.} as the full subcategory $\partial \mathsf{fBin}_k \subseteq \mathsf{fBin}_k$ spanned by the all objects except for the forest of $k$ edges.

\begin{lemma}
    Consider a smooth manifold $M^d$. For each $k \geq 2$, the latching map $\mathsf{Latch}_k(\tilde{\mathbb{E}}_M) \to \tilde{\mathbb{E}}_M(k)$ is weakly equivalent to the inclusion of the boundary of the Fulton--MacPherson compactification
    $$ \partial \mathsf{FM}_M(k) \longrightarrow \mathsf{FM}_M(k)$$
    for all $k \geq 2$. 
    \label{FM_identification}
\end{lemma}

\begin{proof}
    The proof of this result for Euclidean spaces is done in \citep[Ex. 3.1.1.3]{GopplWeiss}, and for a general manifold the same argument also holds without much change. By the discussion above, we have identifications
    $$ \mathsf{Latch}_k(\tilde{\mathbb{E}}_M) \simeq \hocolim_{\mathsf{Fact}_k^{\op}} \tilde{\mathbb{E}}_M(G) \simeq \hocolim_{\partial \mathsf{fBin}_k} \tilde{\mathbb{E}}_M(G).$$
    From the fact that $\tilde{\mathbb{E}}_M(G)$ coincides with the space $\mathsf{FM}_M(G)$ that we defined in \cref{examples_right_modules} \ref{example_FM_manifolds}, as well as all the maps $\mathsf{FM}_M(G) \to \mathsf{FM}_M(G')$ being cofibrations -- they are the inclusions of strata -- we conclude that we can replace the homotopy colimit above with the strict colimit. This then describes the boundary $\partial \mathsf{FM}_M(k)$, as we wanted to show.
\end{proof}

\begin{remark}
    By our discussion at the end of the previous section, we want to understand the cofiber of
    \begin{equation}
     \partial \mathsf{FM}_M(k) \longrightarrow \mathsf{FM}_M(k)
     \label{FM_M}
    \end{equation}
     
    and therefore one is interested in finding small dimensional models for this space. A reasonable approach to this is to find another smooth manifold $M'$ which is homotopy equivalent to $M$, but is obtained by attaching handles of index strictly smaller that $\dim M$. 
    
    In order to address this, Goodwillie and Klein \cite{GoodwillieKleinMultipleSmooth} define the \textit{handle dimension} of $M$ to be the smallest index of handles we need to attach in order to construct any model of $M$ (more formally, the handle dimension of $M$ is $\leq n$ if $M$ is the interior of a manifold built out of handles with index $\leq n$). Denoting it by $\mathrm{hdim}\: M$, we note that the inequality $\mathrm{hdim} \: M \leq \dim M$ holds, and if $M$ is compact then it is an equality. However, for general manifolds these number will differ. 
    
    The conclusion from the discussion in this remark is that, if $M'$ is a manifold realizing the handle dimension of $M$, then we can model $\mathsf{Latch}_k(\tilde{\mathbb{E}}_M) \to \tilde{\mathbb{E}}_M(k)$ via 
    $$ \partial \mathsf{FM}_{M'}(k) \longrightarrow \mathsf{FM}_{M'}(k)$$
    instead, which has homotopical dimension $k \cdot \hdim M$.
    \label{homotopical_dimension}
\end{remark}

We can now finally collect our results and give the desired connectivity estimates for the layers of the Goodwillie--Weiss tower.

\begin{corollary}
    Let $M^d$ and $N^{d+n}$ be smooth manifolds. For $k \geq 2$, the $k^{th}$ layer of the Goodwillie--Weiss tower on $\mathrm{Emb}(M,N)$, that is, the homotopy fiber of 
    $$ T_{k} \mathrm{Emb}(M,N) \longrightarrow T_{k-1}\mathrm{Emb}(M,N)$$
    has connectivity $k \left( \dim N - \hdim M -2 \right)-(\dim N -2)$. If one additionally assumes that $\dim N - \hdim M \geq 3$, then the map
    $$ T_{\infty} \mathrm{Emb}(M,N) \longrightarrow T_{k}\mathrm{Emb}(M,N)$$
    has connectivity given by at least $k \left( \dim N - \hdim M -2 \right)-(\dim N -2)$.
    \label{layers_estimates_connectivity}
\end{corollary}

\begin{proof}
    From \cref{general_connectivity} and \cref{homotopical_dimension}, the connectivity of the layer in question admits the lower bound
$$ \mathrm{con} \left( \tilde{\mathbb{E}}_N(k) \to \mathsf{Match}_k\left( \tilde{\mathbb{E}}_N \right) \right)- k \cdot \hdim M,$$
The result now follows by the connectivity estimate for the cubical diagram of configuration spaces in \cref{cube_identification}.

If $\dim N - \hdim M \geq 3$, then the slope for the linear estimate is positive, and therefore grows arbitrarily large when $k \to \infty$. Consequently, the connectivity estimates we just proved still hold for $T_{\infty} \mathrm{Emb}(M,N) \to T_k \mathrm{Emb}(M,N)$.
\end{proof}

For the next corollary, we will need to appeal to the convergence of the Goodwillie--Weiss tower, which is assured by the Goodwillie--Klein convergence theorem stated below.

\begin{theorem}
    Let $M^d$ and $N^{d+n}$ be smooth manifolds such that $\dim N - \hdim M \geq 3$ holds. Then the Goodwillie--Weiss tower converges, that is, the map 
    $$ \mathrm{Emb}(M,N) \longrightarrow T_{\infty} \mathrm{Emb}(M,N)$$
    is a weak homotopy equivalence.
\end{theorem}

\begin{remark}
    Although the inequality $\dim N - \hdim M \geq 3$ appears both in the Goodwillie--Klein theorem and in \cref{layers_estimates_connectivity}, we should stress that our results don't address any convergence properties of the Goodwillie--Weiss tower. This is a problem of a different nature altogether.
\end{remark}

\begin{corollary}
    Let $M^d$ and $N^{d+n}$ be smooth manifolds such that $\dim N - \hdim M \geq 3$ holds. Then the inclusion of immersions into embeddings
    $$ \mathrm{Emb}(M,N) \longrightarrow \mathrm{Imm}(M,N)$$
    is $(\dim N - 2 \cdot \hdim M -2)$-connected.
    \label{conn_imm_emb}
\end{corollary}

\begin{proof}
    From \cref{layers_estimates_connectivity} we know that 
    \begin{equation}
        T_{\infty} \mathrm{Emb}(M,N) \longrightarrow T_1 \mathrm{Emb}(M,N)
        \label{Tinfinity_T1}
    \end{equation}
    has connectivity given by $2(\dim N - \hdim M -2) - (\dim N -2) = \dim N - 2 \hdim M -2$. As we are assuming $\dim N - \hdim M \geq 3$, we are in the range for the convergence of the embedding tower, which therefore ensures that we have a weak homotopy equivalence
    $$ \mathrm{Emb}(M,N) \simeq T_{\infty}\mathrm{Emb}(M,N).$$
    Moreover \cref{T1_example} together with the Hirsch--Smale theorem identify the right hand side of \eqref{Tinfinity_T1} with the space of immersions $\mathrm{Imm}(M,N)$, which yields the desired result.
\end{proof}

\appendix

\section{Simplicial diagrams of model categories}

The aim of this section is to give a brief overview of \textit{simplicial diagrams of model categories}, which is crucially used in Section 4 of this work. The material in this part of the work is contained in Section 13.2 of \cite{HeutsMoerdijkDendroidal}, and our objective in recalling the constructions below is merely to make this text more self-contained.

\begin{definition}
    A \textit{simplicial diagram of model categories} is a $\mathsf{\Delta}^{\mathsf{op}}$-diagram\footnote{We didn't draw the degeneracies out of convenience.} $\mathcal{M}_{\bullet}$
    \[\begin{tikzcd}[cramped]
	{\mathcal{M}_0} & {\mathcal{M}_1} & {\mathcal{M}_2} & \cdots
	\arrow[shift right, from=1-2, to=1-1]
	\arrow[shift left, from=1-2, to=1-1]
	\arrow[shift right=2, from=1-3, to=1-2]
	\arrow[shift left=2, from=1-3, to=1-2]
	\arrow[from=1-3, to=1-2]
	\arrow[shift left=3, from=1-4, to=1-3]
	\arrow[shift left, from=1-4, to=1-3]
	\arrow[shift right=3, from=1-4, to=1-3]
	\arrow[shift right, from=1-4, to=1-3]
\end{tikzcd}\]
    of model categories $\mathcal{M}_i$, such that all the face and degeneracy maps define left Quillen functors, and the simplicial identities are satisfied up to natural isomorphism in a coherent fashion.

    Given a simplicial morphism $f \colon [n] \to [m]$ in $\mathsf{\Delta}$, we will write
    \[\begin{tikzcd}
	{\mathcal{M}_m} & {\mathcal{M}_n}.
	\arrow["{\mathcal{M}(f)_!}", shift left, from=1-1, to=1-2]
	\arrow["{\mathcal{M}(f)^*}", shift left, from=1-2, to=1-1]
\end{tikzcd}\]
    for the corresponding Quillen pair.
\end{definition}

\begin{example}
    Let $V$ be a forest space. Then an example of a simplicial diagram of model categories which is of interest to us is given by setting $\mathcal{M}_i = \left( \fSets_{/V_i} \right)_{\ctv}$ for each $i \geq 0$.
    \label{example_1}
\end{example}

Having a simplicial diagram of model categories $\mathcal{M}_{\bullet}$ at our disposal, we would like to extract from it a certain model category that in some sense captures the information contained in $\mathcal{M}_{\bullet}$. The first step in making this idea more concrete is the definition below.

\begin{definition}
    Let $\mathcal{M}_{\bullet}$ be a simplicial diagram of model categories. We define the \textit{totalization} of $\mathcal{M}_{\bullet}$ to be the category $\Tot(\mathcal{M})$ given by the following data:
    \begin{itemize}[label=$\diamond$]
        \item The objects of $\Tot(\mathcal{M})$ are pairs $(X, \alpha)$, where $X = \{ X_i : i \geq 0\}$ is a sequence of objects with $X_i \in \mathcal{M}_i$, and to each simplicial map $f \colon [m] \to [n]$ we associate a morphism 
        $$ \alpha(f) : X_n \longrightarrow \mathcal{M}(f)^*(X_m),$$
        which should be compatible with the coherence structure of $\mathcal{M}_{\bullet}$ in the following sense: given $f \colon [n] \to [m]$ and $g \colon [m] \to [\ell]$, the following diagram 
        \[\begin{tikzcd}[cramped]
	{X_{\ell}} & {\mathcal{M}(g)^*(X_m)} \\
	{\mathcal{M}(gf)^*(X_n)} & {\mathcal{M}(g)^*\mathcal{M}(f)^*(X_n)}
	\arrow["{\alpha(g)}", from=1-1, to=1-2]
	\arrow["{\alpha(gf)}"', from=1-1, to=2-1]
	\arrow["{\mathcal{M}(g)^*(\alpha(f))}", from=1-2, to=2-2]
	\arrow["\cong", from=2-1, to=2-2]
\end{tikzcd}\]
commutes, where the bottom map comes from the coherent information of $\mathcal{M}_{\bullet}$.
    \item A morphism $\varphi \colon (X, \alpha) \to (Y, \beta)$ is given by a collection of morphism $\left\{ \varphi_i \colon X_i \to Y_i : i \geq 0 \right\}$, each in $\mathcal{M}_i$, such that, for each $f\colon [n] \to [m]$ in $\mathsf{\Delta}$, the diagram below
\[\begin{tikzcd}[cramped]
	{X_n} && {Y_n} \\
	{\mathcal{M}(f)^*(X_m)} && {\mathcal{M}(f)^*(Y_m)}
	\arrow["{{\varphi_n}}", from=1-1, to=1-3]
	\arrow["{{\alpha(f)}}"', from=1-1, to=2-1]
	\arrow["{\beta(f)}", from=1-3, to=2-3]
	\arrow["{{\mathcal{M}(f)^*(\varphi_m)}}", from=2-1, to=2-3]
\end{tikzcd}\]
    commutes.
    \end{itemize}
\end{definition}

\begin{remark}
    In the previous definition we assumed that $\mathcal{M}_{\bullet}$ is a simplicial diagram of model categories, but it is clear from the definition that the fact that $\mathcal{M}_i$ is a model category doesn't play a role.
\end{remark}

\begin{example}
Continuing with \cref{example_1}, there is an equivalence of categories
 \[\begin{tikzcd}
	{\Tot(\mathcal{M})} & {\fSpaces_{/V} }.
	\arrow["{L}", shift left, from=1-1, to=1-2]
	\arrow["{R}", shift left, from=1-2, to=1-1]
\end{tikzcd}\]

Here $L$ sends a pair $(X, \alpha)$ to the forest map $X \to V$ given by gathering all the $X_i \to V_i$ together, and the naturality in the simplicial degree is encoded in $\alpha$. As for $R$, it sends $X \to V$ to the pair $(X, \alpha)$ where the first component registers each map $X_i \to V_i$ and, for each $f\colon [m] \to [n]$ in $\mathsf{\Delta}$, the structure map $\alpha(f)$ is the pullback map
$$ \alpha(f) : X_n \longrightarrow X_m \times_{V_m} V_n.$$
\label{example_2}
\end{example}

The next result tells us that the totalization of a simplicial diagram of model categories also inherits a model structure, which we will refer to as the \textit{projective model structure} and write $\Tot(\mathcal{M})_{\mathsf{P}}$.

\begin{proposition}
    Let $\mathcal{M}_{\bullet}$ be a simplicial diagram of model categories. Then $\Tot(\mathcal{M})$ carries a model structure where a morphism $\varphi \colon (X, \alpha) \to (Y, \beta)$ is a weak equivalence/fibration if each morphism $\varphi_i$ is a weak equivalence/fibration in $\mathcal{M}_i$, for each $i \geq 0$. This model category is left proper if each $\mathcal{M}_i$ is left proper.
\end{proposition}

\begin{proof}
    This model structure is obtained via right transfer of the model structure on the product $\prod_{i \geq 0} \mathcal{M}_i$ equipped with the termwise cofibrations/fibrations/weak equivalences. For more details see Proposition 13.13 in \cite{HeutsMoerdijkDendroidal}.
\end{proof}

Similarly to what happens in the situation of simplicial presheaves, one can define a version of a Reedy model structure for the totalization $\Tot(\mathcal{M})$. In order to do so, we first define, for each object $(X, \alpha) \in \Tot(\mathcal{M})$ and $k \geq 0$, a version of the \textit{matching} and \textit{latching} objects of $(X, \alpha)$, respectively defined as
$$ M_k (X, \alpha) = \lim_{f \colon [n] \to [k]} \mathcal{M}(f)^*(X_n) \hspace{2em} \text{and} \hspace{2em} L_k (X, \alpha) = \colim_{f \colon [k] \to [n]} \mathcal{M}(f)_!(X_n),$$
where the limit is taken along the proper surjections and the colimit along the proper injections. We can then define \textit{Reedy cofibrations} and \textit{Reedy fibrations} in a manner entirely analogous to the classical Reedy set up.

\begin{theorem}
    Let $\mathcal{M}_{\bullet}$ be a simplicial diagram of model categories. Then the following hold:
    \begin{enumerate}[label=(\alph*)]
        \item The Reedy fibrations and Reedy cofibrations defined above are respectively the classes of fibrations and cofibrations of a model structure $\Tot(\mathcal{M})_{\mathsf{R}}$, which is left proper if each $\mathcal{M}_i$ is also left proper. Moreover, the weak equivalences coincide with the projective weak equivalences.
        \item There is a Quillen equivalence
        \[\begin{tikzcd}
	{\Tot(\mathcal{M})_{\mathsf{P}}} & {\Tot(\mathcal{M})_{\mathsf{R}}}
	\arrow["{\mathrm{id}_!}", shift left, from=1-1, to=1-2]
	\arrow["{\mathrm{id}^*}", shift left, from=1-2, to=1-1]
\end{tikzcd}\]
between the projective model structure $\Tot(\mathcal{M})_{\mathsf{P}}$ and the Reedy model structure $\Tot(\mathcal{M})_{\mathsf{R}}$.
    \end{enumerate}
    \label{Reedy_simplicial_diagrams}
\end{theorem}

\begin{proof}
    This is the content of Theorem 13.14 in \cite{HeutsMoerdijkDendroidal}.
\end{proof}

It will be convenient to define an extra model structure on $\Tot(\mathcal{M})$.

\begin{proposition}
    Suppose $\mathcal{M}_{\bullet}$ is a simplicial diagram of model categories such that each component $\mathcal{M}_i$ is cofibrantly generated and left proper, with cofibrations being monomorphisms. Then there exists a model structure on $\Tot(\mathcal{M})$ such that:
    \begin{enumerate}[label = (\alph*)]
        \item The cofibrations coincide with the Reedy cofibrations.
        \item An object $(X, \alpha)$ is fibrant if it is Reedy fibrant and, for each morphism $f \colon [n] \to [m]$ in $\mathsf{\Delta}$, the map $\alpha(f)$ is a weak equivalence in $\mathcal{M}_m$.
    \end{enumerate}
    We will denote this model category by $\Tot(\mathcal{M})_{\mathsf{w}}$.
    \label{RSC_Totalization}
\end{proposition}

\begin{proof}
    This model structure is a left Bousfield localization of the Reedy model structure $\Tot(\mathcal{M})_{\mathsf{R}}$, and the morphisms we localize at are easily obtained by unravelling what $\alpha(f)$ being a weak equivalence means. See Proposition 13.20 \cite{HeutsMoerdijkDendroidal} for more details. 
\end{proof}

We will end this appendix with some results concerning the comparison between distinct simplicial diagrams of model categories, so that it will in turn induce a Quillen pair for some of the model structures we have thus far. The required notion is that of a Quillen pair of simplicial diagrams.

\begin{definition}
    Let $\mathcal{M}_{\bullet}$ and $\mathcal{N}_{\bullet}$ be simplicial diagrams of model categories. We define a \textit{Quillen adjunction} between $\mathcal{M}_{\bullet}$ and $\mathcal{N}_{\bullet}$ to be a collection of Quillen pairs
    \[\begin{tikzcd}
	{\mathcal{M}_i} & {\mathcal{N}_i}
	\arrow["{F_i}", shift left, from=1-1, to=1-2]
	\arrow["{G_i}", shift left, from=1-2, to=1-1]
\end{tikzcd}\]
for each $i \geq 0$, such that these are compatible the the simplicial structure of $\mathcal{M}_{\bullet}$ and $\mathcal{N}_{\bullet}$ up to isomorphism. We call such such a Quillen adjunction a \textit{Quillen equivalence} if each pair $(F_i, G_i)$ is a Quillen equivalence.

\end{definition}

The following result explains the interaction between Quillen adjunction of simplicial diagrams and Quillen adjunctions of their totalizations.

\begin{theorem}
    Let $\mathcal{M}_{\bullet}$ and $\mathcal{N}_{\bullet}$ be simplicial diagrams of model categories, together with a Quillen pair $F \colon \mathcal{M}_{\bullet} \rightleftarrows \mathcal{N}_{\bullet} \colon G$. Then:
    \begin{enumerate}[label=(\alph*)]
        \item There are Quillen adjunctions
    \[\begin{tikzcd}
	{\Tot(\mathcal{M})_{\mathsf{P}}} & {\Tot(\mathcal{N})_{\mathsf{P}}}
	\arrow["{\Tot(F)}", shift left, from=1-1, to=1-2]
	\arrow["{\Tot(G)}", shift left, from=1-2, to=1-1]
\end{tikzcd} \hspace{1em} and \hspace{1em} \begin{tikzcd}
	{\Tot(\mathcal{M})_{\mathsf{R}}} & {\Tot(\mathcal{N})_{\mathsf{R}}}
	\arrow["{\Tot(F)}", shift left, from=1-1, to=1-2]
	\arrow["{\Tot(G)}", shift left, from=1-2, to=1-1]
\end{tikzcd}\]
    which are Quillen equivalences if each pair $(F,G)$ is a Quillen equivalence of simplicial diagrams of model categories.
        \item If the conditions of \cref{RSC_Totalization} are met by both $\mathcal{M}_{\bullet}$ and $\mathcal{N}_{\bullet}$, then there is a Quillen adjunction
        \[\begin{tikzcd}
	{\Tot(\mathcal{M})_{\mathsf{w}}} & {\Tot(\mathcal{N})_{\mathsf{w}}}
	\arrow["{\Tot(F)}", shift left, from=1-1, to=1-2]
	\arrow["{\Tot(G)}", shift left, from=1-2, to=1-1]
\end{tikzcd}\]
which is a Quillen equivalence if $(F,G)$ is a Quillen equivalence of simplicial diagrams of model categories.
    \end{enumerate}
    \label{comparison_totalizations}
\end{theorem}
 
\begin{proof}
    This is Proposition 13.21 in \cite{HeutsMoerdijkDendroidal}.
\end{proof}

 \bibliographystyle{alpha}
\bibliography{bibliography.bib}
\end{document}